\documentstyle{amltd}
\begin{document}
\annalsline{158}{2003}
\received{December 23, 1999}
\startingpage{1}
\def\bye{\end{document}}
 \font\tenrm=cmr10
\def\ritem#1{\item[{\rm #1}]}
\input amssym.def
\input amssym.tex
\def\joinrel{\mathrel{\mkern-4mu}}
\def\relbar{\mathrel{\smash-}}
\def\lrar{\relbar\joinrel\relbar\joinrel\relbar\joinrel\relbar\joinrel\relbar\joinrel\rightarrow}
\def\srar{\relbar\joinrel\relbar\joinrel\relbar\joinrel\rightarrow}
\def\rsrar{\relbar\joinrel\relbar\joinrel\rightarrow}
\def\slar{\leftarrow\joinrel\relbar\joinrel\relbar\joinrel\relbar}
\def\lline{\relbar\joinrel\relbar\joinrel\relbar\joinrel\relbar}
\def\vvlar{\relbar\joinrel\relbar\joinrel\relbar\joinrel\relbar\joinrel
\relbar\joinrel\relbar\joinrel\relbar\joinrel\relbar\joinrel
\relbar\joinrel\relbar\joinrel\relbar\joinrel\relbar\joinrel\relbar\joinrel
\relbar\joinrel\relbar\joinrel\relbar\joinrel\relbar\joinrel\relbar\joinrel\relbar\joinrel\relbar\joinrel\relbar\joinrel\relbar\joinrel\relbar\joinrel\relbar\joinrel\relbar\joinrel\rightarrow}
\def\dline{\begin{array}{c}\lline\\[-11pt]\lline\end{array}}
\def\eqref#1{(\ref{#1})}
\def\ddownarrow{\begin{array}{c}\big\downarrow\\[-16pt]\big\downarrow\end{array}}
\def\stck#1{\stackrel{#1}{\srar}}
\def\shtck#1{\stackrel{#1}{\longrightarrow}}
\def\dlar{\srar\hskip-13pt\to}
\def\scs#1{\big\downarrow{\scriptstyle{#1}}}
\def\scr#1{\scriptstyle{#1}}
\def\hsm{\hskip-8pt}
\def\speqnua#1{\speqnu{{\rm A.{#1}}}}
\catcode`\@=11
\font\twelvemsb=msbm10 scaled 1100
\font\tenmsb=msbm10
\font\ninemsb=msbm10 scaled 800
\newfam\msbfam
\textfont\msbfam=\twelvemsb  \scriptfont\msbfam=\ninemsb
  \scriptscriptfont\msbfam=\ninemsb
\def\msb@{\hexnumber@\msbfam}
\def\Bbb{\relax\ifmmode\let\next\Bbb@\else
 \def\next{\errmessage{Use \string\Bbb\space only in math
mode}}\fi\next}
\def\Bbb@#1{{\Bbb@@{#1}}}
\def\Bbb@@#1{\fam\msbfam#1}
\catcode`\@=12

 \catcode`\@=11
\font\twelveeuf=eufm10 scaled 1100
\font\teneuf=eufm10
\font\nineeuf=eufm7 scaled 1100
\newfam\euffam
\textfont\euffam=\twelveeuf  \scriptfont\euffam=\teneuf
  \scriptscriptfont\euffam=\nineeuf
\def\euf@{\hexnumber@\euffam}
\def\frak{\relax\ifmmode\let\next\frak@\else
 \def\next{\errmessage{Use \string\frak\space only in math
mode}}\fi\next}
\def\frak@#1{{\frak@@{#1}}}
\def\frak@@#1{\fam\euffam#1}
\catcode`\@=12
\def\operatorname#1{\mathop{\rm #1}\nolimits} 
\newcommand{\Fil}{\operatorname{Fil}}
\newcommand{\Tor}{\operatorname{Tor}}
\newcommand{\Aut}{\operatorname{Aut}}
\newcommand{\Hom}{\operatorname{Hom}}
\newcommand{\End}{\operatorname{End}}
\newcommand{\END}{\mathbold{\operatorname{End}}}
\newcommand{\Ext}{\operatorname{Ext}}
\newcommand{\Gal}{\operatorname{Gal}}
\newcommand{\Spec}{\operatorname{Spec}}
\newcommand{\p}{^{^\wedge}_p}
\newcommand{\f}{\hat{\phantom{i}}}
\newcommand{\s}{_{\mathbold{\cdot}}}
\newcommand{\cs}{^{\mathbold{\cdot}}}
\newcommand{\cd}{\operatorname{cd}}
\newcommand{\im}{\operatorname{im}}
\newcommand{\coim}{\operatorname{coim}}
\newcommand{\coker}{\operatorname{coker}}
\newcommand{\gr}{\operatorname{gr}}
\newcommand{\Br}{\operatorname{Br}}
\newcommand{\TF}{\operatorname{TF}}
\newcommand{\TC}{\operatorname{TC}}
\newcommand{\TR}{\operatorname{TR}}
\newcommand{\THH}{\operatorname{THH}}
\renewcommand{\TH}{\operatorname{TH}}
\newcommand{\HH}{\operatorname{HH}}
\newcommand{\N}{{\Bbb N}}
\newcommand{\Z}{{\Bbb Z}}
\newcommand{\Q}{{\Bbb Q}}
\newcommand{\R}{{\Bbb R}}
\newcommand{\C}{{\Bbb C}}
\newcommand{\Zp}{{\Bbb Z}_p}
\newcommand{\Qp}{{\Bbb Q}_p}
\newcommand{\F}{{\Bbb F}}
\newcommand{\Fp}{{\Bbb F}_p}
\newcommand{\Sm}{{\rm   Sm}}
\newcommand{\Sch}{{\rm   Sch}}
\newcommand{\Aff}{{\rm   Aff}}
\newcommand{\et}{{{\rm  \'et}}}
\newcommand{\Zar}{{{\rm   Zar}}}
\newcommand{\Nis}{{{\rm   Nis}}}
\newcommand{\cont}{{{\rm   cont}}}
\newcommand{\op}{{{\rm   op}}}
\newcommand{\tw}{\operatorname{tw}}
\def\operatornamewithlimits#1{{\displaystyle\mathop{\rm #1}\limits}}
\def\underset#1#2{\begin{array}{c}\noalign{\vskip-6pt}{#2}\\[-8pt]{#1}\end{array}}
\newcommand{\rlim}[3]{\mathchoice
        {{\underset{\scriptstyle#1}{\varinjlim}}^{#2}#3}
        {\varinjlim^{#2}#3}
        {}{}}
\def\varprojlim{{\displaystyle\lim_{\longleftarrow}}}
 \newcommand{\llim}[3]
{\mathchoice
       {{\underset{\scriptstyle#1}{\varprojlim}}^{#2}#3}
         {\varprojlim^{#2}#3}
       {}{}}

\newcommand{\hcl}[2]{\mathchoice{{\displaystyle\mathop{\rm holim}_{\underset{ 
        \scriptstyle#1}{\longrightarrow}} #2}}
        {{\displaystyle\operatornamewithlimits{holim}_{\longrightarrow}#2}}
        {}{}}

\newcommand{\hl}[2]{\mathchoice{{\displaystyle\mathop{\rm holim}_{\underset{ 
        \scriptstyle#1}{\longleftarrow}} #2}}
         {{\displaystyle\operatornamewithlimits{holim}_{\longleftarrow}#2}}
        {}{}}

\newcommand{\ul}[1]{\underline{#1}}
\newcommand{\cy}{{{\rm   cy}}}
\newcommand{\ob}{\operatorname{ob}}
\newcommand{\id}{\operatorname{id}}
\newcommand{\pr}{\operatorname{pr}}
\newcommand{\can}{\operatorname{can}}
\newcommand{\tr}{\operatorname{tr}}
\newcommand{\ev}{\operatorname{ev}}
\newcommand{\sd}{\operatorname{sd}}
\newcommand{\res}{\operatorname{res}}
\renewcommand{\:}{\colon}

\newcommand{\Cyl}{\operatorname{Cyl}}
\def\xto#1{\stackrel{#1}{\longrightarrow}}
\def\xleftarrow#1{\stackrel{#1}{\longleftarrow}}
\newcommand{\e}{\varepsilon}
\newcommand{\tate}{\displaystyle{\hat{{\Bbb H}}}}
\newcommand{\borel}{\displaystyle{{\Bbb H}\hskip1pt_{\mathbold{\cdot}}}}
\newcommand{\coborel}{\displaystyle{{\Bbb H}\,{}^{\mathbold{\cdot}}}}
\newcommand{\red}{\operatorname{red}}
\newcommand{\ann}{\operatorname{ann}}
\newcommand{\supp}{\operatorname{supp}}
\newcommand{\gp}{\operatorname{gp}}
\newcommand{\susp}{\operatorname{susp}}


 \font\emi= cmmi10 scaled 1700 
\title{On the {\emi K}-theory of local fields} 
 \def\titleheadline#1{\def\one{#1}\ifx\one\empty\else
\gdef\thetitle{{\frenchspacing%
\let\\ \relax
{#1}}}\fi}
\newif\ifshort
\def\shortname#1{\global\shorttrue\xdef
\theauthors{{\eightsc\uppercase{#1}}}}
\let\shorttitle\titleheadline
\shorttitle{ \eightsc\uppercase{On the} {\eightpoint \it K}\eightsc\uppercase{-theory of local fields} }

  \acknowledgements{The first named author was supported in part by NSF Grant and the
Alfred P. Sloan Foundation.  The second named author was supported in part by The American Institute of
Mathematics.}
 \twoauthors{Lars Hesselholt}{Ib Madsen}
 \institutions{Massachusetts Institute of Technology, Cambridge, Massachusetts\\
{\eightpoint {\it E-mail address\/}: larsh@math.mit.edu}\\
\vglue6pt 
Matematisk Institut, Aarhus Universitet, Denmark\\
{\eightpoint {\it E-mail address\/}: imadsen@imf.au.dk}}
 
\def\sni#1{\smallbreak\noindent{#1}. }
\def\ssni#1{\vglue-1pt\noindent\hskip18pt {#1}.}
\centerline{\bf Contents}
\smallbreak\noindent Introduction
\sni{1} Topological Hochschild homology and localization
\sni{2} The homotopy groups of $T(A|K)$
\sni{3} The de Rham-Witt complex and ${\rm TR}^{\mathbold{\cdot}}_\ast(A|K;p)$
\sni{4} Tate cohomology and the Tate spectrum
\sni{5} The Tate spectral sequence for $T(A|K)$
\sni{6} The pro-system ${\rm TR}^{\mathbold{\cdot}}_\ast(A | K;p,\Bbb Z/p^v)$
\smallbreak\noindent Appendix A. Truncated polynomial algebras
\smallbreak\noindent References
\vglue12pt

\intro

In this paper we establish a connection between the Quillen $K$-theory
of certain local fields and the de Rham-Witt complex of their rings of
integers with logarithmic poles at the maximal ideal. The fields $K$
we consider are complete discrete valuation fields of characteristic
zero with perfect residue field $k$ of characteristic $p>2$. When $K$
contains the $p^v$-th roots of unity, the relationship between the
$K$-theory with $\Z/p^v$-coefficients and the de Rham-Witt complex can
be described by a sequence
$$\cdots\to K_*(K,\Z/p^v)\to
W\,\omega_{(A,M)}^*\otimes S_{\Z/p^v}(\mu_{p^v})\xto{1-F}
W\,\omega_{(A,M)}^*\otimes S_{\Z/p^v}(\mu_{p^v})\xto{\partial}\cdots
$$
which is exact in degrees $\geq 1$. Here $A={\cal O}_K$ is the
valuation ring and $W\,\omega_{(A,M)}^*$ is the de Rham-Witt complex of
$A$ with log poles at the maximal ideal. The factor
$S_{\Z/p^v}(\mu_{p^v})$ is the symmetric algebra of $\mu_{p^v}$ 
considered as a $\Z/p^v$-module located in degree two. Using this
sequence, we evaluate the $K$-theory with $\Z/p^v$-coefficients of
$K$. The result, which is valid also if $K$ does not contain the
$p^v$-th roots of unity, verifies the Lichtenbaum-Quillen conjecture
for $K$,~\cite{lichtenbaum},~\cite{quillen2}: 

\specialnumber{A}
\proclaim{Theorem}\label{thmA} There are natural isomorphisms
for $s\geq 1${\rm ,} 
\begin{eqnarray*}
K_{2s}(K,\Z/p^v)&=&H^0(K,\mu_{p^v}^{\otimes s})\oplus
H^2(K,\mu_{p^v}^{\otimes(s+1)}),\\
K_{2s-1}(K,\Z/p^v)&=&H^1(K,\mu_{p^v}^{\otimes s}).\\
\end{eqnarray*}
\endproclaim

The Galois cohomology on the right can be effectively calculated when
$k$ is finite, or equivalently, when $K$ is a finite extension of
$\Qp$,~\cite{soule}. For $m$ prime to $p$,
$$K_i(K,\Z/m)=K_i(k,\Z/m)\oplus K_{i-1}(k,\Z/m)$$
by Gabber-Suslin,~\cite{suslin1}, and for $k$ finite, the $K$-groups
on the right are known by Quillen,~\cite{Quillen1}.

For any linear category with cofibrations and weak equivalences in the
sense of~\cite{w}, one has the cyclotomic trace
$$\tr\:K({\cal C})\to\TC({\cal C};p)$$
from $K$-theory to topological cyclic homology,~\cite{dm}. It
coincides in the case of the exact category of finitely generated
projective modules over a ring with the original definition in
\cite{bhm}. The exact sequence above and Theorem \ref{thmA} are based
upon calculations of $\TC_*({\cal C};p,\Z/p^v)$ for certain
categories associated with the field $K$. Let $A={\cal O}_K$ be the
valuation ring in $K$, and let ${\cal P}_A$ be the category of
finitely generated projective $A$-modules. We consider three
categories with cofibrations and weak equivalences: the category
$C^b_z({\cal P}_A)$ of bounded complexes in ${\cal P}_A$ with
homology isomorphisms as weak equivalences, the subcategory with
cofibrations and weak equivalences $C_z^b({\cal P}_A)^q$ of
complexes whose homology is torsion, and the category
$C^b_q({\cal P}_A)$ of bounded complexes in ${\cal P}_A$ with
rational homology isomorphisms as weak equivalences. One then has a
cofibration sequence of $K$-theory spectra
$$K(C_z^b({\cal P}_A)^q)\xto{i^!}
K(C_z^b({\cal P}_A))\xto{j}
K(C_q^b({\cal P}_A))\xto{\partial}
\Sigma K(C_z^b({\cal P}_A)^q),$$
and by Waldhausen's approximation theorem, the terms in this
sequence may be identified with the $K$-theory of the exact categories
${\cal P}_k$, ${\cal P}_A$ and ${\cal P}_K$. The associated
long-exact sequence of homotopy groups is the localization sequence of
\cite{quillen},
$$\dots\to K_i(k)\xto{i^!}K_i(A)\xto{j_*}K_i(K)
\xto{\partial}K_{i-1}(k)\to\dots$$
The map $\partial$ is a split surjection by~\cite{harrissegal}. We
show in Section~1.5 below that one has a similar
cofibration sequence of topological cyclic homology spectra
$$\TC(C_z^b({\cal P}_A)^q;p)\xto{i^!}
\TC(C_z^b({\cal P}_A);p)\xto{j}
\TC(C_q^b({\cal P}_A);p)\xto{\partial}
\Sigma\TC(C_z^b({\cal P}_A)^q;p),$$
and again Waldhausen's approximation theorem allows us to identify the
first two terms on the left with the topological cyclic homology of the
exact categories ${\cal P}_k$ and ${\cal P}_A$. But the third
term is different from the topological cyclic homology of
${\cal P}_K$. We write
$$\TC(A|K;p)=\TC(C_q^b({\cal P}_A);p),$$
and we then have a map of cofibration sequences
\def\sto#1{\stackrel{#1}{\longrightarrow}}
$$\begin{array}{ccccccc}
{K(k)} &\sto{i!}&
{K(A)}&\sto{j_*} &
{K(K)}&\sto{\partial} &\Sigma K(k)
 \\
\big\downarrow{\scriptstyle{\rm tr}}&&\big\downarrow{\scriptstyle{\rm tr}}&&\big\downarrow{\scriptstyle{\rm
tr}}&&\downarrow{\scriptstyle{\rm tr}}\\
 {\TC(k;p)}&\sto{i!} &
{\TC(A;p)} &\sto{j_*} &
{\TC(A|K;p)} &\sto{\partial} &
{\Sigma\TC(k;p).} \end{array}$$
By~\cite[Th.\ D]{hm}, the first two vertical maps from the left
induce isomorphisms of homotopy groups with $\Z/p^v$-coefficients in
degrees $\geq 0$. It follows that the remaining two vertical maps
induce isomorphisms of homotopy groups with $\Z/p^v$-coefficients in
degrees $\geq 1$,
$$\tr\:K_i(K,\Z/p^v)\xto{\sim}\TC_i(A|K;p,\Z/p^v),\hskip4mm i\geq 1.$$
It is the right-hand side we evaluate.  

The spectrum $\TC({\cal C};p)$ is defined as the homotopy fixed
points of an operator called Frobenius on another spectrum
$\TR({\cal C};p)$; so there is a natural cofibration sequence
$$\TC({\cal C};p)\to\TR({\cal C};p)\xto{1-F}\TR({\cal C};p)
\to\Sigma\TC({\cal C};p).$$
The spectrum $\TR({\cal C};p)$, in turn, is the homotopy limit of a
pro-spectrum $\TR\cs({\cal C};p)$, its homotopy groups given by the
Milnor sequence
$$0\to{\displaystyle{\lim_{{\longleftarrow}\atop R}}}^1  {\TR_{s+1}\cs({\cal C};p)}
\to\TR_s({\cal C};p)
\to{\displaystyle{\lim_{{\longleftarrow}\atop R}}}{\TR_s\cs({\cal C};p)}\to 0,$$
and there are maps of pro-spectra
\begin{eqnarray*}
&&F\:\TR^n({\cal C};p)\to\TR^{n-1}({\cal C};p),\\
&&V\:\TR^{n-1}({\cal C};p)\to\TR^n({\cal C};p).\\
\end{eqnarray*}
The spectrum $\TR^1({\cal C};p)$ is the topological Hochschild
homology $T({\cal C})$. It has an action by the circle group
${\Bbb T}$ and the higher levels in the pro-system by definition are
the fixed sets of the cyclic subgroups of ${\Bbb T}$ of $p$-power
order,
$$\TR^n({\cal C};p)=T({\cal C})^{C_{p^{n-1}}}.$$
The map $F$ is the obvious inclusion and $V$ is the accompanying
transfer. The structure map $R$ in the pro-system is harder to define
and uses the so-called cyclotomic structure of $T({\cal C})$; see
Section~1.1 below.

The homotopy groups $\TR_*\cs(A|K;p)$ of this pro-spectrum
with its various operators have a rich algebraic structure which we
now describe. The description involves the notion of a log differential
graded ring from~\cite{k}. A log ring $(R,M)$ is a ring $R$ with a
pre-log structure, defined as a map of monoids
$$\alpha\:M\to (R,\,\cdot\,),$$
and a log differential  graded ring $(E^*,M)$ is a differential graded
ring $E^*$, a pre-log structure $\alpha\:M\to E^0$ and a map of monoids
$d\log\:M\to (E^1,+)$ which satisfies $d\circ d\log=0$ and
$d\alpha(a)=\alpha(a)d\log a$ for all $a\in M$. There is a universal
log differential graded ring with underlying log ring $(R,M)$: the de
Rham complex with log poles $\omega_{(R,M)}^*$. 

The groups $\TR_*^1(A|K;p)$ form a log differential graded ring whose
underlying log ring is $A={\cal O}_K$ with the {\it canonical}
pre-log structure given by the inclusion
$$\alpha\:M=A\cap K^\times\to A.$$
We show that the canonical map
$$\omega_{(A,M)}^*\to\TR_*^1(A|K;p)$$
is an isomorphism in degrees $\leq 2$ and that the left-hand side is
uniquely divisible in degrees $\geq 2$. We do not know a
{\it natural} description of the higher homotopy groups, but we do
for the homotopy groups with $\Z/p$-coefficients. The Bockstein
$$\TR_2^1(A|K;p,\Z/p)\xto{\sim}{}_p\!\TR_1^1(A|K;p)$$
is an isomorphism, and we let $\kappa$ be the element on the left
which corresponds to the class $d\log(-p)$ on the right. The abstract
structure of the groups $\TR_*^1(A;p)$ was determined in~\cite{lm}. We
use this calculation in Section~2 below to show:

\specialnumber{B}\proclaim{Theorem} \label{thmB} There is a natural isomorphism of
log differential graded rings
$$\omega_{(A,M)}^*\otimes_{\Z}S_{\Fp}\{\kappa\}
\xto{\sim}\TR^1_*(A|K;p,\Z/p),$$
where $d\kappa=\kappa d\log(-p)$.
\endproclaim

The higher levels $\TR_*^n(A|K;p)$ are also log differential graded
rings. The underlying log ring is the ring of Witt vectors $W_n(A)$
with the pre-log structure
$$M\xto{\alpha}A\to W_n(A),$$
where the right-hand map is the multiplicative section
$\ul{a}_n=(a,0,\dots,0)$. The maps $R$, $F$ and $V$ extend the
restriction, Frobenius and Verschiebung of Witt vectors. Moreover,
$$F\:\TR_*^n(A|K;p)\to\TR_*^{n-1}(A|K;p)$$
is a map of pro-log graded rings, which satisfies
\begin{eqnarray*}
Fd\log_na&\hsm=\hsm&d\log_{n-1}a,\hskip6.6mm
\hbox{for all $a\in M=A\cap K^\times$,}\\
Fd\ul{a}_n&\hsm=\hsm&\ul{a}_{n-1}^{p-1}d\ul{a}_{n-1},\hskip5mm
\hbox{for all $a\in A$,} 
\end{eqnarray*}
and $V$ is a map of pro-graded modules over the pro-graded ring
$\TR_*\cs(A|K;p)$,
$$V\:F^*\TR_*^{n-1}(A|K;p)\to\TR_*^n(A|K;p).$$
Finally,
$$
FdV=d,\qquad
FV=p. 
$$
The algebraic structure described here makes sense for any log ring
$(R,M)$, and we show that there exists a universal example: the de
Rham-Witt pro-complex with log poles $W\s\,\omega_{(R,M)}^*$. For log
rings of characteristic $p>0$, a different construction has been given
by Hyodo-Kato,~\cite{hyodokato}.

We show in Section~3 below that the canonical map
$$W\s\,\omega_{(A,M)}^*\to\TR_*\cs(A|K;p)$$
is an isomorphism in degrees $\leq 2$ and that the left-hand side is
uniquely divisible in degrees $\geq 2$. Suppose that $\mu_{p^v}\subset
K$. We then have a map
$$S_{\Z/p^v}(\mu_{p^v})\to\TR_*\cs(A|K;p,\Z/p^v)$$
which takes $\zeta\in\mu_{p^v}$ to the associated Bott element defined
as the unique element with image $d\log\s\zeta$ under the Bockstein
$$\TR_2\cs(A|K;p,\Z/p^v)\xto{\sim}{}_{p^v}\!\TR_1\cs(A|K;p).$$
The following is the main theorem of this paper.

\specialnumber{C}\proclaim{Theorem}\label{thmC} Suppose that $\mu_{p^v}\subset
K$. Then the canonical map
$$W\s\,\omega_{(A,M)}^*\otimes_{\Z}S_{\Z/p^v}(\mu_{p^v})
\xto{\sim}\TR_*\cs(A|K;p,\Z/p^v)$$
is a pro\/{\rm -}\/isomorphism.
\endproclaim

We explain the structure of the groups in the theorem for $v=1$; the
structure for $v>1$ is unknown. Let $E\s^*$ stand for either side of
the statement above. The group $E_n^i$ has a natural descending
filtration of length $n$ given by
$$\Fil^sE_n^i=V^sE_{n-s}^i+dV^sE_{n-s}^{i-1}\subset E_n^i,
\hskip4mm 0\leq s<n.$$
There is a natural $k$-vector space structure on $E_n^i$, and for all
$0\leq s<n$ and all $i\geq 0$,
$$\dim_k\gr^sE_n^i=e_K,$$
the absolute ramification index of $K$. In particular, the domain and
range of the map in the statement are abstractly isomorphic.

The main theorem implies that for $s\geq 0$,
\begin{eqnarray*}
\TC_{2s}(A|K;p,\Z/p^v)&=&H^0(K,\mu_{p^v}^{\otimes s})
\oplus H^2(K,\mu_{p^v}^{\otimes(s+1)}),\\
\TC_{2s+1}(A|K;p,\Z/p^v)&=&H^1(K,\mu_{p^v}^{\otimes(s+1)}),\\
\end{eqnarray*}
and thus, in turn, Theorem \ref{thmA}.

It is also easy to see that the canonical map
$$K_*(K,\Z/p^v)\to K_*^{\hbox{{\eightpoint \'et}}}(K,\Z/p^v)$$
is an isomorphism in degrees $\geq 1$. Here the
right-hand side is the Dwyer-Friedlander \'etale $K$-theory of $K$ with
$\Z/p^v$-coefficients. This may be defined as the homotopy
groups with $\Z/p^v$-coefficients of the spectrum
$$K^{\hbox{{\eightpoint \'et}}}(K)= \hcl{L/K} {\coborel(G_{L/K},K(L))},$$
where the homotopy colimit runs over the finite Galois extensions $L/K$
contained in an algebraic closure $\bar K/K$, and where the spectrum
$\coborel(G_{L/K},K(L))$ is the group cohomology spectrum or
homotopy fixed point spectrum of $G_{L/K}$ acting on $K(L)$. There is
a spectral sequence
$$E^2_{s,t}=H^{-s}(K,\mu_{p^v}^{\otimes(t/2)})
\Rightarrow K_{s+t}^{\hbox{{\eightpoint \'et}}}(K,\Z/p^v),$$
where the identification of the $E^2$-term is a consequence of the
celebrated theorem of Suslin,~\cite{suslin}, that
$$K_t(\bar K,\Z/p^v)=\mu_{p^v}^{\otimes(t/2)}.$$
For $K$ a finite extension of $\Qp$, the $p$-adic homotopy type of
the $K^{\hbox{{\eightpoint \'et}}}(K)$ is known by~\cite{thomason} and
\cite{dwyermitchell}. Let $F\Psi^r$ be the homotopy fiber
$$F\Psi^r\to\Z\times BU\stck{\Psi^r-1}BU.$$
It follows from this calculation and from the isomorphism above that:

\specialnumber{D}\proclaim{Theorem}\label{thmD} If $K$ is a finite extension of $\Qp${\rm ,} then
after $p$\/{\rm -}\/completion
$$\Z\times BGL(K)^+\simeq\;F\Psi^{g^{p^{a-1}d}}\times
BF\Psi^{g^{p^{a-1}d}}\times U^{|K\,\:\Qp|},$$
where $d=(p-1)/|K(\mu_p):K|${\rm ,} $a=\max\{v\;|\;\mu_{p^v}\subset
K(\mu_p)\}${\rm ,} and where $g\in\Zp^\times$ is a topological generator.
\endproclaim

The proof of theorem \ref{thmC} is given in Section~6 below. It is
based on the calculation in Section~5 of the Tate spectra for the cyclic
groups $C_{p^n}$ acting on the topological Hochschild spectrum
$T(A|K)$: Given a finite group $G$ and $G$-spectrum $X$, one has the
Tate spectrum $\tate(G,X)$ of~\cite{greenlees},~\cite{greenleesmay}.
Its homotopy groups are approximated by a spectral sequence
$$E^2_{s,t}=\hat H^{-s}(G,\pi_tX)\Rightarrow\pi_{s+t}\tate(G,X),$$
which converges conditionally in the sense of~\cite{boardman}. In Section~4
below we give a slightly different construction of this spectral
sequence which is better suited for studying multiplicative
properties. The cyclotomic structure of $T(A|K)$ gives rise to a map
$$\hat\Gamma_K\:\TR^n(A|K;p)\to\tate(C_{p^n},T(A|K)),$$
and we show in Section~5 that this map induces an isomorphism of homotopy
groups with $\Z/p^v$-coefficients in degrees $\geq 0$. We then
evaluate the Tate spectral sequence for the right-hand side.

Throughout this paper, $A$ will be a complete discrete valuation
ring with field of fractions $K$ of characteristic zero and perfect
residue field $k$ of characteristic $p>2$. All rings are assumed
commutative and unital without further notice. Occasionally, we
will write $\bar\pi_*(-)$ for homotopy groups with
$\Z/p$-coefficients.

This paper has been long underway, and we would like to acknowledge
the financial support and hospitality of the many institutions we have
visited while working on this project: Max Planck Institut f\"ur
Mathematik in Bonn, The American Institute of Mathematics at Stanford,
Princeton University, The University of Chicago, Stanford University,
the SFB 478 at Universit\"at M\"unster, and the SFB 343 at
Universit\"at Bielefeld. It is also a pleasure to thank Mike Hopkins
and Marcel B\"okstedt for valuable help and comments. We are
particularly indebted to Mike Mandell for a conversation which was
instrumental in arriving at the definition of the spectrum $T(A|K)$ as
well as for help at various other points. Finally, we thank an unnamed
referee for valuable suggestions on improving the exposition.

\section{Topological Hochschild homology and localization}\label{s}

1.1.\quad  This section contains the construction of
$\TR^n(A|K;p)$. The main result is the localization sequence of
Theorem~\ref{localizationtr}, which relates this spectrum to
$\TR^n(A;p)$ and $\TR^n(k;p)$. We make extensive use of the machinery
developed by Waldhausen  in~\cite{w} and some familiarity with this
material is assumed.

The stable homotopy category is a triangulated category and a closed
symmetric monoidal category, and the two structures are compatible;
see e.g.~\cite[Appendix]{hoveypalmieristrickland}. By a spectrum we
will mean an object in this category, and by a ring spectrum we will
mean a monoid in this category. The purpose of this section is to
produce the following. Let ${\cal C}$ be a {\it linear} category
with cofibrations and weak equivalences in the sense of \cite[\S
1.2]{w}. We define a pro-spectrum $\TR\cs({\cal C};p)$ together
with maps of pro-spectra
\begin{eqnarray*}
  F\:\TR^n({\cal C};p) &\hskip-6pt \to \hskip-6pt& \TR^{n-1}({\cal C};p), \\
  V\:\TR^{n-1}({\cal C};p) &\hskip-6pt \to \hskip-6pt& \TR^n({\cal C};p), \\
 \mu\:S^1_+\wedge\TR^n({\cal C};p)&\hskip-6pt \to \hskip-6pt&\TR^n({\cal C};p).\\
\end{eqnarray*}
The spectrum $\TR^1({\cal C};p)$ is the topological Hochschild
spectrum of ${\cal C}$. The cyclotomic trace is a map of
pro-spectra
$$\tr\:K({\cal C})\to\TR\cs({\cal C};p),$$
where the algebraic $K$-theory spectrum on the left is regarded as a
constant pro-spectrum.  

Suppose that the category ${\cal C}$ has a strict symmetric
monoidal structure such that the tensor product is bi-exact. Then
there is a natural product on $\TR\cs({\cal C};p)$ which makes it a
commutative pro-ring spectrum. Similarly, $K({\cal C})$ is
naturally a commutative ring spectrum and the maps $F$ and $\tr$ are
maps of ring-spectra.

The pro-spectrum $\TR\cs({\cal C};p)$ has a preferred homotopy
limit $\TR({\cal C};p)$, and there are preferred lifts to the
homotopy limit of the maps $F$, $V$ and $\mu$. Its homotopy groups are
related to those of the pro-system by the Milnor sequence
$$0\to{\displaystyle{\lim_{{\longleftarrow}\atop R}}}^1{\TR_{s+1}\cs({\cal C};p)}
\to\TR_s({\cal C};p)
\to{\displaystyle{\lim_{{\longleftarrow}\atop R}}}{\TR_s\cs({\cal C};p)}\to 0.$$
There is a natural cofibration sequence
$$\TC({\cal C};p)\to\TR({\cal C};p)\stck{R-F}\TR({\cal C};p)
\to\Sigma\TC({\cal C};p),$$
where $\TC({\cal C};p)$ is the topological cyclic homology
spectrum of ${\cal C}$. The cyclotomic trace has a preferred lift
to a map
$$\tr\:K({\cal C})\to\TC({\cal C};p),$$
and in the case where ${\cal C}$ has a bi-exact strict symmetric
monoidal product, the natural product on $\TR\cs({\cal C};p)$ have
preferred lifts to natural products on $\TR({\cal C};p)$ and
$\TC({\cal C};p)$, and the maps $F$ and $\tr$ are ring maps.

Let $G$ be a compact Lie group. One then has the $G$-stable category
which is a triangulated category with a compatible closed symmetric
monoidal structure. The objects of this category are called $G$-spectra,
and the monoids for the smash product are called ring $G$-spectra.
Let $H\subset G$ be a closed subgroup and let $W_HG=N_GH/H$ be the
Weyl group. There is a forgetful functor which to a $G$-spectrum $X$
assigns the underlying $H$-spectrum $U_HX$. We also write $|X|$ for
$U_{\{1\}}X$. It comes with a natural map of spectra
$$\mu_X\:G_+\wedge |X|\to |X|.$$
One also has the $H$-fixed point functor which to a $G$-spectrum $X$
assigns the $W_HG$-spectrum $X^H$. If $H\subset K\subset G$ are two
closed subgroups, there is a map of spectra
$$\iota_H^K\:|X^K|\to |X^H|,$$
and if $|K\!:\!H|$ is finite, a map in the opposite direction
$$\tau_H^K\:|X^H|\to |X^K|.$$
If $X$ is a ring $G$-spectrum then $U_HX$ is a  ring $H$-spectrum and
$X^H$ is a ring $W_GH$-spectrum.

Let ${\Bbb T}$ be the circle group, and let $C_r\subset{\Bbb T}$
be the cyclic subgroup of order~$r$. We then have the canonical
isomorphism of groups
$$\rho_r\:{\Bbb T}\xto{\sim}{\Bbb T}/C_r=W_{{\Bbb T}}C_r$$
given by the $r$-th root. It induces an isomorphism of the
${\Bbb T}/C_r$-stable category and of the ${\Bbb T}$-stable category
by assigning to a ${\Bbb T}/C_r$-spectrum $Y$ the
${\Bbb T}$-spectrum $\rho_r^*Y$. Moreover, there is a transitive
system of natural isomorphisms of spectra
$$\varphi_r\:|\rho_r^*Y|\xto{\sim}|Y|,$$
and the following diagram commutes
$$\begin{array}{ccc}
{{\Bbb T}_+\wedge |\rho_r^*Y|}&\stackrel{\mu}{\longrightarrow}&{|\rho_r^*Y|}\\[4pt]
\phantom{\scriptstyle\rho\wedge\varphi_r} \big\downarrow{\scriptstyle\rho\wedge\varphi_r} & &
\quad \big\downarrow{\scriptstyle\varphi_r} \\[4pt]
{{\Bbb T}/C_{r+}\wedge |Y|}&\stackrel{\mu}{\longrightarrow} &
{|Y|.} \end{array}
$$
 
We will define a ${\Bbb T}$-spectrum $T({\cal C})$ such that
$$\TR^n({\cal C};p)=|\rho_{p^{n-1}}^*T({\cal C})^{C_{p^{n-1}}}|$$
with the maps $F$ and $V$ given by the composites
\begin{eqnarray*}
F=\varphi_{p^{n-2}}^{-1}\iota_{C_{p^{n-2}}}^{C_{p^{n-1}}}\varphi_{p^{n-1}}
\:&&\hskip-20pt |\rho_{p^{n-1}}^*T({\cal C})^{C_{p^{n-1}}}|
\to|\rho_{p^{n-2}}^*T({\cal C})^{C_{p^{n-2}}}|, \\
V=\varphi_{p^{n-1}}^{-1}\tau_{C_{p^{n-2}}}^{C_{p^{n-1}}}\varphi_{p^{n-2}}
\:&&\hskip-20pt |\rho_{p^{n-2}}^*T({\cal C})^{C_{p^{n-2}}}|
\to|\rho_{p^{n-1}}^*T({\cal C})^{C_{p^{n-1}}}|, \\
\end{eqnarray*}
and the map $\mu$ given by
$$\mu=\mu_{\rho_{p^{n-1}}^*T({\cal C})^{C_{p^{n-1}}}}\:
{\Bbb T}_+\wedge|\rho_{p^{n-1}}^*T({\cal C})^{C_{p^{n-1}}}|
\to|\rho_{p^{n-1}}^*T({\cal C})^{C_{p^{n-1}}}|.$$
There is a natural map
$$K({\cal C})\to T({\cal C})^{{\Bbb T}},$$
and the cyclotomic trace is then the composite of this map and
$\varphi_{p^{n-1}}^{-1}\iota_{C_{p^{n-1}}}^{{\Bbb T}}$. The
definition of the structure maps in the pro-system
$\TR\cs({\cal C};p)$ is more complicated and uses the 
{\it cyclotomic} structure on $T({\cal C})$ which we now explain.

There is a cofibration sequence of ${\Bbb T}$-CW-complexes
$$E_+\to S^0\to\tilde E\to\Sigma E_+,$$
where $E$ is a free contractible ${\Bbb T}$-space, and where the left-hand map collapses $E$ to the nonbase point of $S^0$. It
induces, upon smashing with a ${\Bbb T}$-spectrum~$T$, a cofibration sequence of
${\Bbb T}$-spectra
$$E_+\wedge T\to T\to\tilde E\wedge T\to\Sigma E_+\wedge T,$$
and hence the following basic cofibration sequence of spectra
$$|\rho_{p^n}^*(E_+\wedge T)^{C_{p^n}}|
\to |\rho_{p^n}^*T^{C_{p^n}}|
\to |\rho_{p^n}^*(\tilde E\wedge T)^{C_{p^n}}|
\to\Sigma |\rho_{p^n}^*(E_+\wedge T)^{C_{p^n}}|,$$
natural in $T$. The left-hand term is written
$\borel(C_{p^n},T)$ and called the group homology spectrum or
Borel spectrum. Its homotopy groups are approximated by a strongly
convergent first quadrant homology type spectral sequence
$$E^2_{s,t}=H_s(C_{p^n},\pi_tT)
\Rightarrow\pi_{s+t}\borel(C_{p^n},T).$$
The cyclotomic structure on $T({\cal C})$ means that there is a
natural map of\break ${\Bbb T}$-spectra
$$r\:\rho_p^*(\tilde E\wedge T({\cal C}))^{C_p}
\to T({\cal C})$$
such that $U_{C_{p^s}}r$ is an isomorphism of $C_{p^s}$-spectra, for
all $s\geq 0$. More generally, since
$$\rho_{p^n}^*(\tilde E\wedge T({\cal C}))^{C_{p^n}}
= \rho_{p^{n-1}}^*(\rho_p^*(\tilde E\wedge 
T({\cal C}))^{C_p})^{C_{p^{n-1}}},$$
the map $r$ induces a map of ${\Bbb T}$-spectra
$$r_{n+1}\:\rho_{p^n}^*(\tilde E\wedge T({\cal C}))^{C_{p^n}}
\to\rho_{p^{n-1}}^*T({\cal C})^{C_{p^{n-1}}}$$
such that $U_{C_{p^s}}r_{n+1}$ is an isomorphism of $C_{p^s}$-spectra,
for all $s\geq 0$. The map
$$R\:\TR^n({\cal C};p)\to\TR^{n-1}({\cal C};p)$$
is then defined as the composite
$$|\rho_{p^{n-1}}^*T({\cal C})^{C_{p^{n-1}}}|
\to |\rho_{p^{n-1}}^*(\tilde E\wedge T({\cal C}))^{C_{p^{n-1}}}|
\begin{array}{c} {\scriptstyle r_{\scriptscriptstyle n}}\\[-8pt] \longrightarrow\\[-10pt]
{\scriptstyle\sim}
\end{array} |\rho_{p^{n-2}}^*T({\cal C})^{C_{p^{n-2}}}|,$$ where the left-hand map is the middle map in
the cofibration sequence above. We thus have a natural cofibration sequence of spectra
$$\borel(C_{p^{n-1}},T({\cal C}))\xto{N}
\TR^n({\cal C};p)\xto{R}
\TR^{n-1}({\cal C};p)\xto{\partial}
\Sigma\borel(C_{p^{n-1}},T({\cal C})).$$
When ${\cal C}$ has a bi-exact strict symmetric monoidal product,
the map $r$ is a map of ring ${\Bbb T}$-spectra, and hence $R$ is
a map of ring spectra. The cofibration sequence above is a sequence of
$\TR^n({\cal C};p)$-module spectra and maps.

For any ${\Bbb T}$-spectrum $X$, one has the function spectrum
$F(E_+,X)$, and the projection $E_+\to S^0$ defines a natural map
$$\gamma\:X\to F(E_+,X).$$
This map induces an isomorphism of group homology spectra. One defines
the group cohomology spectrum and the Tate spectrum,
\begin{eqnarray*}
\coborel(C_{p^n},X) &
=&|\rho_{p^n}^*F(E_+,X)^{C_{p^n}}|, \\
\tate\;(C_{p^n},X) &
=&|\rho_{p^n}^*(\tilde E\wedge F(E_+,X))^{C_{p^n}}|.  
\end{eqnarray*}
Their homotopy groups are approximated by homology type spectral
sequences
\begin{eqnarray*}
E^2_{s,t}&\hskip-8pt =\hskip-8pt&H^{-s}(C_{p^n},\pi_tX)
\Rightarrow\pi_{s+t}\coborel(C_{p^n},X),\\
\hat E^2_{s,t}&\hskip-8pt=\hskip-8pt&\hat H^{-s}(C_{p^n},\pi_tX)
\Rightarrow\pi_{s+t}\tate\;(C_{p^n},X), 
\end{eqnarray*}
both of which converge  conditionally in the sense of
\cite[Def.\ 5.10]{boardman}. The latter sequence, called the Tate
spectral sequence, will be considered in great detail in
Section~\ref{tatecohomology} below. Taking $T=F(E_+,X)$ in the basic
cofibration sequence above, we get the Tate cofibration sequence of
spectra 
$$\borel(C_{p^n},X)\xto{N^h}
\coborel(C_{p^n},X)\xto{R^h}
\tate(C_{p^n},X)\xto{\partial^h}
\Sigma\borel(C_{p^n},X).$$

Finally, if $X=T({\cal C})$, the map
$$\gamma\:T({\cal C})\to F(E_+,T({\cal C}))$$
induces a map of cofibration sequences
$${\ninepoint \begin{array}{ccccccc}
{\borel(C_{p^n},T({\cal C}))}&\sto{N}  &
{\TR^{n+1}({\cal C};p)} &\sto{R}&{\TR^n({\cal C};p)}&\sto{\partial}&{\Sigma\borel(C_{p^n},T({\cal
C}))}\\[4pt]
\big\Vert&&\big\downarrow{\scriptstyle\Gamma} &&
 \big\downarrow{\scriptstyle\hat\Gamma} &&\big\Vert\\[4pt]
{\borel(C_{p^n},T({\cal C}))}&\sto{N^h} &
{\coborel(C_{p^n},T({\cal C}))}&\sto{R^h} &
{\tate(C_{p^n},T({\cal C}))} &\sto{\partial^h} &
{\Sigma\borel(C_{p^n},T({\cal C})),} \end{array}}$$
in which all maps commute with the action maps $\mu$. Moreover, if
${\cal C}$ is strict symmetric monoidal with bi-exact tensor
product, the four spectra in the middle square are all ring spectra
and $R$, $R^h$, $\Gamma$ and $\hat\Gamma$ are maps of ring
spectra. In this case, the diagram is a diagram of
$\TR^{n+1}({\cal C};p)$-module spectra,~\cite[pp.\ 71--72]{hm}. 

\vglue12pt 1.2.\quad In order to construct the ${\Bbb T}$-spectrum
$T({\cal C})$ we need a model category for the ${\Bbb T}$-stable
category. The model category we use is the category of symmetric
spectra of orthogonal ${\Bbb T}$-spectra, see \cite{mandellmay} and
\cite[Th.\ 5.10]{hovey}. We first recall the topological Hochschild
space $\THH({\cal C})$. See \cite{dm}, \cite{gh} and \cite{hm} for
more details.

A linear category ${\cal C}$ is naturally enriched over the
symmetric monoidal category of symmetric spectra. The symmetric
spectrum of maps from $c$ to $d$, $\ul{\Hom}_{\cal C}(c,d)$, is the
Eilenberg-MacLane spectrum for the abelian group $\Hom_{\cal C}(c,d)$
concentrated in degree zero. In more detail, if $X$ is a pointed
simplicial set, then
$$\Z(X)=\Z\{X\}/\Z\{x_0\}$$
is a simplicial abelian group whose homology is the reduced singular
homology of $X$. Here $\Z\{X\}$ denotes the degree-wise free abelian
group generated by~$X$. Let $S^i$ be the $i$-fold smash product of the
standard simplicial circle
$S^1=\Delta[1]/\partial\Delta[1]$. Then the spaces
$\{|\Z(S^i)|\}_{i\geq 0}$ is a symmetric ring spectrum with the
homotopy type of an Eilenberg-MacLane spectrum for $\Z$ concentrated
in degree zero, and we define
$$\ul{\Hom}_{\cal C}(c,d)_i=|\Hom_{\cal C}(c,d)\otimes\Z(S^i)|.$$

Let $I$ be the category with objects the finite sets
$$\ul{i}=\{1,2,\dots,i\},\hskip6mm i\geq 1,$$
and the empty set $\ul 0$, and morphisms all injective maps. It is a
strict monoidal category under concatenation of sets and maps. There
is a functor $V_k({\cal C};X)$ from $I^{k+1}$ to the category of
pointed spaces which on objects is given by
$$V_k({\cal C};X)(\ul{i_0},\dots,\ul{i_k})=
\bigvee_{c_0,\dots,c_k\in\ob{\cal C}}
\ul{\Hom}_{\cal C}(c_0,c_k)_{i_0}\wedge\dots\wedge
\ul{\Hom}_{\cal C}(c_k,c_{k-1})_{i_k}\wedge X.$$
It induces a functor $G_k({\cal C};X)$ from $I^{k+1}$ to pointed
spaces with
$$G_k(C;X)(\ul{i_0},\dots,\ul{i_k})=F(S^{i_0}\wedge\dots\wedge
S^{i_k},V_k({\cal C};X)(\ul{i_0},\dots,\ul{i_k})),$$
and we define
$$\THH_k({\cal C})=\hcl{I^{k+1}}{G_k({\cal C};S^0)}.$$
This is naturally the space of $k$-simplices in a cyclic
space and, by definition,
$$\THH({\cal C})=|[k]\mapsto\THH_k({\cal C})|.$$
It is a ${\Bbb T}$-space by Connes' theory of cyclic
spaces,~\cite[7.1.9]{loday}.

More generally, let $(n)$ be the finite ordered set
$\{1,2,\dots,n\}$ and let $(0)$ be the empty set. The product category
$I^{(n)}$ is a strict monoidal category under component-wise
concatenation of sets and maps. Concatenation of sets and maps
according to the ordering of $(n)$ also defines a functor
$$\sqcup_n\:I^{(n)}\to I,$$
but this does {\it not} preserve the monoidal structure. By
convention $I^{(0)}$ is the category with one object and one morphism,
and $\sqcup_0$ includes this category as the full subcategory on the
object $\ul{0}$. We let $G^{(n)}_k({\cal C};X)$ be the functor from
$(I^{(n)})^{k+1}$ to the category of pointed spaces given by
$$G^{(n)}_k({\cal C};X)=G_k({\cal C};X)\circ(\sqcup_n)^{k+1},$$
and define
$$\THH^{(n)}_k({\cal C};X) =
\hcl{(I^{(n)})^{k+1}}G^{(n)}_k({\cal C};X).$$
In particular, $\THH^{(0)}_k({\cal C};X) =
N_k^\cy({\cal C})\wedge X$, where
$$N_k^\cy({\cal C}) =
\bigvee_{c_0,\dots,c_k\in\ob{\cal C}}
\Hom_{{\cal C}}(c_0,c_k)\wedge\dots\wedge
\Hom_{{\cal C}}(c_k,c_{k-1})$$
is the cyclic bar construction of ${\cal C}$. Again this is the
space of $k$-simplices in a cyclic space, and hence we have the
$\Sigma_n\times{\Bbb T}$-space
$$\THH^{(n)}({\cal C};X)=|[k]\mapsto\THH^{(n)}_k({\cal C};X)|.$$
There is a natural product
$$\THH^{(m)}({\cal C};X)\wedge\THH^{(n)}({\cal D};Y)
\to\THH^{(m+n)}({\cal C}\otimes{\cal D};X\wedge Y),$$
which is $\Sigma_m\times\Sigma_n\times{\Bbb T}$-equivariant if
${\Bbb T}$ acts diagonally on the left. Here the category
${\cal C}\otimes{\cal D}$ has as objects all pairs $(c,d)$ with
$c\in\ob{\cal C}$ and $d\in{\cal D}$, and
$$\Hom_{{\cal C}\otimes{\cal D}}((c,d),(c',d'))
=\Hom_{{\cal C}}(c,c')\otimes\Hom_{{\cal D}}(d,d').$$

For any category ${\cal C}$, the nerve category
${\bf N}\s{\cal C}$ is the simplicial category with
$k$-simplicies the functor category
$${\bf N}_k{\cal C}={\cal C}^{[k]},$$
where the partially ordered set $[k]=\{0,1,\ldots,k\}$ is viewed as
a category. An order-preserving map $\theta\:[k]\to [l]$ may be
viewed as a functor and hence induces a functor
$$\theta^*\:{\bf N}_l{\cal C}\to{\bf N}_k{\cal C}.$$
The objects of ${\bf N}\s{\cal C}$ comprise the nerve of
${\cal C}$, $N\s{\cal C}$. Clearly, the nerve category is a
functor from categories to simplicial categories.

Suppose now that ${\cal C}$ is a category with cofibrations and
weak equivalences in the sense of \cite[\S 1.2]{w}. We then
define
$${\bf N}\s^w{\cal C}\subset{\bf N}\s{\cal C}$$
to be the {\it full} simplicial subcategory with
$$\ob{\bf N}\s^w{\cal C}=N\s w{\cal C}.$$
There is a natural structure of simplicial categories with cofibrations
and weak equivalences on ${\bf N}\s^w{\cal C}$:
$\operatorname{co}{\bf N}\s^w{\cal C}$ and
$w{\bf N}\s^w{\cal C}$ are the simplicial subcategories which
contain all objects but where morphisms are natural transformations
through cofibrations and weak equivalences in ${\cal C}$,
respectively. With these definitions there is a natural isomorphism of
bi-simplicial categories \pagebreak with cofibrations and weak equivalences
\begin{equation}\label{twist}
{\bf N}\s S\s{\cal C}\cong S\s{\bf N}\s{\cal C}, \speqnu{1.2.1}
\end{equation}
where $S\s{\cal C}$ is Waldhausen's construction,
\cite[\S 1.3]{w}.

Let $V$ be a finite-dimensional orthogonal ${\Bbb T}$-representation.
We define the $(n,V)$-th space in the symmetric orthogonal
${\Bbb T}$-spectrum $T({\cal C})$ by
\begin{equation}\label{tophoch}
T({\cal C})_{n,V}=|\THH^{(n)}({\bf N}\s^wS\s^{(n)}{\cal C};S^V)|. \speqnu{1.2.2}
\end{equation}
There are two ${\Bbb T}$-actions on this space: one which comes
from the topological Hochschild space, and another induced from the
${\Bbb T}$-action on $S^V$. We give $T({\cal C})_{n,V}$ the
diagonal ${\Bbb T}$-action. There are also two $\Sigma_n$-actions:
one which comes from the $\Sigma_n$-action on the topological
Hochschild space, and another induced from the permutation of the
simplicial directions in the $n$-simplicial category
$S^{(n)}\s{\cal C}$; compare \cite[6.1]{gh}. We also give
$T({\cal C})_{n,V}$ the diagonal $\Sigma_n$-action. In particular,
the $(0,0)$-th space is the cyclic bar construction
$$T({\cal C})_{0,0}=|N\s^\cy({\bf N}^w\s{\cal C})|.$$
In general, the ${\Bbb T}$-fixed set of the realization of a cyclic
space $X\s$ is given by
$$|X\s|^{{\Bbb T}}=\{x\in X_0\,|\,s_0(x)=t_1s_0(x)\},$$
and hence, we have a canonical map
$$|\ob{\bf N}\s^wS\s^{(n)}{\cal C}\wedge S^{V^{{\Bbb T}}}|
\to (T({\cal C})_{n,V})^{{\Bbb T}}.$$
The space on the left is the $(n,V^{{\Bbb T}})$-th space of a
symmetric orthogonal spectrum, which represents the spectrum
$K({\cal C})$ in the stable homotopy category, and the map above
defines the cyclotomic trace. Moreover, by a construction similar to
that of \cite[\S2]{hm}, there are ${\Bbb T}$-equivariant maps
$$\rho_p^*(T({\cal C})_{n,V})^{C_p}
\to T({\cal C})_{n,\rho_p^*V^{C_p}},$$
and one can prove that for fixed $n$, the object of the
${\Bbb T}$-stable category defined by the orthogonal spectrum
$V\mapsto T({\cal C})_{n,V}$ has a cyclotomic structure.

Suppose that ${\cal C}$ is a strict symmetric monoidal category and
that the tensor product is bi-exact. There is then an induced
$\Sigma_m\times\Sigma_n$-equivariant product
$$S^{(m)}\s{\cal C}\otimes S^{(n)}\s{\cal C}\to
S^{(m+n)}\s{\cal C},$$
and hence
$$T({\cal C})_{m,V}\wedge T({\cal C})_{n,W}
\to T({\cal C})_{m+n,V\oplus W}.$$
This product makes $T({\cal C})$ a monoid in the symmetric monoidal
category of symmetric orthogonal ${\Bbb T}$-spectra.

\vglue12pt 1.3.\quad We need to recall some of the properties of this
construction. It is convenient to work in a more general setting.

Let $\Phi$ be a functor from a category of categories with
cofibrations and weak equivalences to the category of pointed
spaces. If ${\cal C}\s$ is a simplicial category with cofibrations
and weak equivalences, we define
$$\Phi({\cal C}\s)=|[n]\mapsto\Phi({\cal C}_n)|.$$
We shall assume that $\Phi$ satisfies the following axioms:
\begin{itemize}
\item[(i)] The trivial category with cofibrations and weak equivalences is
mapped to a one-point space.

\item[(ii)] For any pair ${\cal C}$ and ${\cal D}$ of
categories with cofibrations and weak equivalences, the canonical map
$$\Phi({\cal C}\times{\cal D})\xto{\sim}
\Phi({\cal C})\times\Phi({\cal D})$$
is a weak equivalence.

\item[(iii)] If $f\s\:{\cal C}\s\to{\cal D}\s$ is a map of simplicial
categories with cofibrations and weak equivalences, and if for all
$n$,
$\Phi(f_n)\:\Phi({\cal C}_n)\to\Phi({\cal D}_n)$ is a
weak equivalence, then
$$\Phi(f\s)\:\Phi({\cal C}\s)\to\Phi({\cal D}\s)$$
is a weak equivalence.
\end{itemize}

In \cite{w}, $\Phi$ is the functor which to a category assigns the set
of objects. Here our main concern is the functor $\THH$ and variations
thereof.

We next recall some generalities. Let 
$$f,g\:{\cal C}\s\to{\cal D}\s$$
be two exact simplicial functors. An {\it exact} simplicial
homotopy from $f$ to $g$ is an exact simplicial functor
$$h\:\Delta[1]\s\times{\cal C}\s\to{\cal D}\s$$
such that $h\circ(d^1\times\id)=f$ and $h\circ(d^0\times\id)=g$. Here
$\Delta[n]\s$ is viewed\break as a discrete simplicial category with its
unique structure of a simplicial category with cofibrations and weak
equivalences. An exact simplicial functor\break
$f\:{\cal C}\s\to{\cal D}\s$ is an exact simplicial homotopy
equivalence if there exists an exact simplicial functor
$g\:{\cal D}\s\to{\cal C}\s$ and exact simplicial homotopies of
the two composites to the respective identity simplicial functors. 

\specialnumber{1.3.1} \proclaim{Lemma}\label{sim} An exact simplicial homotopy
$\Delta[1]\s\times{\cal C}\s\to{\cal D}\s$ induces a
homotopy
$$\Delta[1]\times\Phi({\cal C}\s)\to\Phi({\cal D}\s).$$
Hence $\Phi$ takes exact simplicial homotopy equivalences to homotopy
equivalences.
\endproclaim

\demo{Proof} There is a natural transformation
$$\Delta[1]_k\times\Phi({\cal C}_k)
\to\Phi(\Delta[1]_k\times{\cal C}_k).$$
Indeed, $\Delta[1]_k\times\Phi({\cal C}_k)$ and
$\Delta[1]_k\times{\cal C}_k$ are coproducts in the category of
spaces and the category of categories with cofibrations and weak
equivalences, respectively, indexed by the set $\Delta[1]_k$. The map
exists by the universal property of coproducts.
\enddemo

\specialnumber{1.3.2} \proclaim{Lemma}\label{nerve} An exact functor of categories with
cofibrations and weak equivalences $f\:{\cal C}\to{\cal D}$
induces an exact simplicial functor ${\bf N}\s^w
f\:{\bf N}\s^w{\cal C}\to{\bf N}\s^w{\cal D}$. A natural
transformation through weak equivalences of ${\cal D}$ between two
such functors $f$ and $g$ induces an exact simplicial homotopy between
${\bf N}\s^w f$ and ${\bf N}\s^w g$.
\endproclaim

\demo{Proof} The first statement is clear. We view the partially
ordered set $[1]$ as a category with cofibrations and weak
equivalences where the nonidentity map is a weak equivalence but not
a cofibration. Then the natural transformation defines an exact
functor $[1]\times{\cal C}\to{\cal D}$, and the required exact
simplicial homotopy is given by the composite
$$\Delta[1]\s\times{\bf N}\s^w{\cal C}\to{\bf N}\s^w[1]\times
{\bf N}\s^w{\cal C}\to{\bf N}\s^w([1]\times{\cal C})
\to{\bf N}\s^w{\cal D},
$$
where the first and the middle arrow are the canonical simplicial
functors, and the last is induced from the natural transformation.
(Note that ${\bf N}\s^w[n]$ is not a discrete category.)
\enddemo

\proclaimtitle{[48, Lemma 1.4.1]}
\specialnumber{1.3.3} \proclaim{Lemma}\label{Speriod}  Let
$f,g\:{\cal C}\to{\cal D}$ be a pair of exact functors of
categories with cofibrations. A natural isomorphism from $f$ to $g$
induces an exact simplicial homotopy
$$\Delta[1]\s\times S\s{\cal C}\to S\s{\cal D}$$
from $S\s f$ to $S\s g$.\ 
\endproclaim

\specialnumber{1.3.4}\proclaim{{C}orollary}\label{i} Let ${\cal C}$ be a category with
cofibrations{\rm ,} and let $i{\cal C}$ be the subcategory of
isomorphisms. Then the map induced from the degeneracies in the nerve
direction induces a weak equivalence
$$\Phi(S\s{\cal C})\xto{\sim}\Phi({\bf N}\s^iS\s{\cal C}).$$
\endproclaim

\demo{Proof} For each $k$, the iterated degeneracy functor
$$s\:{\cal C}={\bf N}^i_0{\cal C}\to{\bf N}^i_k{\cal C},$$
has the retraction
$$\theta^*\:{\bf N}^i_k{\cal C}\to{\cal C},$$
where $\theta\:[0]\to[k]$ is given by $\theta(0)=0$. Moreover, there
is a natural isomorphism $\id\xto{\sim}\theta^*$, and hence \pagebreak by
Lemma~\ref{Speriod},
$$S\s s\:S\s{\cal C}\to S\s{\bf N}^i_k{\cal C}
={\bf N}^i_kS\s{\cal C}$$
is an exact simplicial homotopy equivalence. The corollary follows
from\break Lemma~\ref{sim} and from property (iii) above.
\enddemo

Let ${\cal A}$, ${\cal B}$ and ${\cal C}$ be categories with
cofibrations and weak equivalences and suppose that ${\cal A}$ and
${\cal B}$ are subcategories of ${\cal C}$ and that the
inclusion functors are exact. Following \cite[p.~335]{w}, let
$E({\cal A},{\cal C},{\cal B})$ be the category with
cofibrations and weak equivalences given by the pull-back diagram
$$\begin{array}{ccc}
{E({\cal A},{\cal C},{\cal B})}&\stackrel{(s,t,q)}{\srar} &
{{\cal A}\times{\cal C}\times{\cal B}}  \\[4pt]
\big\downarrow&&\big\downarrow\\[4pt]
{S_2{\cal C}}&\stackrel{(d_2,d_1,d_0)}{\lrar} &
{{\cal C}\times{\cal C}\times{\cal C}}.\end{array}
$$
In other words, $E({\cal A},{\cal C},{\cal B})$ is the
category of cofibration sequences in ${\cal C}$ of the form
$$A\rightarrowtail C\twoheadrightarrow B,\hskip5mm A\in{\cal A},\;
B\in{\cal B}.$$
The exact functors $s$, $t$ and $q$ take this sequence to $A$, $C$ and
$B$, respectively. The extension of the additivity theorem to the
present situation is due to McCarthy, \cite{mc}. Indeed, the proof
given   there for $\Phi$ the cyclic nerve functor
generalizes mutatis mutandis to prove the statement (1) below. The
equivalence of the four statements follows from \cite[Prop.\
1.3.2]{w}. 

\proclaimtitle{Additivity theorem}
\specialnumber{1.3.5}\proclaim{Theorem} \label{additivity} The
following equivalent assertions hold\/{\rm :}\/
\begin{itemize}
\item[{\rm (1)}] The exact functors $s$ and $q$ induce a weak equivalence
$$\Phi({\bf N}\s^wS\s E({\cal A},{\cal C},{\cal B}))
\xto{\sim}\Phi({\bf N}\s^wS\s{\cal A})\times
\Phi({\bf N}\s^wS\s{\cal B}).$$

\item[{\rm (2)}] The exact functors $s$ and $q$ induce a weak equivalence
$$\Phi({\bf N}\s^wS\s E({\cal C},{\cal C},{\cal C}))
\xto{\sim}\Phi({\bf N}\s^wS\s{\cal C})\times
\Phi({\bf N}\s^wS\s{\cal C}).$$

\item[{\rm (3)}] The functors $t$ and $s\vee q$ induce homotopic maps
$$\Phi({\bf N}\s^wS\s E({\cal C},{\cal C},{\cal C}))
\to\Phi({\bf N}\s^wS\s{\cal C}).$$

\item[{\rm (4)}] Let $F'\rightarrowtail F\twoheadrightarrow F''$ be
a cofibration sequence of exact functors
${\cal C}\to{\cal D}$. Then the exact functors $F$ and $F'\vee
F''$ induce homotopic maps
$$\Phi({\bf N}\s^wS\s{\cal C})\to\Phi({\bf N}\s^wS\s{\cal D}).
$$
\end{itemize}

\endproclaim

Let $f\:{\cal C}\to{\cal D}$ be an exact functor and let
$S\s(f\:{\cal C}\to{\cal D})$ be Waldhausen's relative
construction, \cite[Def.\ 1.5.4]{w}. Then the commutative square
\begin{equation}\label{loc}
\begin{array}{ccc}
{\Phi({\bf N}\s^wS\s{\cal C})}&\srar&
{\Phi({\bf N}\s^wS\s S\s(\id\:{\cal C}\to{\cal C}))}\\[4pt]
\big\downarrow&&\big\downarrow\\[4pt]
{\Phi({\bf N}\s^wS\s{\cal D})}&\srar&
{\Phi({\bf N}\s^wS\s S\s(f\:{\cal C}\to{\cal D}))}\end{array} \speqnu{1.3.6}
\end{equation}
is homotopy cartesian, and there is a canonical contraction of the
upper right-hand term. In particular, if we let ${\cal D}$ be the
category with one object and one morphism, this shows that the
canonical map
$$\Phi({\bf N}\s^wS\s{\cal C})\xto{\sim}
\Omega\Phi({\bf N}\s^wS\s S\s{\cal C})$$
is a weak equivalence.  

\specialnumber{1.3.7}\numbereddemo{Definition} \label{F-equivalence} A map $f\:X\to Y$ of
${\Bbb T}$-spaces is called an ${\cal F}$-{\it equivalence} if for all
$r\geq 1$ the induced map of $C_r$-fixed points is a weak equivalence
of spaces.
\enddemo

\specialnumber{1.3.8}\proclaim{Proposition} \label{fibrant} Let ${\cal C}$ be a linear category with
cofibrations and weak equivalences{\rm ,} and let $T({\cal C})$ be the
topological Hochschild spectrum. Then for all orthogonal
${\Bbb T}$\/{\rm -}\/representations $W$ and $V${\rm ,} the spectrum structure maps
$$T({\cal C})_{n,V}\xto{\sim}
F(S^m\wedge S^W,T({\cal C})_{m+n,W\oplus V})$$
are ${\cal F}$\/{\rm -}\/equivalences{\rm ,} provided that $n\geq 1$.
\endproclaim

\demo{Proof} We factor the map in the statement as
$$T({\cal C})_{n,V}\to F(S^m,T({\cal C})_{m+n,V})
\to F(S^m,F(S^W,T({\cal C})_{m+n,W\oplus V})).$$
Since $S^m$ is $C_r$-fixed the map of $C_r$-fixed sets induced from
the first map may be identified with the map
$$(T({\cal C})_{n,V})^{C_r}\to
\Omega^m(T({\cal C})_{m+n,V})^{C_r},$$
and by definition, this is the map
$$\THH^{(n)}({\bf N}\s^wS\s^{(n)}{\cal C};S^V)^{C_r}
\to\Omega^m\THH^{(m+n)}({\bf N}\s^wS\s^{(m+n)}{\cal C};S^V)^{C_r}.$$
By the approximation lemma, \cite[Th.\ 1.6]{bokstedt} or
\cite[Lemma 2.3.7]{m}, we can replace the functor $\THH^{(k)}(-;-)$ by
the common functor $\THH(-;-)$, and the claim now follows
from~\eqref{loc} applied to the functor
$$\Phi({\cal C})=\THH({\cal C};S^V)^{C_r}.$$
Finally, it follows from the proof of \cite[Prop.\ 2.4]{hm} that
$$(T({\cal C})_{m+n,V})^{C_r}
\to F(S^W,T({\cal C})_{m+n,W\oplus V}))^{C_r}$$
is a weak equivalence.
\enddemo

We next extend Waldhausen's fibration theorem to the present
situation. We follow the original proof in \cite[\S 1.6]{w},
where also the notion of a cylinder functor is defined.

\specialnumber{1.3.9} \proclaim{Lemma}\label{barww} Suppose that ${\cal C}$ has a cylinder
functor{\rm ,} and that $w{\cal C}$ satisfies the cylinder axiom and the
saturation axiom. Then
$$\Phi({\bf N}\s^{\bar{w}}{\cal C})
\xto{\sim}\Phi({\bf N}\s^w{\cal C})$$
is a weak equivalence. Here $\bar
w{\cal C}=w{\cal C}\cap\operatorname{co}{\cal C}$.
\endproclaim

\demo{Proof} The proof is analogous to the proof of \cite[Lemma
1.6.3]{w}, but we need the proof of \cite[Th.\ A]{quillen} and not
just the statement. We consider the bi-simplicial category
${\bf T}({\cal C})$ whose category of $(p,q)$-simplices has, as 
objects, pairs of diagrams in ${\cal C}$ of the form
$$(A_q\to\cdots\to A_0,A_0\to B_0\to\cdots\to B_p),$$
and morphisms, all natural transformations of such pairs of
diagrams. We let
$${\bf T}^{\bar{w},w}({\cal C})\subset{\bf T}({\cal C})$$
be the {\it full} subcategory with objects the pairs of diagrams
with the left-hand diagram in $\bar{w}{\cal C}$ and the right-hand
diagram in $w{\cal C}$. There are bi-simplicial functors
$${\bf N}^{\bar{w}}({\cal C}^\op)R\stackrel{p_1}{\longleftarrow}
{\bf T}^{\bar{w},w}({\cal C})\xto{p_2}{\bf N}^w({\cal C})L,$$
where for a simplicial object $X$, the bi-simplicial objects $XL$
and $XR$ are obtained by precomposing $X$ with projections $\pr_1$
and $\pr_2$ from $\mathbold{\Delta}\times\mathbold{\Delta}$ to
$\mathbold{\Delta}$, respectively. Applying $\Phi$ in each
bi-simplicial degree, we get corresponding maps of bi-simplicial
spaces. We show that both maps induce weak equivalences after realization. 

For fixed $q$, the simplicial functor
$$p_1\:{\bf T}_{\mathbold{\cdot},q}^{\bar{w},w}({\cal C})
\to{\bf N}_q^{\bar{w}}({\cal C}^\op)$$
is a simplicial homotopy equivalence, and hence induces a homotopy
equivalence upon realization. It follows that
$$\Phi(p_1)\:\Phi({\bf T}^{\bar{w},w}({\cal C}))\xto{\sim}
\Phi({\bf N}\s^{\bar{w}}({\cal C}^\op))$$
is a weak equivalence of spaces.

Similarly, we claim that for fixed $p$, the simplicial functor
$$p_2\:{\bf T}_{p,\mathbold{\cdot}}^{\bar{w},w}({\cal C})
\to{\bf N}_p^w({\cal C})$$
is a simplicial homotopy equivalence. The homotopy inverse $\sigma$
maps
$$(B_0\to\cdots\to B_p)\mapsto
(B_0\xto{\id}\dots\xto{\id}B_0,B_0\xto{\id}B_0\to\cdots\to B_p).$$
Following the proof of \cite[Lemma 1.6.3]{w} we consider the
simplicial functor
$$t\:{\bf T}_{p,\mathbold{\cdot}}^{\bar{w},w}({\cal C})
\to{\bf T}_{p,\mathbold{\cdot}}^{\bar{w},w}({\cal C})$$
which maps
\begin{eqnarray*}
{}&(A_q\to\cdots\to A_0,A_0\to B_0\to\dots B_p)\\
{}&\mapsto(T(A_q\to B_0)\to\cdots\to T(A_0\to B_0),T(A_0\to
B_0)\xto{p}B_0\to\cdots\to B_p),\\
\end{eqnarray*}
where $T$ is the cylinder functor. There are exact simplicial
homotopies from $\sigma\circ p_2$ to $t$ and from the identity functor
to $t$. Hence
$$\Phi(p_2)\:\Phi({\bf T}^{\bar{w},w}({\cal C}))
\xto{\sim}\Phi({\bf N}^w({\cal C}))$$
is a weak equivalence of spaces.

Finally, consider the diagram of bi-simplicial categories
$$\begin{array}{ccccc}
{{\bf N}^{\bar{w}}({\cal C}^\op)R}  &\stackrel{p_1}{\slar}&{{\bf T}^{\bar{w},w}({\cal C})}
& \stackrel{p_2}{\srar}&{{\bf N}^w({\cal C})L} \\[4pt]
\big\downarrow{\scriptstyle i} &
&\big\downarrow{\scriptstyle i'}&&\big\Vert\\[4pt]
{{\bf N}^{w}({\cal C}^\op)R} &\stackrel{p_1}{\slar}&
{{\bf T}^{w,w}({\cal C})}&\stackrel{p_2}{\srar} &
{{\bf N}^w({\cal C})L}, \end{array}$$
where $i'$ is the obvious inclusion functor. Applying $\Phi$, we see that the
horizontal functors all induce weak equivalences. The lemma follows.
\enddemo 

Let ${\cal C}$ be a category with cofibrations and two categories
of weak equivalences $v{\cal C}$ and $w{\cal C}$, and write
$${\bf N}^{v,w}{\cal C}={\bf N}\s^v({\bf N}\s^w{\cal C})
\cong{\bf N}\s^w({\bf N}\s^v{\cal C}).$$
This is a bi-simplicial category with cofibrations which again has two
categories of weak equivalences.

\proclaimtitle{Swallowing lemma}
\specialnumber{1.3.10} \proclaim{Lemma}   If $v{\cal C}\subset
w{\cal C}$ then
$$\Phi({\bf N}\s^w{\cal C})=\Phi(({\bf N}^w{\cal C})R)
\xto{\sim}\Phi({\bf N}^{v,w}{\cal C})$$
is a homotopy equivalence with a canonical homotopy inverse.
\endproclaim

\demo{Proof} We claim that for fixed $m$, the iterated degeneracy in
the $v$-direction,
$${\bf N}\s^w{\cal C}
\to{\bf N}^w\s({\bf N}^v_m{\cal C}),$$
is an exact simplicial homotopy equivalence. Given this, the lemma
follows from Lemma~\ref{sim} and from property (iii). The
iterated degeneracy above is induced from the (exact) iterated
degeneracy map ${\cal C}\to{\bf N}_m^v{\cal C}$ in the
simplicial category ${\bf N}\s^v{\cal C}$. This map has a
retraction given by the (exact) iterated\break face map which takes
$c_0\to\cdots\to c_m$ to $c_0$. The other composite takes
$c_0\to\cdots\to c_m$ to the appropriate sequence of identity maps on
$c_0$. There is a natural transformation from this functor to the
identity functor, given by
$$\begin{array}{ccccc}
c_0 & \dline&c_0&\dline\cdots\dline& c_0\phantom{{\scriptstyle
{f_m\circ\dots\circ f_1}}}\\[4pt]
\phantom{\scr{\id}}\big\downarrow{\scriptstyle
\id}&&\phantom{\scr{f_1}}\big\downarrow{\scriptstyle f_1}&&\big\downarrow{\scriptstyle
{f_m\circ\dots\circ f_1}}\\[4pt]  c_0&\stackrel{f_1}{\srar} &
c_1&\stackrel{f_2}{\srar}{\dots}\stackrel{f_m}{\srar} & c_m.\phantom{{\scriptstyle
{f_m\circ\dots\circ f_1}}}\end{array}
$$ 
The natural transformation is through arrows in $v{\cal C}$, and
hence in $w{\cal C}$. The claim now follows from Lemma~\ref{nerve}.
\enddemo 

The proof of \cite[Th.\ 1.6.4]{w} now gives:

\proclaimtitle{Fibration theorem}
\specialnumber{1.3.11}\proclaim{Theorem} \label{fibration} 
Let ${\cal C}$ be a category with cofibrations equipped and two
categories of weak equivalences $v{\cal C}\subset w{\cal C}${\rm ,}
and let ${\cal C}^w$ be the subcategory with cofibrations of
${\cal C}$ given by the objects $A$ such that $*\to A$ is in
$w{\cal C}$. Suppose that ${\cal C}$ has a cylinder functor{\rm ,}
and that $w{\cal C}$ satisfies the cylinder axiom{\rm ,} the saturation
axiom{\rm ,} and the extension axiom. Then
$$\begin{array}{ccc}
{\Phi({\bf N}\s^vS\s{\cal C}^w)}&\srar&
{\Phi({\bf N}\s^wS\s{\cal C}^w)}\\[3pt]
\big\downarrow&&\big\downarrow\\[3pt]
{\Phi({\bf N}\s^vS\s{\cal C})}&\srar&
{\Phi({\bf N}\s^wS\s{\cal C})}\end{array}$$
is a homotopy cartesian square of pointed spaces, and there is a
canonical contraction of the upper right-hand term.
\endproclaim

1.4.\quad Let ${\cal A}$ be an abelian category. We view
${\cal A}$ as a category with cofibrations and weak equivalences
by choosing a null-object and taking the monomorphisms as the
cofibrations and the isomorphisms as the weak equivalences. Let
${\cal E}$ be an additive category embedded as a full subcategory
of ${\cal A}$, and assume that for every exact sequence in
${\cal A}$,
$$0\to A'\to A\to A''\to 0,$$
if $A'$ and $A''$ are in ${\cal E}$ then $A$ is in ${\cal E}$,
and if $A$ and $A''$ are in ${\cal E}$ then $A'$ is in
${\cal E}$. We then view ${\cal E}$ as a subcategory with
cofibrations and weak equivalences of ${\cal A}$ in the sense of
\cite[\S 1.1]{w}.

The category $C^b({\cal A})$ of bounded complexes in ${\cal A}$
is a category with cofibrations and weak equivalences, where the
cofibrations are the degree-wise monomorphisms and the weak
equivalences $zC^b({\cal A})$ are the quasi-isomorphisms. We view
the category $C^b({\cal E})$ of bounded complexes in ${\cal E}$
as a subcategory with cofibrations and weak equivalences of
$C^b({\cal A})$. The inclusion ${\cal E}\to C^b({\cal E})$
of ${\cal E}$ as the subcategory of complexes concentrated in
degree zero, is an exact functor. The assumptions of the fibration
Theorem~\ref{fibration} are satisfied for $C^b({\cal E})$.

\specialnumber{1.4.1}\proclaim{Theorem} \label{resolution} With ${\cal E}$ as above{\rm ,} the
inclusion induces an equivalence
$$\Phi({\bf N}^i\s S\s{\cal E})\xto{\sim}
\Phi({\bf N}^z\s S\s C^b({\cal E})).$$
\endproclaim

\demo{Proof} We follow the proof of \cite[Th.\ 1.11.7]{tt}. Since
the category $C^b({\cal E})$ has a cylinder functor which satisfies
the cylinder axiom with respect to quasi-isomorphisms, the fibration
theorem shows that the right-hand square in the diagram
$$\begin{array}{ccccc}
{\Phi({\bf N}\s^iS\s{\cal E}^i)}&\srar &
{\Phi({\bf N}\s^iS\s C^b({\cal E})^z)}&\srar&
{\Phi({\bf N}\s^zS\s C^b({\cal E})^z)} \\[3pt]
\big\downarrow&&\big\downarrow&&\big\downarrow\\[3pt]
{\Phi({\bf N}\s^iS\s{\cal E})} &\srar &
{\Phi({\bf N}\s^iS\s C^b({\cal E}))} &\srar&
{\Phi({\bf N}\s^zS\s C^b({\cal E}))} \end{array}
 $$
is homotopy cartesian. Moreover, the composite of the maps in the
lower row is equal to the map of the statement, and the upper left-hand and upper right-hand terms are contractible. Hence the theorem is
equivalent to the statement that the left-hand square, and thus the outer
square, are homotopy cartesian. 

Let ${\cal C}_a^b$ be the full subcategory of $C^b({\cal E})$
consisting of the complexes $E_*$ with $E_i=0$ for $i>b$ and
$i<a$. Then $C^b({\cal E})$ is the colimit of the categories
${\cal C}_a^b$ as $a$ and $b$ tend  to $-\infty$ and $+\infty$,
respectively. We consider ${\cal C}_a^b$ as a subcategory with
cofibrations of $C^b({\cal E})$. We first show that there is a weak
equivalence
$$\Phi({\bf N}^i\s S\s{\cal C}^b_a)
\to\prod_{a\leq s\leq b}\Phi({\bf N}^i\s S\s{\cal E}),\hskip6mm
E_*\mapsto (E_b,E_{b-1},\dots,E_a).$$
The map is an isomorphism for $b=a$. If $b>a$, the functor
$$e\:{\cal C}_a^b\to
E({\cal C}_a^a,{\cal C}_a^b,{\cal C}_{a+1}^b),$$
which takes $E_*$ to the extension
$$\sigma_{\leq a}E_*\rightarrowtail
E_*\twoheadrightarrow\sigma_{>a}E_*,$$
is an exact equivalence of categories. Here $\sigma_{\leq n}E_*$ is
the brutal truncation,~\cite[1.2.7]{weibel}. The inverse, given by the
total-object functor, is also exact. Hence, the induced map
$$\Phi({\bf N}^i\s S\s{\cal C}_a^b)\xto{\sim}
\Phi({\bf N}^i\s S\s
E({\cal C}_a^a,{\cal C}_a^b,{\cal C}_{a+1}^b)),
$$
is a homotopy equivalence by Lemma~\ref{nerve}. The additivity
Theorem~\ref{additivity} then shows that
$$(s,q)\:\Phi({\bf N}^i\s S\s
E({\cal C}_a^a,{\cal C}_a^b,{\cal C}_{a+1}^b))\xto{\sim}
\Phi({\bf N}\s^iS\s{\cal C}_a^a)\times
\Phi({\bf N}\s^iS\s{\cal C}_{a+1}^b);$$
thus, we have a weak equivalence
$$\Phi({\bf N}\s^iS\s{\cal C}_a^b)\xto{\sim}
\Phi({\bf N}\s^iS\s{\cal E})\times
\Phi({\bf N}\s^iS\s{\cal C}_{a+1}^b),\hskip6mm
E_*\mapsto (E_a,\sigma_{>a}E_*).$$
It now follows by easy induction that the map in question is a weak
equivalence.

Next, we claim that the map
$$\Phi({\bf N}\s^iS\s{\cal C}_a^{bz})\to
\prod_{a\leq s<b}\Phi({\bf N}\s^iS\s{\cal E}),\hskip6mm
E_*\mapsto(B_{b-1},B_{b-2},\dots,B_a),$$
where $B_i\subset E_i$ are the boundaries, is a weak equivalence. Note
that the exactness of the functors $E_*\mapsto B_i$ uses the fact  that the
complex $E_*$ is acyclic. If $a=b-1$ the functor $E_*\mapsto B_{b-1}$
is an equivalence of categories with exact inverse
functor. Therefore, in this case, the claim follows from Lemma
\ref{nerve}. If $b-1>a$, we consider the functor
$${\cal C}_a^{bz}\to
E({\cal C}_{b-1}^{bz},{\cal C}_a^{bz},{\cal C}_a^{(b-1)z}),$$
which takes the acyclic complex $E_*$ to the extension
$$\tau_{\geq b-1}E_*\rightarrowtail
E_*\twoheadrightarrow\tau_{<b-1}E_*,$$
where $\tau_{\geq n}E_*$ is the good
truncation,~\cite[1.2.7]{weibel}. The functor is exact, since we only
consider acyclic complexes, and it is an equivalence of categories
with exact inverse given by the total-object functor. Hence the
induced map
$$\Phi({\bf N}\s^iS\s{\cal C}_a^{bz})\xto{\sim}
\Phi({\bf N}\s^iS\s
E({\cal C}_{b-1}^{bz},{\cal C}_a^{bz},{\cal C}_a^{(b-1)z}))$$
is a homotopy equivalence by Lemma~\ref{nerve}. The additivity theorem
now shows that
$$\Phi({\bf N}\s^iS\s{\cal C}_a^{bz})\xto{\sim}
\Phi({\bf N}\s^iS\s{\cal E})\times
\Phi({\bf N}\s^iS\s{\cal C}_a^{b-1}),\hskip6mm
E_*\mapsto(B_{b-1},\tau_{<b-1}E_*),$$
is a weak equivalence, and the claim follows by induction.

Statement (4) of the additivity theorem shows that there is a homotopy
commutative diagram
$$\begin{array}{ccc}
{\Phi({\bf N}\s^iS\s{\cal C}_a^{bz})}&\stackrel{\sim}{\srar}  &
{\;\prod_{a\leq s<b}\Phi({\bf N}\s^iS\s{\cal E})}\\[3pt]
\big\downarrow&&\big\downarrow\\[3pt]
{\Phi({\bf N}\s^iS\s{\cal C}_a^b)}&\stackrel{\sim}{\srar} &
{\;\prod_{a\leq s\leq b}\Phi({\bf N}\s^iS\s{\cal E})} \end{array}
 $$
where the horizontal maps are the equivalences established above,
and where the right-hand vertical map takes $(x_s)$ to
$(x_s+x_{s-1})$. It follows that the diagram
$$\begin{array}{ccc}
{\Phi({\bf N}\s^iS\s{\cal C}_0^{0z})} &\srar &
{\Phi({\bf N}\s^iS\s{\cal C}_a^{bz})}\\[3pt]
\big\downarrow&&\big\downarrow\\[3pt]
{\Phi({\bf N}\s^iS\s{\cal C}_0^0)}&\srar &
{\Phi({\bf N}\s^iS\s{\cal C}_a^b))}, \end{array}
 $$
where the maps are induced by the canonical inclusions, is homotopy
cartesian. Indeed, the map of horizontal homotopy fibers may be
identified with the map
$$\prod_{a\leq s<b}\Omega\Phi({\bf N}\s^iS\s{\cal E})\to
\prod_{a\leq s\leq b, s\not=0}\Omega\Phi({\bf N}\s^iS\s{\cal E}),$$
which takes $(x_s)$ to $(x_s+x_{s-1})$, and this, clearly, is a
homotopy equivalence. Taking the homotopy colimit over $a$ and $b$, we
see that the left-hand square in the diagram at the beginning of the
proof is homotopy cartesian.
\enddemo

1.5.\quad  In the remainder of this
section, $A$ will be a discrete valuation ring with quotient field
$K$ and residue field $k$. The main result is
Theorem~\ref{localization} below. It seems unlikely that this result
is valid in the generality of the previous section. Indeed, the proof
of the corresponding result for $K$-theory uses the approximation
theorem \cite[Th.\ 1.6.7]{w}, and this fails for general $\Phi$,
topological Hochschild homology included. Our proof of
Theorem~\ref{localization} uses the equivalence criterion of 
Dundas-McCarthy for topological Hochschild homology, which we now
recall.

If ${\cal C}$ is a category and $n\geq 0$ an integer, we let
$\END_n({\cal C})$ be the category where an object is a tuple
$(c;v_1,\dots,v_n)$ with $c$ an object of ${\cal C}$ and
$v_1,\dots,v_n$ endomorphisms of $c$, and where a morphism from
$(c;v_1,\dots,v_n)$ to $(d;w_1,\dots,w_n)$ is a morphism $f\:c\to d$
in ${\cal C}$ such that $fv_i=w_if$, for $1\leq i\leq n$. We note
that $\END_0({\cal C})={\cal C}$.

\proclaimtitle{[7, Prop.\ 2.3.3]}
\specialnumber{1.5.1}\proclaim{Proposition} \label{eqcriterion}
Let $F\:{\cal C}\to{\cal D}$ be an exact functor of linear
categories with cofibrations and weak equivalences{\rm ,} and suppose that
for all $n\geq 0${\rm ,} the map $|\ob{\bf N}^w\s S\s\END_n(F)|$ 
is an equivalence. Then
$$F_*\:\THH({\bf N}^w\s S\s{\cal C}) \xto{\sim}
\THH({\bf N}^w\s S\s{\cal D})$$
is an ${\cal F}$\/{\rm -}\/equivalence {\rm (}see
Def.~{\rm \ref{F-equivalence})}. 
\endproclaim

Let ${\cal M}_A$ be the category of finitely generated
$A$-modules. We consider two categories with cofibrations and weak
equivalences, $C_z^b({\cal M}_A)$ and $C_q^b({\cal M}_A)$, both
of which have the category of bounded complexes in ${\cal M}_A$
with degree-wise monomorphisms as their underlying category with
cofibrations. The weak equivalences are the categories
$zC^b({\cal M}_A)$ of quasi-isomorphisms and $qC^b({\cal M}_A)$
of chain maps which become quasi-isomorphisms in $C^b({\cal M}_K)$,
respectively. We note that $C^b({\cal M}_A^q)$ and
$C^b({\cal M}_A)^q$ are the categories of bounded complexes of
finitely generated torsion $A$-modules and bounded complexes of
finitely generated $A$-modules with torsion homology, respectively.  

\specialnumber{1.5.2}\proclaim{Theorem} \label{localization} The inclusion functor induces an
${\cal F}$\/{\rm -}\/equivalence
$$\THH({\bf N}\s^zS\s{\cal C}^b({\cal M}_A^q))\xto{\sim}
\THH({\bf N}\s^zS\s{\cal C}^b({\cal M}_A)^q).$$
\endproclaim

\demo{Proof} We show that the assumptions of Proposition
\ref{eqcriterion} are satisfied. The proof relies on Waldhausen's
approximation theorem, \cite[Th.\ 1.6.7]{w}, but in a formulation
due to Thomason, \cite[Th.\ 1.9.8]{tt}, which is particularly
well suited to the situation at hand.

For $n\geq 0$, let $A_n$ be the ring of polynomials in $n$
noncommuting variables with coefficients in $A$, and let
${\cal M}_{A,n}\subset{\cal M}_{A_n}$ be the category of 
$A_n$-modules which are finitely generated as $A$-modules. Then the
category $\END_n(C^b({\cal M}_A))$
(resp.\ $\END_n(C^b({\cal M}_A))^q$,
resp.\ $\END_n(C^b({\cal M}_A^q))$) is canonically isomorphic to the
category $C^b({\cal M}_{A,n})$ (resp.\ $C^b({\cal M}_{A,n})^q$,
resp.\ $C^b({\cal M}_{A,n}^q)$). Here
$C^b({\cal M}_{A,n})^q\subset C^b({\cal M}_{A,n})$ is the full
subcategory of complexes whose image under the forgetful functor
$C^b({\cal M}_{A,n})\to C^b({\cal M}_A)$
lies in $C^b({\cal M}_A)^q$, and similarly for
${\cal M}_{A,n}^q$. We must show that the inclusion functor induces
a weak equivalence
$$|\ob{\bf N}\s^zS\s C^b({\cal M}_{A,n}^q)|\xto{\sim}
|\ob{\bf N}\s^zS\s C^b({\cal M}_{A,n})^q|,$$
for which we use \cite[Th.\ 1.9.8]{tt}. The categories
$C^b({\cal M}_{A,n}^q)$ and $C^b({\cal M}_{A,n})^q$ are both 
complicial bi-Waldhausen categories in the sense of \cite[1.2.4]{tt},
which are closed under the formation of canonical homotopy pushouts
and homotopy pullbacks in the sense of \cite[1.9.6]{tt}. The inclusion
functor
$$F\:C^b({\cal M}_{A,n}^q)\to C^b({\cal M}_{A,n})^q$$
is a complicial exact functor in the sense of \cite[1.2.16]{tt}. We
must verify the conditions \cite[1.9.7.0--1.9.7.3]{tt}. These
conditions are easily verified with the exception of condition 1.9.7.1
which reads: for every object $B$ of $C^b({\cal M}_{A,n})^q$, there
exist an object $A$ of $C^b({\cal M}_{A,n}^q)$ and a map
$FA\xto{\sim}B$ in $zC^b({\cal M}_{A,n})^q$. This follows from
Lemma~\ref{qis} below.
\enddemo

\specialnumber{1.5.3} \proclaim{Lemma}\label{qis}
Let $A$ be a commutative noetherian ring{\rm ,} and
let $B$ be a not necessarily commutative $A$\/{\rm -}\/algebra. Let $C_*$ be a
bounded complex of left $B$\/{\rm -}\/modules which as $A$\/{\rm -}\/modules are finitely
generated and suppose that the homology of $C_*$ is annihilated by some
power of an ideal $I\subset A$. Then there exists a quasi\/{\rm -}\/isomorphism
$$C_*\xto{\sim}D_*$$
with $D_*$ a bounded complex of left $B$\/{\rm -}\/modules which as $A$\/{\rm -}\/modules
are finitely generated and annihilated by some power of $I$.
\endproclaim

\demo{Proof} Let $n$ be an integer such that for all $i\geq n$, $C_i$
is annihilated by some power of $I$. We construct a quasi-isomorphism
$C\xto{\sim}C''$ to a bounded complex $C''$ of left $B$-modules which
as $A$-modules are finitely generated and such that for all $i\geq
n-1$, $C_i''$ is annihilated by some power of $I$. The lemma follows
by easy induction. To begin we note that the exact sequences
\begin{eqnarray*}
{}& 0\to Z_n \to C_n \xto{d} B_{n-1} \to 0, \\
{}& 0\to B_{n-1} \to Z_{n-1} \to H_{n-1} \to 0, \\
\end{eqnarray*}
show that $Z_{n-1}$ is annihilated by some power of $I$, say, by
$I^r$. As an $A$-module $Z_{n-1}$ is finitely generated because
$C_{n-1}$ is a finitely generated $A$-module and because $A$ is
noetherian. Hence, by the Artin-Rees lemma, \cite[Th.\ 8.5]{mat},
we can find $s\geq 1$ such that $Z_{n-1}\cap I^sC_{n-1}\subset
I^rZ_{n-1}=0$. We now define $C''$ to be the complex with $C_i''=C_i$,
if $\neq n-1,n-2$, with $C_{n-1}''=C_{n-1}/I^sC_{n-1}$, and with
$C_{n-2}''$ given by the pushout square
$$\begin{array}{ccc}
{C_{n-1}} &\stackrel{d}{\srar} &{C_{n-2}}\\[3pt]
 \ddownarrow{\scriptstyle\pr} &&\ddownarrow\\[4pt]
 {C_{n-1}''}&\srar&
{C_{n-2}''.} \end{array}
 $$ 
There is a unique differential on $C''$ such that the canonical
projection $C\to C''$ is a map of complexes. The kernel complex $C'$
is concentrated in degrees $n-1$ and $n-2$. The differential
$C_{n-1}'\to C_{n-2}'$ is injective, since $Z_{n-1}\cap I^sC_{n-1}$ is
zero, and surjective, since the square is a pushout. Hence, the
homology sequence associated with the short exact sequence of
complexes
$$0\to C'\to C\to C''\to 0$$
shows that $C\to C''$ is a quasi-isomorphism. And by construction,
some power of $I$ annihilates $C_i''$, if $i\geq n-1$.
\enddemo

We thank Thomas Geisser and Stefan Schwede for help with the argument
above.

Let $C_z^b({\cal P}_A)$ and $C_q^b({\cal P}_A)$ be the category
of bounded complexes of finitely generated projective $A$-modules
considered as a subcategory with cofibrations and weak equivalences of
$C_z^b({\cal M}_A)$ and $C_q^b({\cal M}_A)$, respectively.

\specialnumber{1.5.4}\proclaim{Proposition} \label{fgproj} The inclusion functor induces an
${\cal F}$\/{\rm -}\/equivalence
$$\THH({\bf N}\s^zS\s C^b({\cal P}_A)^q)
\xto{\sim}\THH({\bf N}\s^zS\s C^b({\cal M}_A)^q).$$
\endproclaim

\demo{Proof} Let $A_n$ and ${\cal M}_{A,n}$ be as in the proof of
Theorem~\ref{localization}, and let ${\cal P}_{A,n}$ be the full
subcategory of ${\cal M}_{A,n}$ consisting of the $A_n$-modules
which as $A$-modules are finitely generated projective. Then 
$\END_n(C^b({\cal M}_A))^q$ and $\END_n(C^b({\cal P}_A))^q$ are
canonically isomorphic to $C^b({\cal M}_{A,n})^q$ and
$C^b({\cal P}_{A,n})^q$, respectively, and we must show that the
inclusion functor induces a weak equivalence
$$|\ob{\bf N}\s^zS\s C^b({\cal P}_{A,n})^q|
\xto{\sim}|\ob{\bf N}\s^zS\s C^b({\cal M}_{A,n})^q|.$$
Again, we use \cite[Th.\ 1.9.8]{tt}, where the nontrivial thing to
check is condition 1.9.7.1: for every object $C_*$ of
$C^b({\cal M}_{A,n})^q$, there exists an object $P_*$ of
$C^b({\cal P}_{A,n})^q$ and a map $P_*\xto{\sim}C_*$ in
$zC^b({\cal M}_{A,n})^q$. But this follows from \cite[Chap.\ XVII,
Prop.\ 1.2]{ce}. Indeed, let $\e\:P_{*,*}\to C_*$ be a projective
resolution of $C_*$ regarded as a complex of $A$-modules. We may
assume that each $P_{i,j}$ is a finitely generated $A$-module, and
since $A$ is regular, that $P_{i,j}$ is zero for all but finitely many
$(i,j)$. Furthermore, it is proved in {\it loc.cit.} that there
exists an $A_n$-module structure on $P_{*,*}$ such that $\e$ is
$A_n$-linear. Hence, the total complex
$P_*=\operatorname{Tot}(P_{*,*})$ is in $C^b({\cal P}_{A,n})$ and
$\operatorname{Tot}(\e)\:P_*\xto{\sim}C_*$ is in
$zC^b({\cal M}_{A,n})$. It follows that $P_*$ is in
$C^b({\cal P}_{A,n})^q$ as desired.
\enddemo

\specialnumber{1.5.5}\numbereddemo{Definition} \label{def}We define ring ${\Bbb T}$-spectra
$$T(A|K)=T(C_q^b({\cal P}_A)),\hskip4mm
T(A)=T(C_z^b({\cal P}_A)),\hskip4mm
T(k)=T(C_z^b({\cal P}_A)^q)$$
and let $\TR^n(A|K;p)$, $\TR^n(A;p)$, and $\TR^n(k;p)$ be the
associated $C_{p^{n-1}}$-fixed point ring spectra.
\enddemo

We show that the definition of the spectra $\TR^n(A;p)$ and $\TR^n(k;p)$
given here agrees with the usual definition. By Morita
invariance,~\cite[Prop.\ 2.1.5]{dm}, it suffices to show
that there are canonical isomorphisms of spectra
$$\TR^n(A;p)\simeq\TR^n({\cal P}_A;p),\hskip4mm
\TR^n(k;p)\simeq\TR^n({\cal P}_k;p),$$
compatible with the maps $R$, $F$, $V$, and $\mu$. Here the exact
category ${\cal P}_R$ is considered a category with cofibrations
and weak equivalences in the usual way. It follows from
Theorem~\ref{resolution}, applied to the functor
$\Phi({\cal C})=\THH({\cal C})^{C_r}$, and Proposition~\ref{fibrant} that the map induced by the
inclusion functor
$$T({\cal P}_A)\to T(C_z^b({\cal P}_A))=T(A)$$
is an ${\cal F}$-equivalence. This gives the first of the stated
isomorphisms of spectra. A similar argument shows that the
inclusion functor induces an ${\cal F}$-equivalence
$$T({\cal P}_k)=T({\cal M}_k)\to T(C_z^b({\cal M}_k)).$$
By devisage, \cite[Th.\ 1]{dundas}, the same is true for
$$T(C_z^b({\cal M}_k))\to T(C_z^b({\cal M}_A^q)).$$
Finally, Theorem~\ref{localization} and Proposition~\ref{fgproj}
show that the maps induced from the inclusion functors
$$T(C_z^b({\cal M}_A^q)) \xto{\sim}
T(C_z^b({\cal M}_A)^q) \stackrel{\sim}{\longleftarrow}
T(C_z^b({\cal P}_A)^q)=T(k)$$
are both ${\cal F}$-equivalences. This establishes the second of
the stated isomorphisms of spectra. Let
$$i_*\:\TR^n(A;p)\to\TR^n(k;p)$$
be the map induced from the reduction.

\specialnumber{1.5.6}\proclaim{Theorem} \label{localizationtr}
For all $n\geq 1${\rm ,} there is a
natural cofibration sequence of spectra
$$\TR^n(k;p)\xto{i^!}\TR^n(A;p)\xto{j_*}\TR^n(A|K;p)
\xto{\partial}\Sigma\TR^n(k;p),$$
and all maps in the sequence commute with the maps $R${\rm ,} $F${\rm ,}
$V${\rm ,} and $\mu$. The map $j_*$ is a map of ring spectra{\rm ,} and the maps
$i^!$ and $\partial$ are maps of $\TR^n(A;p)$\/{\rm -}\/module spectra. Here
$\TR^n(k;p)$ is considered a $\TR^n(A;p)$\/{\rm -}\/module spectrum via the map
$i_*$. Moreover{\rm ,} the preferred homotopy limits form a cofibration
sequence of spectra.
\endproclaim

\demo{Proof} We have a commutative square of symmetric orthogonal
${\Bbb T}$-spectra
$$\begin{array}{ccc}
{T(C_z^b({\cal P}_A)^q)}&\lrar&
{T(C_q^b({\cal P}_A)^q)} \\[3pt]
\big\downarrow&&\big\downarrow\\[4pt]
{T(C_z^b({\cal P}_A))} &\lrar &
{T(C_q^b({\cal P}_A)),} \end{array}
$$
and the fibration Theorem~\ref{fibration} applied to the functor
$\Phi({\cal C})=\THH({\cal C})^{C_r}$ shows that the corresponding
square of $C_r$-fixed point spectra is homotopy cartesian. It follows
that there is natural cofibration sequence of spectra
$$\TR^n(k;p)\xto{i^!}\TR^n(A;p)\xto{j_*}\TR^n(A|K;p)
\xto{\partial}\Sigma\TR^n(k;p),$$
compatible with $R$, $F$, $V$ and $\mu$. It is clear that this is a
sequence of $\TR^n(A;p)$-module spectra.
\enddemo

\specialnumber{1.5.7}\proclaim{Addendum} \label{localizationtc}
There is a natural map of
cofibration sequences
$$\begin{array}{ccccccc}
{K(k)}&\hskip-6pt\stackrel{i^!}{\srar}\hskip-6pt&  
{K(A)} &\hskip-6pt\stackrel{j_*}{\srar}\hskip-6pt&  
{K(K)}&\hskip-6pt\stackrel{\partial}{\srar}\hskip-6pt&   
{\Sigma K(k)}  \\[4pt]
\big\downarrow{\scriptstyle\tr} &&\big\downarrow{\scriptstyle\tr}
&&\big\downarrow{\scriptstyle\tr} &&\big\downarrow{\scriptstyle\tr} \\[4pt]
 {\TC(k;p)} &\hskip-6pt\stackrel{i^!}{\srar}\hskip-6pt&  
{\TC(A;p)} &\hskip-6pt\stackrel{j_*}{\srar}\hskip-6pt&  
{\TC(A|K;p)} &\hskip-6pt\stackrel{\partial}{\srar}\hskip-6pt&   
{\Sigma\TC(k;p)} \end{array}
$$
and the vertical maps are all maps of ring spectra.
\endproclaim

\specialnumber{1.5.8}\numbereddemo{{R}emark} Let $X$ be a regular affine scheme and let
$i\:Y\hookrightarrow X$ be a closed subscheme with open complement
$j\:U\hookrightarrow X$. Then, more generally, the proof of
Theorem~\ref{localizationtr} gives a cofibration sequence of spectra
$$\TR'\!{}^n(Y;p) \xto{i^!}
\TR^n(X;p) \xto{j_*}
\TR^n(X|U;p) \xto{\partial}
\Sigma\TR'\!{}^n(Y;p),$$
where the three terms are   as in Definition~\ref{def} with
${\cal P}_A$ replaced by the category ${\cal P}_X$ of locally
free ${\cal O}_X$-modules of finite rank. The weak equivalences are
the quasi-isomorphisms, $zC^b({\cal P}_X)$, and the chain maps
which become quasi-isomorphisms after restriction to $U$,
$qC^b({\cal P}_X)$, respectively. Similarly, the argument following
Definition~\ref{def} gives canonical isomorphisms of spectra
$$\TR^n(X;p)\simeq\TR^n({\cal P}_X;p),\hskip4mm
\TR'\!{}^n(Y;p)\simeq\TR^n({\cal M}_Y;p),$$
where ${\cal M}_Y$ is the category of coherent
${\cal O}_Y$-modules. Moreover, if $Y$ is regular, the resolution
theorem, \cite[prop. 2.2.3]{dm}, shows that $\TR^n({\cal M}_Y;p)$
is canonically isomorphic to $\TR^n({\cal P}_Y;p)$.
\enddemo

\section{The homotopy groups of $T(A|K)$} \label{log}

2.1.\quad In this section we evaluate the
homotopy groups with $\Z/p$-coefficients of the topological Hochschild
spectrum $T(A|K)$. We first fix some conventions.

Let $G$ be a finite group and let $k$ be a commutative ring. The
category of chain complexes of left $kG$-modules and chain homotopy
classes of chain maps is a triangulated category and a closed
symmetric monoidal category, and the two structures are compatible. The
same is true for the category of $G$-CW-spectra and homotopy classes
of cellular maps. We  fix our choices for the triangulated and
closed symmetric monoidal structures in such a way that the cellular chain functor
preserves our choices.

We first consider complexes. If $f\:X\to Y$ is a chain map, we define
the mapping cone $C_f$ to be the complex
$$(C_f)_n=Y_n\oplus X_{n-1},\hskip5mm d(y,x)=(dy-f(x),-dx),$$
and the suspension $\Sigma X$ to be the cokernel of the inclusion
$\iota\:Y\to C_f$ of the first summand. More explicitly,
$$(\Sigma X)_n=X_{n-1},\hskip5mm d_{\Sigma X}(x)=-d_X(x).$$
Then, by definition, a sequence $X\xto{f} Y\xto{g} Z\xto{h} \Sigma X$
is a triangle or a cofibration sequence if it isomorphic to the
{\it distinguished} triangle
$$X\xto{f}Y\xto{i}C_f\xto{\partial}\Sigma X,$$
where $\partial$ is the canonical projection. If $X\xto{f}Y\xto{g}Z$
is a short exact sequence of complexes then the projection $p\:C_f\to
Z$, $p(y,x)=g(y)$, is a quasi-isomorphism and the composite 
$$
H_nZ
\begin{array}{c}  \scr{p_\ast}\\[-8pt]{\longleftarrow}\\[-10pt]\scr{\sim}\end{array}
H_nC_f\xto{\partial_\ast}H_n\Sigma X=H_{n-1}X$$ is equal to the connecting homomorphism.

Let $X$ and $Y$ be two complexes. We define the tensor product
complex by
$$(X\otimes Y)_n=\bigoplus_{s+t=n}X_s\otimes Y_t;\hskip5mm
d(x\otimes y)=dx\otimes y+(-1)^{|x|}x\otimes dy,$$
and the complex of ($k$-linear) homomorphisms by
$$\Hom(X,Y)_n=\prod_{s\in\Z}\Hom(X_s,Y_{n+s});\hskip5mm
d(f(x))=(df)(x)+(-1)^{|f|}f(dx).$$
We note that $Z_0\Hom(X,Y)$ is equal to the set of chain maps from $X$
to $Y$ and that $H_0\Hom(X,Y)$ is equal to the set of chain homotopy
classes of chain maps from $X$ to $Y$. The adjunction and
twist isomorphisms are given by
\begin{eqnarray*}
&&\phi\:\Hom(X\otimes Y,Z)\to\Hom(X,\Hom(Y,Z)),\hskip5mm
\phi(f)(x)(y)=f(x\otimes y),\\
&&\gamma\: X\otimes Y\to Y\otimes X,\hskip5mm\gamma(x\otimes
y)=(-1)^{|x||y|}y\otimes x. 
\end{eqnarray*}
The triangulated and closed symmetric monoidal structures are
compatible in the sense that 
$$\Sigma(X\otimes Y)=(\Sigma X)\otimes Y$$
and that if $W$ is a complex and $X\xto{f}Y\xto{g}Z\xto{h}\Sigma X$ is
a triangle, then so is
$$X\otimes W\xto{f\otimes 1}Y\otimes W\xto{g\otimes 1}Z\otimes W
\xto{h\otimes 1}\Sigma X\otimes W.$$
Indeed, the isomorphism
$$\rho\:C_f\otimes W\xto{\sim}C_{f\otimes W},\hskip5mm
\rho((y,x)\otimes w)=(y\otimes w,x\otimes w),$$
and the identity map of $X\otimes W$, $Y\otimes W$, and $\Sigma
X\otimes W$ define an isomorphism of the appropriate distinguished
triangles.

We define the homology of $X$ with $\Z/m$-coefficients by
$$H_*(X,\Z/m) = H_*(M_m\otimes X),$$
where $M_m$ is the Moore complex given by the distinguished triangle
$$k \xto{m} k \xto{\iota} M_m \xto{\beta} \Sigma k.$$
Suppose that $X$ is $m$-torsion free such that $X \xto{m} X \xto{\pr}
X/mX$ is a short-exact sequence of complexes. Then the composite
$$H_n(X/mX) \begin{array}{c} \scr{p_*}\\[-8pt]{\longleftarrow}\\[-10pt] \scr{\sim}\end{array}
H_n(C_{m}) \begin{array}{c} \scr{\rho_*}\\[-8pt]{\longleftarrow}\\[-10pt] \scr{\sim}\end{array}
H_n(M_m\otimes X) \xto{\beta}
H_n(\Sigma X)=H_{n-1}(X)$$
is equal to the connecting homomorphism.

We next consider the category of $G$-CW-spectra and homotopy classes
of cellular maps, see \cite[Chap.\ I, \S5]{lms}. This category, we
recall, is equivalent to the $G$-stable category. In one direction,
the equivalence associates to a $G$-CW-spectrum $X$ the underlying
$G$-spectrum $UX$. In the other direction, we choose a functorial
$G$-CW-replacement $\Gamma X$ such that $U\Gamma X\xto{\sim}X$.

If $X$ and $Y$ are two $G$-CW-spectra, the smash product
$UX\wedge UY$ has a canonical $G$-CW-structure. But the function
spectrum $F(UX,UY)$ usually does not. Instead we consider $\Gamma
F(UX,UY)$. This defines the closed symmetric monoidal structure. 

The mapping cone of a celluar map $f\:X\to Y$ is defined by
$$C_f=Y\cup_X([0,1]\wedge X),$$
where we use $1$ as the base point for the smash product. The interval
is given the usual CW-structure with a single $1$-cell oriented from
$0$ to $1$, and the mapping cone is given the induced $G$-CW-structure.
Collapsing the image of the canonical inclusion $i\:Y\to C_f$ to the
base point defines the map
$$\partial\:C_f\to S^1\wedge X=\Sigma X,$$
where $S^1=[0,1]/\partial[0,1]$ with the induced CW-structure. We then
define the distinguished triangles to be   sequences of the form
$$X\xto{f}Y\xto{i}C_f\xto{\partial}\Sigma X.$$
Again, the triangulated and the closed symmetric monoidal structures
are compatible. Indeed, the associativity isomorphism, which is part
of the monoical structure, gives rise to canonical isomorphisms
$$\alpha\:\Sigma(X\wedge W) \xto{\sim} (\Sigma X)\wedge W, \hskip6mm
\rho\:C_f\wedge W\xto{\sim}C_{f\wedge W}.$$

The choices made above are preserved by the cellular chain
functor. To be more precise, if $X$ (resp.\ $f\:X\to Y$) is a
$G$-CW-spectrum (resp.\ a cellular map), then the suspension
isomorphism gives rise to a canonical isomorphism of complexes $\Sigma
C_*(X;k)\xto{\sim}C_*(\Sigma X;k)$
(resp.\ $C_*(C_f;k)\xto{\sim}C_{f_*}$).  Under these identifications,
the cellular chain functor carries the distinguished triangles of
$G$-CW-spectra to the distinguished triangles of complexes of left
$kG$-modules. Similarly, if $X$ and $Y$ are two $G$-CW-complexes, then
the K\"unneth isomorphism gives a canonical isomorphism
$C_*(X;k)\otimes C_*(Y;k)\xto{\sim} C_*(X\wedge Y;k)$.

We define the homotopy groups of $X$ with $\Z/m$-coefficients by
$$\pi_*(X,\Z/m) = \pi_*(M_m\wedge X),$$
where $M_m$ is the Moore spectrum given by the distinguished triangle
$$S^0 \xto{m} S^0 \xto{\iota} M_m \xto{\beta} S^1,$$
and the homotopy groups with $\Zp$-coefficients by
$$\pi_*(X,\Zp) = \pi_*(\hl{v}{(M_{p^v}\wedge X)}).$$
The latter are related to the former by the Milnor sequence
$$0 \to \lim_{\underset{\scr{v}}{\longleftarrow}}{\hskip-6pt}^{1} {\pi_{q+1} (X,\Z/p^v)} \to
\pi_q(X,\Zp) \to \lim_{\underset{\scr{v}}{\longleftarrow}} {\pi_q(X,\Z/p^v)} \to 0.$$
We shall often abbreviate $\pi_q(X,\Z/p)$ and write $\bar\pi_q(X)$. Let
$H\Z/m$ be the Eilenberg-MacLane spectrum for $\Z/m$. It is a ring
spectrum, and we let $\e\in\pi_1(H\Z/m,\Z/m)$ be the unique element
such that $\beta(\e)=1$. Then for left $H\Z/m$-module spectra $X$, we
have a natural sum-diagram  
\begin{equation}\label{sumdiagram}
{X}\begin{array}{c}{\scriptstyle\iota\wedge\id} \\[-8pt]\srar\\[-10pt]\slar\\[-8pt] {\scriptstyle
r}\end{array}  {M_m\wedge X}
\begin{array}{c}{\scriptstyle \beta\wedge\id}\\[-8pt]\srar\\[-10pt]\slar\\[-8pt] {\scriptstyle
s}\end{array} {\Sigma X,}  \speqnu{2.1.1}
\end{equation} 
where $s$ is the composite
$$S^1\wedge X \stackrel{\e\wedge\id}{\srar}
M_m\wedge H\Z/m\wedge X \stackrel{\id\wedge\mu}{\srar}
M_m\wedge X,$$
and where $r$ is determined by the requirement that $r\circ\iota=\id$
and $r\circ s=0$. 

We recall Connes' operator. Let ${\Bbb T}$ be the space
$S(\C)$ of complex numbers of length $1$ considered as a group under
multiplication. We give ${\Bbb T}$ the orientation induced from the
standard orientation of the complex plane, and let $[{\Bbb T}]\in
H_1({\Bbb T})$ be the corresponding fundamental class. The reduced
homology of a ${\Bbb T}$-space $X$ has a natural differential given
by the composite
$$d\:\tilde H_q(X) \xto{[{\Bbb T}]}
\tilde H_{q+1}({\Bbb T}_+\wedge X) \xto{\mu}
\tilde H_{q+1}(X),$$
where the left-hand map is given by the K\"unneth isomorphism and the
right-hand map is induced by the action map. There is a sum-diagram
$$ 
{\Z/2\cdot\eta=\pi_1^S(S^0)}\begin{array}{c} {\scriptstyle e} \\[-8pt]\srar\\[-10pt]\slar\\[-8pt] 
{\scriptstyle c} \end{array}{\pi_1^S({\Bbb T}_+)}\begin{array}{c} {\scriptstyle
h}\\[-8pt]\srar\\[-10pt]\slar\\[-8pt]  {\scriptstyle \sigma} \end{array} H_1({\Bbb T})=\Z\cdot[{\Bbb
T}] $$
where $h$ is the Hurewitz homomorphism, $e$ is induced from
the map $S^0\to{\Bbb T}_+$ which takes the nonbase-point of
$S^0$ to $1\in{\Bbb T}$, $c$ is induced from the map
${\Bbb T}_+\to S^0$ which collapses ${\Bbb T}$ to the
nonbase-point of $S^0$, and $\sigma$ is determined by $h\sigma=\id$
and $c\sigma=0$. Let $T$ be a ${\Bbb T}$-spectrum. Then Connes'
operator is the map
\begin{equation}\label{connes}
d\:\pi_q(T) \stackrel{[{\Bbb T}]\wedge-}{\srar}
\pi_{q+1}({\Bbb T}_+\wedge T) \xto{\mu_T}
\pi_{q+1}(T). \speqnu{2.1.2}
\end{equation}
If $T=\HH(A)$ is the Hochschild spectrum of a ring $A$, then this
definition agrees with Connes' original definition, \cite[Prop.\ 1.4.6]{h}. We recall from {\it op.~cit.},
Lemma~1.4.2, that, in general, $dd=d\eta=\eta d$. Hence, $d$ is a differential, provided
that multiplication by $\eta$ is trivial on $\pi_*(T)$. This is the
case, for instance, if multiplication by $2$ on $\pi_*(T)$ is an
isomorphism. 

\vglue12pt 2.2. \quad We next recall the notion of differentials with
logarithmic poles. The standard reference for this material
is~\cite{k}. A pre-log structure on a ring $R$ is a map of monoids
$$\alpha\:M\to R,$$
where $R$ is considered a monoid under multiplication. By a log ring
we mean a ring with a pre-log structure. A derivation of a log ring
$(R,M)$ into an $R$-module $E$ is a pair of maps
$$(D,D\log)\:(R,M)\to E,$$
where $D\:R\to E$ is a derivation and $D\log\:M\to E$ a map of
monoids, such that for all $a\in M$,
$$\alpha(a)D\log a=D\alpha(a).$$
A log differential graded ring $(E^*,M)$ consists of a differential
graded ring $E^*$, a pre-log structure $\alpha\:M\to E^0$, and a
derivation $(D,D\log)\:(E^0,M)\to E^1$ such that $D$ is equal to the
differential $d\:E^0\to E^1$ and such that $d\circ D\log=0$. 

There is a universal example of a derivation of a log ring $(R,M)$
given by the $R$-module
$$\omega_{(R,M)}^1=(\Omega_R^1\oplus(R\otimes_{\Z}M^{\gp}))/
\langle d\alpha(a)-\alpha(a)\otimes a\mid a\in M\rangle,$$
where $M^{\gp}$ is the group completion (or Grothendieck group) of
$M$ and $\langle\dots\rangle$ denotes the submodule generated by the
indicated elements. The structure maps are 
$$
\begin{array}{rlrl}
d\: R &\hsm\to\omega_{(R,M)}^1,  &\quad da&\hsm =da\oplus 0,\\
d\log\: M&\hsm \to\omega_{(R,M)}^1,&\quad  d\log a&\hsm =0\oplus(1\otimes a). 
\end{array}
$$
The exterior algebra
$$\omega_{(R,M)}^*=\Lambda_R^*(\omega_{(R,M)}^1)$$
endowed with the usual differential is the universal log differential
graded ring whose underlying log ring is $(R,M)$. We stress that here
and throughout we use $\Omega_R^1$ to mean the {\it absolute}
differentials.

Let $A$ be a complete discrete valuation ring with quotient field $K$
and {\it perfect} residue field $k$ of mixed characteristic
$(0,p)$. We recall the structure of $A$ from \cite[\S5, Th.\
4]{serre}. Let $W(k)$ be the ring of Witt vectors in $k$, and let
$K_0$ be the quotient field of $W(k)$. There is a unique ring
homomorphism
$$f\:W(k)\to A$$
such that the induced map of residue fields is the identity
homomorphism. We will always view $A$ as an algebra over $W(k)$ via
the map $f$. Moreover, if $\pi_K$ is a generator of the maximal ideal
${\frak m}_K\subset A$, then
\begin{equation}\label{structure}
A=W(k)[\pi_K]/(\phi_K(\pi_K)), \speqnu{2.2.1}
\end{equation}
and the minimal polynomial takes the form
$$\phi_K(x)=x^{e_K}+p\theta_K(x),$$
where $e_K=|K\!:\!K_0|$ is the ramification index and where
$\theta_K(x)$ is a polynomial of degree less that $e_K$ such that
$\theta_K(0)$ is a unit in $W(k)$. It follows that $\theta_K(\pi_K)$
is a unit and that
$$-p=\pi_K^{e_K}\theta_K(\pi_K)^{-1}.$$
We will use this formula on numerous occasions in the following. The
valuation ring $A$ has a canonical pre-log structure given by the
inclusion
$$\alpha\:M=A\cap K^\times\hookrightarrow A.$$
Let $v_K\:K^\times\to\Z$ be the valuation.

\specialnumber{2.2.2}\proclaim{Proposition} \label{poincare}
There is a natural short exact sequence
$$0\to\Omega_A^1\to\omega_{(A,M)}^1\xto{\res}k\to 0,$$
where \pagebreak $\res(ad\log b)=av_K(b)+{\frak m}_K$.
\endproclaim

\demo{Proof} If $a\in A\cap K^\times$ then $av_K(a)\in{\frak m}_K$,
and hence, the composition of the two maps in the statement is zero.
Only the exactness in the middle needs proof. Let $ad\log b$ be an
element of $\omega_{(A,M)}^1$ and write $b=\pi_K^iu$ with $u\in
A^\times$. Then
$$ad\log b=iad\log\pi_K+au^{-1}du.$$
Suppose that $\res(ad\log b)=ia+{\frak m}_K$ is trivial. Then
$ia\in{\frak m}_K$, which implies that $ia\pi_K^{-1}\in A$, and
hence, $iad\log\pi_K=ia\pi_K^{-1}d\pi_K$.
\enddemo

We define the module of relative differentials
$$\omega_{(A,M)/W(k)}^1
=\big(\Omega_{A/W(k)}^1\oplus(A\otimes_{\Z}K^\times)\big)
\big/\big\langle da-a\otimes a\;\big|\;a\in A\cap
K^\times\big\rangle.$$
Again, there is a natural exact sequence
$$0\to\Omega_{A/W(k)}^1\to\omega_{(A,M)/W(k)}^1\xto{\res}k\to 0.$$

\specialnumber{2.2.3} \proclaim{Lemma}\label{logdiff}
Let $\pi_K\in A$ be a uniformizer with
minimal polynomial $\phi_K(x)$. Then the element $d\log\pi_K$
generates the $A$\/{\rm -}\/module $\omega_{(A,M)/W(k)}^1${\rm ,} and its
annihilator is the ideal generated by $\phi_K'(\pi_K)\pi_K$. This
ideal contains $p$.
\endproclaim

\demo{Proof} Since every element of $K^\times$ can be written as a
product $\pi_K^iu$ with $i\in\Z$ and $u\in A^\times$, the formula
$$d\log(\pi_K^iu)=id\log\pi_K+u^{-1}du$$
shows that $\omega_{(A,M)/W(k)}^1$ is generated by
$d\log\pi_K$. The relation identifies 
$$\phi_K'(\pi_K)\pi_Kd\log\pi_K=d(\phi_K(\pi_K))=0,$$
so the annihilator ideal is generated $\phi_K'(\pi_K)\pi_K$.
\enddemo

\specialnumber{2.2.4} \proclaim{Lemma}\label{divisible}
For all $i>0${\rm ,} there is a natural exact
sequence
$$A\otimes_{W(k)}\Omega_{W(k)}^i\to\omega_{(A,M)}^i
\to\omega_{(A,M)/W(k)}^i\to 0,$$
and the left\/{\rm -}\/hand group is uniquely divisible.
\endproclaim

\demo{Proof} The stated sequence for $i=1$ follows from the diagram
$$\begin{array}{ccccccccc}
{0}&\srar& 
{\Omega_A^1} &\srar &
{\omega_{(A,M)}^1}&\srar &
{k} &\srar &
{0} \\[4pt]
&&\big\downarrow&&\big\downarrow&&\big\Vert\\[4pt]
{0}&\srar&
{\Omega_{A/W(k)}^1} &\srar&
{\omega_{(A,M)/W(k)}^1}&\srar&
{k} &\srar&
{0} \end{array}
$$
with horizontal exact sequences and from the standard exact sequence
$$A\otimes_{W(k)}\otimes\,\Omega_{W(k)}^1\to\Omega_A^1\to
\Omega_{A/W(k)}^1\to 0.$$
We show that the group $\Omega_{W(k)}^1\xto{\sim}\HH_1(W(k))$ is a \pagebreak
uniquely divisible group or, more generally, that $\HH_i(W(k))$ is
uniquely divisible, for all $i>0$. Since $W(k)$ is torsion-free and
since $W(k)/p=k$, the coefficient sequence takes the form
$$\cdots\to\HH_{i+1}(k)\to\HH_i(W(k))\xto{p}\HH_i(W(k))
\to\HH_i(k)\to\cdots\; .$$
But $\HH_i(k)=0$, for $i>0$, since $k$ is perfect, \cite[Lemma 
5.5]{hm}. This proves the lemma for $i=1$. In particular, the maximal
divisible sub-$A$-module of $\omega_{(A,M)}^1$ is equal to the image
of $A\otimes_{W(k)}\Omega_{W(k)}^1$, and $\omega_{(A,M)}^1$ is the sum
of this divisible module $D$ and the cyclic torsion $A$-module
$\omega_{(A,M)/W(k)}^1$. It follows that for $i>1$,
$\omega_{(A,M)}^i=\Lambda^i_AD$, and this in turn is the image of the left-hand map of the statement.
\enddemo

\specialnumber{2.2.5}\proclaim{{C}orollary}\label{ptorsion}
The $p$\/{\rm -}\/torsion submodule of
$\omega_{(A,M)}^1$ is
$${}_p\omega_{(A,M)}^1=A/p\cdot d\log(-p).$$
\endproclaim

\demo{Proof} It follows from Lemma \ref{divisible} that the canonical
map
$${}_p\omega_{(A,M)}^1\xto{\sim}{}_p\omega_{(A,M)/W(k)}^1$$
is an isomorphism. Let $\pi_K$ be a uniformizer with minimal
polynomial $\phi_K(x)$. Then by Lemma~\ref{logdiff},
$$\omega_{(A,M)/W(k)}^1
=A/(\pi_K\phi_K'(\pi_K))\cdot d\log\pi_K.$$
We write $\phi_K(x)=x^{e_K}+p\theta_K(x)$ such that
$-p=\pi_K^{e_K}\theta_K(\pi_K)^{-1}$. Hence, on the one 
hand, we have
$$\pi_K\phi_K'(\pi_K)=e_K\pi_K^{e_K}+p\pi_K\theta_K'(\pi_K)
=(e_K-\pi_K\theta_K'(\pi_K)\theta_K(\pi_K)^{-1})\pi_K^{e_K},$$
and on the other hand,
$$d\log(-p)=d\log(\pi_K^{e_K}\theta_K(\pi_K)^{-1})
=(e_K-\pi_K\theta_K'(\pi_K)\theta_K(\pi_K)^{-1})d\log\pi_K.$$
The claim follows.
\enddemo

Let $L$ be a finite extension of $K$, let $B$ be the integral
closure of $A$ in $L$, and let $e_{L/K}=e_L/e_K$ be the ramification
index of $L/K$. Then the following diagram commutes
$$\begin{array}{ccc}
\omega_{(A,M_A)/W(k)}^1&\stackrel{\res_A}{\srar}&
A/{\frak m}_K \\[4pt]
\big\downarrow{\scriptstyle i_*} &&\big\downarrow {\scriptstyle e_{L/K}\cdot i}\\[4pt]
\omega_{(B,M_B)/W(k)}^1&\stackrel{\res_B}{\srar} &
B/{\frak m}_L. \end{array}
$$
Recall that 
$B\otimes_A\Omega_{A/W(k)}^1\to\Omega_{B/W(k)}^1$
is an \pagebreak isomorphism if and only if $e_{L/K}=1$.

\specialnumber{2.2.6} \proclaim{Lemma}\label{tamelemma}
The canonical map
$$B\otimes_A\omega_{(A,M_A)/W(k)}^1\to\omega_{(B,M_B)/W(k)}^1$$
is an isomorphism if and only if $p$ does not divide $e_{L/K}$.
\endproclaim

\demo{Proof} Suppose that $p$ does not divide $e_{L/K}$. If
$e_{L/K}=1$ the lemma follows from the natural exact sequence 
$$0\to\Omega_{A/W(k)}^1\to\omega_{(A,M)/W(k)}^1\to A/{\frak m}_K
\to 0$$ 
and from the isomorphism mentioned before the lemma. Thus, replacing $K$
by the maximal subfield of $L$ which is unramified over $K$, we may
assume that the extension is totally ramified. Then there exists
$\pi_K\in A$ such that
$$L=K\left(\pi_K^{1/e_{L/K}}\right).$$
Indeed, if $\pi_K$ and $\pi_L$ are uniformizers of $A$ and $B$ over
$W(k)$, then $\pi_K=u\pi_L^{e_{L/K}}$, where $u\in B^\times$ is a
unit. But the sequence
$$1\to U^1_B\to B^\times\xto{r}k^\times\to 1$$
is split by the composition of the Teichm\"uller character
$\tau\:k^\times\to W(k)^\times$ and the inclusion
$W(k)^\times\hookrightarrow B^\times$. Therefore, replacing $\pi_K$ by
$\tau(r(u))^{-1}\pi_K$, we can assume that the unit $u$ lies in the
subgroup $U^1_B$ of units in $B$ which are congruent to $1$ mod
${\frak m}_L$.  But every element of $U^1_B$ has an $e_{L/K}$-th
root, so replacing $\pi_L$ by $u^{1/e_{L/K}}\pi_L$ we may assume that
$u=1$.

Let $\pi_K$ and $\pi_L$ be uniformizers of $A$ and $B$ over $W(k)$
such that $\pi_K=\pi_L^{e_{L/K}}$, and let $\phi_K(x)$ be the
minimal polynomial of $\pi_K$. Then
$$\phi_L(x)=\phi_K(x^{e_{L/K}})$$
is the minimal polynomial of $\pi_L$. The $A$-module
$\omega_{(A,M_A)/W(k)}^1$ is generated by $d\log\pi_K$ with
annihilator $(\phi_K'(\pi_K)\pi_K)$, and similarly, the $B$-module
$\omega_{(B,M_B)/W(k)}^1$ is generated by $d\log\pi_L$ with
annihilator $(\phi_L'(\pi_L)\pi_L)$. But
$$d\log\pi_K=d\log(\pi_L^{e_{L/K}})=e_{L/K}d\log\pi_L$$
and
$$\phi_L'(\pi_L)\pi_L=\phi_K'(\pi_L^{e_{L/K}})\cdot
e_{L/K}\pi_L^{e_{L/K}}=e_{L/K}\phi_K'(\pi_K)\pi_K;$$
so the claim follows since $e_{L/K}$ is a unit. It is also clear from
this argument that the map of the statement cannot be an isomorphism
if the extension $L/K$ is wildly ramified.
\enddemo

2.3.\quad In this section we show that the homotopy groups
$(\pi_*T(A|K),M)$ form a log differential graded ring. In effect, we
prove the more general statement:

\specialnumber{2.3.1}\proclaim{Proposition} \label{diflogring}
The homotopy groups $(\TR_*^n(A|K;p),M)$
form a log differential graded ring{\rm ,} if $p$ is odd or $n=1$.
\endproclaim

The homotopy groups $\TR_*^n(A|K;p)$ form a graded-commutative
differential graded ring with the differential given by Connes'
operator~(\ref{connes}), \cite[\S1]{h}. It remains to define the maps
\begin{equation}\label{alpha_n}
\alpha_n\:M\to\TR_0^n(A|K;p),\hskip6mm
d\log_n\:M\to\TR_1^n(A|K;p) \hskip.5in \speqnu{2.3.2}
\end{equation}
and to verify the relation $\alpha_n(a)d\log_na=d\alpha_n(a)$. We
define $\alpha_n$ as the composite of the inclusion $M=A\cap
K^\times\hookrightarrow A$ and the multiplicative map
$$\ul{\phantom{x}}_n\:A\to\TR_0^n(A|K;p).$$
This, we recall, is the map of components induced from the composite
$$A \xto{\;i\;} |N\s^\cy({\bf N}\s^q{\cal C})|
\begin{array}{c} \scr{D_r\circ\Delta_r}\\[-8pt]\lrar\\[-10pt]\scr{\sim}\end{array} |N\s^\cy({\bf
N}\s^q{\cal C})|^{C_r} = \TR^n(A|K;p)_{0,0},$$
where ${\cal C}=C_q^b({\cal P}_A)$ of Definition~\ref{def},
$i(a)$ is the $0$-simplex $A\xto{a}A$, and $r=p^{n-1}$. We refer
the reader to \cite[\S1]{bhm} for the definition of the maps
$\Delta_r$ and~$D_r$.

In general, if ${\cal C}$ is a category with cofibrations and weak
equivalences and if $X$ is an object of ${\cal C}$, there is a
natural map in the stable category
$$\widetilde\det\:\Sigma^\infty B\Aut(X)\to K({\cal C}),$$
where $\Aut(X)$ is the monoid of endomorphisms of $X$ in the category
$w{\cal C}$ of weak equivalences. The inclusion of $\Aut(X)$ as a
full subcategory of $w{\cal C}$ induces
$$B\Aut(X)=|N\s\Aut(X)|\to |N\s w{\cal C}|=K({\cal C})_0,$$
but this map does not preserve the basepoint (unless $X$ is the chosen
null object). However, we still get a map of symmetric spectra
$$\det\:\Sigma^\infty B\Aut(X)_+\to K({\cal C}).$$
To get the map $\widetilde\det$, we use the fact that for every pointed space $B$,
there is a natural isomorphism $S^0\vee\Sigma^\infty
B\xto{\sim}\Sigma^\infty B_+$ in the stable category. The inverse is
induced from the map which collapses $B$ to the nonbase point in
$S^0$ and the map which identifies the extra base point with the base
point in $B$.

We again let ${\cal C}=C_q^b({\cal P}_A)$ and view $A$ as a
complex concentrated in degree zero. Then $\Aut(A)=A\cap K^\times=M$
such that we have a map of monoids
$$M \to \pi_1BM \xto{\widetilde\det_*} \pi_1K({\cal C}),$$
and we define $d\log_n$ to be the composite of this map and the
cyclotomic trace. Spelling out the definition, we see that  $d\log_n$  is given by
the composite
\begin{eqnarray*}
S^{l+1}\wedge M_+ &\hskip-12pt&
\stck{\sigma\wedge\id} S^l\wedge{\Bbb T}_+\wedge M_+
\stck{\id\wedge j} S^l\wedge |N\s^{\cy}({\bf N}\s^q{\cal C})| \\
{} &\hskip-12pt& \begin{array}{c} \scr{D_r\circ\Delta_r}\\[-8pt]\lrar\\[-10pt]\scr{\sim}\end{array}
S^l\wedge |N\s^\cy({\bf N}\s^q{\cal C})|^{C_r} 
\stck{\lambda_{l,0}} \TR^n(A|K;p)_{l,0}, 
\end{eqnarray*}
where the map $j$, when restricted to ${\Bbb T}\times\{a\}$, traces
out the loop in the realization given by the $1$-simplex (in the
diagonal simplicial set):
$$\begin{array}{ccccc}
{A} &\hskip-8pt \stackrel{1}{\rsrar}\hskip-8pt&
{A} &\hskip-8pt \stackrel{1}{\rsrar}\hskip-8pt&
{A} \\[4pt]
\phantom{\scriptstyle a}\big\downarrow{\scriptstyle a} &&\phantom{\scriptstyle
a}\big\downarrow{\scriptstyle a} &&\phantom{\scriptstyle a}\big\downarrow{\scriptstyle a}
\\[4pt]
 {A} &\hskip-8pt \stackrel{1}{\rsrar}\hskip-8pt&
{A} &\hskip-8pt \stackrel{1}{\rsrar}\hskip-8pt&
{A}.\end{array}
$$

\specialnumber{2.3.3} \proclaim{Lemma}\label{da=adloga}
For all $a\in M${\rm ,}
$d\alpha_n(a)=\alpha_n(a)d\log_na$.
\endproclaim

\demo{Proof} Spelling out the definitions, one readily recognizes that
it will suffice to show that the following diagram homotopy-commutes:
$$\begin{array}{c} 
{{\Bbb T}_+\wedge M_+} \stackrel{\id\wedge\Delta}{\lrar}{{\Bbb T}_+\wedge M_+\wedge M_+}
\stackrel{j\wedge\id}{\srar}{|N\s^\cy({\bf N}\s^q{\cal C})| \wedge
|N\s^\cy({\bf N}\s^q{\cal C})|} 
\\[4pt]
  \big\downarrow{\scriptstyle \id\wedge i}  \hskip3in
\big\downarrow{\scriptstyle {\mu_{0,0}}}\qquad\
\\[4pt]
{{\Bbb T}_+\wedge |N\s^\cy({\bf N}\s^q{\cal C})|} 
\stackrel{\mu}{\vvlar}
{|N\s^\cy({\bf N}\s^q{\cal C})|.}\quad\qquad\end{array}$$ 
Since $M$ is discrete, we may check this separately for each $a\in
M$. The composite of the upper horizontal maps and the right-hand
vertical map, when restricted to ${\Bbb T}\times\{a\}$, traces out
the loop in the realization given by the $1$-simplex (in the diagonal
simplicial set) on the left below. Similarly, the composite of the
left-hand vertical map and the lower horizontal map, when restricted
to ${\Bbb T}\times\{a\}$, traces out the loop given by the
$1$-simplex on the right below:
$$
\begin{array}{ccccc}
{A} &\hskip-8pt \stackrel{a}{\rsrar}\hskip-8pt&
{A} &\hskip-8pt \stackrel{1}{\rsrar}\hskip-8pt&
{A} \\[4pt]
\phantom{\scriptstyle a}\big\downarrow{\scriptstyle a} &&\phantom{\scriptstyle
a}\big\downarrow{\scriptstyle a} &&\phantom{\scriptstyle a}\big\downarrow{\scriptstyle a}
\\[4pt]
 {A} &\hskip-8pt \stackrel{a}{\rsrar}\hskip-8pt&
{A} &\hskip-8pt \stackrel{1}{\rsrar}\hskip-8pt&
{A}  \end{array} \quad,\quad
\begin{array}{ccccc}
{A} &\hskip-8pt \stackrel{1}{\rsrar}\hskip-8pt&
{A} &\hskip-8pt \stackrel{a}{\rsrar}\hskip-8pt&
{A} \\[4pt]
\phantom{\scriptstyle 1}\big\downarrow{\scriptstyle 1} &&\phantom{\scriptstyle
1}\big\downarrow{\scriptstyle 1} &&\phantom{\scriptstyle 1}\big\downarrow{\scriptstyle 1}
\\[4pt]
 {A} &\hskip-8pt \stackrel{1}{\rsrar}\hskip-8pt&
{A} &\hskip-8pt \stackrel{a}{\rsrar}\hskip-8pt&
{A}.\end{array}
$$ 
Note that both loops are based at the vertex $A\xto{a}A$. We must show
that the two loops are homotopic through loops based at $A\xto{a}A$. To
this end, we consider the $2$-simplices 
$$\begin{array}{ccccccc}
{A} &\hskip-8pt \stackrel{a}{\rsrar}\hskip-8pt&
{A} &\hskip-8pt \stackrel{1}{\rsrar}\hskip-8pt&
{A} &\hskip-8pt \stackrel{1}{\rsrar}\hskip-8pt&
{A} \\[4pt]
\phantom{\scriptstyle 1}\big\downarrow{\scriptstyle 1} &&\phantom{\scriptstyle
1}\big\downarrow{\scriptstyle 1} &&\phantom{\scriptstyle 1}\big\downarrow{\scriptstyle
1}&&\phantom{\scriptstyle 1}\big\downarrow{\scriptstyle 1}
\\[4pt]
{A} &\hskip-8pt \stackrel{a}{\rsrar}\hskip-8pt&
{A} &\hskip-8pt \stackrel{1}{\rsrar}\hskip-8pt&
{A} &\hskip-8pt \stackrel{1}{\rsrar}\hskip-8pt&
{A} \\[4pt]
\phantom{\scriptstyle a}\big\downarrow{\scriptstyle a} &&\phantom{\scriptstyle
1}\big\downarrow{\scriptstyle 1} &&\phantom{\scriptstyle a}\big\downarrow{\scriptstyle
a}&&\phantom{\scriptstyle a}\big\downarrow{\scriptstyle a}
\\[4pt]
 {A} &\hskip-8pt \stackrel{1}{\rsrar}\hskip-8pt&
{A} &\hskip-8pt \stackrel{a}{\rsrar}\hskip-8pt&A&\hskip-8pt \stackrel{1}{\rsrar}\hskip-8pt&
{A}  \end{array}\quad ,\quad
\begin{array}{ccccccc}
{A} &\hskip-8pt \stackrel{1}{\rsrar}\hskip-8pt&
{A} &\hskip-8pt \stackrel{1}{\rsrar}\hskip-8pt&
{A} &\hskip-8pt \stackrel{a}{\rsrar}\hskip-8pt&
{A} \\[4pt]
\phantom{\scriptstyle 1}\big\downarrow{\scriptstyle 1} &&\phantom{\scriptstyle
1}\big\downarrow{\scriptstyle 1} &&\phantom{\scriptstyle 1a}\big\downarrow{\scriptstyle
a}&&\phantom{\scriptstyle 1}\big\downarrow{\scriptstyle 1}
\\[4pt]
{A} &\hskip-8pt \stackrel{1}{\rsrar}\hskip-8pt&
{A} &\hskip-8pt \stackrel{a}{\rsrar}\hskip-8pt&
{A} &\hskip-8pt \stackrel{1}{\rsrar}\hskip-8pt&
{A} \\[4pt]
\phantom{\scriptstyle 1}\big\downarrow{\scriptstyle 1} &&\phantom{\scriptstyle
1}\big\downarrow{\scriptstyle 1} &&\phantom{\scriptstyle 1}\big\downarrow{\scriptstyle
1}&&\phantom{\scriptstyle 1}\big\downarrow{\scriptstyle 1}
\\[4pt]
 {A} &\hskip-8pt \stackrel{1}{\rsrar}\hskip-8pt&
{A} &\hskip-8pt \stackrel{a}{\rsrar}\hskip-8pt&A&\hskip-8pt \stackrel{1}{\rsrar}\hskip-8pt&
\phantom{.}{A} . \end{array}
$$  
The $2$-simplex on the left gives a homotopy through loops based at
$A\xto{a}A$ between the loop given by the left-hand $1$-simplex above
and the loop given by the $1$-simplex
$$ 
\begin{array}{ccccc}
{A} &\hskip-8pt \stackrel{a}{\rsrar}\hskip-8pt&
{A} &\hskip-8pt \stackrel{1}{\rsrar}\hskip-8pt&
{A} \\[4pt]
\phantom{\scriptstyle a}\big\downarrow{\scriptstyle a} &&\phantom{\scriptstyle
1}\big\downarrow{\scriptstyle 1} &&\phantom{\scriptstyle a}\big\downarrow{\scriptstyle a}
\\[4pt]
 {A} &\hskip-8pt \stackrel{1}{\rsrar}\hskip-8pt&
{A} &\hskip-8pt \stackrel{a}{\rsrar}\hskip-8pt&
\phantom{.}{A} .  \end{array} 
$$
Similarly, the $2$-simplex on the right gives a homotopy through loops
based at $A\xto{a}A$ between this loop and the loop given by the right-hand $1$-simplex above.
\enddemo 

\specialnumber{2.3.4}\proclaim{Proposition} \label{pi1}
The canonical map
$$\omega_{(A,M)}^q\to\pi_qT(A|K)$$
is an isomorphism{\rm ,} for $q\leq 2${\rm ,} and a rational isomorphism{\rm ,} for all
$q\geq 0$.
\endproclaim

\demo{Proof} We consider the long exact sequence of homotopy groups
associated with the sequence of Theorem~\ref{localizationtr},
$$T(k) \xto{i^!} T(A) \xto{j_*} T(A|K) \xto{\partial} \Sigma T(k),$$
and note that $i^!\:\pi_qT(k)\to\pi_qT(A)$ is zero, if
$q=0,1$. Indeed, for $q=0$ this is a map from a torsion group to a
torsion-free group, and for $q=1$ the domain is isomorphic to the
group $\Omega_k^1$ which vanishes since $k$ is a perfect, \cite[Lemma 
5.5]{hm}. This proves the statement for $q=0$. It also shows that the
top sequence in the following diagram of $A$-modules and $A$-linear
maps,
$$\begin{array}{ccccccccc}
{0} &\hskip-8pt \srar \hskip-8pt &
{\pi_1T(A)} &\hskip-8pt \stackrel{j_*}{\srar}\hskip-8pt &
{\pi_1T(A|K)}&\hskip-8pt \stackrel{\partial_*}{\srar} \hskip-8pt &
{\pi_0T(k)} &\hskip-8pt \srar \hskip-8pt &
{0} \\[4pt]
&&\big\uparrow&&\big\uparrow&&\big\uparrow\\[3pt]
0&\hskip-8pt \srar\hskip-8pt &{\Omega_A^1} &\hskip-8pt \srar\hskip-8pt &
{\omega_{(A,M)}^1}&\hskip-8pt \stackrel{\res}{\srar}\hskip-8pt &
{k} &\hskip-8pt \srar\hskip-8pt &
{0,} \end{array}
$$
is exact. The lower sequence is the exact sequence of
Proposition~\ref{poincare} and the vertical maps are the canonical
maps. The left-hand square commutes since $j_*$ preserves the
differential. The commutativity of the right-hand square is equivalent
to the statement that $\partial_*(d\log x)=v_K(x)$, for all $x\in
M$. But this follows from the definition of the map $d\log$ in
(\ref{alpha_n}) and from the commutativity of the right-hand square in
Addendum~\ref{localizationtc}. Since the left- and right-hand vertical
maps in the diagram are isomorphisms, so is the middle vertical map.
This proves the statement for $q=1$.  

We next argue that the map of the statement is a rational isomorphism,
for all $q\geq 0$. Since $\pi_*T(k)$ is torsion the long exact
sequence associated with the cofibration sequence above shows that
$$j_*\:\pi_*T(A)\otimes\Q \xto{\sim} \pi_*T(A|K)\otimes\Q$$
is an isomorphism. Moreover, the linearization map induces an
isomorphism
$$l\:\pi_*T(A)\otimes\Q \xto{\sim} \HH_*(A)\otimes\Q,$$
and the right-hand side is canonically isomorphic to
$\HH_*(K)$. It thus remains to prove that the canonical map
$\Omega_K^*\to\HH_*(K)$ is an isomorphism. This in turn follows from
\cite{hkr} and the fact that $K$ can be written as a filtered
colimit of smooth $\Q$-algebras,~\cite[IV.17.5.1]{ega4}.

It remains to show that $\pi_2T(A|K)$ is uniquely divisible. 
The structure of the $p$-adic homotopy groups $\pi_*(T(A),\Zp)$ is
known from \cite[Th.\ 5.1]{lm}. (The assumption that the
residue field is finite is not needed. For {\it op.~cit.},
Propositions 5.3 and 5.4 and \cite[Th.\ 7.1]{boardman} show  that
the Bockstein spectral sequence converges strongly.) The result is
that for $m>0$, $\pi_{2m}(T(A),\Zp)$ vanishes and
$\pi_{2m-1}(T(A),\Zp)$ is isomorphic to $A/(m\phi_K'(\pi_K))$. The
latter is a torsion group of bounded exponent. It follows that for
$m>0$, $\pi_{2m}T(A)$ is a uniquely divisible group and
$\pi_{2m-1}T(A)$ is the sum of a uniquely divisible group and the
torsion group $\pi_{2m-1}(T(A),\Zp)$. Since $\pi_1T(k)$ is trivial, we
see that $\pi_2T(A|K)$ is uniquely divisible as stated.
\enddemo

2.4.\quad It follows from Proposition~\ref{diflogring} that the
homotopy groups with $\Z/p$-coefficients $\bar\pi_*T(A|K)$ form a log
differential graded $k$-algebra. We now evaluate this log differential
graded $k$-algebra and prove Theorem~\ref{thmB} of the introduction.

The proof of Theorem~\ref{thmB} is based on the calculation in
\cite[Ths.\ 4.4, 4.6]{lm} of the graded $k$-algebra
$\bar\pi_*T(A)=\pi_*(T(A),\Z/p)$. The result, which we now recall,
depends on whether $p$ divides $e_K$ or not. We consider the graded
$k$-algebra
$$B=A/p\otimes\Lambda\{\alpha_1\}\otimes S\{\alpha_2\}$$
with the generators in the indicated degrees. Let $C\subset B$ be the
subalgebra generated by all elements $a\alpha_1^{\e}\alpha_2^m$ for
which $a\in{\frak m}_K/pA$ or $\e=1$ or $p$ divides $m$, and let
$I\subset C$ be the ideal generated by all elements
$a\alpha_1\alpha_2^{m-1}$ for which $a\in{\frak m}_K^{e_K-1}/pA$
and $m$ is prime to $p$. Then as graded $k$-algebras,
$$\bar\pi_*T(A)\cong\left\{\begin{array}{ll}
B, & \hbox{if $p$ divides $e_K$,} \\
C/I, & \hbox{if $p$ does not divide $e_K$.} \end{array}
\right.$$
We note that, in the former case, the dimension of the $k$-vector
space $\bar\pi_qT(A)$ is equal to $e_K$, for all $q\geq 0$. In the 
latter case, this dimension is equal to $e_K$, if $q$ is congruent to
either $-1$ or $0$ modulo $2p$, and is equal to $e_K-1$, otherwise.

We also recall from \cite[Th.\ 5.2, Cor.\ 5.5]{hm} that as a
graded $k$-algebra,
$$\bar\pi_*T(k)=\Lambda\{\e\}\otimes S\{\sigma\},$$
with the generators $\varepsilon$ and $\sigma$ characterized as follows:
the Bockstein takes $\varepsilon$ to~$1$ and Connes'
operator~\eqref{connes} takes $\e$ to $\sigma$. It follows from the
proof of \cite[Ths.\ 4.4, 4.6]{lm} that the reduction map
$$i_*\:\bar\pi_*T(A) \to \bar\pi_*T(k)$$
is induced from a $k$-algebra map $B \to \bar\pi_*T(k)$, which in
degree zero is given by the reduction $A/p\to k$, and which takes the
generators $\alpha_1$ and $\alpha_2$ to zero and a unit times
$\sigma$, respectively.

Since the group $\pi_2T(A|K)$ is uniquely divisible, by
Proposition~\ref{pi1}, the integral Bockstein induces an isomorphism
$$\beta\:\bar\pi_2T(A|K)\xto{\sim}{}_p\pi_1T(A|K).$$
We define $\kappa\in\bar\pi_2T(A|K)$ to be the class which corresponds
to the generator $d\log(-p)$ on the right. (Note that
$\kappa\in\bar\pi_2T(\Zp|\Qp)$.) We now prove Theorem~\ref{thmB} of
the introduction:

\specialnumber{2.4.1} \proclaim{Theorem} \label{pisubast}
There is a natural isomorphism of
log differential graded rings
$$\omega_{(A,M)}^*\otimes_{\Z}S_{\Z/p}\{\kappa\}
\xto{\sim} \bar\pi_*T(A|K),$$
where $d\kappa=\kappa d\log(-p)$.
\endproclaim

\demo{Proof} It is clear that there is a map of graded $k$-algebras as
stated. We show that this is an isomorphism.

Suppose first that $p$ divides $e_K$. We know from
Proposition~\ref{pi1} that the map of the statement is an isomorphism
in degrees $q\leq 1$. So it suffices to show that multiplication by
$\kappa$ induces an isomorphism
$$\kappa\:\bar\pi_qT(A|K) \xto{\sim} \bar\pi_{q+2}T(A|K).$$
To this end, we consider the long-exact sequence associated with the
cofibration sequence of Theorem~\ref{localizationtr},
$$\cdots\to
\bar\pi_qT(k) \xto{i^!}
\bar\pi_qT(A) \xto{j_*}
\bar\pi_qT(A|K) \xto{\partial}
\bar\pi_{q-1}T(k) \to \cdots.$$
This is a sequence of graded $\bar\pi_*T(A)$-modules, where
$\bar\pi_*T(A|K)$ (resp.\ $\bar\pi_*T(k)$) is viewed as a graded
$\bar\pi_*T(A)$-module via the map $j_*$ (resp.\ $i_*$). We claim that
the map $j_*$ is an isomorphism for $q=2$. Granting this for the
moment, there exists $\tilde\kappa\in\bar\pi_2T(A)$ such that
$\kappa=j_*(\tilde\kappa)$. And since $\bar\pi_2T(A)$ and
$\bar\pi_2T(A|K)$ are both free $A/p$-modules of rank one, the
class $\tilde\kappa$ is necessarily a generator. It follows that in
the diagram  
$$\begin{array}{ccccccccc}
{\cdots} \longrightarrow \hsm& 
{\bar\pi_qT(k)}&\hsm\shtck{i^!} \hsm &  
{\bar\pi_qT(A)} &\hsm\shtck{j_*} \hsm & 
{\bar\pi_qT(A|K)} &\hsm\shtck{\partial} \hsm & 
{\bar\pi_{q-1}T(k)}  & \hsm
\longrightarrow{\cdots} \\[4pt]
&{\scriptstyle\sim}\big\downarrow\scr{\tilde\kappa}&\hsm\hsm&
{\scriptstyle\sim}\big\downarrow\scr{\tilde\kappa}&\hsm\hsm&
\big\downarrow\scr{\tilde\kappa}&\hsm\hsm&{\scriptstyle
\sim}\big\downarrow\scr{\tilde\kappa}\\[4pt] {\cdots}\longrightarrow&
{\bar\pi_{q+2}T(k)} &\hsm\shtck{i^!} \hsm &
{\bar\pi_{q+2}T(A)} &\hsm\shtck{j_*} \hsm &
{\bar\pi_{q+2}T(A|K)}&\hsm\shtck{\partial}\hsm  &
{\bar\pi_{q+1}T(k)} &\hsm
\longrightarrow{\cdots,} \end{array}
$$
two out of three of the vertical maps are isomorphisms. Hence, so are
the remaining vertical maps. To prove the claim, we consider the
diagram of\break $A/p$-modules
$$\begin{array}{ccccc}
{\bar\pi_2T(A)} &\begin{array}{c}\scr{\beta}\\[-10pt]\srar\\[-8pt]\scr{\sim}\end{array}&
{{}_p\pi_1T(A)} &\stackrel{\sim}{\slar}&
{{}_p\Omega_A^1} \\[3pt]
\big\downarrow{\scriptstyle j_*} &&\big\downarrow{\scriptstyle j_*}
&&\big\downarrow{\scriptstyle j_*}\\[4pt]
 {\bar\pi_2T(A|K)} &\begin{array}{c}\scr{\beta}\\[-10pt]\srar\\[-8pt] \scr{\sim}\end{array} &
{{}_p\pi_1T(A|K)} & \stackrel{\sim}{\slar}&
{{}_p\omega_{(A,M)}^1}.\end{array}
$$
The left-hand horizontal maps are isomorphisms since $\pi_2T(A)$
and $\pi_2T(A|K)$ are (uniquely) divisible. It follows from
Proposition~\ref{pi1} that the right-hand horizontal maps are 
isomorphisms and that the right-hand vertical map is a monomorphism.
But the domain and range of the latter are both $k$-vector spaces of
dimension $e_K$. Hence, this map is an isomorphism. This proves the
claim.

Suppose now that $p$ does not divide $e_K$. Let $L/K$ be a totally
ramified extension such that $p$ divides $e_{L/K}$, and let $B$ be the
integral closure of $A$ in~$L$. Then we have a commutative diagram
$$\begin{array}{ccc}
{\omega_{(A,M_A)}^*\otimes S\{\kappa\}}&\srar&
{\bar\pi_*T(A|K)}\\[4pt]
\big\downarrow&&\big\downarrow\\[4pt]
{\omega_{(B,M_B)}^*\otimes S\{\kappa\}}&\stackrel{\sim}{\srar} &
{\bar\pi_*T(B|\,L),}\end{array}
$$
and the lower horizontal map is an isomorphism. It is easy to see
that there exists $L/K$ for which the left-hand vertical map is a
monomorphism. For example, one can take
$L=K[\pi_L]/(\pi_L^{e_{L/K}}+\pi_K(\pi_L+1))$. Hence, the upper
horizontal map is a monomorphism. The domain and range of this map are
graded $k$-vector spaces concentrated in nonnegative degrees. The
dimension of the domain is equal to $e_K$ in each degree. Hence the
dimension of the range is at least $e_K$ in each degree. We can
estimate the dimension of the range further by means of the exact
sequence of $k$-vector spaces
$$\cdots\to
\bar\pi_qT(k) \xto{i^!}
\bar\pi_qT(A) \xto{j_*}
\bar\pi_qT(A|K) \xto{\partial}
\bar\pi_{q-1}T(k) \to \cdots\; .$$
The dimension of $\bar\pi_qT(A)$ is equal to $e_K$, if 
$q\equiv -1,0$ (mod $2p$), and is equal to $e_K-1$, otherwise. The
dimension of $\bar\pi_qT(k)$ is equal to one, for all $q\geq 0$. It
follows that the dimension of $\bar\pi_qT(A|K)$ is equal to either
$e_K$ or $e_K+1$, if $q\equiv -1,0$ (mod $2p$), and is equal to $e_K$
otherwise. We  argue that for $q\equiv -1,0$, the dimension
of $\bar\pi_qT(A|K)$ is equal to $e_K$. This happens if and only if
for all $s\geq 0$, the map
$$i^{!}\:\bar\pi_{2ps-1}T(k)\to\bar\pi_{2ps-1}T(A)$$
is nonzero. We show that the class
$\pi_K^{e_K-1}\alpha_1\alpha_2^{ps-1}$ on the right is in the image of
$i^!$, or equivalently, that it maps to zero under $j_*$. If $e_K>1$,
we can write
$$\pi_K^{e_K-1}\alpha_1\alpha_2^{ps-1}=
\pi_K^{e_K-2}\alpha_1\cdot\pi_K\alpha_2^{ps-1}.$$
The image of this class under $j_*$ is equal to a unit in $A/p$ times
the class
$$\pi_K^{e_K-1}d\log\pi_K\cdot\pi_K\kappa^{ps-1}.$$
But this class is in the image of the ring homomorphism
$$\omega_{(A,M)}^*\otimes_{\Z}S_{\Z/p}\{\kappa\}
\to\bar\pi_*T(A|K)$$
and the product is equal to zero on the left. Hence
$j_*(\pi_K^{e_K-1}\alpha_1\alpha_2^{ps-1})$ is equal to zero. Finally,
in the unramified case we choose a totally ramified extension $K/K_0$
such that $p$ does not divide $e_K$ and consider the diagram 
$$\begin{array}{ccc}
{\bar\pi_{2ps-1}T(k)} &\stackrel{i^!}{\srar}& 
{\bar\pi_{2ps-1}T(W(k))} \\[3pt]
\big\downarrow{\scriptstyle e_K} &&\big\downarrow \\[4pt]
{\bar\pi_{2ps-1}T(k)} &\stackrel{i^!}{\srar} &
{\bar\pi_{2ps-1}T(A).} \end{array}
 $$
We have just proved that the lower horizontal map is a monomorphism,
for all $s\geq 0$. And the left-hand vertical map is an isomorphism
since $e_K$ is prime to $p$. Hence the top horizontal map is a
monomorphism. We have proved that the map of the statement is an
isomorphism of graded $k$-algebras for all $K$. In particular, the
class $d\kappa$ on the right is the image of an element on the
left. To determine this element, we may assume that $K=\Qp$. In the
diagram
$$\begin{array}{ccc}
{\bar\pi_3T(\Zp|\Qp)} &\stackrel{\partial}{\srar} &
{\bar\pi_2T(\Fp)} \\[4pt]
\big\uparrow{\scriptstyle d}&&\big\uparrow{\scriptstyle d}\\[4pt]
{\bar\pi_2T(\Zp|\Qp)}&\stackrel{\partial}{\srar} &
{\bar\pi_1T(\Fp)}
\end{array}
$$
the horizontal maps and the right-hand vertical map are
isomorphisms. Hence also the left-hand vertical map is an
isomorphism. This shows that $d\kappa=u\kappa d\log(-p)$ with
$u\in\Fp^\times$. We show in remark~\ref{unitisone} below that in fact
$u=1$.
\enddemo

\specialnumber{2.4.2}\numbereddemo{{R}emark} \label{concodd}
An argument similar to \cite[\S5]{lm}
shows that for $m>0$, there exists a noncanonical isomorphism
$$\pi_{2m-1}(T(A|K),\Zp)\cong A/(m\pi_K\phi_K'(\pi_K))$$
and that $\pi_{2m}(T(A|K),\Zp)$ vanishes. It would be interesting to
give a functorial description of the left-hand group analogous to
Proposition~\ref{pi1}.
\enddemo

Let $L/K$ be a Galois extension with Galois group $G_{L/K}$. The
descent problem for topological Hochschild homology asks under what
conditions the canonical map
$$T(A|K)\to\coborel(G_{L/K},T(B|\,L))$$
is a weak equivalence. It is not hard to see from
Theorem~\ref{pisubast} that this is false in general, e.g. for a
cyclotomic extension $\Qp(\mu_{p^n})/\Qp$ with $n>1$. However:

\specialnumber{2.4.3}\proclaim{Theorem} \label{tamedescent}
Let $L/K$ be a finite and tamely
ramified Galois extension. Then the canonical map induces an
isomorphism
$$\bar\pi_*T(A|K) \xto{\sim}
\bar\pi_*\coborel(G_{L/K},T(B|\,L)).$$
\endproclaim

\demo{Proof} It will suffice to show that for all $t\geq 0$, the
$G_{L/K}$-module $\bar\pi_tT(B|\,L)$ is isomorphic to $B/p$. Indeed,
a classical theorem of Noether,~\cite[I.3, Th.\ 3]{frohlich},
states that $B$ is isomorphic to $A[G_{L/K}]$ as a $G_{L/K}$-module,
if and only if $L/K$ is tamely ramified. Hence, the spectral sequence
$$E^2_{s,t}=H^{-s}(G_{L/K},\bar\pi_tT(B|\,L))
\Rightarrow\bar\pi_{s+t}\coborel(G_{L/K},T(B|\,L))$$
collapses to yield the isomorphism of the statement.

We use Theorem~\ref{pisubast} to get the natural isomorphisms
$$\kappa^i\:\bar\pi_{\e}T(B|\,L) \xto{\sim}
\bar\pi_{2i+\e}T(B|\,L).$$
Hence, we only need to consider $\bar\pi_0T(B|\,L)$ and
$\bar\pi_1T(B|\,L)$. The former is naturally isomorphic to $B/p$
regardless of whether $L/K$ is tamely ramified or not, and the latter
is naturally isomorphic to $\omega_{(B,M_B)}^1/p$. We have from
Lemma~\ref{logdiff} that
$$\omega_{(A,M_A)}^1/p=A/p\cdot d\log\pi_K,$$
and since $L/K$ is tamely ramified, Lemma~\ref{tamelemma} shows that
$$\omega_{(B,M_B)}^1/p=B/p\cdot d\log\pi_K.$$
Hence, also $\omega_{(B,M_B)}^1/p$ is is isomorphic to $B/p$ as a
$G_{L/K}$-module.
\enddemo

\section{The de Rham-Witt complex and $\TR_*\cs(A|K;p)$} \label{pione}
3.1.\quad In this paragraph, we evaluate the integral homotopy
groups $\TR_i\cs(A|K;p)$, for $i\leq 2$. We first consider Witt
vectors, see e.g.~\cite[Appendix]{mumford1}.

The ring $W_n(R)$ of Witt vectors of length $n$ in $R$ is the set of
$n$-tuples in $R$ but with a new ring
structure characterized by the requirement that the ``ghost'' map
$$w\:W_n(R)\to R^n,$$ 
which to the vector $(a_0,,\dots,a_{n-1})$ associates the sequence
$(w_0,\dots,w_{n-1})$ with
$$w_s=a_0^{p^s}+pa_1^{p^{s-1}}+\cdots+p^sa_s,$$
be a natural transformation of functors from rings to rings. If $R$
has no $p$-torsion then the ghost map is injective. If, in addition,
there exists a ring endomorphism $\phi\:R\to R$ such that
$a^p\equiv\phi(a)$ (mod $pR$), then a sequence $(w_0,\dots,w_{n-1})$
is in the image if and only if
$w_s\equiv\phi(w_{s-1})$ (mod $p^sR$), for all $0<s<n$. If
$R=\Z[X_\alpha]$, the ring homomorphism which maps $X_\alpha$ to
$X_\alpha^p$ is such an endomorphism. Let
$$\ul{\phantom{x}}_n\:R\to W_n(R)$$
be the multiplicative section given by $\ul{a}_n=(a,0,\dots,0)$.

\specialnumber{3.1.1} \proclaim{Lemma}\label{p>2}If $p$ is odd then $V(1)\equiv\ul{-p}_n$ 
and $\ul{-1}_n\equiv -1$ modulo $pW_n(R)$.
\endproclaim

\demo{Proof} By naturality, we may assume that $R=\Z$. Now,
$$w(\ul{p}_n+V(\ul1))=p(1,1+p^{p-1},1+p^{p^2-1},1+p^{p^3-1},\dots),$$
and therefore it is enough to show that the sequence
$$(1,1+p^{p-1},1+p^{p^2-1},\dots,1+p^{p^{n-1}-1})$$
is in the image of the ghost map. The identity $\phi\:\Z\to\Z$ has the
property that $a^p\equiv\phi(a)$ (mod $p\Z$). Hence, this sequence
is in the image of the ghost map if and only if for all $1<s<n$,
$$1+p^{p^s-1}\equiv 1+p^{p^{s-1}-1}\hskip6mm\hbox{(mod $p^s$)}.$$
This is true, if $p$ is odd, but fails for $p=2$ and $s=2$. The second
congruence of the statement is proved in a similar manner.
\enddemo

In general, $\ul{(x+y)}_n$ and $\ul{x}_n+\ul{y}_n$ are not equivalent
modulo $pW_n(A)$. However, we have the following:

\specialnumber{3.1.2} \proclaim{Lemma}\label{ppower}
For all $x,y\in R${\rm ,}
$$\ul{(x+y)}_n^p\equiv(\ul{x}_n+\ul{y}_n)^p
\equiv\ul{x}_n^p+\ul{y}_n^p$$
modulo $pW_n(R)$.
\endproclaim

\demo{Proof} The right-hand congruence is valid in any ring. To prove
the left-hand congruence, we place ourselves in the universal case
$R=\Z[x,y]$. The ghost map 
$$w\:W_n(R)\to R^n$$
is an injection and maps the vector
$\ul{x}_n^p+\ul{y}_n^p-(\ul{x+y}_n)^p$ to the sequence
$$(x^p+y^p-(x+y)^p,\dots,x^{p^n}+y^{p^n}-(x+y)^{p^n}).$$
As an element of $R^n$ this is divisible by $p$. We must show
that the quotient is in the image of the ghost map. By the criterion
recalled above, we must show that
$$(x^{p^{s+1}}+y^{p^{s+1}}-(x+y)^{p^{s+1}})/p\equiv
(x^{p^{s+1}}+y^{p^{s+1}}-(x^p+y^p)^{p^s})/p\hskip4mm\hbox{(mod
$p^s$)},$$
or equivalently, that
$$(x+y)^{p^{s+1}}\equiv(x^p+y^p)^{p^s}\hskip4mm\hbox{(mod
$p^{s+1}$)}.$$
But this follows from
$$(x+y)^p\equiv x^p+y^p\hskip4mm\hbox{(mod $p$)}$$
and from the fact, valid in any commutative ring, that $a\equiv b$
(mod $p$) implies $a^{p^s}\equiv b^{p^s}$ (mod $p^{s+1}$). Indeed, one
easily sees that $a\equiv b$ (mod $p^k$) implies that $a^p\equiv b^p$
(mod $p^{k+1}$), and the desired formula then follows by simple induction.
\enddemo

Now, from Lemma~\ref{ppower}, for every ring $R$, the map
$$R\to\bar{W}_n(R)=W_n(R)/p,$$
which takes $x$ to the class of $\ul{x}_n^p$, is a ring
homomorphism. Let $A$ be a complete discrete valuation ring with
quotient field $K$ and perfect residue field $k$ of mixed
characteristic $(0,p)$. We recall from~\eqref{structure} that there is
a unique ring homomorphism $f\:W(k)\to A$ such that the induced map of
residue fields is the identity homomorphism. Hence, we have a ring
homomorphism
\begin{equation}\label{rho}
\rho_n\: k \to \bar{W}_n(A) \speqnu{3.1.3}
\end{equation}
which to $x$ assigns
$\ul{f(\smash{\widetilde{x^{1/p}}})\,}_n^p+pW_n(A)$. Here 
$\widetilde{x^{1/p}}\in W(k)$ is any element whose residue class
modulo $p$ is the unique $p$-th root of $x$. We will always view
$\bar{W}_n(A)$ as a $k$-algebra via the map $\rho_n$. We note that
$$R(\rho_n(x))=\rho_{n-1}(x),\hskip5mm
F(\rho_n(x))=\rho_{n-1}(x^p).$$
Let $\pi=\pi_K$ be a uniformizer with minimal polynomial
$x^{e_K}+p\theta_K(x)$. We introduce the modified Verschiebung
\begin{equation}\label{modified}
V_\pi\:\bar{W}_{n-1}(A)\to \bar{W}_n(A),\hskip4mm
V_\pi(a)=\theta_K(\ul{\pi}_n)V(a), \speqnu{3.1.4}
\end{equation}
where $\theta_K(\ul{\pi}_n)$ is the image of 
$\theta_K(x)$ under the $k$-algebra map $k[x]\to\bar{W}_n(A)$
which to $x$ assigns the class of $\ul{\pi}_n$. The composite
$FV_{\pi}$ is zero modulo $p$.

\specialnumber{3.1.5}\proclaim{Proposition} \label{wittvectors}
Suppose that $p$ is odd. Then the
$k$\/{\rm -}\/algebra $\bar{W}_n(A)$ is generated by the elements
$V_\pi^s(\ul\pi^i)$ with $0\leq s<n$ and $i\geq 0$ subject to the
relations
\begin{eqnarray*}
V_\pi^s(\ul\pi^i)\cdot V_\pi^t(\ul\pi^j)&=&\left\{\begin{array}{ll}
V_\pi^t(\ul\pi^{p^ti+j}), & \hbox{if $0=s\leq t<n$,} \\
0, & \hbox{if $0<s\leq t<n$,} \\
\end{array}\right.\\
V_\pi^s(\ul\pi^{e_K+i})&=&V_\pi^{s+1}(\ul\pi^{pi}). 
\end{eqnarray*}
\endproclaim

\demo{Proof} The $k$-vector space $\bar{W}_n(A)$ is generated by
$V^s(\ul\pi^i)$ with $0\leq s<n$ and $i\geq 0$. Indeed, write
$a\in A$ as $a=x_d\pi^d+\cdots+x_0$ with $x_i\in W(k)$. Then
$$V^s(\ul a)\equiv V^s(\ul{x_d\pi^d})+\cdots+V^s(\ul{x_0})
\equiv V^s(\rho_{n-s}(\bar{x}_d)\ul{\pi}^d)+\cdots+V^s(\rho_{n-s}(\bar{x}_0))$$
modulo $V^{s+1}\bar{W}_n(A)$, and
$$V^s(\rho_{n-s}(\bar{x}_i)\ul\pi^i)=\rho_n(\bar{x}^{1/p^s})V^s(\ul\pi^i).$$
Since $\theta_K(\ul\pi)$ is a unit, we instead can use the elements
$V_\pi^s(\ul\pi^i)$ as generators. In general, for $s\leq t$, 
$$V_{\pi}^s(\ul{\pi}^i)\cdot V_{\pi}^t(\ul{\pi}^j)
=V_{\pi}^t(F^tV_{\pi}^s(\ul{\pi}^i)\cdot\ul{\pi}^j)$$
from which the first relation follows. Next, Lemmas~\ref{p>2} and
\ref{ppower} show that
\begin{eqnarray*}
\ul\pi^{e_K}
{} & =& \ul{-p}\cdot\ul{\theta_K(\pi)}
\equiv V(1)\ul{\theta_K(\pi)}
= V((\ul{\theta_K(\pi)})^p) \\
{} & \equiv &V(\theta_K^{(1)}(\ul\pi^p))
= V(1)\theta_K(\ul\pi)
= V_\pi(1), 
\end{eqnarray*}
where $\theta_K^{(1)}(x)$ denotes the image of $\theta_K(x)$ under the
automorphism of $k[x]$ induced by the Frobenius of $k$. The second
relation is an immediate consequence. It remains to prove that there
are no further relations. The sequences
$$0\to A/p\stck{V^{n-1}}\bar{W}_n(A)\xto{R}\bar{W}_{n-1}(A)\to 0$$
are exact, since $W_n(A)$ is torsion-free, and show that
$\bar{W}_n(A)$ is an $ne_K$-dimensional $k$-vector space. 
The relations of the statement imply that
$$\gr_V^s\bar{W}_n(A)=k\big\{V_\pi^s(\ul\pi^i)\;\big|\;0\leq
i<e_K\big\},$$
which is an $e_K$-dimensional $k$-vector space. Thus there can be no
further relations among the $V_\pi^s(\ul\pi^i)$.
\enddemo

3.2.\quad A pre-log structure $\alpha\:M\to R$ on a ring $R$
induces one on $W_n(R)$ upon composition with the multiplicative
section $\ul{\phantom{a}}_n\:R\to W_n(R)$. We write $(W_n(R),M)$ for
this log ring. We now assume that $p$ is odd and that $R$ is a
$\Z_{(p)}$-algebra.

\specialnumber{3.2.1} \numbereddemo{Definition} \label{logwittcx}A log Witt complex over $(R,M)$
consists of:
\medbreak
(i) a pro-log differential graded ring $(E\s^*,M_E)$ together with a
map of pro-log rings $\lambda\:(W\s(R),M)\to (E\s^0,M_E)$;
\medbreak
(ii) a map of pro-log graded rings
$$F\:E_n^*\to E_{n-1}^*,$$
such that $\lambda F=F\lambda$ and such that
\begin{eqnarray*}
Fd\log_na&=&d\log_{n-1}a,\hskip6.6mm\hbox{for all $a\in M$,}\\
Fd\ul{a}_n&=&\ul{a}_{n-1}^{p-1}d\ul{a}_{n-1},\hskip5mm\hbox{for all
$a\in R$;} 
\end{eqnarray*}
\medbreak
(iii) a map of pro-graded modules over the pro-graded ring $E^*\s$,
$$V\:F^*E_n^*\to E_{n+1}^*,$$
such that $\lambda V=V\lambda$, $FV=p$ and $FdV=d$.

A map of log Witt complexes over $(R,M)$ is a map of pro-log
differential graded rings which commutes with the maps $\lambda$, $F$
and $V$.
\enddemo

The following relations are valid in any log Witt complex:
$$dF=pFd,\hskip4mm Vd=pdV,\hskip4mm V(xdy)=V(x)dV(y).$$
Indeed, $V(xdy)=V(xFdV(y))=V(x)dV(y)$, and
\begin{eqnarray*}
dF(x)&=&FdVF(x)=Fd(V(1)x)=FdV(1)F(x)+FV(1)F(dx)\\
{}&=&d(1)F(x)+pFd(x)=pFd(x),\\
Vd(x)&=&V(1)dV(x)=d(V(1)V(x))-dV(1)V(x)\\
{}&=&dV(xFV(1)-V(xd(1)))=pdV(x). 
\end{eqnarray*}

\specialnumber{3.2.2}\proclaim{Proposition} \label{derhamwitt}
The category of log Witt complexes over
$(R,M)$ has an initial object $W\s\,\omega_{(R,M)}^*$. Moreover{\rm ,} the
canonical map is surjective\/{\rm :}
$$\lambda\:\omega_{(W\s(R),M)}^*\twoheadrightarrow
W\s\,\omega_{(R,M)}^*.$$
\endproclaim

\demo{Proof} This is a fairly straightforward application of the Freyd
adjoint functor theorem,~\cite[p.\ 116]{maclane}. For a detailed proof,
we refer the reader to\break \cite[\S1]{hm3}.
\enddemo

We note that $W\s\,\omega_{(R,M)}^0=W\s(R)$. For we may consider
$(W\s(R),M)$ a log Witt complex concentrated in degree
zero. Moreover, from~\cite[Th.\ D]{hm3} we have:

\specialnumber{3.2.3}\proclaim{Addendum} The canonical map is an isomorphism\/{\rm :}\/
$$\lambda\:\omega_{(R,M)}^* \xto{\sim}
W_1\,\omega_{(R,M)}^*.
$$
\endproclaim

The filtration of a log Witt complex by the differential graded ideals
$$\Fil^sE_n^i=V^sE_{n-s}^i+dV^sE_{n-s}^{i-1}\subset E_n^i$$
is called the standard filtration. It satisfies
\begin{eqnarray*}
F(\Fil^sE_n^i)&\subset&\Fil^{s-1}E_{n-1}^i,\\
V(\Fil^sE_n^i)&\subset&\Fil^{s+1}E_{n+1}^i, 
\end{eqnarray*}
but in general is not multiplicative.

\specialnumber{3.2.4} \proclaim{Lemma}\label{filtration}
The restriction map induces an isomorphism
$$W_n\,\omega_{(R,M)}^i/\Fil^sW_n\,\omega_{(R,M)}^i
\xto{\sim}W_s\,\omega_{(R,M)}^i.$$
\endproclaim

\demo{Proof} For a fixed value of $n-s$, the filtration quotients
$${}'W_s\,\omega_{(R,M)}^i
=W_n\,\omega_{(R,M)}^i/\Fil^sW_n\,\omega_{(R,M)}^i$$
form a log Witt complex over $(R,M)$. We show that it has the
universal property. Let $(E\s^*,M_E)$ be a log Witt complex over
$(R,M)$. Then there exists a map of log Witt complexes over $(R,M)$:
$${}'W\s\,\omega_{(R,M)}^*\to E\s^*.$$
Indeed, the standard filtration is natural, so we have maps
$$W_n\,\omega_{(R,M)}^i/\Fil^sW_n\,\omega_{(R,M)}^i
\to E_n^i/\Fil^sE_n^i\to E_s^i,$$
where the right-hand map is induced from the restriction maps in
$E\s^*$. We must show that this map of log Witt complexes is unique. To
prove this, it will suffice to show that the canonical map
$$\omega_{(W_s(R),M)}^i\to{}'W_n\,\omega_{(R,M)}^i$$
is surjective. But this follows from the commutativity of the diagram
$$\begin{array}{ccc}
{\omega_{(W_n(R),M)}^i} &\srar\hskip-13pt\to &
{W_n\,\omega_{(R,M)}^i}\\[8pt]
\ddownarrow{}&&\ddownarrow\\[4pt]
{\omega_{(W_s(R),M)}^i} &\srar &
{{}'W_s\,\omega_{(R,M)}^i}\end{array}$$
since the top horizontal and right-hand vertical maps are surjective.
\enddemo 

We define a map $F^{n-1}d\:W_n(R)\to\omega_{(R,M)}^1$ by the formula
$$F^{n-1}d(a)=a_0^{p^{n-1}-1}da_0+a_1^{p^{n-2}-1}da_1+\cdots+da_{n-1},$$
where $a=(a_0,\dots,a_{n-1})$. One easily verifies that $F^{n-1}d$ is
a derivation of $W_n(R)$ into the $W_n(R)$-module
$(F^{n-1})^*\omega_{(R,M)}^1$  and that the following relation holds:
$$dF^{n-1}=p^{n-1}F^{n-1}d.$$
It follows immediately from the derivation property that the formula
$$a\cdot(\omega_1,\omega_2)
=(F^{n-1}(a)\omega_1,F^{n-1}(a)\omega_2-F^{n-1}da\cdot\omega_1)$$
defines a $W_n(R)$-module structure on
$\omega_{(R,M)}^{i-1}\oplus\omega_{(R,M)}^i$. And the relation shows
that
$$(F^{n-1})^*\omega_{(R,M)}^{i-1} \to
\omega_{(R,M)}^{i-1}\oplus\omega_{(R,M)}^i,
\hskip4mm \omega \mapsto (p^{n-1}\omega,-d\omega),$$
is a map of $W_n(R)$-modules. We let ${}_hW_n\,\omega_{(R,M)}^i$ be
the quotient $W_n(R)$-module. This definition is motivated by
Lemma~\ref{hTR} below.

\specialnumber{3.2.5} \proclaim{Lemma}\label{hW}
There is a natural exact sequence of
$W_n(R)$\/{\rm -}\/modules
\begin{eqnarray*}
{}&&(F^{n-1})^*{}_{p^{n-1}}\omega_{(R,M)}^{i-1}
\xto{d} (F^{n-1})^*\omega_{(R,M)}^i \\
{} && \xto{\iota_2} {}_hW_n\,\omega_{(R,M)}^i \xto{\pr_1}
(F^{n-1})^*(\omega_{(R,M)}^{i-1}/p^{n-1})
\to 0. 
\end{eqnarray*}
\endproclaim

\demo{Proof} Indeed, as an abelian group, ${}_hW_n\,\omega_{(R,M)}^i$ is
equal to the push out
$$\begin{array}{ccc}
{\omega_{(R,M)}^{i-1}}&\stackrel{d}{\srar} &
{\omega_{(R,M)}^i} \\[4pt]
\big\downarrow{\scriptstyle p^{n-1}} &&\big\downarrow {\scriptstyle\iota_2}\\[4pt]
{\omega_{(R,M)}^{i-1}}&\stackrel{\iota_1}{\srar} &
{{}_hW_n\,\omega_{(R,M)}^i} \end{array}
$$
so the underlying sequence of abelian groups is exact. One readily
verifies that the various maps are $W_n(R)$-linear.
\enddemo

\specialnumber{3.2.6}\proclaim{Proposition} \label{kernelR}
For any log ring $(R,M)${\rm ,} there is a
natural exact sequence of $W_n(R)$\/{\rm -}\/modules{\rm ,}
$${}_hW_n\,\omega_{(R,M)}^i\xto{N}
W_n\,\omega_{(R,M)}^i\xto{R}W_{n-1}\,\omega_{(R,M)}^i\to 0,$$
where $N(\omega_1,\omega_2)
=dV^{n-1}\lambda(\omega_1)+V^{n-1}\lambda(\omega_2)$.
\endproclaim

\demo{Proof} It follows immediately from Definition~\ref{logwittcx}
that for all $a\in W_n(R)$,
$$\lambda(F^{n-1}da)=F^{n-1}d\lambda(a),$$
and hence $N$ is $W_n(R)$-linear. Since the image of $N$ is equal to
$\Fil^{n-1}W_n\,\omega_{(R,M)}^i$, the statement follows from Lemma~\ref{filtration}.
$\phantom{\sum^|}$
\enddemo

\specialnumber{3.2.7}\proclaim{{C}orollary}\label{drwdivisible}
Let $A$ be a complete discrete
valuation ring of mixed characteristic $(0,p)$ with perfect residue
field{\rm ,} and let $\alpha\:M\to A$ be the canonical log structure. Then
for all $n\geq 1$ and $i\geq 2${\rm ,} $W_n\,\omega_{(A,M)}^i$ is a uniquely
divisible group.
\endproclaim

\demo{Proof} Lemma~\ref{divisible} shows that $\omega_{(A,M)}^i$
is uniquely divisible, if $i\geq 2$. It follows that
${}_hW_n\,\omega_{(A,M)}^i$ is uniquely divisible, if $i\geq 3$, and
an induction argument based on Proposition~\ref{kernelR} shows that so is
  $W_n\,\omega_{(A,M)}^i$. The group ${}_hW_n\,\omega_{(A,M)}^2$
is a direct sum of a uniquely divisible group and the group
$\omega_{(A,M)}^1/p^{n-1}$. Hence $W_n\,\omega_{(A,M)}^2$ is a direct
sum of a uniquely divisible group and a finitely generated torsion
$W(k)$-module. It is therefore enough to show that the modulo $p$
reduction $\bar W_n\,\omega_{(A,M)}^2$ is trivial. Inductively, it
suffices to show that the map
$$dV^{n-1}\:\bar\omega_{(A,M)}^1\to\bar W_n\omega_{(A,M)}^2 \pagebreak $$
is trivial. The map is $k$-linear, and the domain is generated as a
$k$-vector space by the elements $\pi_K^id\log\pi_K$ with $0\leq
i<e_K$. Now the relation
$$\ul{\pi}_n^{e_K}+\theta_K(\ul{\pi}_n)V(1),$$
valid in $\bar W_n(A)$, shows that
$V^{n-1}(\pi_K^id\log\pi_K)=V^{n-1}(\pi_K^i)d\log_n\pi_K$ is either
trivial or contained in the span of elements of the form
$\ul{\pi_K^j}_nd\log_n\pi_K$. But these elements have vanishing
differential.
\enddemo

3.3.\quad We refer the reader to~\cite[\S2]{hm3} for a fuller
discussion of the following result.

\specialnumber{3.3.1}\proclaim{Proposition}
 The homotopy groups $\TR_*\cs(A|K;p)$ form a log Witt
complex over $(A,M)${\rm ,} provided that $p$ is odd. In particular{\rm ,} there
is a canonical map
$$W\s\,\omega_{(A,M)}^*\to\TR\cs_*(A|K;p).$$
\endproclaim

\demo{Proof}
 We recall from Proposition~\ref{diflogring} above that
for all $n\geq 1$, the homotopy groups $\TR_*^n(A|K;p)$ form a log
differential graded ring whose underlying log ring is $(W_n(R),M)$.
The relation that for all $a\in M$,
$$Fd\log_na=d\log_{n-1}a,$$
is immediate from the definition of the maps $F$ and $d\log_n$, and
the remaining relations are proved in~\cite[Lemma  3.3]{hm} and
\cite[Lemmas 1.5.1 and 1.5.6]{h}.
\enddemo

The homotopy groups of the homotopy orbit spectra,
$${}_h\!\TR^n_*(A|K;p)=\pi_*(\borel(C_{p^{n-1}},T(A|K))),$$
are differential graded modules over $\TR^n_*(A|K;p)$, and there are
$\TR^n_*(A|K;p)$-linear maps
\begin{eqnarray*}
&&F\:{}_h\!\TR^n_*(A|K;p)\to F^*({}_h\!\TR^{n-1}_*(A|K;p)),\\
&&V\: F^*({}_h\!\TR^{n-1}_*(A|K;p))\to{}_h\!\TR^n_*(A|K;p),
\end{eqnarray*}
which satisfy  that $FdV=d$ and $FV=p$. Moreover, there is a natural
spectral sequence of $W_n(A)$-modules,
\begin{equation}\label{homotopyorbitss}
E^2_{s,t}=H_s(C_{p^{n-1}},(F^{n-1})^*\,\pi_tT(A|K))\Rightarrow
{}_h\!\TR^n_{s+t}(A|K;p). \hskip.25in\speqnu{3.3.2}
\end{equation}
The reader is referred to~\cite[\S1]{h} and~\cite[\S5]{hm} for   proofs
of these statements.

\specialnumber{3.3.3} \proclaim{Lemma}\label{hTR}
Let $\iota\:\omega_{(A,M)}^i\to\pi_iT(A|K)$ be
the canonical map. Then the map 
\begin{eqnarray*}
{}_hW_n\,\omega_{(A,M)}^i&\to&{}_h\!\TR^n_i(A|K;p),\\
(\omega_1,\omega_2) &\mapsto&
dV^{n-1}\iota(\omega_1)+V^{n-1}\iota(\omega_2), 
\end{eqnarray*}
is a map of $W_n(A)$-modules. It is an isomorphism  for $i\leq 1${\rm ,} and
for $i=2${\rm ,} there is an exact sequence
$$(F^{n-1})^*(A/p^{n-1}) \to {}_hW_n\,\omega_{(A,M)}^2
\to {}_h\!\TR_2^n(A|K;p)\to 0,$$
where the map on the left takes $a$ to $(da,0)$.
\endproclaim

\demo{Proof} If $a\in W_n(A)$, $\omega_1\in\omega_{(A,M)}^{i-1}$
and $\omega_2\in\omega_{(A,M)}^i$, then
\begin{eqnarray*}
a \cdot dV^{n-1}\iota(\omega_1)
& = &d(a\cdot V^{n-1}\iota(\omega_1))
- da\cdot V^{n-1}\iota(\omega_1) \\
{} & =& dV^{n-1}(F^{n-1}a\cdot \iota(\omega_1))
- V^{n-1}(F^{n-1}da\cdot \iota(\omega_1)) \\
{} & =& dV^{n-1}\iota(F^{n-1}a\cdot\omega_1)
- V^{n-1}\iota(F^{n-1}da\cdot\omega_1)),\\
a \cdot V^{n-1}\iota(\omega_2)
& = &V^{n-1}(F^{n-1}a\cdot \iota(\omega_2))\\
{} & =& V^{n-1}\iota(F^{n-1}a\cdot\omega_2), 
\end{eqnarray*}
which shows that the map of the statement is indeed a map of
$W_n(A)$-modules. The map $\iota$ is an isomorphism for $i\leq 2$. So
the spectral sequence gives an isomorphism of $W_n(A)$-modules
$$\iota_0\:(F^{n-1})^*A\xto{\sim}{}_h\!\TR_0^n(A|K;p)$$
and a natural exact sequence of $W_n(A)$-modules
$$0\to(F^{n-1})^*\omega_{(A,M)}^1
\xto{\iota_1}{}_h\!\TR_1^n(A|K;p)
\to(F^{n-1})^*(A/p^{n-1})\to 0.$$
The sequence of Lemma~\ref{hW} maps to the sequence above, and the map
of the left-hand terms is an isomorphism. It remains to show that the
same holds for the map of the right-hand terms. This map is induced
from the composite
$$A\to {}_hW_n\,\omega_{(A,M)}^1\to{}_h\TR_1^n(A|K;p)
\to A/p^{n-1}$$
which in turn may be identified with the map
$$H_0(C_{p^{n-1}},A)\to H_1(C_{p^{n-1}},A)$$
given by multiplication by the fundamental class
$[{\Bbb T}/C_{p^{n-1}}]$. This map is an epimorphism with kernel
$p^{n-1}A$, and the lemma follows for $i=1$. The statement for $i=2$
is proved in a similar manner, using the spectral sequence
in total degree $\leq 3$ and Proposition~\ref{d2=connes} below.
\enddemo

\specialnumber{3.3.4}\numbereddemo{{R}emark} For $i\leq 1$, the proof above does not use the fact that
$A$ is a discrete valuation ring beyond the definition of $T(A|K)$. In
effect, the same proof gives an isomorphism
$${}_hW_n\,\Omega_R^1\xto{\sim}\pi_1\borel(C_{p^{n-1}},T(R)),$$
for any $\Z_{(p)}$-algebra $R$.
\enddemo

Since $\omega_{(A,M)}^2$ is a uniquely divisible group, by
Lemma~\ref{divisible}, the spectral sequence \eqref{homotopyorbitss}
gives an exact sequence of $W_n(A)$-modules
$$(F^{n-1})^*(A/p^{n-1})\xto{d}(F^{n-1})^*(\omega_{(A,M)}^1/p^{n-1})\to
{}_h\!\TR_2^n(A|K;p,\Zp)\to 0,$$
and $d$ is $W_n(A)$-linear since $dF^{n-1}=p^{n-1}F^{n-1}d$. If $\pi_K$
is a uniformizer, then $d\log\pi_K$ represents a class in the
cokernel. We denote this class by $[d\log\pi_K]_n$.

\specialnumber{3.3.5} \proclaim{Lemma}
The map of $W_n(A)$\/{\rm -}\/modules
$$F\:{}_h\!\TR_2^n(A|K;p,\Zp)\to{}_h\!\TR_2^{n-1}(A|K;p,\Zp)$$
is a surjection whose kernel is generated by $p^{n-2}[d\log\pi_K]_n$.
\endproclaim

\demo{Proof} The exact sequence above shows that the map of the
statement is a surjection and that the kernel is a quotient of the
cokernel of the following map:
$$(F^{n-1})^*(p^{n-2}A/p^{n-1}A) \xto{d}
(F^{n-1})^*(p^{n-2}\omega_{(A,M)}^1/p^{n-1}\omega_{(A,M)}^1).$$
Hence, it suffices to show that this cokernel is generated by
$p^{n-2}[d\log\pi_K]_n$. We consider the polynomial ring $P=W(k)[x]$
with the pre-log structure $\alpha\:\N_0\to P$ given by
$\alpha(i)=x^i$. The map of $W(k)$-algebras $\e\:P\to A$, $\e(x)=\pi_K$,
preserves the pre-log structure and induces a surjection
$\omega_{(P,\,\N_0)}^1\twoheadrightarrow\omega_{(A,M)}^1$. It follows
that the map $p^i\omega_{(P,\,\N_0)}^1\twoheadrightarrow p^i\omega_{(A,M)}^1$
is a surjection for $i\geq 0$, and therefore  it will be enough to
show that the cokernel of the map
$$(F^{n-1})^*(p^{n-2}P/p^{n-1}P) \xto{d}
(F^{n-1})^*(p^{n-2}\omega_{(P,\,\N_0)}^1/p^{n-1}\omega_{(P,\,\N_0)}^1)$$
is generated as a $W_n(P)$-module by the canonical image of
$p^{n-2}d\log x$. Now as a $P$-module, the quotient
$p^{n-2}\omega_{(P,\,\N_0)}^1/p^{n-1}\omega_{(P,\,\N_0)}^1$ 
is generated by $p^{n-2}d\log x$, and hence the $W_n(P)$-module
$(F^{n-1})^*(p^{n-2}\omega_{(P,\,\N_0)}^1/p^{n-1}\omega_{(P,\,\N_0)}^1$
is generated by the elements $p^{n-2}d\log x$ and
$p^{n-2}x^{p^i}d\log x$, $0\leq i<n-1$. But the last $n-1$ generators
are all in the image of the map $d$:
$$p^{n-2}x^{p^i}d\log x=p^{n-2-i}d(x^{p^i}).$$
Hence the cokernel of $d$ is generated by $p^{n-2}d\log x$, and the
lemma follows.
\enddemo

\specialnumber{3.3.6}\proclaim{Proposition} \label{pi_1exact}
The sequences
$$0 \to {}_h\!\TR_i^n(A|K;p)
\xto{N} \TR_i^n(A|K;p)
\xto{R} \TR_i^{n-1}(A|K;p) \to 0$$
are exact for $i\leq 1$, and \,$\TR_2^n(A|K;p)$ is uniquely
divisible.
\endproclaim

\demo{Proof} The statement for $i=0$ is~\cite[Prop.\ 3.3]{hm}
and for $i=1$ is equivalent to the statement that the norm map
is injective. The corresponding sequence of maximal uniquely
divisible subgroups is exact, since $F^{n-1}\circ N$ is injective on
this part. Hence, it suffices to show that $\TR_2^{n-1}(A|K;p)$ is
uniquely divisible.

We show by induction on $m\geq 1$ that $\TR_2^m(A|K;p)$ is uniquely
divisible, or equivalently, that $\TR_2^m(A|K;p,\Zp)$ vanishes. The
basic case $m=1$ follows from Proposition~\ref{pi1} and
Lemma~\ref{divisible}. In the induction step, we show that
$$\partial_{K,m}\:\TR_3^{m-1}(A|K;p,\Zp)
\to{}_h\!\TR_2^m(A|K;p,\Zp)$$
is surjective. We first consider the case $m=2$. In the diagram of
$W_2(A)$-modules
$$\begin{array}{ccc}
{\TR_3^1(A|K;p,\Zp)} &\stackrel{\partial_{K,2}}{\srar}&{{}_h\!\TR_2^2(A|K;p,\Zp)} \\[4pt]
\ddownarrow{\scriptstyle\delta} &&\ddownarrow {\scriptstyle\delta}\\[4pt]
{\TR_2^1(k;p)} &\stackrel{\partial_k}{\srar\hskip-13pt\to} &
{{}_h\!\TR_1^2(k;p)},\end{array}$$
the lower horizontal map and the left-hand vertical map are both
surjections. Indeed, for the former, this was proved in~\cite[Th.\ 5.5]{hm}, and for the latter, it follows from the fact, proved
in~\cite{lm}, that $\TR_2^1(A;p,\Zp)$ is trivial. The upper right-hand
group $Q$ is a quotient of the $W_2(A)$-module
$M=F^*(\omega_{(A,M)}^1/p)$. We claim that $M$ is annihilated by the
ideal $I=VW_2(A)+pW_2(A)$. Indeed, as an abelian group $M$ is
$p$-torsion and $FV=p$. It follows that also $Q$ is annihilated by
$I$, and we can therefore view it as a module over the quotient ring
$W_2(A)/I$. This ring is isomorphic to $A/p$, the isomorphism given by
$$W_2(A)/I\xto{\sim}A/p,\hskip5mm a+I\mapsto R(a)+pA,$$
and we let $g\:A/p\to W_2(A)/I$ denote the inverse. The $A/p$-module
$g^*Q$ is generated by the class $[d\log\pi_K]_2$. The image of this
class under the right-hand vertical map is a generator $\iota_1$ of
the $W_2(A)$-module ${}_h\!\TR_1^2(k;p)$, which is isomorphic to
$k$. We now pick $\alpha\in\TR_3^1(A|K;p,\Zp)$ such that
$\delta(\partial_{K,2}(\alpha))=\iota_1$. The difference
$\beta=\partial_{K,2}(\alpha)-[d\log\pi_K]_2$ is in the kernel of
$\delta$, and therefore, 
$$\beta=g(x\pi_K)\cdot[d\log\pi_K]_2,$$
for some $x\in A/p$. We then have
$$g(1+x\pi_K)\cdot[d\log\pi_K]_2=\partial_{K,2}(\alpha),$$
and since $(1+x\pi_K)\in(A/p)^\times$,
$$[d\log\pi_K]_2=(g(1+x\pi_K)^{-1})\cdot\partial_{K,2}(\alpha).$$
We would like to know that the map of units
$$W_2(A)^\times\to(W_2(A)/I)^\times$$
is a surjection. This will follow if we know that the $I$-adic
topology on $W_2(A)$ is complete and separated. But the formula
$$V(x)\cdot V(y)=V(FV(x)y)=V(pxy)=pV(xy)$$
implies that the $I$-adic and $p$-adic topologies on $W_2(A)$
coincide, and the $p$-adic topology is complete and separated. So we
can find a unit $u\in W_2(A)^\times$ such that $u+I=g(1+x\pi_K)$. Since
$\partial_{K,2}$ is $W_2(A)$-linear, 
$$[d\log\pi_K]_2=u^{-1}\partial_{K,2}(\alpha)=\partial_{K,2}(u^{-1}\alpha),$$
which concludes the proof for $m=2$.

We now proceed inductively, and consider the diagram
$$\begin{array}{ccccc}
{\TR_3^{m-1}(A|K;p,\Zp)}&\hskip-4pt\stck{\partial_{K,m}}\hskip-4pt& 
{{}_h\!\TR_2^m(A|K;p,\Zp)} &\hskip-4pt\shtck{N}\hskip-4pt&
{\TR_2^m(A|K;p,\Zp)} \\[4pt]
\ddownarrow{\scriptstyle F} &&\ddownarrow{\scriptstyle F} &&\ddownarrow{\scriptstyle F} \\[4pt]
{\TR_3^{m-2}(A|K;p,\Zp)} &\hskip-4pt\stck{\partial_{K,m-1}}\hskip-4pt &
{{}_h\!\TR_2^{m-1}(A|K;p,\Zp)}  &\hskip-4pt\shtck{N}\hskip-4pt&
{\TR_2^{m-1}(A|K;p,\Zp).} \end{array}$$
Inductively, the map $\partial_{K,m-1}$ is surjective, and the left-hand vertical map $F$ is surjective by Lemma~\ref{frobsurj}. Moreover,
the kernel of the middle vertical map is generated as a
$W_m(A)$-module by the class $p^{m-2}[d\log\pi_K]_m$. It therefore
suffices to show that this class is in the image of $\partial_{K,m}$
in the top row, and this in turn will follow if we show that the
class $[d\log\pi_K]_m$ is in the image of $\partial_{K,m}$.  To see
this, we pick $\alpha\in\TR_3^{m-1}(A|K;p),\Zp)$ such that
$\partial_{K,m-1}(F(\alpha))=[d\log\pi_K]_{m-1}$. Then
$\beta=\partial_{K,m}(\alpha)-[d\log\pi_K]_m$ is in the kernel of the
middle vertical map, so we can write $\beta=x\cdot p^{m-2}d\log\pi_K$,
for some $x\in W_m(A)$. But then
$$(1+p^{m-2}x)[d\log\pi_K]_m=\partial_{K,m}(\alpha),$$
and hence
$$[d\log\pi_K]_m=(1+p^{m-2}x)^{-1}\partial_{K,m}(\alpha)
=\partial_{K,m}((1+p^{m-2}x)^{-1}\alpha),$$
where the inverse exists since the $p$-adic topology on $W_m(A)$ is
complete and separated.
\enddemo 

\specialnumber{3.3.7}\proclaim{Addendum} \label{trAdivisible}
The group $\TR^n_2(A;p)$ is
uniquely divisible for all $n$.
\endproclaim

\demo{Proof} It suffices to show that $\TR_2^n(A;p,\Zp)$ is
trivial. We prove this by induction, and refer to the proof of
Proposition~\ref{pi1} for the case $n=1$. Since
$\TR_2^n(A|K;p,\Zp)$ vanishes, there is an exact sequence
$$\TR_3^n(A|K;p,\Zp)\xto{\delta_n}\TR_2^n(k;p)
\to\TR_2^n(A;p,\Zp)\to 0,$$
and we must prove that the map $\delta_n$ is surjective. We
consider the diagram
$$\begin{array}{ccc}
{\TR_3^n(A|K;p,\Zp)} &\stck{\delta_n} &
{\TR_2^n(k;p)}  \\[4pt]
 \ddownarrow{\scriptstyle F} &&\ddownarrow{\scriptstyle F}\\[4pt]
{\TR_3^{n-1}(A|K;p,\Zp)}&\stackrel{\delta_{n-1}}{\srar\hskip-13pt\to} &
{\TR_2^{n-1}(k;p).}\end{array}
$$
The map $\delta_{n-1}$ is surjective by induction, and the left-hand
vertical map is surjective by Lemma~\ref{frobsurj}. Moreover, it was
proved in~\cite[Th.\ 5.5]{hm} that the right-hand vertical map $F$
is a surjection whose kernel is equal to the image of the map
$$V\:\TR^{n-1}_2(k;p)\to\TR^n_2(k;p).$$
Since the square
$$\begin{array}{ccc}
{\TR_3^{n-1}(A|K;p,\Zp)} &\stackrel{\delta_{n-1}}{\dlar}  &
{\TR^{n-1}_2(k;p)} \\[4pt]
\scs{V}&&\scs{V}\\[4pt]
{\TR_3^n(A|K;p,\Zp)} &\stck{\delta_n} &
{\TR^n_2(k;p)} \end{array}
$$
commutes and the top horizontal map is a surjection, the proof of the
induction step is complete.
\enddemo  

\specialnumber{3.3.8}\proclaim{Theorem} \label{pionetr}
The canonical map
$$W_n\,\omega_{(A,M)}^q\to\TR_q^n(A|K;p)$$
is an isomorphism{\rm ,} for $q\leq 2${\rm ,} and a rational isomorphism{\rm ,} for all
$q\geq 0$.
\endproclaim

\demo{Proof} The proof is by induction on $n$ starting from
Proposition~\ref{pi1}. In the induction step, we use the exact
sequences of Lemma~\ref{kernelR} and Proposition~\ref{pi_1exact},
$$\begin{array}{ccccccccc}
&&
{{}_hW_n\,\omega_{(A,M)}^q} &\longrightarrow &
{W_n\,\omega_{(A,M)}^q} &\shtck{R} &
{W_{n-1}\,\omega_{(A,M)}^q} &\longrightarrow&
{0} \\[4pt]
&&\big\downarrow&&\big\downarrow&&\scs{\sim}\\[4pt]
{0} &\longrightarrow&
{{}_h\!\TR_q^n(A|K;p)} &\shtck{N} &
{\TR_q^n(A|K;p)} &\shtck{R} &
{\TR_q^{n-1}(A|K;p)}&\longrightarrow &
{0,}\end{array}$$
where the lower sequence is exact, for $q\leq 1$, and exact modulo
torsion, for all~$q$. If $q\leq 1$, the left-hand vertical map is an
isomorphism by Lemma~\ref{hTR}, and hence the statement follows in
this case. When $q=2$, the left-hand vertical map is an epimorphism
with torsion kernel. Since the domain and range of the middle and
right-hand vertical maps are both divisible groups, the statement
follows.
\enddemo

In the proof of   Proposition~\ref{pi_1exact},
Addendum~\ref{trAdivisible}, and Theorem~\ref{pionetr} above for $n>3$
we have used Lemma~\ref{frobsurj} below. However, the lemma is not
needed to prove these statements for $n\leq 3$. In particular, the
proof of the following result does not use Lemma~\ref{frobsurj}.

\specialnumber{3.3.9}\proclaim{Addendum} \label{connecting}
The connecting homomorphism
$$\partial\:\TR^1_2(A|K;p,\Z/p)\to{}_h\!\TR^2_1(A|K;p,\Z/p)$$
maps $\kappa$ to $dV(1)-V(d\log(-p))$.
\endproclaim

\demo{Proof} To prove the statement, we apply Lemma~\ref{3x3}
below to the $3\times 3$-diagram obtained from the smash product of
the coefficient cofibration sequence
$$S^0\xto{p}S^0\xto{i}M_p\xto{\beta}S^1$$
and the fundamental cofibration sequence
$${}_h\!\TR^n(A|K;p)\xto{N}\TR^n(A|K;p)\xto{R}\TR^{n-1}(A|K;p)
\xto{\partial}\Sigma({}_h\!\TR^n(A|K;p)).$$
Since $\TR_2(A|K;p)$ is uniquely divisible and $\TR_0(A|K;p)$ torsion-free, the lemma shows that the
connecting homomorphism of the statement is equal to the opposite of the connecting homomorphism
associated with the diagram
$$\begin{array}{ccccccccc}
{0}&\longrightarrow &
{{}_h\TR^2_1(A|K;p)} &\shtck{N}  &
{\TR^2_1(A|K;p)} &\shtck{R}   &
{\TR^1_1(A|K;p)}&\longrightarrow  &
{0} \\[4pt]
&&\scs{p}&&\scs{p}&&\scs{p}\\[4pt]
{0}&\longrightarrow &
{{}_h\TR^2_1(A|K;p)} &\shtck{N} &
{\TR^2_1(A|K;p)} &\shtck{R} &
{\TR^1_1(A|K;p)}&\longrightarrow &
{0.} \end{array}
$$
And by Theorem~\ref{pionetr}, this diagram is canonically isomorphic
to the diagram
$$\begin{array}{ccccccccc}
{0}&\longrightarrow &
{{}_hW_2\omega_{(A,M)}^1} &\shtck{N}  &
{W_2\omega_{(A,M)}^1}&\shtck{R}  &
{W_1\omega_{(A,M)}^1} &\longrightarrow &
{0} \\[4pt]
&&\scs{p}&&\scs{p}&&\scs{p}\\[4pt]
{0}& \longrightarrow &
{{}_hW_2\,\omega_{(A,M)}^1} &\shtck{N} &
{W_2\,\omega_{(A,M)}^1} &\shtck{R} &
{W_1\,\omega_{(A,M)}^1}  &\longrightarrow &
{0.} \end{array}
$$
The Bockstein maps $\kappa$ to $d\log(-p)\in W_1\,\omega_{(A,M)}^1$,
which is the image by the restriction of $d\log_2(-p)\in 
W_2\,\omega_{(A,M)}^1$. To evaluate $pd\log_2(-p)$ we use the formula
$$-\ul{(-p)}_2+V(1)=p(1+p^{p-2}V(1)),$$
which one readily verifies using the ghost map. Differentiating,
we find
$$-d\ul{(-p)}_2+dV(1)=p^{p-2}dV(1)=0,$$
and if we multiply by $d\log_2(-p)$, we get
$$-d\ul{(-p)}_2+V(d\log(-p))=pd\log_2(-p)+p^{p-2}V(d\log(-p))=pd\log_2(-p).$$
This shows that $pd\log_2(-p)=V(d\log(-p))-dV(1)$ as desired.
\enddemo

\specialnumber{3.3.10} \proclaim{Lemma}\label{3x3}
Given a $3\times 3$\/{\rm -}\/diagram of cofibration
sequences
$$\begin{array}{ccccccc}
{E_{11}} &\stck{f_{11}} & 
{E_{12}}&\stck{f_{12}}  &
{E_{13}}&\stck{f_{13}}  &
{\Sigma E_{11}}  \\[4pt]
\phantom{\scr{g_{11}}}\scs{g_{11}}&&\phantom{\scr{g_{11}}}\scs{g_{12}}&&\phantom{\scr{g_{11}}}\scs{g_{13}}&&\phantom{\scr{g_{11}}}\scs{\Sigma g_{11}}\\[4pt]
{E_{21}} &\stck{f_{21}} &
{E_{22}} &\stck{f_{22}} &
{E_{23}} &\stck{f_{23}} &
{\Sigma E_{21}} \\[4pt]
 \phantom{\scr{g_{11}}}\scs{g_{21}}&& \phantom{\scr{g_{11}}}\scs{g_{22}}&& \phantom{\scr{g_{11}}}\scs{g_{23}}&& \phantom{\scr{g_{11}}}\scs{\Sigma g_{21}}\\[4pt]
{E_{31}}& \stck{f_{31}}  &
{E_{32}}& \stck{f_{32}}  &
{E_{33}}& \stck{f_{33}} &
{\Sigma E_{31}} \\[4pt]
\phantom{\scr{g_{11}}}\scs{g_{31}}&&\phantom{\scr{g_{11}}}\scs{g_{32}}&&\phantom{\scr{g_{11}}}\scs{g_{33}}& {\scriptstyle (-1)}  & \phantom{\scr{g_{11}}}\scs{-\Sigma g_{11}}\\[4pt]
{\Sigma E_{11}} &\stck{\Sigma f_{11}} &
{\Sigma E_{12}} &\stck{\Sigma f_{12}} &
{\Sigma E_{13}} &\stck{-\Sigma f_{13}} &
{\Sigma^2 E_{11}} \end{array}$$
and classes $e_{ij}\in\pi_*E_{ij}$ such that
$g_{33}(e_{33})=\Sigma f_{12}(e_{12})$ and
$f_{33}(e_{33})=\Sigma g_{21}(e_{21})${\rm .} Then the sum
$f_{21}(e_{21})+g_{12}(e_{12})$ is in the image of
$\pi_*E_{11}\to\pi_*E_{22}$. 
\endproclaim

3.4.\quad The $k$-algebra $\bar{W}_n(A)$ was evaluated in
Proposition~\ref{wittvectors} above. We now evaluate the differential
graded $k$-algebra $\bar{W}_n\,\omega_{(A,M)}^*$. Let $\pi=\pi_K$ be a
uniformizer. Then the modified Verschiebung from~(\ref{modified})
satisfies
$$FdV_\pi(a)=\theta_K(\ul{\pi})^pda.$$
Let $r=r(i,e_K)=v_p(i-pe_K/(p-1))$.

\specialnumber{3.4.1}\proclaim{Proposition} \label{drwp}
The differential graded $k$\/{\rm -}\/algebra
$E^*=\bar{W}_n\,\omega_{(A,M)}^*$ is concentrated in degrees
$0$ and $1$ and satisfies\/{\rm :}\/
\medbreak
{\rm(i)} A $k$-basis for $E_n^1$ is given by the elements
$V_\pi^s(\ul\pi^id\log\pi)${\rm ,} where $0\leq i<e_K$ and $0\leq s\leq r${\rm ,}
and $dV_\pi^s(\ul\pi^i)${\rm ,} where $0\leq i<e_K$ and $r<s<n$.
Moreover{\rm ,} $V_\pi^s(\ul\pi^id\log\pi)$ vanishes{\rm ,} if $s>r${\rm ,}
$dV_\pi^s(\ul{\pi}^i)$ vanishes{\rm ,} if $s<r${\rm ,} and
$$dV_\pi^r(\ul\pi^i) = p^{-r}(i-pe_K/(p-1))\cdot
V_\pi^r(\ul\pi^id\log\pi).$$
\medbreak
{\rm(ii)} The $E_n^0$\/{\rm -}\/module structure on $E_n^1$ is given by
\begin{eqnarray*}
&&\hskip-24pt V_\pi^s(\ul\pi^i)dV_\pi^t(\ul\pi^j)  =\left\{ \begin{array}{ll}
dV_\pi^t(\ul{\pi}^{p^ti+j}) - iV_\pi^t(\ul{\pi}^{p^ti+j}d\log\pi)
& \hbox{if $0=s\leq t$}, \\
-iV_\pi^t(\theta_K(\ul{\pi})^{p^{t-s}(\frac{p^{s+1}-1}{p-1}-1)}
\ul{\pi}^{p^{t-s}i+j}d\log\pi)  & \hbox{if $0<s\leq t$},\\
jV_\pi^s(\theta_K(\ul\pi)^{p^{s-t}(\frac{p^{t+1}-1}{p-1}-1)}
\ul\pi^{i+p^{s-t}j}d\log\pi) 
&\hbox{if $s\geq t$}, \\
\end{array}\right. \\[5pt]
&&\hskip-24pt V_\pi^s(\ul\pi^i)V_\pi^t(\ul\pi^jd\log\pi) =\left\{ \begin{array}{ll}
V_\pi^t(\ul{\pi}^{p^ti+j}d\log\pi)  & \hbox{if $s=0$},\\
V_\pi^s(\ul{\pi}^{i+p^sj}d\log\pi)   & \hbox{if $t=0$},\\
0  &\hbox{otherwise.} \\
\end{array}\right. 
\end{eqnarray*}
\endproclaim

\demo{Proof} It follows from Propositions~\ref{wittvectors}
and~\ref{derhamwitt} that $E_n^*$ is generated, as a graded
$k$-vector space, by the monomials in the variables
$V_{\pi}^s(\ul\pi^i)$, $dV_{\pi}^s(\ul\pi^i)$, $V_{\pi}^s(\ul\pi^i
d\log\pi)$, and $dV_{\pi}^s(\ul\pi^i d\log\pi)$ with $0\leq s<n$ and
$i\geq 0$.  Theorem~\ref{pionetr} and Corollary~\ref{drwdivisible}
show that $E_n^q$ vanishes, for $q\geq 2$. In particular, the latter
generators, which are of degree two, must vanish. 

We verify the relations in (i). If $s\leq r$ then
$p^{-s}(i+pe_K(p^s-1)/(p-1))$ is an integer, and iterated use of
the second relation in Proposition~\ref{wittvectors} shows that
$$V_\pi^s(\ul\pi^i)=\ul\pi^{p^{-s}(i+pe_K\frac{p^s-1}{p-1})}.$$
It follows that $dV_{\pi}^s(\ul{\pi}^i)$ vanishes, if $s<r$, and that
$dV_{\pi}^r(\ul{\pi}^i)$ and $V_{\pi}^r(\ul{\pi}^id\log\pi)$ are
related as stated. And since $V_{\pi}d$ is the zero homomorphism, this
also shows that for $s>r$,
$V_{\pi}^s(\ul{\pi}^id\log\pi)=V_{\pi}^{s-r}V_{\pi}^r(\ul{\pi}^id\log\pi)$
vanishes.

The formulas in (ii) are readily obtained by differentiating the first
set of relations in Proposition~\ref{wittvectors}. If, for instance,
$0<s\leq t<n$, we find that
\begin{eqnarray*}
V_{\pi}^s(\ul{\pi}^i)dV_{\pi}^t(\ul{\pi}^j) & =&
-dV_{\pi}^s(\ul{\pi}^i)V_{\pi}^t(\ul{\pi}^j) 
= - V_{\pi}^t(F^tdV_{\pi}^s(\ul{\pi}^i)\ul{\pi}^j) \\[4pt]
{} & =& - i V_{\pi}^t(
\theta_K(\ul{\pi})^{p^{t-s}(\frac{p^{s+1}-1}{p-1}-1)}
\ul{\pi}^{p^{t-s}i+j}d\log\pi), 
\end{eqnarray*}
and the remaining formulas are verified in a similar manner. It
remains to prove that this gives all relations in $E_n^1$. This is
the case if and only if $E_n^1$ is an $ne_K$-dimensional $k$-vector
space. We prove in Proposition~\ref{dimension} below that this is
indeed the case, and hence there can be no further relations.
\enddemo 

\section{Tate cohomology and the Tate spectrum} \label{tatecohomology}

4.1.\quad Let $G$ be a finite group and let $k$ be a commutative
ring. The norm element $N_G\in kG$ is defined as the sum of all the
elements of $G$. If $M$ is a left $kG$-module, multiplication by $N_G$
defines a map
$$N_G\:M_G\to M^G$$
from the coinvariants $M_G=k\otimes_{kG}M$ to the invariants
$M^G=\Hom_{kG}(k,M)$. We note that for left $kG$-modules $M$ and $N$,
there are canonical isomorphisms
$$(M\otimes N)_G  \cong c^*M\otimes_{kG}N,\hskip 5mm
\Hom(M,N)^G  \cong \Hom_{kG}(M,N), $$
where $c^*M$ denotes the right $kG$-module with $m\cdot g=g^{-1}m$.

Let $\e\:P\to k$ be a projective resolution and let $\tilde P$ be the
mapping cone of $\e$ such that there is a distinguished triangle
(see \S 2.1 above)
$$P\xto{\e} k \xto{\iota} \tilde P \xto{\partial} \Sigma P.$$

\specialnumber{4.1.1}\numbereddemo{Definition} \label{tate} Let $M$ be a left $kG$-module. The
{\it Tate cohomology of $G$ with coefficients} in $M$ is given by
$$\hat H^*(G,M)=H_{-*}((\tilde P\otimes\Hom(P,M))^G).$$
\enddemo

It is clear that the Tate cohomology groups are well-defined up to
canonical isomorphism. We show that the definition given here agrees
with the usual definition in terms of complete resolutions,
\cite[Chap.\ XII, \S3]{ce}.

\specialnumber{4.1.2} \proclaim{Lemma}\label{grphomology} The following maps are
quasi\/{\rm -}\/isomorphisms\/{\rm :}\/
$$(P\otimes M)_G\begin{array}{c}  {\scriptstyle{N}}\\[-8pt]\longrightarrow\\[-8pt]
{\scriptstyle\sim}\end{array}
 (P\otimes M)^G 
\begin{array}{c} 
{\scriptstyle\id\otimes\e^*}\\[-8pt]\srar\\[-8pt]{\scriptstyle\sim}\end{array}
(P\otimes\Hom(P,M))^G.$$
\endproclaim  

\demo{Proof} We first show that the norm map is an isomorphism of
complexes. It will suffice to show that the norm map
$$(kG\otimes M)_G \xto{N} (kG\otimes M)^G$$
is an isomorphism, for both sides commute with the formation of
arbitrary direct sums. Let $\eta\:k\to kG$ and $\e\:kG\to k$ be the
unit and co-unit of the Hopf algebra $kG$, respectively. Then we have
an isomorphism of left $kG$-modules
$$\xi\:kG\otimes\e^*\eta^*M\xto{\sim}kG\otimes M,\hskip4mm
\xi(g\otimes x)=g\otimes gx.$$
The left-hand side is isomorphic to a direct sum indexed by the
elements of $G$ of copies of $M$, and $G$ acts by permuting the
summands. Hence $N$ is an isomorphism.

In order to show that the right-hand map of the statement is a
quasi-isomorphism, we filter the double complex on the right after
the first tensor factor. This gives, by~\cite[Th.\ 6.1]{boardman},
a strongly convergent fourth quadrant spectral sequence
$$E_{s,t}^1=H_t((P_s\otimes\Hom(P,M))^G) \Rightarrow
H_{s+t}((P\otimes\Hom(P,M))^G),$$
and hence, it suffices to show that for all $s\geq 0$, the map
$$(P_s\otimes M)^G \stck{\id\otimes\e^*} (P_s\otimes\Hom(P,M))^G$$
is a quasi-isomorphism. Since both sides commute with filtered
colimits in the first tensor factor, we can further assume that the
projective $kG$-module $P_s$ is finitely generated. In this case, the
dual $DP_s=\Hom(P_s,k)$ again is a (finitely generated) projective
$kG$-module, and there is a commutative diagram
$$\begin{array}{ccc}
{(P_s\otimes M)^G} &\stck{\id\otimes\e^*} &
{(P_s\otimes\Hom(P,M))^G}   \\[4pt]
\scs{\sim} &&\scs{\sim} \\[4pt]
{\Hom(DP_s,M)^G} &\stck{(\e\otimes\id)^*} &
{\Hom(P\otimes DP_s,M)^G,} \end{array}
$$
with the vertical maps isomorphisms. The map
$$\e\otimes\id\:P\otimes DP_s \xto{\sim} DP_s$$
is a quasi-isomorphism between bounded below complexes of projective
$kG$-modules. Therefore, it is a chain homotopy equivalence, and
hence, so is the lower horizontal map in the diagram above. The
lemma follows.
\enddemo 

\specialnumber{4.1.4}\numbereddemo{{R}emark} \label{completeres} The triangle preceding
Definition~\ref{tate} and Lemma~\ref{grphomology} gives rise to natural
isomorphisms
$$\hat H^i(G,M)\cong\left\{ \begin{array}{ll}
H^i(G,M)&\hbox{if $i\geq 1$}\\
H_{-i-1}(G,M)&\hbox{if $i\leq-1$} 
\end{array}\right.$$
and to a natural exact sequence
$$0\to\hat H^{-1}(G,M)\xto{\partial}
H_0(G,M)\xto{N}H^0(G,M)\xto{i}\hat H^0(G,M)\to 0.$$
Hence, the definition of Tate cohomology given here agrees with the
original one in terms of complete resolutions,~\cite[Chap.\ XII,
\S3]{ce}. This can also be seen more directly as follows. Let
$\e\:\hat P\to k$ be a complete resolution in the sense of
{\it loc.~cit.}, and let $P$ and $P^-$ be the complexes whose nonzero
terms are $P_i=\hat P_i$, if $i\geq 0$, and $P^-_i=\hat P_i$, if
$i<0$, respectively. Then $\e\:P\to k$ is a resolution of $k$ by
finitely generated projective left $kG$-modules and there is a
canonical triangle

\centerline{$P^-\to\hat P\to P\to\Sigma P^-. $}
\smallbreak\noindent 
An argument similar to the proof of Lemma~\ref{grphomology} shows
that the canonical maps
$$\Hom(\hat P,M)^G \xto{\sim}
(\tilde P\otimes\Hom(\hat P,M))^G \xto{\sim}
(\tilde P\otimes\Hom(P,M))^G$$
are quasi-isomorphisms.
\enddemo

{\it Definition} 4.1.5. The {\it cup product}
\smallbreak
\centerline{$\hat H^*(G,M)\otimes\hat H^*(G,M')\to\hat H^*(G,M\otimes M')$}
\smallbreak\noindent 
is the map on homology induced by the composite
\begin{eqnarray*} \noalign{\vskip-2pt}
(\tilde P\otimes\Hom(P,M))^G &\hsm\otimes\hsm&(\tilde P\otimes\Hom(P,M'))^G
\\
&\hsm\to\hsm&(\tilde P\otimes\tilde P\otimes\Hom(P\otimes P,M\otimes M'))^G\\
{}&\hsm\to\hsm&(\tilde P\otimes\Hom(P,M\otimes M'))^G,\\
\noalign{\vskip-20pt} 
\end{eqnarray*}
where the first map is the canonical map, and the second map is
induced from a choice of chain maps $P\to P\otimes P$ and 
$\tilde P\otimes\tilde P\to\tilde P$ compatible with the canonical
isomorphisms $k\to k\otimes k$ and $k\otimes k\to k$,
respectively.
\vglue6pt

It is well-known that the chain map $P\to P\otimes P$ exists and is
unique up to chain homotopy. The analogous statement for the map
$\tilde P\otimes\tilde P\to\tilde P$ is proved in a  similar
manner. Hence, the cup product is well-defined. It makes $\hat
H^*(G,k)$ a graded commutative graded ring and $\hat H^*(G,M)$ a
graded module over this ring.\pagebreak 

  4.2.\quad Let $C$ be a cyclic group of order
$r$ and let $g\in C$ be a generator. We let $\e\:W\to k$ be the
standard resolution which in degree $s\geq 0$ is a free $kC$-module on
a single generator $x_s$ with differential
$$dx_s=\left\{ \begin{array}{ll}
Nx_{s-1},&\hbox{$s$ even,}\\
(g-1)x_{s-1},&\hbox{$s$ odd,}\\
\end{array}\right.$$
and with augmentation $\e(x_0)=1$. Then $\tilde W$ is the complex
which in degree $s>0$ is a free $kC$-module on the generator
$y_s=(0,x_{s-1})$ and in degree $s=0$ is a trivial $kC$-module on the
generator $y_0=(1,0)$. The differential is
$$dy_s=\left\{ \begin{array}{ll}
-(g-1)y_{s-1},&\hbox{$s$ even,}\\
-Ny_{s-1},&\hbox{$s>1$ odd,}\\
-y_0&s=1.\\
\end{array}\right.$$
The dual of $x_s$ is the element $x_s^*\in DW_s=\Hom(W_s,k)$ given by
$x_s^*(g^ix_s)=\delta_{i,0}$. Note that $g^i\cdot x_n^*=(g^ix_n)^*$
and that the map $(g^i)^*\:DW_s\to DW_s$ maps 
$x_s^*\mapsto g^{-i}x_s^*$. Thus
$$dx_s^*=\left\{ \begin{array}{ll}
(g^{-1}-1)x_{s+1}^*,&\hbox{$s$ even,}\\
Nx_{s+1}^*,&\hbox{$s$ odd.}\\
\end{array}\right.$$

\specialnumber{4.2.1} \proclaim{Lemma}\label{tateparc} Suppose that the order of $C$ is odd and
congruent to zero in $k$. Then as a graded $k$\/{\rm -}\/algebra
$$\hat H^*(C,k)=\Lambda\{u\}\otimes S\{t^{\pm1}\}$$
where $t$ and $u$ are the classes of $y_0\otimes Nx_2^*$ and
$y_0\otimes Nx_1^*${\rm ,} respectively. Moreover{\rm ,} the classes $1${\rm ,} $ut^{-1}$
and $t^{-1}$ are represented by the elements $y_0\otimes Nx_0^*${\rm ,}
$-Ny_1\otimes Nx_0^*$ and $Ny_2\otimes Nx_0^*${\rm ,} respectively.
\endproclaim

\demo{Proof} We first evaluate the homology of the complex
$$(\tilde W\otimes\Hom(W,k))^C=(\tilde W\otimes DW)^C.$$
This is the total complex of a double complex, and the filtration
after the first tensor factor gives rise to a fourth quadrant homology
type spectral sequence which converges strongly to the homology of the
total complex,~\cite[Th.\ 6.1]{boardman}. We have
$$E^1_{s,t}=H_{s+t}(\tilde W_s\otimes DW)^C
\xto{\sim}H_{s+t}(\Hom(W,\tilde W_s)^C),$$
which vanishes if both $s$ and $t$ are nonzero. Hence
$E^2_{s,t}=E^\infty_{s,t}$ and it is easy to see that if either $s$
or $t$ is zero, this is a free $k$-module of rank one generated by the
classes of $y_0\otimes Nx_{-t}^*$ and $Ny_s\otimes Nx_0^*$,
respectively. We note that these elements are also cycles in the total
complex.

To evaluate the multiplicative structure, we choose liftings
\begin{eqnarray*}
&&\Psi\:W\to W\otimes W,\\
&&\Phi \:\tilde W\otimes\tilde W\to\tilde W 
\end{eqnarray*}
of the canonical maps $k\to k\otimes k$ and $k\otimes k\to k$,
respectively:
$$\Psi_{m,n}(g^sx_{m+n})=\left\{ \begin{array}{ll}
{\displaystyle{\sum_{s\leq p<q<s}g^px_m\otimes g^qx_n}},&
\hbox{$m$ and $n$ odd}\\
g^sx_m\otimes g^{s+1}x_n,&\hbox{$m$ odd, $n$ even}\\
g^sx_m\otimes g^sx_n,&\hbox{$m$ even} 
\end{array}\right.$$
and
$$\Phi_{m,n}(g^py_m\otimes g^qy_n)=\left\{ \begin{array}{ll}
{\displaystyle{\sum_{p\leq s<q<p}g^sy_{m+n}}},&\hbox{$m$ and $n$ odd}\\
\delta_{p,q+1}g^py_{m+n},&\hbox{$m$ odd, $n$ even}\\
\delta_{p,q}g^py_{m+n},&\hbox{$m$ even}, 
\end{array}\right.$$
where in the first line the sum ranges over the $g^s$ between $g^p$
and $g^{q-1}$, both included, in the cyclic ordering of $C$ specified
by the generator $g$. The sum is zero if and only if $p=q$. The map
$\Psi$ induces a product map on the dual $DW$ given by the composite
$$\Psi^*\:DW\otimes DW\xto{\nu}D(W\otimes W)\xto{D\Psi}DW,$$
or
$$\Psi^*_{m,n}(g^{-p}x_m^*\otimes g^{-q}x_n^*)=\left\{ \begin{array}{ll}
\displaystyle{-\sum_{p\leq s<q<p}g^{-s}x_{m+n}^*},&\hbox{$m$ and $n$ odd}\\
\delta_{p,q+1}g^{-p}x_{m+n}^*,&\hbox{$m$ odd, $n$ even}\\
\delta_{p,q}g^{-p}x_{m+n}^*,&\hbox{$m$ even}. 
\end{array}\right.$$
We find that
$$(y_0\otimes Nx_m^*)\cdot(y_0\otimes Nx_n^*)=\left\{ \begin{array}{ll}
\displaystyle{-\frac{r(r-1)}{2}\;y_0\otimes Nx_{m+n}^*},&
\hbox{$m$ and $n$ odd}\\
y_0\otimes Nx_{m+n}^*&\hbox{otherwise} 
\end{array}\right.$$
and
$$(Ny_m\otimes Nx_0^*)\cdot(Ny_n\otimes Nx_0^*)=\left\{ \begin{array}{ll}
\displaystyle{\frac{r(r-1)}{2}\;Ny_{m+n}\otimes Nx_0^*},&
\hbox{$m$ and $n$ odd}\\
Ny_{m+n}\otimes Nx_0^*&\hbox{otherwise}. 
\end{array}\right.$$
Moreover, the product
$$(y_0\otimes Nx_2^*)\cdot(Ny_2\otimes Nx_0^*)=Ny_2\otimes Nx_2^*$$
is homologous to $y_0\otimes Nx_0^*$, which represents the
multiplicative unit in the cohomology ring. Indeed, with $\Delta(N)=\sum_{0\le s<n}g^s\otimes g^s$
$$d(\Delta(N)(y_1\otimes x_0^*)+\Delta(N)(y_2\otimes x_1^*))
=-y_0\otimes Nx_0^*+Ny_2\otimes Nx_2^*.$$
Hence $Ny_2\otimes Nx_0^*$ represents the class $t^{-1}$. Finally, for
any element $\alpha\in kC$,
$$(1\otimes\alpha)\Delta(N)=(\bar\alpha\otimes 1)\Delta(N),$$
where $\bar\alpha=c(\alpha)$ is the antipode. Therefore, if $\alpha\in
kC$ is such that $(g-1)\alpha=r-N$ (for example
$\alpha=1+2g+\cdots+rg^{r-1}$ is such an element), then
\begin{eqnarray*}
&&\hskip-59pt d((\alpha \otimes 1)\Delta(N)(y_2\otimes x_0^*))\\
{}&=&-((g-1)\otimes 1)(\alpha\otimes1)
\Delta(N)(y_1\otimes x_0^*)\\
{}& &
-(1\otimes(\bar g-1))(1\otimes\bar\alpha)
\Delta(N)(y_2\otimes x_1^*)\\
{}&=&Ny_1\otimes Nx_0^*+Ny_2\otimes Nx_1^*
-r\Delta(N)(y_1\otimes x_0^*+y_2\otimes x_1^*), 
\end{eqnarray*}
and hence, $Ny_1\otimes Nx_0^*$ represents the class $-ut^{-1}$ in the
cohomology ring.
\enddemo  

\specialnumber{4.2.2}\proclaim{Addendum} \label{boundarymap}The boundary map $\partial\:\hat
H^{-1}(C,k)\to H_0(C,k)$ takes $ut^{-1}$ to the class of $-1$.
\endproclaim

\demo{Proof} The boundary map, by definition, is induced by the 
composite
\begin{eqnarray*}
&&(\tilde W\otimes\Hom(W,k))^C   \stck{\partial\otimes\id}
(\Sigma W\otimes\Hom(W,k))^C \begin{array}{c} {\scriptstyle\id\otimes\e^*}\\[-8pt] \slar\\[-8pt]
{\scriptstyle 
\sim}\end{array}  (\Sigma W\otimes k)^C \\[4pt]
&&\phantom{(\tilde W\otimes\Hom(W,k))^C}
 \begin{array}{c} {\scriptstyle N}\\[-8pt] \longleftarrow\\[-8pt]{\scriptstyle\sim} \end{array}
(\Sigma W\otimes k)_C \xto{\e\otimes 1}
\Sigma k_C.  
\end{eqnarray*}
The class $ut^{-1}$ is represented by the element $-Ny_1\otimes
Nx_0^*$ whose image under $\partial\otimes\id$ is $-Nx_0\otimes
Nx_0^*$. This element is equal to $(\id\otimes\e^*)(-Nx_0\otimes 1)$
and $-Nx_0\otimes 1=N(-x_0\otimes 1)$. Finally
$(\e\otimes\id)(-x_0\otimes 1)$ is equal to the class of~$-1$.
\enddemo 

4.3.\quad  We recall that for spectra $X$ and $Y$,
there are natural maps
\begin{equation}\label{wedge}
\begin{array}{l}
 \wedge  \:\pi_sX\otimes\pi_tY \to \pi_{s+t}(X\wedge Y), \\
 \vee  \:\pi_{s+t}F(X,Y) \to \Hom(\pi_{-s}X,\pi_tY),  
\end{array} \speqnu{4.3.1}
\end{equation}
where $\wedge$ is the external product and $\vee$ is the adjoint of
the composite
$$\pi_{s+t}F(X,Y)\otimes\pi_{-s}X\xto{\wedge}\pi_t(F(X,Y)\wedge
X)\xto{\ev}\pi_tY.$$
Let $X$ be a $G$-CW-spectrum with an increasing filtration
$\{X_s\}$ by sub-$G$-CW-spectra. Then the exact couple
$$D_{s-1,t+1}\xto{i}D_{s,t}\xto{j}E_{s,t}\xto{\partial}D_{s-1,t}$$
with
\begin{equation}\label{exactcouple}
\begin{array}{l}
 D_{s,t}(X)=\pi_{s+t}((X_s)^G)\\
 E_{s,t}(X) =\pi_{s+t}((X_s/X_{s-1})^G) 
\end{array} \speqnu{4.3.2}
\end{equation}
gives rise to a spectral sequence which abuts the homotopy groups of
$X^G$. The spectral sequence converges conditionally in the sense of
\cite[Def.\ 5.10]{boardman}, provided that $\bigcup X_s=X$ and
$\hl{}{(X_s)^G}$ is contractible.

If $X$ and $X'$ are two $G$-CW-spectra with such filtrations, we give
the smash product $X\wedge X'$ the usual product filtration
$$(X\wedge X')_n=\bigcup_{s+s'=n}X_s\wedge X'_{s'}.$$
with filtration quotients
$$(X\wedge X')_n/(X\wedge X')_{n-1}
=\bigvee_{s+s'=n}X_s/X_{s-1}\wedge X_{s'}/X_{s'-1}.$$
The external product \eqref{wedge} and the inclusions
\begin{eqnarray*}
X_s\wedge X'_{s'}&\to&(X\wedge X')_{s+s'},\\
X_s/X_{s-1}\wedge X'_{s'}/X'_{s'-1}&\to &
(X\wedge X')_{s+s'}/(X\wedge X')_{s+s'-1} 
\end{eqnarray*}
then give rise to pairings
\begin{eqnarray*}
&&D_{s,t}(X)\otimes D_{s',t'}(X')\to D_{s+s',t+t'}(X\wedge X'),\\
&&E_{s,t}(X)\otimes E_{s',t'}(X')\to E_{s+s',t+t'}(X\wedge X'). 
\end{eqnarray*}
These, in turn, give rise to an external pairing of the associated
spectral sequences, that is, pairings
$$E^r_{s,t}(X)\otimes E^r_{s',t'}(X') \to
E^r_{s+s',t+t'}(X\wedge X'),$$
for all $r\geq1$, which satisfies the Leibnitz rule
$$d^r(xx')=d^rxx'+(-1)^{|x|}xd^rx'.$$
Here $|x|$ is the total degree of $x$. A filtration-preserving product
map $X\wedge X\to X$ induces a map of the associated spectral
sequences which, pre-composed by the external product, give  an
internal product on the spectral sequence $E^*(X)$. The differentials
act as derivations for this product, and if the product on $X$ is
associative, commutative or unital, the same holds for the internal
product in the spectral sequence. Commutativity in the spectral
sequence is up to the usual sign.

Let $G$ be a finite group and let $E$ be a free
contractible $ G$-CW-complex. Let $\tilde E$ be the mapping cone of the
projection $\pr\:E_+\to S^0$ which collapses $E$ to the nonbase point
of $S^0$. The associated suspension-$G$-CW-spectra (we make
no change in notation) form a distinguished triangle 
$$E_+\xto{\pr}S^0\to\tilde E\xto{\partial}\Sigma E_+.$$
Let $P$ and $\tilde P$ be the cellular complexes of $E_+$ and
$\tilde E$ with coefficients in a commutative ring $k$. We then have a
distinguished triangle
$$P\xto{\pr_*}k\to\tilde P\to\Sigma P$$
in the category of chain complexes.

The Tate spectrum of a $G$-spectrum $T$ is defined by
$$\tate(G,T)=(\tilde E\wedge\Gamma F(E_+,T))^G,$$
where $\Gamma X\xto{\sim}X$ is a functorial $G$-CW-substitute. If
$T$ and $T'$ are two\break $G$-spectra, we define a pairing
\begin{equation}\label{tatespprod}
\tate(G,T)\wedge\tate(G,T')\to\tate(G,T\wedge T') \speqnu{4.3.3}
\end{equation}
as follows. By elementary obstruction theory, there are cellular
$G$-homotopy equivalences $E_+\to E_+\wedge E_+$ and $\tilde
E\wedge\tilde E\to\tilde E$ compatible with the canonical isomorphisms
$S^0\to S^0\wedge S^0$ and $S^0\wedge S^0\to S^0$, respectively, and
any two such equivalences are $G$-homotopic. The pairing then is given
by
\begin{eqnarray*}
(\tilde E\wedge\Gamma F(E_+,T))^G&\hsm
\wedge\hsm&(\tilde E\wedge\Gamma F(E_+,T))^G
\to(\tilde E\wedge\tilde E\wedge\Gamma F(E_+\wedge E_+,T\wedge T'))^G\\
&\hsm\to\hsm&(\tilde E\wedge\Gamma F(E_+,T\wedge T'))^G, 
\end{eqnarray*}
where the first map is the canonical map and the second is induced
from the chosen $G$-equivalences. If $T$ is a $G$-ring spectrum, the
composition of the external product with the map of Tate spectra
induced from the product map on $T$, makes $\tate(G,T)$ a ring
spectrum. This ring spectrum is   associative, commutative or unital
if the $G$-ring spectrum $T$ is associative, commutative or unital,
respectively.

The CW-filtrations of $E$ and $\tilde E$ give rise to a double
filtration of the Tate spectrum. In more detail, we define
\begin{eqnarray*}
X_{r,s}&\hsm =\hsm &\tilde E_r\wedge\Gamma F(E/E_{-s-1},T),\\
Y_{r,s}&\hsm =\hsm &\tilde E_r/\tilde E_{r-1}\wedge\Gamma F(E/E_{-s-1},T), \\
Z_{r,s}&\hsm =\hsm &\tilde E_r\wedge\Gamma F(E_{-s}/E_{-s-1},T)), \\
W_{r,s}&\hsm =\hsm &\tilde E_r/\tilde E_{r-1}\wedge\Gamma F(E_{-s}/E_{-s-1},T)). 
\end{eqnarray*}
To get an honest filtration by sub-$G$-CW-spectra, we let
$$\bar X_{r,s}=\hcl{}{X_{r',s'}},$$
where the homotopy colimit runs over all $0\leq r'\leq r$ and $s'\leq
s\leq 0$. There is a canonical homotopy equivalence
$\bar X_{r,s}\xto{\sim}X_{r,s}$ and $\bar X_{r,s}$ is a
sub-$G$-CW-spectrum of the $G$-CW-spectrum $\bar X=\bar X_{\infty,0}$.
We also let
\begin{eqnarray*}
\bar Y_{r,s}&\hsm =\hsm &\bar X_{r,s}/\bar X_{r-1,s}\\
\bar Z_{r,s}&\hsm =\hsm &\bar X_{r,s}/\bar X_{r,s-1}\\
\bar W_{r,s}&\hsm =\hsm &\bar X_{r,s}/\bar X_{r-1,s}\cup\bar X_{r,s-1} 
\end{eqnarray*}
and define
$$\bar X_n=\bigcup_{r+s\leq n}\bar X_{r,s}\subset\bar X.$$
The exact couple~\eqref{exactcouple} associated with the filtration
$\{\bar X_n\}$ gives rise to a conditionally convergent spectral
sequence
$$\hat E^*(G,T)=E^*(\bar X) \Rightarrow \pi_*(\tate(G,T)).$$

\specialnumber{4.3.4} \proclaim{Lemma}\label{E^1}
There is a canonical isomorphism of complexes
$$\hat E^1_{*,t}(G,T)\cong(\tilde P\otimes\Hom(P,\pi_tT))^G$$
and hence $\hat E^2_{s,t}(G,T)\cong\hat H^s(G,\pi_tT)$.
\endproclaim 

\demo{Proof} The inclusions $\bar X_{r,s}\to\bar X_{r+s}$ induce an
isomorphism
$$\bigvee_{r+s=n}\bar W_{r,s}\xto{\sim}\bar X_n/\bar X_{n-1}$$
such that the boundary map
$$\bar X_n/\bar X_{n-1}\to\Sigma\bar X_{n-1}
\to\Sigma(\bar X_{n-1}/\bar X_{n-2})$$
maps the summand $\bar W_{r,s}$ to the summands $\Sigma\bar W_{r-1,s}$
and $\Sigma\bar W_{r,s-1}$ by the maps
\begin{eqnarray*}
\partial'&\hsm\:\hsm&\bar W_{r,s}\to\Sigma\bar Y_{r,s-1}\to\Sigma\bar W_{r,s-1},\\
\partial''&\hsm\:\hsm&\bar W_{r,s}\to\Sigma\bar Z_{r-1,s}\to\Sigma\bar W_{r-1,s}, 
\end{eqnarray*}
respectively. We identify
$$\pi_{r+s+t}((\bar W_{r,s})^G) \cong
(\tilde P_r\otimes\Hom(P_{-s},\pi_tT))^G$$
as follows: If $X$ and $Y$ are two $G$-spectra, we have the canonical
map
$$\pi_*((X\wedge Y)^G)\to(\pi_*(X\wedge Y))^G.$$
This is an isomorphism, for example, if $X$ is a wedge of free
$G$-cells. The desired isomorphism is the composition of the inverse
of this map with $X=\tilde E_r/\tilde E_{r-1}$ and
$Y=\Gamma F(E_{-s}/E_{-s-1},T)$ and the map of $G$-fixed sets
induced by
\begin{eqnarray*}
&&\hskip-36pt \pi_{r+s+t}(\tilde E_r/\tilde E_{r-1}\wedge\Gamma
F(E_{-s}/E_{-s-1},T)) \\
&&\begin{array}{c}\scr{\wedge}\\[-8pt]
\longleftarrow\\[-8pt]{\scriptstyle\sim}\end{array}\pi_r(\tilde E_r/\tilde E_{r-1})\otimes
\pi_{s+t}\Gamma F(E_{-s}/E_{-s-1},T) \\
&&  \hskip6pt  \xto{\sim}\pi_r  
(\tilde E_r/\tilde E_{r-1})\otimes
\pi_{s+t}F(E_{-s}/E_{-s-1},T) \\
&&\begin{array}{c}\scr{h\otimes\vee}\\[-8pt] \srar\\[-8pt]\scr{\sim}\end{array} H_r(\tilde
E_r/\tilde E_{r-1})\otimes\Hom(\pi_{-s}(E_{-s}/E_{-s-1}),\pi_tT)\\
&& \begin{array}{c} \scr{1\otimes h^*}\\[-8pt] \slar\\[-8pt] \scr{\sim}\end{array} H_r(\tilde
E_r/\tilde E_{r-1})
\otimes\Hom(H_{-s}(E_{-s}/E_{-s-1}),\pi_tT). 
\end{eqnarray*}
Here $h$ is the Hurewitz homomorphism. One readily shows that
under this identification, $\pi_*(\partial')$ and $\pi_*(\partial'')$
correspond to the differentials in the algebraic double complex.
\enddemo  

The pairing \eqref{tatespprod} induces a pairing $\bar X(T)\wedge\bar
X(T')\to\bar X(T\wedge T')$, and since the equivalences $E_+\to
E_+\wedge E_+$ and $\tilde E\wedge\tilde E\to\tilde E$ were chosen to be
cellular, this pairing preserves the filtration by the sub-CW-spectra
$\{\bar X_n\}$. Hence, we get an induced pairing of the associated
spectral sequences.

\specialnumber{4.3.5}\proclaim{Proposition}
 Let $T$ and $T'$ be two $G$\/{\rm -}\/spectra. Then the pairing of
Tate spectra {\rm \eqref{tatespprod}} induces a pairing of the associated
spectral sequences. On $E^2$\/{\rm -}\/terms{\rm ,} this pairing corresponds to the
pairing on Tate cohomology
$$\hat H^*(G,\pi_*T)\otimes\hat H^*(G,\pi_*T')
\to \hat H^*(G,\pi_*(T\wedge T'))$$
under the isomorphism of Lemma~{\rm \ref{E^1}}. In particular{\rm ,} if $T$ is an
associative $G$\/{\rm -}\/ring spectrum{\rm ,} then $E^2\cong\hat H^*(G,\pi_*T)$ as a
bi\/{\rm -}\/graded ring.
\endproclaim

\demo{Proof} The equivalences $E_+\to E_+\wedge E_+$ and $\tilde
E\wedge\tilde E\to\tilde E$ induces chain maps $P\to P\otimes P$ and
$\tilde P\otimes\tilde P\to\tilde P$ which lift the canonical maps
$k\to k\otimes k$ and $k\otimes k\to k$, respectively. Now suppose $T$
and $T'$ are two $G$-spectra and consider the spectral sequences
corresponding to the filtrations $\{(\bar X(T)\wedge\bar X(T'))_n\}$
and $\{\bar X(T\wedge T')_n\}$. An argument analogous to the proof of
  Lemma~\ref{E^1} identifies the $E^1$-terms of the associated
spectral sequences with the complexes
$$(\tilde P\otimes\Hom(P,\pi_*T)\otimes\tilde
P\otimes\Hom(P,\pi_*T'))^G$$
and
$$(\tilde P\otimes\Hom(P,\pi_*(T\wedge T')))^G,$$
respectively. We claim that under these identifications, the pairing
$$\bar X(T)\wedge\bar X(T')\to\bar X(T\wedge T')$$
corresponds to the composition
\begin{eqnarray*}
(\tilde P&\hsm\otimes\hsm&\Hom(P,\pi_*T))^G
\otimes (\tilde P\otimes\Hom(P,\pi_*T'))^G \\[2pt]
{}&\hsm \to\hsm& (\tilde P\otimes\tilde P\otimes
\Hom(P\otimes P,\pi_*T\otimes\pi_*T'))^G
\to (\tilde P\otimes\Hom(P,\pi_*(T\otimes T')))^G,  
\end{eqnarray*}
where the first map is the canonical map of chain complexes (which involves
sign changes) and the second map is induced from the maps $P\to
P\otimes P$ and $\tilde P\otimes\tilde P\to\tilde P$ and from the
exterior product \eqref{wedge}. This is straightforward to
check. Similarly, under the isomorphism of Lemma~\ref{E^1} and the
analogous isomorphism above, the external pairing corresponds to the
canonical map (no sign changes)
\begin{eqnarray*}
(\tilde P\otimes\Hom(P,\pi_*T))^G&\hsm\otimes\hsm&(\tilde
P\otimes\Hom(P,\pi_*T'))^G\\[2pt]
{}&\hsm\to\hsm&(\tilde P\otimes\Hom(P,\pi_*T)\otimes\tilde
P\otimes\Hom(P,\pi_*T'))^G. 
\end{eqnarray*}
But this was our definition of the pairing in Tate cohomology;
see~(4.1.5).
\enddemo

\specialnumber{4.3.6}\numbereddemo{{R}emark} We show that the spectral sequence $\hat E^*(G,T)$
considered here is canonically isomorphic to the spectral sequence
obtained from Greenlees' $\Z$-graded `filtration' of $\tilde E$,
\cite{greenlees}, \cite{greenleesmay}. This is the sequence of
$G$-CW-spectra,
$$\cdots\to\tilde E_{r-1}\to\tilde E_r\to\tilde E_{r+1}\to\cdots,$$
where, for $r\geq 0$, $\tilde E_r$ is the suspension $G$-spectrum of
the $r$-skeleton of $\tilde E$, and for $r<0$, 
$\tilde E_r$ is the dual $D(\tilde E_{-r})=\Gamma F(\tilde
E_{-r},S^0)$. In particular, $\tilde E_0=S^0$ is the sphere
$G$-spectrum. The maps $\tilde E_{r-1}\to\tilde E_r$ are induced from
the canonical inclusions, and for $r=0$, from the canonical map
$D(S^0)\xto{\sim}S^0$. In the definition of the $G$-CW-spectra
$\bar X_{r,s}$ and $\bar X_n$, we now may vary $r$ over all
integers. Let $\bar X_{r,s}'$ and $\bar X_n'$ denote the
$G$-CW-spectra so obtained. Then, for $r\geq 0$, the canonical
inclusion $\bar X_{r,s}\xto{\sim}\bar X_{r,s}'$ is a homotopy
equivalence. We have maps of filtrations
$$\{\bar X_n\}_{n\in\Z}\to\{\bar X_n'\}_{n\in\Z}\leftarrow
\{\bar X_{r,0}'\}_{r\in\Z},$$
and the filtration on the right is Greenlees' filtration. We show that
both maps induce isomorphisms of the $E^2$-terms of the associated
spectral sequences. In order to identify the $E^1$-terms, let
$\e\:\hat P\to k$ be the complete resolution, where
$$(\Sigma\hat P)_s=H_s(\tilde E_s\cup C\tilde E_{s-1};k)$$
with differential
$$H_s(\tilde E_s\cup C\tilde E_{s-1}) \xto{\partial_*}
H_s(\Sigma\tilde E_{s-1}) \begin{array}{c}\scr{\susp}\\[-8pt] \slar\\[-8pt]  \scr{\sim}\end{array}
H_{s-1}(\tilde E_{s-1}) \xto{i_*}
H_{s-1}(E_{s-1}\cup C\tilde E_{s-2})$$
and with structure map
$$\e\:\hat P_0=H_1(\tilde E_1\cup C\tilde E_0) \xto{\partial_*}
H_1(\Sigma E_0) \begin{array}{c}\scr{\susp}\\[-8pt] \slar\\[-8pt]  \scr{\sim}\end{array} H_0(E_0)=k.$$
The map of distinguished triangles
$$\begin{array}{ccccccc}
{P} &\shtck{\e} & 
{k}\quad& \longrightarrow &
{\tilde P} &\longrightarrow&
{\Sigma P} \\[4pt]
 \big\Vert&&\scs{\e^*}&&\big\downarrow&&\big\Vert\\[4pt]
{P} &\longrightarrow&
{\Sigma P^-} &\longrightarrow &
{\Sigma\hat P}&\longrightarrow &
{\Sigma P} \end{array}
$$
defines a quasi-isomorphism of the mapping cones of the two middle
vertical maps. (See remark~\ref{completeres} for the definition of the
lower triangle.) Now an argument similar to the proof of
Lemma~\ref{E^1} identifies the maps of $E^1$-terms induced from the
above maps of filtrations with the canonical maps
$$(\tilde P\otimes\Hom(P,M))^G\to(\Sigma\hat P\otimes\Hom(P,M))^G
\leftarrow(\Sigma\hat P\otimes\Hom(k,M))^G.$$
Finally, an argument similar to the proof of Lemma~\ref{grphomology}
shows that both maps are quasi-isomorphisms.
\enddemo

4.4.\quad Again let $C$ be a cyclic group of order $r$
and let $g$ be a generator. As our model for $E$, we choose the unit
sphere
$$E=S(\C^\infty),$$
where the generator $g$ acts on $\C$ by multiplication by $e^{2\pi
i/r}$. We give $E$ the usual $C$-CW-structure with one free cell in
each dimension. The skeletons are
$$E_n=\left\{ \begin{array}{ll}
S(\C^d)&\hbox{$n=2d-1$}\\
S(\C^d)*(C\cdot1)&\hbox{$n=2d$},\\
\end{array}\right.$$
where in the latter case, we identify the join with its image under
the canonical homeomorphism $S(\C^n)*S(\C)\cong S(\C^n\oplus\C)$. The
attaching maps
$$\alpha_n\:D^n\times C\to E_n$$
are defined in even dimensions by the composite
$$D^{2d}\times C\xto{\xi}D(\C^d)\times C\xto{\pi}S(\C^d)*(C\cdot1),$$
where $\xi(z,g^s)=(g^s\cdot z,g^s)$ and $\pi$ is the canonical
projection. We define 
$$\alpha_1(x,g^s)=g^s\cdot e^{\pi i(x+1)/r}$$
and let $\alpha_{2d+1}$ be the composite
$$D^{2d}\times D^1\times C\xto{\xi}D(\C^d)\times D^1\times
C\xto{1\times\alpha_1}D(\C^d)\times S(\C)\xto{\pi}S(\C^d)*S(\C).$$

We give $D(\C^d)$ the complex orientation and $D^1=D(\R)=[-1,1]$ the
standard orientation from $-1$ to $1$. We may then identify the
cellular complex of $E$ with the standard complex $W$ by the
isomorphism
$$W\xto{\sim}C_*(E;k)$$
which maps the generator $x_n\in W_n$ to the image of the fundamental
class under the composite
$$H_n(D^n,S^{n-1})\xto{\iota_0}H_n(D^n\times C,S^{n-1}\times C)
\xto{\alpha_n}H_n(E_n,E_{n-1}).$$
Here $\iota_0\:D^n\to D^n\times C$ maps $z$ to $(z,1)$.

The $C$-CW-structure on $E$ induces one on $\tilde E$ and the
isomorphism above induces an isomorphism of chain complexes
$$\tilde W\xto{\sim}\tilde C_*(\tilde E;k).$$
We identify $\tilde E$ with $S^{\C^\infty}$ by the homeomorphism
$$CS(\C^\infty)_+\mathbold{\cup} S^0\xto{\sim}D(\C^\infty)/S(\C^\infty)$$
which maps $t\wedge z\mapsto tz$. Note that under this homeomorphism,
the orientation of the cells in $\tilde E$ corresponds \pagebreak to the complex
orientation of $S^{\C^\infty}$. In particular, the composite
$$H_2(S^\C)\stackrel{\sim}{\longleftarrow}H_2(\tilde E_2)
\xto{\pr_*}H_2(\tilde E_2,\tilde E_1)\stackrel{\sim}{\longleftarrow}\tilde W_2$$
maps the fundamental class $[S^\C]$ to the class $Ny_2$.

Let $C\subset{\Bbb T}$ be the subgroup of order $r$. We give
${\Bbb T}$ the $C$-CW-structure of $S(\C)=E_1$. Then the
multiplication is cellular, and hence, the cellular complex
$$\Lambda=C_*({\Bbb T};k)$$
is naturally a differential graded Hopf algebra with unit $1=x_0$. The
differential maps $x_1$ to $(g-1)\cdot x_0$, $x_1$ is primitive, the
coproduct on $g$ is $g\otimes g$, and the antipode is given by
$c(x_1)=-x_1$. We note that $x_1$ represents the fundamental class
$[{\Bbb T}]$. The $C$-action on $E=S(\C^{\infty})$ naturally extends
to a ${\Bbb T}$-action, and the action map
$$\mu\:{\Bbb T}\times E \to E$$
is cellular. The induced action on $\tilde E$,
$$\tilde\mu\:{\Bbb T}_+\wedge\tilde E={\Bbb T}_+\wedge
C_{\pr}\xto{\rho}C_{{\Bbb T}_+\wedge\pr}\xto{C_\mu}C_{\pr}=\tilde
E,$$
again is cellular. The induced left $\Lambda$-module structures on
the cellular complexes $W$ and $\tilde W$ are given by
$$x_1\cdot x_s=\left\{ \begin{array}{ll}
x_{s+1}&\hbox{$s$ even}\\
0&\hbox{$s$ odd}\\
\end{array}\right.,\hskip6mm
x_1\cdot y_s=\left\{ \begin{array}{ll}
0&\hbox{$s$ even}\\
-y_{s+1}&\hbox{$s$ odd.}\\
\end{array}\right.$$

Let $T$ be a ${\Bbb T}$-spectrum and let $\bar X=\bar X(T)$ be the
filtered ${\Bbb T}$-CW-spectrum, which gives rise to the spectral
sequence $\hat E^*(C,T)$. We give ${\Bbb T}/C$ the skeleton
filtration such that $\Lambda_C=C_*({\Bbb T}/C;k)$. Then the
${\Bbb T}$-actions on $E$, $\tilde E$, and $T$ induce a filtration-preserving map
$$\omega\:{\Bbb T}/C_+\wedge\bar X^C\to\bar X^C.$$
An argument similar to the proof of Lemma~\ref{E^1} identifies the
induced map of $E^1$-terms of the associated spectral sequences with the map
$$\Lambda_C\otimes(\tilde W\otimes\Hom(W,\pi_*T))^C
\to(\tilde W\otimes\Hom(W,\pi_*T))^C$$
given by the composite
\begin{eqnarray*}
\Lambda_C & \hsm\otimes\hsm&(\tilde W\otimes\Hom(W,\pi_*T))^C
\begin{array}{c}\scr{N\otimes\id}\\[-8pt] \srar\\[-8pt]\scr{\sim}\end{array}
\Lambda^C\otimes(\tilde W\otimes\Hom(W,\pi_*T))^C \\
{} &\hsm \to\hsm& (\Lambda\otimes\tilde W\otimes\Hom(W,\pi_*T))^C
\xto{\omega_*} (\tilde W\otimes\Hom(W,\pi_*T))^C.
\end{eqnarray*}

\specialnumber{4.4.1}\proclaim{Proposition} \label{connesandtate}
Let $T$ be a $\,{\Bbb T}$\/{\rm -}\/spectrum.
Then $\hat E^*(C,T)$ is a spectral sequence of left
$\Lambda_C$-modules. Moreover{\rm ,} if the class $a\in\pi_*(\tate(C,T))$ is 
represented by the infinite cycle $z\in E^1_{s,t}${\rm ,} and if $x_1\cdot
z\in E^1_{s+1,t}$ is nonzero{\rm ,} then $x_1\cdot z$ is an infinite cycle
and represents the class of $da\in\pi_*(\tate(C,T))$. 
\endproclaim

Let $k$ be a perfect field of odd characteristic $p$ and let $T(k)$ be
the topological Hochschild spectrum of $k$. Then as a differential
graded $k$-algebra,
$$\pi_*(T(k),\Z/p)=\Lambda\{\e\}\otimes S\{\sigma\}$$
with the classes $\e\in\pi_1(T(k),\Z/p)$ and
$\sigma\in\pi_2(T(k),\Z/p)$ characterized by $\beta(\e)=1$ and
$d(\e)=\sigma$. The Tate spectral sequence takes the form
$$\hat E^2(C_p, M_p\wedge T(k))=\Lambda\{u_1,\e\}\otimes S\{t^{\pm 1},\sigma\}
\Rightarrow\pi_*(\hat{\Bbb H}(C_p, T(k)),\Z/p),$$
where $u_1=u$ and $t$ are the generators of $\hat H^*(C_p,k)$ from
Lemma~\ref{tateparc}. The nonzero differentials are multiplicatively
generated from $d^2(\e)=t\sigma$.

\specialnumber{4.4.2}\proclaim{{C}orollary}\label{gammahatk}
The image of the classes $\e$ and
$\sigma$ under the map induced from
$$\hat\Gamma_k\:T(k)\to\tate(C_p,T(k))$$
are represented by the infinite cycles $ut^{-1}$ and $t^{-1}${\rm ,}
respectively.
\endproclaim

\demo{Proof} We recall from Section~1.1 that $\hat\Gamma_k$ is
defined as the composite
$$T(k)\begin{array}{c} \scr{r}\\[-8pt] \leftarrow\\[-8pt] \scr{\sim}\end{array}\rho_{C_p}^*(\tilde
E\wedge T)^{C_p}
\to\rho_{C_p}^*(\tilde E\wedge F(E_+,T))^{C_p}.$$
Both maps are ${\Bbb T}$-equivariant, so $\hat\Gamma$ commutes with
Connes' operator. It also commutes with the Bockstein operator. Hence,
it suffices to show that $ut^{-1}\otimes 1$ represents the unique class whose
Bockstein is the multiplicative unit $1$, and that $t^{-1}\otimes 1$
represents the image under Connes' operator of this class. To this
end, we recall from Lemma~\ref{tateparc} that the classes $1$,
$ut^{-1}$ and $t^{-1}$ in $\hat H^*(C_p,\Fp)$ are represented by the
elements $y_0\otimes Nx_0^*$, $-Ny_1\otimes Nx_0^*$ and $Ny_2\otimes
Nx_0^*$, respectively. We recall from Section~2.1 above
that the Bockstein
$$\beta\:\hat H^*(C_p,\Fp) \to \hat H^{*+1}(C_p,\Z)$$
is equal to the connecting homomorphism associated with the exact
sequence
$$0 \to (\tilde W\otimes\Hom(W,\Z))^{C_p} \xto{p}
(\tilde W\otimes\Hom(W,\Z))^{C_p} \xto{\pr}
(\tilde W\otimes\Hom(W,\Fp))^{C_p} \to 0.$$
This takes $-Ny_1\otimes Nx_0^*$ to $y_0\otimes Nx_0^*$, and hence
$\beta(ut^{-1})=1$. Next,
$$x_1\cdot(-Ny_1\otimes Nx_0^*)=-N(x_1\cdot y_1)\otimes Nx_0^*
+Ny_1\otimes N(x_1\cdot x_0^*)=Ny_2\otimes Nx_0^*,$$
and so by Proposition~\ref{connesandtate}, the image under Connes'
operator of the class represented by $ut^{-1}\otimes 1$ is represented
by $t^{-1}\otimes 1$.
\enddemo

Finally, for a ${\Bbb T}$-spectrum $T$, we will also consider the
${\Bbb T}$-Tate spectrum
$$\tate({\Bbb T},T)=(\tilde E\wedge\Gamma F(E_+,T))^{{\Bbb T}},$$
where again $E=S(\C^{\infty})$. The filtration of $E$ by the odd
skeletons $E_{2d-1}$, $d\geq 1$, and the associated filtration of
$\tilde E$ both are preserved by the ${\Bbb T}$-action. The
induced filtration of the Tate spectrum gives a conditionally
convergent spectral sequence
$$\hat E^2({\Bbb T},T)=S\{t^{\pm 1}\}\otimes\pi_*(T)
\Rightarrow\pi_*(\tate({\Bbb T},T)),$$
with the generator $t$ in bi-degree $(-2,0)$. Let
$C\subset{\Bbb T}$ be a subgroup. Then the canonical inclusion
$\tate({\Bbb T},T)\to\tate(C,T)$ induces a map of spectral sequences
$$\hat E^*({\Bbb T},T)\to\hat E^*(C,T).$$
If the order of $C$ is odd and annihilates $\pi_*(T)$ then, on
$E^2$-terms, this map is the canonical inclusion which maps $t$ 
to the generator $t$ of Lemma~\ref{tateparc}. This is the case, for
instance, if $T=M_p\wedge T(A|K)$ and $C=C_{p^n}$.

\specialnumber{4.4.3}\proclaim{Proposition} \label{d2=connes}
Let $T$ be a ${\Bbb T}$\/{\rm -}\/spectrum and
let $C\subset{\Bbb T}$ be a subgroup whose order $r$ is odd and
annihilates $\pi_*(T)$. Then the $d^2$\/{\rm -}\/differential in
$$\hat E^2(C,T)=\hat H^*(C,\Z/r)\otimes\pi_*(T)
\Rightarrow\pi_*(\tate(C,T))$$
is given by $d^2(\gamma\otimes\tau)=\gamma t\otimes d\tau${\rm ,} where
$d$ is Connes\/{\rm '} operator.
\endproclaim

\demo{Proof} It was proved in \cite[Lemma  1.4.2]{h} that  
in the ${\Bbb T}$-Tate spectral sequence, the $d^2$-differential is
given by the formula of the statement. Moreover, every $C$-spectrum
$T$ is a module $C$-spectrum over the sphere $C$-spectrum
$S^0$. Hence, it suffices to show that the class $u$ is a $d^2$-cycle
in the spectral sequence $\hat E^*(C,S^0)$. But $\pi_1(S^0,\Z/r)$
vanishes since $r$ is odd.
\enddemo

\section{The Tate spectral sequence for $T(A|K)$} \label{cartier}

5.1.\quad The Tate spectral sequence $\hat E^*(C_{p^n},M_p\wedge
T(A|K))$ is a spectral sequence of bi-graded $k$-algebras
in a canonical way, which we now explain. (We will abuse notation and
write $\hat E^*(C_{p^n},T(A|K))$ for this spectral sequence.) For every $C_{p^n}$-ring spectrum
$T$, $\tate(C_{p^n},T)$ is a $T^{C_{p^n}}$-algebra spectrum, and the
Tate spectral sequence is one of bi-graded
$\bar{\pi}_*(T^{C_{p^n}})$-algebras,
$$\hat E^2(C_{p^n},T) = \hat H^{-*}(C_{p^n},F^{n*}\bar\pi_*(T))
\Rightarrow \bar\pi_*(\tate(C_{p^n},T)),$$
where $F^n\:T^{C_{p^n}} \to T$ is the natural inclusion. Here 
$F^{n*}\bar\pi_*T$ denotes the graded ring $\bar\pi_*T=\pi_*(T,\Z/p)$
considered as a $\bar\pi_*(T^{C_{p^n}})$-algebra via the ring
homomorphism induced by $F^n$. In the case at hand, we consider this a
spectral sequence of bi-graded $k$-algebras via the ring
homomorphism~\eqref{rho},
$$\rho_{n+1}\: k \to \bar W_{n+1}(A) = \bar\pi_0(T(A|K)^{C_{p^n}}).$$
We recall that $F^n\circ\rho_{n+1}=\rho_1\circ\varphi^n$, where
$\varphi\:k\to k$ is the Frobenius. The latter is an automorphism, by
our assumption that $k$ is perfect, and hence
$$\hat E^2(C_{p^n},T(A|K))=\Lambda\{u_n,d\log\pi_K\}\otimes
S\{\pi_K,\kappa,t^{\pm1}\}/(\pi_K^{e_K}),$$
where $u_n$ and $t$ are the canonical generators from Lemma
\ref{tateparc}. The differential structure of this spectral sequence is
evaluated in this section. We briefly outline the
argument.

The $d^2$-differential in $\hat E^*(C_{p^n},T(A|K))$ is given by
Proposition~\ref{d2=connes} in terms of Connes'
operator~\eqref{connes} on $\bar\pi_*T(A|K)$. Hence, by
Theorem~\ref{pisubast},
$$d^2\pi_K  = td\log\pi_K\cdot\pi_K, \hskip6mm
d^2\kappa  = td\log(-p)\cdot\kappa,$$
and we can use the equation $-p=\pi_K^{e_K}\theta_K(\pi_K)^{-1}$ to
express $d\log(-p)$ as a polynomial in $\pi_K$ times $d\log\pi_K$.
In Section~5.2, we replace $\kappa$ by a new generator
$\alpha_K$, defined as a certain linear combination of the elements
$\pi_K^i\kappa$ with $0\leq i<e_K$, which satisfies that
$d\alpha_K=e_Kd\log\pi_K\cdot\alpha_K$. In particular, $\alpha_K$ is a
$d^2$-cycle, if $p$ divides $e_K$. We also replace $t$ by a new
generator $\tau_K$ defined in a similar manner.

The key results that make it possible to completely evaluate the
spectral sequence are consequences of the map
$$\hat\Gamma_{A|K}\:  T(A|K)^{C_{p^{n-1}}} \to \tate(C_{p^n},T(A|K)),$$
and of the unit map of the ring spectrum on the right,
$$\ell\: S^0 \to \tate(C_{p^n},T(A|K)).$$
We show in Section~5.3 that for $n<v_p(e_K)$, 
$\pi_K^{p^n}$ and $-\tau_K\alpha_K$ are infinite cycles which
represent the classes $\hat\Gamma_{A|K}(\ul{\pi_K}_n)$ and
$\hat\Gamma_{A|K}(\ul{\pi_K}_n^{e_K/p^n})$, respectively. We also show
that $-\tau_K\alpha_K^p$ is always an infinite cycle which represents
the image by the unit map of the canonical generator
$v_1\in\bar\pi_{2p-2}(S^0)$. Given these infinite cycles together with
the value of the differentials on the $p$-powers of $\pi_K$, which we
examine by a universal example in Section~5.5, one can
evaluate the spectral sequence, if $n<v_p(e_K)$. The final part of the
argument consists of a somewhat complicated induction argument, which
we present at the end of Section~5.5. The key for this part
is naturality, going back and forth between $T(A|K)$ and $T(B|L)$ for
suitable ramified extensions $L/K$.

The handling of the spectral sequences is algebraically somewhat
complex. To ease the presentation we first consider in
Section~5.4 the case of $\hat E^*(C_p,T(A|K))$. This section
also contains the proof that the map $\hat\Gamma_{A|K}$ induces an
isomorphism of homotopy groups with $\Z/p$-coefficients in
nonnegative degrees.

\vglue12pt 5.2.\quad Let $L$ be a finite and totally ramified
extension of $K$, and let $B$ be the integral closure of $A$ in
$L$. Then $B$ is a complete discrete valuation ring with quotient
field $L$ and  residue field $k$. Let $\pi_K$ and $\pi_L$ be
uniformizers of $A$ and $B$, respectively. The minimal polynomial of
$\pi_L$ over $K$ has the form
$$\phi_{L/K}(x)=x^{e_{L/K}}+\pi_K\theta_{L/K}(x),$$
where $\theta_{L/K}(x)$ is a polynomial over $A$ of degree $<e_{L/K}$
and $\theta_{L/K}(0)\in A^\times$. Moreover, the canonical map
$$A[\pi_L]/(\phi_{L/K}(\pi_L))\xto{\sim}B$$
is an isomorphism. When $K=K_0$ is the quotient field of $W(k)$, we
will always use $\pi_{K_0}=p$ and  write $\theta_L(x)$
instead of $\theta_{L/K_0}(x)$.

\specialnumber{5.2.1} \proclaim{Lemma}\label{lubintate}
Suppose that $\mu_p\subset K$. Then a
choice of a generator $\zeta\in\mu_p$ and a uniformizer $\pi_K\in A$
determines a polynomial $u_K(x)\in W(k)[x]$ of degree $<e_K$ such that
$u_K(\pi_K)^{p-1}=\theta_K(\pi_K)$. Moreover{\rm ,} in $\omega_{(A,M)}^1${\rm ,}
$$d\log\zeta=-\pi_K^{e_K/(p-1)}u_K(\pi_K)^{-1}d\log(-p).$$
\endproclaim 

\demo{Proof} Consider the power series $f(x) = px+x^p$ and $g(x) =
(1+x)^p-1$ and recall from~\cite[\S3, Prop.\ 3]{cf} that there
exists a {\it unique} power series $\varphi(x)$ such that
$f(\varphi(x)) = \varphi(g(x))$ and  $\varphi(x) \equiv x$ modulo
$(x^2)$. Hence, if $\zeta\in\mu_p$ is a generator then
$\varphi(\zeta-1)$ is a $(p-1)$st root of $-p$. We define $u_K(x)$ to
be the unique polynomial of degree $<e_K$ such that
$$u_K(\pi_K)=\pi_K^{e_K/(p-1)}\varphi(\zeta-1)^{-1}.$$
To prove the second statement, we first note that
\begin{eqnarray*}
d\varphi(\zeta-1) & = &\varphi(\zeta-1)d\log\varphi(\zeta-1) \\
{} & =& \pi_K^{e_K/(p-1)}u_K(\pi_K)^{-1}\cdot (p-1)^{-1}d\log(-p) \\
{} & = &-\pi_K^{e_K/(p-1)}u_K(\pi_K)^{-1}\cdot d\log(-p),  
\end{eqnarray*}
where the last equality uses that $d\log(-p)$ is $p$-torsion. Hence,
it suffices to show that $d\varphi(\zeta-1)=d\log\zeta$. We may assume
that $K=\Qp(\mu_p)$, where as a uniformizer, we take $\pi_K=\zeta-1$.
Then $\omega_{(A,M)}^1$ is annihilated by $\pi_K^{p-1}$, and since
$d\varphi(\zeta-1)=\varphi'(\zeta-1)\zeta d\log\zeta$, it remains to
show that $\varphi'(x)\equiv(1+x)^{-1}$ modulo $(x^{p-1})$, or
equivalently, that $\varphi(x)\equiv\log(1+x)$ modulo $(x^p)$. But
this follows from the uniqueness of $\varphi(x)$ and from the
calculation in $\Zp[x]/(x^p)$:
\vglue12pt
\hfill ${\displaystyle \log(1+g(x))=\log((1+x)^p)=p\log(1+x)=f(\log(1+x)).}$\enddemo
 \vglue7pt

\specialnumber{5.2.2}\proclaim{Addendum} \label{lubintateadd}
Let $L/K$ be a finite and totally
ramified extension. Then the inclusion of valuation rings{\rm ,} 
$\iota\:A\to B$, maps
$$\iota(u_K(\pi_K)) =
(-\theta_{L/K}(\pi_L))^{-e_K/(p-1)}u_L(\pi_L).$$
\endproclaim  

\demo{Proof} We can write the
$\iota(\varphi(\zeta-1))=\varphi(\zeta-1)$ as
$$\iota(\pi_K^{e_K/(p-1)}u_K(\pi_K)^{-1})=
\pi_L^{e_L/(p-1)}u_L(\pi_L)^{-1}.$$
Since $\iota(\pi_K)=\theta_{L/K}(\pi_L)^{-1}\pi_L^{e_{L/K}}$, the left-hand side also is equal to
$$(-\theta_{L/K}(\pi_L)^{-1}\pi_L^{e_{L/K}})^{e_K/(p-1)}
\iota(u_K(\pi_K)^{-1}).$$
The formula follows since $e_{L/K}e_K=e_L$ and since
$\pi_L$ is a nonzero divisor.
\enddemo

Suppose that $\mu_p\subset K$. We choose a generator $\zeta\in\mu_p$
and a uniformizer $\pi_K\in A$ and let $u_K(x)$ be the polynomial from
Lemma~\ref{lubintate}. Let $\kappa\in\bar\pi_2T(A|K)$ be the unique
class with $\beta(\kappa)=d\log(-p)$ and define
$\alpha_K=u_K(\pi_K)^{-1}\kappa$.

\specialnumber{5.2.3}\proclaim{Proposition} \label{pi_*alpha}
As a differential graded $k$\/{\rm -}\/algebra
$$\bar\pi_*T(A|K) = \Lambda\{d\log\pi_K\}\otimes
S\{\alpha_K,\pi_K\}/(\pi_K^{e_K})$$
with $d\pi_K=\pi_K d\log\pi_K$ and 
$d\alpha_K=e_K\alpha_K d\log\pi_K$.
\endproclaim

\demo{Proof} It follows from Theorem~\ref{pisubast} and
Lemma~\ref{logdiff} that as a differential graded $k$-algebra
$$\bar\pi_*T(A|K)=\Lambda\{d\log\pi_K\}\otimes
S\{\kappa,\pi_K\}/(\pi_K^{e_K})$$
with the differential given by $d\pi_K=\pi_K d\log\pi_K$ and
$d\kappa=\kappa d\log(-p)$. Moreover, differentiating the equation
$-p=\pi_K^{e_K}\theta_K(\pi_K)^{-1}$, we find
$$d\log(-p)=(e_Kd\log\pi_K-d\log\theta_K(\pi_K)).$$
Finally, $\theta_K(\pi_K)=u_K(\pi_K)^{p-1}$, and hence
\begin{eqnarray*}
d(\alpha_K)&=&-u_K(\pi_K)^{-1}d\log u_K(\pi_K)\cdot\kappa+
u_K(\pi_K)^{-1}\cdot\kappa d\log(-p) \\
{}&=&-\alpha_K d\log u_K(\pi_K)+\alpha_K(e_Kd\log\pi_K-(p-1)d\log
u_K(\pi_K))\\
{}&=&e_K\alpha_K d\log\pi_K 
\end{eqnarray*}
as stated.
\enddemo

We recall the Bott element. Since $p$ is odd, the Bockstein is an
isomorphism,
$$\bar\pi_2(\Sigma^\infty B\mu_{p+}) \xto{\sim}
{}_p\pi_1(\Sigma^\infty B\mu_{p+}) \stackrel{\sim}{\longleftarrow} \mu_p,$$
and by definition, the Bott element $b=b_{\zeta}$ is the class on the
left which corresponds to the chosen generator $\zeta$ on the
right. The spectrum $\Sigma^\infty B\mu_{p+}$ is a ring spectrum and
the $(p-1)$st power $b^{p-1}$, which is independent of the choice of
generator, is the image by the unit map of a generator $v_1$ in
$\bar\pi_{2p-2}(S^0)$. If $\mu_p\subset K$, we have the maps of ring
spectra
$$\Sigma^\infty B\mu_{p+} \xto{\det} K(K)
\xto{\tr} T(A|K)^{C_{p^{n-1}}},$$
and let $b_n=b_{n,\zeta}$ be the image of the Bott element in
$\bar\pi_2(T(A|K)^{C_{p^{n-1}}})$. We note that
$\beta(b_n)=d\log_n\zeta$ and that since $\pi_2(T(A|K)^{C_{p^{n-1}}})$
is uniquely divisible, this equation characterizes $b_n$. In
particular, the calculation
$$\beta(b_1)=d\log\zeta=-\pi_K^{e_K/(p-1)}u_K(\pi_K)^{-1}d\log(-p)
=\beta(-\pi_K^{e_K/(p-1)}\alpha_K)$$
shows that
\begin{equation}\label{bottelement}
b_1=-\pi_K^{e_K/(p-1)}\alpha_K. \speqnu{5.2.4}
\end{equation}
The elements $b_n$ for $n>1$, however, are not well understood.

Let $L/K$ be a finite and totally ramified extension, and let
$\iota\:A\to B$ be the inclusion of valuation rings. Then the map
\begin{equation}\label{incl}
\iota_*\:\bar\pi_*T(A|K)\to\bar\pi_*T(B|L) \speqnu{5.2.5}
\end{equation}
is given by
\begin{eqnarray*}
\iota_*(\pi_K)&=&-\theta_{L/K}(\pi_L)^{-1}\pi_L^{e_{L/K}},\\
\iota_*(d\log\pi_K)&=&e_{L/K}d\log\pi_K-
d\log\theta_{L/K}(\pi_L),\\
\iota_*(\alpha_K)&=&(-\theta_{L/K}(\pi_L))^{e_K/(p-1)}\alpha_L. 
\end{eqnarray*}
The first two equalities follow  immediately from the definition of
$\theta_{L/K}(\pi_L)$, and the last equality follows from Addendum~\ref{lubintateadd}.

Let $f(x)\in k[x]$ and let $n$ be an
integer. We write $f^{(n)}(x)$ for the image of $f(x)$ under the
automorphism $\varphi^n[x]\:k[x]\to k[x]$, which applies $\varphi^n$
to the coefficients of a polynomial. If $R$ is a $k$-algebra and if
$\pi\in R$\, then, as usual, $f(\pi)$ denotes the image of $f(x)$ by
the unique $k$-algebra homomorphism $k[x]\to R$ which takes $x$ to
$\pi$. We note that $f^{(-n)}(\pi)\in\varphi^{n*}R$ and $f(\pi)\in R$
is the same element.

Suppose either $\mu_p\subset K$ or $K=K_0$. In the former case, let
$\pi_K$ be a uniformizer, let $\zeta\in\mu_p$ be a generator and let
$u_K(x)$ be the polynomial from Lemma~\ref{lubintate}. In the latter
case, let $u_{K_0}(x)=1$. Then as a bi-graded $k$-algebra,
\begin{equation}\label{E2term}
\hat E^2(C_{p^n},T(A|K))=\Lambda\{u_n,d\log\pi_K\}\otimes
S\{\pi_K,\alpha_K,\tau_K^{\pm1}\}/(\pi^{e_K}),\qquad \speqnu{5.2.6}
\end{equation}
with the new generators given by
$$\alpha_K = u_K^{(-n)}(\pi_K)^{-1}\kappa, \hskip6mm
\tau_K = u_K^{(-n)}(\pi_K)^p\,t.$$
We note the relations $\tau_K\alpha_K=\theta_K^{(-n)}(\pi_K)t\kappa$
and $\tau_K\alpha_K^p=t\kappa^p$.

It will be important to know how these new generators behave under
extensions. For integers $a,r,d$ with $0\leq r<e_K$ and $d\geq 0$,
we define
$$\{a,r,d\}_K=(pa-d)e_K/(p-1)+r.$$
If $\mu_p\subset K$ then $p-1$ divides $e_K$ such that $\{a,r,d\}_K$
is an integer. Let $L/K$ be a finite and totally ramified extension,
and let $\iota\:A\to B$ be the inclusion of valuation rings. Then
$\{a,e_{L/K}r,d\}_L=e_{L/K}\{a,r,d\}_K$ and
\begin{equation}\label{E2iota}
\iota_*\:\hat E^2(C_{p^n},T(A|K))\to \hat E^2(C_{p^n},T(B|L)), \speqnu{5.2.7}
\end{equation}
is given by
\begin{eqnarray*}
\iota_*(\tau_K^a\pi_K^r\alpha_K^d) & =&
(-\theta_{L/K}^{(-n)}(\pi_L))^{-\{a,r,d\}_K}
\tau_L^a\pi_L^{e_{L/K}r}\alpha_L^d,\\
\iota_*(d\log\pi_K)&=&(e_{L/K}-
\frac{\theta_{L/K}^{(-n)}{}'(\pi_L)\pi_L}{
\theta_{L/K}^{(-n)}(\pi_L)})d\log\pi_L. 
\end{eqnarray*}

\vglue12pt 5.3.\quad In this section, we produce a number of
infinite cycles in the spectral sequence $\hat
E^*(C_{p^n},T(A|K))$. This uses the maps of differential graded
$k$-algebras
$$\bar\pi_*T(A|K) \xleftarrow{\;j_*}
\bar\pi_*T(A) \xto{\;\rho_*}
\bar\pi_*T(A/p),$$
where the right-hand map is induced from the reduction. We evaluate
these maps assuming that $v_p(e_K)>0$. The left  map may be
identified with the map of graded $k$-algebras
$$j_*\:
\Lambda\{d\pi_K\}\otimes S\{\tilde\kappa,\pi_K\}/(\pi_K^{e_K})
\to \Lambda\{d\log\pi_K\}\otimes S\{\kappa,\pi_K\}/(\pi_K^{e_K}),$$
which takes $\pi_K$ to $\pi_K$, $d\pi_K$ to $\pi_Kd\log\pi_K$ and
$\tilde\kappa$ to $\kappa$. (See the discussion preceding
Theorem~\ref{pisubast}.) The group $\pi_2T(A)$ is uniquely divisible
so the Bockstein induces an isomorphism $\beta\:\bar\pi_2T(A)
\xto{\sim} {}_p\pi_1T(A)$ and the class $\tilde\kappa$ corresponds to
the generator
$$d\log(-p)
=-((e_K/p)\pi_K^{e_K-1}+\theta_K'(\pi_K))\theta_K(\pi_K)^{-1}d\pi_K$$
on the right. The differential graded $k$-algebra $\bar\pi_*T(A/p)$ is
evaluated in Proposition~\ref{TparkparPi} of the appendix. We refer to
{\it loc.~cit.}~for the notation.

\specialnumber{5.3.1}\proclaim{Proposition} \label{reductionmodp}
If $v_p(e_K)>0$ the map
$\rho_*\:\bar\pi_*T(A) \to \bar\pi_*T(A/p)$ may be identified with the
inclusion of differential graded $k$\/{\rm -}\/algebras
$$\rho_*\:\Lambda\{d\pi_K\}\otimes S\{\pi_K,\tilde\kappa\}/(\pi_K^{e_K})
\hookrightarrow
\Lambda\{d\bar\pi_K,\e\}\otimes S\{\sigma,\bar\pi_K\}/(\bar\pi_K^{e_K})
\otimes\Gamma\{\bar c_2\}$$
which takes $\pi_K$ to $\bar\pi_K$ and $\tilde\kappa$ to the class
$$\sigma-\theta_K(\bar\pi_K)^{-1}\bar c_2-
\e\cdot((e_K/p)\bar\pi_K^{e_K-1}+\theta_K'(\bar\pi_K))
\theta_K(\bar\pi_K)^{-1}d\bar\pi_K.$$
\endproclaim

\demo{Proof} Only the formula for $\rho_*(\tilde\kappa)$ requires
proof. Consider the diagram
$$\begin{array}{ccccccc}
{T(A)}&\stck{p} &
{T(A)}&\stck{i} &
{M_p\wedge T(A)}&\stck{\beta}& 
{\Sigma T(A)} \\[4pt]
 \scs{T(\rho)}&&\scs{T(\rho)} &&\scs{M_p\wedge T(\rho)}&& \scs{\Sigma T(\rho)}\\[4pt]
{T(A/p)}&\stck{p=0}  &
{T(A/p)}&\begin{array}{c}\scr{i} \\[-8pt] \srar\\[-10pt]\slar\\[-8pt] \scr{r}\end{array}&
{M_p\wedge T(A/p)}& \begin{array}{c}\scr{\beta} \\[-8pt] \srar\\[-10pt]\slar\\[-8pt]
\scr{s}\end{array}& {\Sigma T(A/p),}  \end{array}
 $$
with horizontal triangles, the lower triangle split by the maps
$r$ and $s$ of Section~2.1 above. It shows that
$$\rho_*(\tilde\kappa) = 
\e\cdot((\Sigma T(\rho))_*\circ\beta_*)(\tilde\kappa)
+ (i_*\circ r_*\circ (M_p\wedge T(\rho))_*)(\tilde\kappa).$$
The value of the first summand is easily determined from the diagram 
$$\begin{array}{ccccccc}
{\pi_2(M_p\wedge T(A))} &\stck{\beta_*}&
{\pi_2(\Sigma T(A))} & \begin{array}{c}\scr{\susp}\\[-8pt]{\slar}\\[-8pt]\scr{\sim}\end{array} &
{\pi_1(T(A))}&  \xleftarrow{\sim} &\Omega^1_A \\[4pt]
\scs{(M_p\wedge T(\rho))_*} &
&\scs{(\Sigma T(\rho))_*} &&\scs{T(\rho)_*}&&  \scs{\rho_*}\\[4pt]
{\pi_2(M_p\wedge T(A/p))} &\stck{\beta_*} &
{\pi_2(\Sigma T(A/p))} & \begin{array}{c}\scr{\susp}\\[-8pt]{\slar}\\[-8pt]\scr{\sim}\end{array}&
{\pi_1(T(A/p))} & \xleftarrow{\sim}&
{\Omega_{A/p}^1} \end{array}$$
and the formula for the Bockstein of $\tilde\kappa$ above: 
$$((\Sigma T(\rho))_*\circ\beta_*)(\tilde\kappa) = 
-((e_K/p)\bar\pi_K^{e_K-1}+\theta_K'(\bar\pi_K))
\theta_K(\bar\pi_K)^{-1}d\bar\pi_K.$$
It remains to show that
$$(r_*\circ(M_p\wedge T(\rho))_*)(\tilde\kappa) =
\sigma-\theta_K(\bar\pi_K)^{-1}\bar c_2.$$
We first show that the linearization
$l_*\:\pi_*T(A/p) \to \pi_*\HH(A/p)$ takes this class to
$-\theta_K(\bar\pi_K)^{-1}\bar c_2$. The following diagram
$$\begin{array}{ccccc}
{\pi_*(M_p\wedge T(A))} &\shtck{\rho_*} &
{\pi_*(M_p\wedge T(A/p))} &\shtck{r_*}  &
{\pi_*T(A/p)}  \\[4pt]
\scs{l_*} &&\scs{l_*}&&\scs{l_*}\\[4pt]
{\pi_*(M_p\wedge\HH(A))} &\shtck{\rho_*} &
{\pi_*(M_p\wedge\HH(A/p))} &\shtck{r_*} &
{\pi_*\HH(A/p)} \end{array}
$$
commutes and the composite of the lower horizontal maps is an
isomorphism. Let $c_1,c_2^{[d]}\in\pi_*(M_p\wedge\HH(A))$ be the
classes which correspond to $\bar c_1,\bar c_2^{[d]}\in\pi_*\HH(A/p)$
under this isomorphism. We claim that
$c_2=-\theta_K(\pi_K)\tilde\kappa$. The diagram
$$\begin{array}{ccc}
{\pi_2(M_p\wedge T(A))} &\begin{array}{c}\scr{\beta}\\[-8pt]\srar\\[-8pt]\scr{\sim}\end{array}&
{{}_p\pi_1T(A)} \\[4pt]
\scr{\sim}\scs{l_*} &&\scr{\sim}\scs{l_*} \\[4pt]
{\pi_2(M_p\wedge\HH(A))} &\begin{array}{c}\scr{\beta}\\[-8pt]\srar\\[-8pt]\scr{\sim}\end{array}&
{{}_p\pi_1\HH(A),} \end{array}
$$
where all maps are isomorphisms, shows that to prove this, it will
suffice to show that the lower Bockstein takes $c_2$ to
$((e/p)\pi_K^{e_K-1}+\theta_K'(\pi_K))d\pi_K$. This Bockstein, in
turn, may be identified with the connecting homomorphism in the
diagram
$$\begin{array}{cclclclcc}
{0}&\longrightarrow &
{A\cdot c_2} &\shtck{p}  &
{A\cdot c_2}& \shtck{\rho}  &
{A/p\cdot\bar c_2}&\longrightarrow  &
{0} \\[4pt]
&&\phantom{A\cdot}\scs{\phi_K'(\pi_K)}
&&\phantom{A\cdot} \scs{\phi_K'(\pi_K)}&&\phantom{A\cdot} \scs{\phi_K'(\bar\pi_K) =0}\\[4pt]
{0}& \longrightarrow &
{A\cdot c_2} &\shtck{p} &
{A\cdot c_2} &\shtck{\rho} &
{A/p\cdot\bar c_2}& \longrightarrow &
{0} \end{array}
$$
and the claim follows; compare Section~A.1 below. We
consider again the diagram from the beginning of the proof. This may
be further refined to a diagram of horizontal triangles  
$$\begin{array}{ccccccc}
{T(A;A)} &\hsm \shtck{p}  \hsm & 
{T(A;A)} &\hsm \shtck{i}  \hsm & 
{M_p\wedge T(A;A)} &\hsm \shtck{\beta} \hsm &  
{\Sigma T(A;A)} \\[4pt]
\scs{T(A;\rho)}&\hsm \hsm & \scs{T(A;\rho)}  & \hsm \hsm & \scs{M_p\wedge T(A;\rho)}  &\hsm
\hsm & \scs{\Sigma T(A;\rho)} 
\\[4pt]
{T(A;A/p)} &\hsm \shtck{p}\hsm &  
{T(A;A/p)} &\hsm \begin{array}{c} \scr{i} \\[-8pt]\longrightarrow\\[-10pt]\longleftarrow\\[-8pt]\scr{r}\end{array} \hsm & 
{M_p\wedge T(A;A/p)} &\hsm \begin{array}{c} \scr{\beta}
\\[-8pt]\longrightarrow\\[-10pt]\longleftarrow\\[-8pt]\scr{s}\end{array}
 \hsm & 
{\Sigma T(A;A/p)}  \\[4pt]
\scs{T(\rho;A/p)} &\hsm \hsm & \scs{T(\rho;A/p)}&\hsm \hsm & \scs{M_p\wedge T(\rho;A/p)} &\hsm \hsm & \scs{\Sigma T(\rho;A/p)}  
\\[4pt]
{T(A/p;A/p)} &\hsm \shtck{p=0}  \hsm & 
{T(A/p;A/p)} &\hsm \begin{array}{c} \scr{i} \\[-8pt]\longrightarrow\\[-10pt]\longleftarrow\\[-8pt]\scr{r}\end{array} \hsm & 
{M_p\wedge T(A/p;A/p)} &\hsm \begin{array}{c} \scr{\beta}
\\[-8pt]\longrightarrow\\[-10pt]\longleftarrow\\[-8pt]\scr{r}\end{array} \hsm &  {\Sigma T(A/p;A/p),}\end{array}
$$
where for an $A$-$A$-bimodule $M$, $T(A;M)$ is the topological
Hochschild spectrum of $A$ with coefficients in $M$. It shows that
$$(r_*\circ(M_p\wedge T(\rho))_*)(\tilde\kappa)=
(T(\rho;A/p)_*\circ r_*\circ (M_p\wedge T(A;\rho))_*)
(\tilde\kappa).$$
The map $T(\rho;A/p)_*$ is equal to the edge homomorphism of the
spectral sequence
$$E_{s,t}^2=\pi_sT(A/p,\Tor_t^A(A/p,A/p))
\Rightarrow\pi_{s+t}T(A;A/p)$$
considered in~\cite{lm}. Hence, {\it loc.cit.},~Proposition~4.3,
shows that there is a unique class in the image of
$$T(\rho;A/p)_*\:\pi_2T(A;A/p) \to \pi_2T(A/p;A/p) = \pi_2T(A/p)$$
whose image under $l_*\:\pi_2T(A/p)\to\pi_2\HH(A/p)$ is
$-\theta_K(\bar\pi_K)^{-1}\bar c_2$ and that this class has the form
$\lambda\cdot\sigma-\theta_K(\bar\pi_K)^{-1}\bar c_2$, where
$\lambda\in(\Z/p)^\times$ is a unit. Finally, the following lemma
shows that $\lambda=1$ (or equivalently, that the class $\sigma$ of
{\it loc.cit.} agrees with our class $\sigma$).
\enddemo

\specialnumber{5.3.2} \proclaim{Lemma}\label{reductionmap}The reduction
$i_*\:\bar\pi_*T(A)\to\bar\pi_*T(k)$ maps $\tilde\kappa$ to
$\sigma$.
\endproclaim

\demo{Proof} We proved in Addendum~\ref{connecting} that in the
diagram
$$\begin{array}{ccccc}
{\bar\pi_2T(A|K)}  &\stackrel{\scriptstyle j_*}{\longleftarrow}&
{\bar\pi_2T(A)} &\shtck{i_*}   &
{\bar\pi_2T(k)}  \\[4pt]
\scs{\partial_{A|K}}&&\scs{\partial_A}&&\scs{\partial_k}\\[4pt]
{\bar\pi_1\borel(C_p,T(A|K))} &\stackrel{\scriptstyle j_*}{\longleftarrow}&
{\bar\pi_1\borel(C_p,T(A))} &\shtck{i_*}   &
{\bar\pi_1\borel(C_p,T(k))} \end{array}
$$
the left-hand vertical map takes $\kappa$ to $dV(1)-V(d\log(-p))$. It
follows that the middle vertical map takes $\tilde\kappa$ to
$dV(1)-V(d\log(-p))+aV(\pi_K^{e_K-1}d\pi_K)$ for some $a\in A$. Since
$\Omega_k^1$ vanishes, we conclude that the right-hand vertical map
takes $i_*(\tilde\kappa)$ to $dV(1)$, and since $\bar\pi_2T(k)$ is
a one-dimensional $k$-vector space, it thus suffices to  show also that 
$\partial_k(\sigma)=dV(1)$. To this end, we consider the diagram 
$$\begin{array}{ccccc}
{\bar\pi_2T(k)} &\stackrel{\scr{d}}{\slar}&  {\bar\pi_1T(k)}&
\stck{\hat\Gamma} &
{\bar\pi_1\tate(C_p,T(k))}  \\[4pt]
\scs{\partial_k}& \scr{(-1)} &  \scs{\partial_k} &&\scs{\partial_k^h}
 \\[4pt]
{\bar\pi_1\borel(C_p,T(k))} &\stackrel{\scr{d}}{\slar}&
{\bar\pi_0\borel(C_p,T(k))} &\dline&
{\bar\pi_0\borel(C_p,T(k)),} \end{array}
$$
where the left-hand square anti-commutes by our conventions from
Section~2.1 above. The class $\sigma$, by definition, is
the image of $\e$ under the top differential, and the bottom
differential takes $V(1)$ to $dV(1)$. Hence, it suffices to show that
$\partial_k(\e)=-V(1)$. We recall from corollary~\ref{gammahatk} that
the class $\hat\Gamma(\e)$ is represented in the spectral sequence
$\hat E^*(C_p,T(k))$ by the infinite cycle $u_1t^{-1}$. Hence,
Addendum~\ref{boundarymap} shows that the image of this class by the
right-hand vertical map is $-V(1)$.
\enddemo  

\specialnumber{5.3.3}\numbereddemo{{R}emark} \label{unitisone}
It follows from
Propositions~\ref{reductionmodp} and~\ref{TparkparPi} that in
$\bar\pi_*T(A)$,
$$d\tilde\kappa = 
-\theta_K'(\pi_K)\theta_K(\pi_K)^{-1}d\pi_K\cdot\tilde\kappa.$$
This implies that $d\kappa=\kappa d\log(-p)$ in $\bar\pi_*T(A|K)$ as
stated in Theorem~\ref{pisubast}.
\enddemo

We construct a number of infinite cycles. Recall the map of ring
spectra
$$\hat\Gamma_{A|K}\:T(A|K)^{C_{p^{n-1}}}\to\tate(C_{p^n},T(A|K)).$$

\specialnumber{5.3.4}\proclaim{Proposition} \label{dlogpi}
For all $K${\rm ,} the element $d\log\pi_K\in
\hat E^2(C_{p^n},T(A|K))$ is an infinite cycle and represents the
homotopy class $\hat\Gamma_{A|K}(d\log_n\pi_K)$.
\endproclaim

\demo{Proof} We consider the diagram
$$\begin{array}{ccc}
{T(A|K)^{C_{p^n}}} &\stck{R} &
{T(A|K)^{C_{p^{n-1}}}}\\[4pt]
\scs{\Gamma_{A|K}} && \scs{\hat\Gamma_{A|K}} \\[4pt]
{\coborel(C_{p^n},T(A|K))}& \stck{R^h} &
{\tate(C_{p^n},T(A|K)).} \end{array}
$$
In the spectral sequence
\begin{eqnarray*}
E^2(C_{p^n},T(A|K)) & =&
\Lambda\{u_n,d\log\pi_K\}\otimes S\{\pi_K,t,\kappa\}/(\pi_K^{e_K}) \\
{} & \Rightarrow&\bar\pi_*(\coborel(C_{p^n},T(A|K))) ,
\end{eqnarray*}
the element $d\log\pi_K$ is an infinite cycle and represents
$\Gamma_{A|K}(d\log_{n+1}\pi_K)$. Indeed, if we compose $\Gamma_{A|K}$
and the edge-homomorphism of this spectral sequence, we get the map
$F^n\:\bar\pi_* T(A|K)^{C_{p^n}} \to \bar\pi_*T(A|K)$ which takes
$d\log_{n+1}\pi_K$ to $d\log\pi_K$. The map $R^h$ induces the obvious
inclusion on $E^2$-terms. Hence, the element $d\log\pi_K$ of $\hat
E^*(C_{p^n},T(A|K))$ is an infinite cycle and, since it is not a boundary,
represents
$R^h(\Gamma_{A|K}(d\log_{n+1}\pi_K))=\hat\Gamma_{A|K}(R(d\log_{n+1}\pi_K))
=\hat\Gamma_{A|K}(d\log_n\pi_K)$.
\enddemo

Let $\alpha_A=u_K^{(-n)}(\pi_K)^{-1}\tilde\kappa$ and
$\tau_A=u_K^{(-n)}(\pi_K)^p\,t$ such that $j_*(\alpha_A)=\alpha_K$ and
$j_*(\tau_A)=\tau_K$.

\specialnumber{5.3.5} \proclaim{Lemma}\label{gammahatpi}
Suppose that $\mu_p\subset K$ and let
$n<v_p(e_K)$. Then the elements $\pi_K^{p^n}$ and $-\tau_A\alpha_A$
of $\hat E^2(C_{p^n},T(A))$ are infinite cycles and represent the
homotopy classes $\hat\Gamma_A(\ul{\pi_K}_n)$ and
$\hat\Gamma_A(\ul{\pi_K}_n^{e_K/p^n})${\rm ,} respectively.
\endproclaim 

\demo{Proof} We consider the diagram
$$\begin{array}{ccc}
{\bar\pi_*(T(A)^{C_{p^{n-1}}})} &\stck{\Gamma_A} &
{\bar\pi_*(\tate(C_{p^n},T(A)))} \\[4pt]
 \scs{\rho_*}&& \scs{\rho_*} \\[4pt]
{\bar\pi_*(T(A/p)^{C_{p^{n-1}}})}& \stck{\Gamma_{A/p}} &
{\bar\pi_*(\tate(C_{p^n},T(A))),} \end{array}
$$
with the vertical maps induced from the reduction $\rho\:A\to
A/p$. The lower horizontal map is studied in the appendix. By
Addendum~\ref{tatepAfpp},
\begin{eqnarray*}
\hat E^2(C_{p^n},T(A/p)) &\hsm = \hsm&\Lambda\{u_n,d\pi_K,\e\}\otimes
S\{t^{\pm 1},\pi_K,\sigma\}/(\pi_K^{e_K})\otimes\Gamma\{\bar c_2\}, \\
\hat E^3(C_{p^n},T(A/p)) & \hsm= \hsm&\Lambda\{u_n,d\pi_K\}\otimes
S\{t^{\pm 1},\pi_K^p\}/(\pi_K^{e_K})\otimes\Gamma\{\bar c_2\}, 
\end{eqnarray*}
and $\hat E^3(C_{p^n},T(A/p))=\hat E^\infty(C_{p^n},T(A/p))$.
We compare this to
\begin{eqnarray*}
\hat E^2(C_{p^n},T(A)) &\hsm =\hsm& \Lambda\{u_n,d\pi_K\}\otimes
S\{\tau_A^{\pm1},\alpha_A,\pi_K\}/(\pi_K^{e_K}), \\
\hat E^3(C_{p^n},T(A))&\hsm =\hsm&\Lambda\{u_n,\pi_K^{p-1}d\pi_K\}\otimes
S\{\tau_A^{\pm1},\alpha_A,\pi_K^p\}/(\pi_K^{e_K}).  
\end{eqnarray*}
The map $\rho_*\:\bar\pi_*T(A)\to\bar\pi_*T(A/p)$ was evaluated
above. The induced map
$$\hat E^3(C_{p^n},T(A)) \hookrightarrow \hat E^3(C_{p^n},T(A/p))$$
is the monomorphism which takes $\tau_A\alpha_A$ to $-t\bar
c_2$. Indeed, the map of $E^2$-terms takes the element
$\tau_A\alpha_A$ to
$-t\bar c_2+\theta_K(\pi_K)t\sigma-t\e\cdot\theta_K'(\pi_K)d\pi_K$, 
and the last two summands are equal to the image by the
$d^2$-differential of $\e\cdot\theta_K(\pi_K)$. For $0\leq s\leq 1$,
we have the diagram  
$$\begin{array}{ccc} 
{\hat E_{-s,s}^3(C_{p^n},T(A))\;}&\hookrightarrow&
{\hat E_{-s,s}^3(C_{p^n},T(A/p))} \\[4pt]
\lower7pt\hbox{${\scriptscriptstyle\cup}$}\hskip-7.75pt\big\uparrow
&&\big\uparrow\scr{\sim}\\[4pt] {\hat E_{-s,s}^{\infty}(C_{p^n},T(A))\;}&\rightarrow&
{\hat E_{-s,s}^{\infty}(C_{p^n},T(A/p))}   \end{array}$$
and we conclude that the lower horizontal map is a monomorphism. We
show in Proposition~\ref{gammahatpilemma} that the classes
$\hat\Gamma_{A/p}(\ul{\pi_K}_n)$ and
$\hat\Gamma_{A/p}(\ul{\pi_K}_n^{e_K/p^n})$ are represented in the
spectral sequence $\hat E^*(C_{p^n},T(A/p))$ by the infinite cycles
$\pi_K^{p^n}$ and $t\bar c_2$, respectively. It follows immediately
that $\hat\Gamma_A(\ul{\pi_K}_n)$ is represented by $\pi_K^{p^n}$ as
stated. To conclude that $\hat\Gamma_A(\ul{\pi_K}_n^{e_K/p^n})$ is
represented by $-\tau_A\alpha_A$ we must rule out that an element of
$\hat E_{-s,s}^2(C_{p^n},T(A))$ with $0\leq s\leq 1$ represent this
class. But this follows from the injectivity of the lower horizontal
map in the diagram above. 
\enddemo

\specialnumber{5.3.6}\proclaim{Proposition} \label{taualpha}
Suppose that $\mu_p\subset K$ and let
$n<v_p(e_K)$. Then the elements $\pi_K^{p^n}$ and $-\tau_K\alpha_K$
of $\hat E^2(C_{p^n},T(A|K))$ are infinite cycles and represent the
homotopy classes $\hat\Gamma_{A|K}(\ul{\pi_K}_n)$ and
$\hat\Gamma_{A|K}(\ul{\pi_K}_n^{e_K/p^n})${\rm ,} respectively.
\endproclaim

\demo{Proof} The map 
$j_*\:\hat E^2(C_{p^n},T(A)) \to \hat E^2(C_{p^n},T(A|K))$
takes the infinite cycle $\pi_K^{p^n}$ (resp.\ $-\tau_A\alpha_A$) to
the element $\pi_K^{p^n}$ (resp.\ $-\tau_K\alpha_K$) which therefore
is an infinite cycle. It follows $\pi_K^{p^n}$
(resp.\ $-\tau_K\alpha_K$) either represents
$\hat\Gamma_{A|K}(\ul{\pi_K}_n)$
(resp.\ $\hat\Gamma_{A|K}(\ul{\pi_K}_n^{e_K/p^n})$) or else it is a 
boundary. The element $\pi_K^{p^n}$ cannot be a boundary, but we must
check that $-\tau_K\alpha_K$ is not a $d^3$-boundary. To this end we
consider the diagram
$$\begin{array}{ccc}
{\hat E_{1,0}^3(C_{p^n},T(A|K))} &\shtck{d^3} &
{\hat E_{-2,2}^3(C_{p^n},T(A|K))} \\[4pt]
\scr{\sim}\big\uparrow\scr{j_*}&&\scr{\sim}\big\uparrow\scr{j_*}\\[4pt]
{\hat E_{1,0}^3(C_{p^n},T(A))}& \shtck{d^3}   &
{\hat E_{-2,2}^3(C_{p^n},T(A))} \end{array}
$$
with vertical isomorphisms. The right-hand vertical map takes
$-\tau_A\alpha_A$ to $-\tau_K\alpha_K$, and since $-\tau_A\alpha_A$ is
not a $d^3$-boundary, neither is $-\tau_K\alpha_K$.
\enddemo

Let $\ell\:S^0\to\tate(C_{p^n},T(A|K))$ be the unit map and let
$v_1\in\bar\pi_{2(p-1)}(S^0)$ be the canonical generator. If
$\mu_p\subset K$, then $\ell_*(v_1)=\hat\Gamma_{A|K}(b_n)^{p-1}$.

\specialnumber{5.3.7}\proclaim{Addendum} \label{v_1}
The elements $-t\kappa^p$ and
$(-t\kappa)^{p^n}$ of $\hat E^2(C_{p^n},T(A|K))$ are infinite cycles
which{\rm ,} if not boundaries{\rm ,} represent the homotopy classes $\ell_*(v_1)$
and $V(1)${\rm ,} respectively.
\endproclaim

\demo{Proof} The elements $-t\kappa^p$ and $(-t\kappa)^{p^n}$ are in
the image of
$$\iota_*\:\hat E^2(C_{p^n},T(W(k)|K_0)) \to
\hat E^2(C_{p^n},T(A|K))$$
and so the statement, if valid for some $K$, is valid for all $K$. Now,
suppose that $\mu_p\subset K$ and that $v_p(e_K)>n$. We may argue as
in the proof of Proposition~\ref{dlogpi} that $\hat\Gamma_{A|K}(b_n)$
is represented by the infinite cycle $-\pi_K^{e_K/(p-1)}\alpha_K$.
Indeed, $b_n=R(b_{n+1})$ and $F^n(b_{n+1})=b_1$ and
from~(\ref{bottelement}) we know that $b_1=-\pi_K^{e_K/(p-1)}\alpha_K$.
Now
$$-\pi_K^{e_K/(p-1)}\alpha_K=-(\pi_K^{p^n})^{e_K/p^n(p-1)}\alpha_K,$$
and it follows from Proposition~\ref{taualpha} that
$\pi_K^{e_K/(p-1)}$ is an infinite cycle and represents
$\hat\Gamma_{A|K}(\ul{\pi_K}_n^{e_K/p^n(p-1)})$. Hence, also
$\alpha_K$ is an infinite cycle, and since it is  not a boundary, represents
a homotopy class, say, $\tilde\alpha_K$. Since the classes
$\hat\Gamma_{A|K}(b_n)$ and
$-\hat\Gamma_{A|K}(\ul{\pi_K}_n^{e_K/p^n(p-1)})\tilde\alpha_K$ are
represented in the spectral sequence by the same element so are their
$(p-1)$st powers. We know from Proposition~\ref{taualpha} that
$$\hat\Gamma_{A|K}(\ul{\pi_K}_n^{e_K/p^n(p-1)})^{p-1}=
\hat\Gamma_{A|K}(\ul{\pi_K}_n^{e_K/p^n})$$
is represented by $-\tau_K\alpha_K$. And $\alpha_K^{p-1}$, if not a
boundary, represents $\tilde\alpha_K^{p-1}$. It follows that
$-\tau_K\alpha_K^p$, if not a boundary, represents
$\hat\Gamma_{A|K}(b_n^{p-1})=\ell_*(v_1)$.

We recall from Lemmas~\ref{p>2} and~\ref{ppower} that in the Witt ring
$\bar W_n(A)$,
$$V(1)=\theta_K(\ul{\pi_K}_n)^{-1}\ul{\pi_K}_n^{e_K},$$
and hence, in $\bar\pi_*\tate(C_{p^n},T(A|K))$,
\begin{eqnarray*}
V(1) & =&V(\hat\Gamma_{A|K}(1)) = \hat\Gamma_{A|K}(V(1))
= \hat\Gamma_{A|K}(\theta_K(\ul{\pi_K}_n)^{-1}\ul{\pi_K}_n^{e_K}) \\[4pt]
{} & =  &\hat\Gamma_{A|K}(\theta_K(\ul{\pi_K}_n))^{-1} \cdot
\hat\Gamma_{A|K}(\ul{\pi_K}_n^{e_K/p^n})^{p^n}. 
\end{eqnarray*}
It follows that the infinite cycle
$$\theta_K(\pi_K^{p^n})^{-1}(-\tau_K\alpha_K)^{p^n}
=(-\theta_K^{(-n)}(\pi_K)^{-1}\tau_K\alpha_K)^{p^n}=(-t\kappa)^{p^n},$$
if not a boundary, represents $V(1)$ as stated.
\enddemo

5.4.\quad In this section we evaluate the spectral
sequence $\hat E^*(C_p,T(A|K))$.

\specialnumber{5.4.1} \proclaim{Lemma}\label{cartierlemma}
Let $k$ be a field of characteristic
$p>0${\rm ,} let $f(x)$ be a power series over $k$ with nonzero constant
term{\rm ,} and let
$$\frac{f'(x)x}{f(x)}=a_1x+a_2x^2+\dots$$
be the logarithmic derivative. Then $a_{pi}=a_i^p$, for all $i\geq 1$.
\endproclaim

\demo{Proof} We may assume that $f(x)$ is a polynomial with
$f(0)\in k^\times$. Moreover, replacing $k$ by a splitting field 
for $f(x)$, we can assume that $f(x)$ splits as a product of
linear factors. And since the logarithmic derivative takes products of
power series to sums, we are reduced to the case of a linear
polyonomial. The result in this case is readily verified by computation. 
\enddemo

\specialnumber{5.4.2}\proclaim{Proposition} \label{nequals1}
Suppose either $\mu_p\subset K$ or
$K=K_0$. Then{\rm ,} up to a unit{\rm ,} the nonzero differentials in the
spectral sequence $\hat E^2(C_p,T(A|K))$ are generated from
\begin{eqnarray*}
d^2(\tau_K^a\pi_K^r\alpha_K^d) & =&
\tau_K\,d\log\pi_K\cdot
\tau_K^a\pi_K^r\alpha_K^d, \hskip6mm \hbox{if $v_p\{a,r,d\}_K=0$}, \\
d^{2p+1}(u_1) & = &(\tau_K\alpha_K)^p\tau_K\\
\end{eqnarray*}
and from the fact that  $d\log\pi_K$ is an infinite cycle.
\endproclaim

\demo{Proof} The $d^2$-differential follows from
Proposition~\ref{d2=connes}. If $K=K_0$, we have
$$\hat E^3(C_p,T(W(k)|K_0)) =
\Lambda\{u_1,d\log(-p)\}\otimes S\{t^{\pm1},\kappa^p\},$$
and for degree reasons, the first possible differential is
$d^{2p+1}$. The canonical map
$$\hat E^{2p+1}({\Bbb T},T(W(k)|K_0)) \hookrightarrow
\hat E^{2p+1}(C_p,T(W(k)|K_0))$$
may be identified with the inclusion of the subalgebra generated by
$t$, $\kappa^p$, and $d\log(-p)$. The $d^{2p+1}$-differential on these
elements in the left-hand spectral sequence are zero for degree
reasons. Hence, the $d^{2p+1}$-differential on these classes in the
right-hand spectral sequence are zero as well. We claim that, up to a
unit,
$$d^{2p+1}u_1=t^{p+1}\kappa^p.$$
For if not, $t\kappa^p$ would survive the spectral sequence and
represent the homotopy class $-v_1\cdot 1$. But $\tate(C_p,T(W|K_0))$
is a module spectrum over the\break generalized Eilenberg-MacLane spectrum
$T(W)$, and therefore, is itself a\break generalized Eilenberg-MacLane
spectrum. Hence, multiplication by $v_1$ on
$\bar\pi_*\tate(C_p,T(W|K_0))$ is identically zero. All further
differentials must vanish for degree reasons.

If $\mu_p\subset K$ and $v_p(e_K)>1$, 
$$\hat E^3(C_p,T(A|K))=\Lambda\{u_1,d\log\pi_K\}\otimes
S\{\pi_K^p,\alpha_K,\tau_K^{\pm 1}\}/(\pi_K^{e_K}),$$
and by Proposition~\ref{taualpha}, $\pi_K^p$ and $\tau_K\alpha_K$ are
infinite cycles. From the previous case, we know that $t$ is an
infinite cycle, and hence, so is $\tau_K=u_K(\pi_K^p)\,t$. It follows
that  $\alpha_K$ also is an infinite cycle. Hence, the remaining
nonzero differentials are generated from the differential on
$u_1$. Again all further differentials vanish for degree reasons.

Finally suppose that $\mu_p\subset K$, but with no restriction on
$v_p(e_K)$. Then
$$\hat E^3(C_p,T(A|K))=\Lambda\{u_1,d\log\pi_K\}\otimes
k\{\tau_K^a\pi_K^r\alpha_K^d\,|\,v_p\{a,r,d\}_K\geq 1\},$$
and we need to show that the elements $\tau_K^a\pi_K^r\alpha_K^d$ with
$v_p\{a,r,d\}_K\geq 1$ are $d^q$-cycles, for $3\leq q\leq 2p+1$. To this
end, we let $L/K$ be a totally ramified extension such that
$v_p(e_L)>1$ and consider
$$\iota_*\:\hat E^q(C_p,T(A|K))
\to \hat E^q(C_p,T(B|L)).$$
We have from~(\ref{E2iota}) and Lemma~\ref{cartierlemma} that
\begin{eqnarray*}
\iota_*(\tau_K^a\pi_K^r\alpha_K^d) & =&
(\theta_{L/K}(\pi_L^p))^{-\{a,r,d\}_K/p}
\tau_L^a\pi_L^{e_{L/K}r}\alpha_L^d,\\
\iota_*(d\log\pi_K)&=&(e_{L/K}-\frac{\theta_{L/K}'(\pi_L^p)\pi_L^p}{
\theta_{L/K}(\pi_L^p)})d\log\pi_L, 
\end{eqnarray*}
where in the first line, $v_p\{a,r,d\}_K\geq 1$. We know from the previous
case that the $d^q$-differential on
$\iota_*(\tau_K^a\pi_K^r\alpha_K^d)$ vanishes, and hence, it will
suffice to show that we can find $L/K$ for which the map $\iota_*$ is
injective.

If $v_p(e_{L/K})>1$ and $\theta_{L/K}(x)=x-1$ then, up to a unit,
$$\iota_*(u_1^{\e}\tau_K^a\pi_K^r\alpha_K^dd\log\pi_K)
=u_1^{\e}\tau_L^a\pi_L^{e_{L/K}r+p}\alpha_L^rd\log\pi_L,$$
and hence, for $\iota_*$ to be injective, we need that
$e_{L/K}r+p<e_L$. Since $r\leq e_K-1$ and $e_L=e_{L/K}e_K$, this is
equivalent to the requirement that $e_{L/K}\geq p$. We also need
$v_p(e_L)>1$. The extension $L$ with $e_{L/K}=p^2$ and
$\theta_{L/K}(x)=x-1$ satisfies both requirements. It follows that the
$d^q$-differentials vanish, if $3\leq q\leq 2p$, and that the nonzero
$d^{2p+1}$-differentials are generated from the differential on $u_1$.
All further differentials vanish for degree reasons.
\enddemo

\specialnumber{5.4.3}\proclaim{Theorem} \label{gammahat}
For all $K${\rm ,} and for $i\geq 0${\rm ,} the map
$$\hat\Gamma_{A|K}\:\bar\pi_iT(A|K)
\xto{\sim} \bar\pi_i\tate(C_p,T(A|K)) $$
is an isomorphism\/{\rm .}
\endproclaim

\demo{Proof} If we let $L=K(\mu_p)$, then in the diagram
$$\begin{array}{ccc}
{\bar\pi_*T(A|K)} &\stck{\hat\Gamma_{A|K}} &
{\bar\pi_*\tate(C_p,T(A|K))}  \\[3pt]
\scs{\sim}&&\scs{\sim}\\[3pt]
{\bar\pi_*T(B|L)^{G_{L/K}}} &\stck{\hat\Gamma_{B|L}} &
{\bar\pi_*\tate(C_p,T(B|L))^{G_{L/K}},}\end{array}
$$
the vertical maps are isomorphisms. Indeed, this follows from
Theorem~\ref{tamedescent} and from the Tate spectral sequence, since
the order of $G_{L/K}$ is prime to $p$. Hence, we can assume that
$\mu_p\subset K$.

If $\mu_p\subset K$ and $v_p(e_K)>1$ or if $K=K_0$, then
$$\hat E^\infty(C_p,T(A|K))=\Lambda\{d\log\pi_K\}\otimes
S\{\pi_K^p,\alpha_K,\tau_K^{\pm1}\}/(\pi_K^e,\alpha_K^p),$$
and by Proposition~\ref{taualpha}, there is a multiplicative extension
$(\pi_K^p)^{e/p}=-\tau_K\alpha_K$ in the passage to the actual homotopy
groups. Hence, as a $k$-algebra 
$$\bar\pi_*\tate(C_p,T(A|K))=\Lambda\{\hat\Gamma_{A|K}(d\log\pi_K)\}
\otimes S\{\hat\Gamma_{A|K}(\pi_K),\tilde\tau_K^{\pm1}\}/
(\hat\Gamma_{A|K}(\pi_K)^{e_K}),$$
where the class $\tilde\tau_K$ is a lifting of $\tau_K$. It follows
that $\bar\pi_*T(A|K)$ and the nonnegatively graded part of
$\bar\pi_*\tate(C_p,T(A|K))$ are abstractly isomorphic $k$-algebras,
and that the map $\hat\Gamma_{A|K}$ is an isomorphism for $i=0$ and
$i=1$. To show that $\hat\Gamma_{A|K}$ is an isomorphism, for $i\geq
0$, it will therefore suffice to show that 
$$\hat\Gamma_{W(k)|K_0}\:\bar\pi_2T(W|K_0)\xto{\sim}
\bar\pi_2\tate(C_p,T(W|K_0))$$
is an isomorphism. To this end, we consider the diagram
$$\begin{array}{ccc}
{\bar\pi_2T(W|K_0)}&
\begin{array}{c}\scr{{\beta_1}}\\[-8pt]\srar\\[-10pt]\scr{\sim}\end{array}&
 {\bar\pi_1T(W|K_0)}\\[4pt]
\hskip24pt\scs{\hat\Gamma_{W(k)|K_0}} &&\quad\scr{\sim} \scs{\hat\Gamma_{W(k)|K_0}}\\[4pt]
{\bar\pi_2\tate(C_p,T(W|K_0))} &\stck{\beta_1} &
{\bar\pi_1\tate(C_p,T(W|K_0)),}\end{array}$$
where the upper horizontal map and right-hand vertical maps are
isomorphisms. Since all groups in the diagram are one-dimensional
$k$-vector spaces, the left-hand vertical map and lower horizontal
map must also be isomorphisms. This shows that the map of the
statement is an isomorphism if $\mu_p\subset K$ and $v_p(e_K)>1$ or if
$K=K_0$.

If $\mu_p\subset K$, but there are no restrictions on $v_p(e_K)$,
$$\hat E^\infty(C_p,T(A|K))=
\Lambda\{d\log\pi_K\}\otimes
k\{\tau_K^a\pi_K^r\alpha_K^d\,|\,v_p\{a,r,d\}_K\geq 1,\,d<p\},$$
where $0\leq r<e_K$, $d\in\N_0$ and $a\in\Z$. Again, the domain and
range of $\hat\Gamma_{A|K}$ are abstractly isomorphic $k$-vector
spaces. We choose an extension $L/K$ such that $v_p(e_L)>1$ and such
that $\iota\:\bar\pi_*T(A|K)\to\bar\pi_*T(B|L)$ is a monomorphism.
Since $\hat\Gamma_{B|L}$ is an isomorphism in nonnegative degrees,
$\hat\Gamma_{A|K}$ is a monomorphism, and hence an isomorphism, in
nonnegative degrees.
\enddemo 

\specialnumber{5.4.4}\proclaim{Addendum} \label{gammahatadd}
For all $K${\rm ,} for all $n,v\geq 1${\rm ,}
and for all $i\geq 0${\rm ,} the maps
\begin{eqnarray*}
\hat\Gamma_{A|K} \: \pi_i(T(A|K)^{C_{p^{n-1}}},\Z/p^v) & \xto{\sim}&
\pi_i(\tate(C_{p^n},T(A|K)),\Z/p^v), \\
\Gamma_{A|K} \: \pi_i(T(A|K)^{C_{p^n}},\Z/p^v) & \xto{\sim}&
\pi_i(\coborel(C_{p^n},T(A|K)),\Z/p^v),  
\end{eqnarray*}
are isomorphisms.
\endproclaim

\demo{Proof} If $v=1$ this follows from Theorem~\ref{gammahat} and
the main theorem of~\cite{tsalidis}, and the general case follows by
easy induction based on the Bockstein sequence.
\enddemo

5.5.\quad We now evaluate the spectral sequences
$\hat E^*(C_{p^n},T(A|K))$.

\specialnumber{5.5.1}\proclaim{Theorem} \label{diff}
Suppose either $\mu_p\subset K$ or
$K=K_0$. Then the nonzero differentials in the spectral sequence
\begin{eqnarray*}
\hat E^2(C_{p^n},T(A|K)) & =& \Lambda\{u_n,d\log\pi_K\}\otimes
S\{\pi_K,\alpha_K,\tau_K^{\pm1}\}/(\pi_K^{e_K}) \\
{} & \Rightarrow&\bar\pi_*(\tate(C_{p^n},T(A|K))) 
\end{eqnarray*}
are multiplicatively generated from
\begin{eqnarray*}
d^{2(\frac{p^{v+1}-1}{p-1})}(\tau_K^a\pi_K^r\alpha_K^d) & =&
\lambda\cdot(\tau_K\alpha_K)^{\frac{p^{v+1}-1}{p-1}-1}\tau_K\,d\log\pi_K
\cdot\tau_K^a\pi_K^r\alpha_K^d, \\
d^{2(\frac{p^{n+1}-1}{p-1})-1}(u_n) & = &
\mu\cdot(\tau_K\alpha_K)^{\frac{p^{n+1}-1}{p-1}-1}\tau_K,  
\end{eqnarray*}
and from $d\log\pi_K$ being an infinite cycle. Here $\lambda$ and
$\mu$ are units in $A/p$ and in the first line
$v=v_p\{a,r,d\}_K${\rm ,} where $\{a,r,d\}_K=(pa-d)e_K/(p-1)+r$.
\endproclaim

\specialnumber{5.5.2}\numbereddemo{{R}emark} 
We show that the units $\lambda$ and $\mu$ above are
given by
$$\lambda = -\lambda_v\cdot p^{-v}\{a,r,d\}_K\cdot
u_K^{(v-n)}(\pi_K^{p^v})^{-p}, \hskip6mm
\mu = \mu_n\cdot u_K(\pi_K^{p^n})^{-p},$$
where $\lambda_v$ and $\mu_n$ are units in $\Fp$ independent of $K$.
\enddemo

The proof of Theorem~\ref{diff} is similar to the proof of
Proposition~\ref{nequals1} above, but the individual steps are more
involved. It will be necessary to know   the structure of the
$E^r$-terms, given the differential structure.

\specialnumber{5.5.3} \proclaim{Lemma}\label{Er}
Suppose $\mu_p\subset K$ or $K=K_0${\rm ,} and assume
that Theorem~{\rm \ref{diff}} is true for $K$. Let $\hat E^q=\hat
E^q(C_{p^n},T(A|K))$. Then for $0\leq s<n$ and $2(p^s-1)/(p-1)<q\le 2(p^{s+1}-1)/(p-1)${\rm ,}
\begin{eqnarray*}
  \hat E^{q} &\hsm =\hsm&\bigoplus_{v=1}^{s-1}
\Lambda\{u_n\}  \otimes\; k\big\{\tau_K^a\pi_K^r
\alpha_K^dd\log\pi_K\mid v_p\{a,r,d\}_K\!=v,\,
d<\textstyle{\frac{p^{v+1}\!-1}{p-1}-1}\big\}\\
&&  \oplus\;\Lambda\{u_n,d\log\pi_K\}\otimes 
k\big\{\tau_K^a\pi_K^r\alpha_K^d\mid
v_p\{a,r,d\}_K\!\geq s\big\},\end{eqnarray*}
\begin{eqnarray*}
 \hat E^\infty&\hsm=\hsm & \bigoplus_{v=1}^{n-1}
\Lambda\{u_n\}\otimes k\{\tau_K^a\pi_K^r
\alpha_K^dd\log\pi_K\mid v_p\{a,r,d\}_K=v,\,
d<\textstyle{\frac{p^{v+1}-1}{p-1}-1}\}\\
&&  \oplus\;\Lambda\{d\log\pi_K\}\otimes
k\{\tau_K^a\pi_K^r\alpha_K^d\mid v_p\{a,r,d\}_K\geq n,\,
d<\textstyle{\frac{p^{n+1}-1}{p-1}-1}\}, 
\end{eqnarray*}
where $0\leq r<e_K${\rm ,} $d\in\N_0$ and $a\in\Z$.
\endproclaim

\demo{Proof} When we assume the result for $s<n-1$, Theorem~\ref{diff}
implies that
$$\hat E^{2(\frac{p^{s+2}-1}{p-1})}=\hat E^{2(\frac{p^{s+1}-1}{p-1})+1},$$
and inductively, $\hat E^{2(\frac{p^s-1}{p-1})}$ is given by the statement
of the lemma. Indeed, this is clear in the basic case $s=0$. The
differential $d^{2(p^{s+1}-1)/(p-1)}$ only affects the last summand
on the right-hand side of the statement and does not involve the
tensor factor $\Lambda\{u_n\}$. If we rewrite
\begin{eqnarray*}
&&\hskip-48pt\Lambda\{d\log\pi_K\}\otimes k\{\tau_K^a\pi_K^r\alpha_K^d\,|\,
v_p\{a,r,d\}_K\geq s\}\\[3pt]
&=\hsm & k\{\tau_K^a\pi_K^r\alpha_K^d\,|\,v_p\{a,r,d\}=s\}\;\\[3pt]
&&\oplus\;  k\{\tau_K^a\pi_K^r\alpha_K^dd\log\pi_K|\,v_p\{a,r,d\}_K=s,\,
d\geq\textstyle{\frac{p^{s+1}-1}{p-1}-1}\}\;\\[3pt]
{}&&\oplus\;  k\{\tau_K^a\pi_K^r\alpha_K^dd\log\pi_K|\,v_p\{a,r,d\}_K=s,\,
d<\textstyle{\frac{p^{s+1}-1}{p-1}-1}\}\;\\[3pt]
{}&&\oplus\;  \Lambda\{d\log\pi_K\}\otimes
k\{\tau_K^a\pi_K^r\alpha_K^d|\,v_p\{a,r,d\}_K\geq s+1\}, 
\end{eqnarray*}
the differential $d^{2(p^{s+1}-1)/(p-1)}$ vanishes on the last two
summands and maps the first summand isomorphically onto the
second. Indeed,
$$\textstyle{\{a+\frac{p^{s+1}-1}{p-1},r,d+\frac{p^{s+1}-1}{p-1}-1\}_K
=\{a,r,d\}_K+\frac{p^{s+1}}{p-1}.}$$
Assuming that Theorem~\ref{diff} holds for $K$, we have
$$\hat E^{2(\frac{p^n-1}{p-1})+1}=\hat E^{2(\frac{p^{n+1}-1}{p-1})-1},$$
and the common value has already been determined. Up to a unit,
$$d^{2(\frac{p^{n+1}-1}{p-1})-1}u_n
=(\tau_K\alpha_K)^{\frac{p^{n+1}-1}{p-1}-1}\tau_K.$$
 If we rewrite
\begin{eqnarray*}
&&\hskip-36pt \Lambda\{u_n  ,d\log\pi_K\} \otimes k\{\tau_K^a\pi_K^r\alpha_K^d\,|\,
v_p\{a,r,d\}_K\geq n\}  \\[3pt]
&=\hsm & \Lambda\{d\log\pi_K\}\otimes 
k\{\tau_K^a\pi_K^r\alpha_K^d\,|\,v_p\{a,r,d\}_K\geq n,\,
d\geq\textstyle{\frac{p^{n+1}-1}{p-1}-1}\}\;\\[3pt]
{} & &\oplus\;  \Lambda\{d\log\pi_K\}\otimes 
k\{u_n\tau_K^a\pi_K^r\alpha_K^d\,|\,
v_p\{a,r,d\}_K\geq n\}\;\\[3pt]
{} &&\oplus\;  \Lambda\{d\log\pi_K\}\otimes 
k\{\tau_K^a\pi_K^r\alpha_K^d\,|\,
v_p\{a,r,d\}_K\geq n,\;d<\textstyle{\frac{p^{n+1}-1}{p-1}-1}\}, 
\end{eqnarray*}
the differential $d^{2(p^{n+1}-1)/(p-1)-1}$ maps the second summand
isomorphically onto the first summand and leaves the last summand
unchanged. 
\enddemo

\specialnumber{5.5.4}\proclaim{Proposition} \label{KeqKsub0}
Let $T=T(W(k)|K_0)$. In the spectral sequence
$$\hat E^2(C_{p^n},T) = \Lambda\{u_n,d\log(-p)\}\otimes
S\{t^{\pm1},\kappa\} 
\Rightarrow \bar\pi_*(\tate(C_{p^n},T)),$$
the higher differentials are multiplicatively generated from
\begin{eqnarray*}
d^{2(\frac{p^{v+1}-1}{p-1})}(t^{p^{v-1}}) & =& \lambda_v\cdot
(t\kappa)^{\frac{p^{v+1}-1}{p-1}-1}td\log(-p)\cdot t^{p^{v-1}},\hskip4mm
1\leq v<n,\\
d^{2(\frac{p^{n+1}-1}{p-1})-1}(u_n)&
=&\mu_n\cdot (t\kappa)^{\frac{p^{n+1}-1}{p-1}-1}t,
\end{eqnarray*}
where $\lambda_v,\mu_n\in\Fp$ are units{\rm ,} and from $t\kappa^p$ and
$d\log(-p)$ being infinite cycles. Moreover{\rm ,} the infinite cycles
$(-t\kappa)^{p^{s+1}}d\log(-p)${\rm ,} $1\leq s<n${\rm ,} represent
$dV^{n-s}(1)$.
\endproclaim

\demo{Proof} The proof is by induction on $n$ and is similar to the
proof in~\cite{bm} of the differential structure of the spectral
sequences $\hat E^*(C_{p^n},T(W(k)))$. The basic case $n=1$ was proved
in Proposition~\ref{nequals1}. So assume the statement for $n-1$.

We first argue that in $\bar\pi_*(\tate(C_{p^n},T))$, the class
$v_1^m$ is nonzero if and only if $m<(p^n-1)/(p-1)$. By
Addendum~\ref{gammahatadd}, the maps
$$\bar\pi_*(\tate(C_{p^n},T)) \xleftarrow{\hat\Gamma}
\bar\pi_*(T^{C_{p^{n-1}}}) \xto{\hat\Gamma}
\bar\pi_*(\coborel(C_{p^{n-1}},T))$$
are isomorphisms in nonnegative degrees,  and hence, we may instead
consider the class $v_1^m$ in $\bar\pi_*(\coborel(C_{p^{n-1}},T))$.
To this end, we use the spectral sequence
$$E^2(C_{p^{n-1}},T)  = \Lambda\{u_{n-1},d\log(-p)\} \otimes
S\{t,\kappa\} \Rightarrow \pi_*(\coborel(C_{p^{n-1}},T))$$
whose differential structure is determined by the statement for $n-1$.
We evaluate the $E^r$-term by an argument similar to the proof of
Lemma~\ref{Er}. To state the result, let $P(a,d,v)$ be the statement
$$\textstyle{\hbox{``$a<\frac{p^{v+1}\!-1}{p-1}$ or
$d<\frac{p^{v+1}\!-1}{p-1}-1$, or both''.}}$$
Then for $0\leq s<n-1$,
\begin{eqnarray*}
E^{2(\frac{p^{s+1}\!-1}{p-1})}&\hsm =\hsm
  & \bigoplus_{v=0}^{s-1}\, \Lambda\{u_{n-1}\}\otimes
k\big\{t^a\kappa^dd\log(-p)\mid \hbox{$v_p(pa-d)=v$ and $P(a,d,v)$}
\big\} \\[4pt]
{}&\hsm\hsm&\oplus\; \Lambda\{u_{n-1},d\log(-p)\}\otimes 
k\big\{t^a\kappa^d\mid v_p(pa-d)\geq s\big\},\end{eqnarray*}
\begin{eqnarray*}
E^\infty&\hsm =\hsm
{} &  \bigoplus_{v=0}^{n-2}\, \Lambda\{u_{n-1}\}\otimes
k\big\{t^a\kappa^dd\log(-p)\mid \hbox{$v_p(pa-d)=v$
and $P(a,d,v)$}
\big\} \\[4pt]
{}&\hsm\hsm&\oplus\; \Lambda\{d\log(-p)\}\otimes
k\{t^a\kappa^d\mid \hbox{$v_p(pa-d)\geq n-1$ and $P(a,d,n-1)$}
\big\}.  
\end{eqnarray*}
We know from Addendum~\ref{v_1} that the infinite cycle
$(-t\kappa^p)^m$, if not a boundary, represents the class $v_1^m$.
The smallest power $m_0$ such that $(-t\kappa^p)^{m_0}$ is a boundary
is $m_0=(p^n-1)/(p-1)$, $(-t\kappa^p)^{m_0} =
d^{2m_0-1}(u_{n-1}\kappa^{p^n})$. Hence, $v_1^m$ is nonzero, if
$m<m_0$, and $v_1^{m_0}$ is represented by an element of
$E_{s,2(p-1)m_0-s}^\infty$ with $s<-2m_0$. But these groups are all
zero, and therefore, so is $v_1^{m_0}$.

We next show that in $\hat E^*(C_{p^n},T)$,
$(-t\kappa)^{p^{s+1}}d\log(-p)$, $1\leq s<n-1$, represents
$dV^{n-s}(1)$, and that $(-t\kappa)^{p^n}d\log(-p)$, if not a
boundary, represents $dV(1)$. The latter follows from
Proposition~\ref{dlogpi} and Addendum~\ref{v_1}, since, by
Lemma~\ref{p>2},
$$dV(1) = d(\ul{-p}_n) = \ul{-p}_nd\log_n(-p) = V(1)d\log_n(-p).$$
To prove the former, we consider the map
$$F\:\bar\pi_1(\tate(C_{p^n},T)) \to \bar\pi_1(\tate(C_{p^{n-1}},T)),$$
which, by Lemma~\ref{hTR} and Proposition~\ref{drwp}, is a surjection
whose kernel is generated by $dV(1)$. Moreover, it takes $dV^{n-s}(1)$
to $dV^{n-1-s}(1)$ and the induced map of spectral sequences
$$F\:\hat E^*(C_{p^{n}},T) \to \hat E^*(C_{p^{n-1}},T)$$
takes $(-t\kappa)^{p^{s+1}}d\log(-p)$ to
$(-t\kappa)^{p^{s+1}}d\log(-p)$. The claim follows, inductively,
since the generator $dV(1)$ of the kernel of $F$ is represented by an
element of $\hat E_{m,1-m}^*(C_{p^n},T)$ with $m\leq-2p^n$.  

We now begin the proof of the statement of the proposition for $n$.
Suppose first that $2<r<2(p^n-1)/(p-1)$. The statement for $n-1$
implies that in the spectral sequence
$$\hat E^2({\Bbb T},T) = \Lambda\{d\log(-p)\}\otimes
S\{t^{\pm1},\kappa\} 
\Rightarrow \bar\pi_*\tate({\Bbb T},T),$$
the $d^r$-differential is multiplicatively generated from the stated
differentials on $t^{p^{v-1}}$ and from $d\log(-p)$ and $t\kappa^p$
being infinite cycles. Indeed, one shows inductively that the
canonical map induces an isomorphism
$$\gamma_{n-1}\:\Lambda\{u_{n-1}\}\otimes \hat E^r({\Bbb T},T) 
\xto{\sim} \hat E^r(C_{p^{n-1}},T).$$ 
We claim that for $r$ in the stated range, $d^r(u_n)$ is zero. To see
this, we consider the map of spectral sequences induced from $V\:
\bar\pi_*(\tate(C_{p^{n-1}},T)) \to \bar\pi_*(\tate(C_{p^n},T))$, 
$$V \: \hat E^r(C_{p^{n-1}},T) \to \hat E^r(C_{p^n},T).$$
The map of $E^2$-terms is given by the transfer map in Tate
cohomology. It follows that $u_n=V(u_{n-1})$, and hence,
$d^r(u_n)=V(d^r(u_{n-1}))$, which is zero for $r$ in the stated
range. We now conclude, by induction on $r$, that
$$\gamma_n\:\Lambda\{u_n\}\otimes\hat E^r({\Bbb T},T)
\xto{\sim}\hat E^r(C_{p^n},T)$$
is an isomorphism and that the $d^r$-differential is as stated. Before
we proceed, we note that in $\hat E^*({\Bbb T},T)$, and hence in
$\hat E^*(C_{p^n},T)$, the elements $t^i\kappa^jd\log(-p)$ are infinite
cycles. This follows, by arguments similar to~\cite[\S 5.3]{hm},
from the fact that the homotopy groups of $T$ with $\Zp$-coefficients
are concentrated in degree zero and in odd positive degree.

If $r=2(p^n-1)/(p-1)$, the possible nonzero differentials are
generated from
\begin{eqnarray*}
d^r(t^{p^{n-2}}) & = &\lambda_{n-1}\cdot
(t\kappa)^{\frac{p^n-1}{p-1}-1}t d\log(-p)\cdot t^{p^{n-2}}, \\
d^r(u_n) & =& \nu_n\cdot
(t\kappa)^{\frac{p^n-1}{p-1}-1}t d\log(-p)\cdot u_n,  
\end{eqnarray*}
where $\lambda_{n-1},\nu_n\in\Fp$. We first show that $\lambda_{n-1}$
is a unit. If $n=2$, the $k$-vector space
$\bar\pi_1(\tate(C_{p^2},T))$ is generated by the classes
$d\log_2(-p)$ and $dV(1)$. The former is represented by $d\log(-p)$
and the latter by an element of $\hat E_{m,1-m}^*(C_{p^2},T)$ with
$m\leq-2p^2$. Hence, the infinite cycle $(t\kappa)^pd\log(-p)$ 
must be hit by a differential, and this can happen only if $\lambda_1$
is a unit. If $n>2$, we consider $dV^2(1)\in\bar\pi_1(\tate(C_{p^n},T))$
which is represented by $(-t\kappa)^{p^{n-1}}d\log(-p)$. We know
that $v_1^{(p^{n-2}-1)/(p-1)}$ annihilates 
$1\in\bar\pi_*(\tate(C_{p^{n-2}},T))$ and hence also
$dV^2(1)\in\bar\pi_1(\tate(C_{p^n},T))$. Therefore,
$$(-t\kappa^p)^{(p^{n-2}-1)/(p-1)} \cdot (-t\kappa)^{p^{n-1}}d\log(-p)
=(-t\kappa)^{\frac{p^{n-1}\!-1}{p-1}-1}td\log(-p)\cdot t^{-p^{n-2}}$$
must be hit by a differential, and this can happen only if
$\lambda_{n-1}$ is a unit.

We next show that $\nu_n$ is zero. If not, then
$d^r(u_nt^{p^{n-2}c})=0$, for some $0<c<p$, and for degree reasons,
the next possible nonzero differential is
$$d^{2(\frac{p^n-1}{p-1}+p^{n-1}c)-1}(u_nt^{p^{n-2}c})
= \xi_n\cdot (t\kappa^p)^{\frac{p^{n-1}\!-1}{p-1}+p^{n-2}c} \cdot
t^{p^{n-1}(c+1)}.$$
But this must be zero, or else $v_1^{(p^{n-1}\!-1)/(p-1)+p^{n-2}c}$
would be zero. For degree reasons, the next possible nonzero
differential is $d^{2(p^{n+1}\!-1)/(p-1)}$. In particular, no
differential can hit $v_1^{(p^n-1)/(p-1)}$. So we must have
$\nu_n=0$.

The next possible differential is the stated one on $u_n$, and since
$v_1^{(p^n-1)/(p-1)}$ is zero, $\mu_n$ must be a unit. For degree
reasons, all further differentials vanish.
\enddemo 

We next prove Theorem~\ref{diff}, when $\mu_p\subset K$ and $n<v_p(e_K)$.

\specialnumber{5.5.5}\proclaim{Proposition} \label{nlessv}
If $\mu_p\subset K$ and if $n<v_p(e_K)${\rm ,} the
nonzero differentials in the spectral sequence $\hat
E^*(C_{p^n},T(A|K))$ are multiplicatively generated from
\begin{eqnarray*}
d^{2(\frac{p^{v+1}\!-1}{p-1})}(\pi_K^{p^v}) & =& -\lambda_v \cdot
(t\kappa)^{\frac{p^{v+1}-1}{p-1}-1}t\,d\log\pi_K\cdot\pi_K^{p^v},
\hskip5mm 0\leq v<n, \\
d^{2(\frac{p^{n+1}\!-1}{p-1})-1}(u_n) & =& \mu_n \cdot
(t\kappa)^{\frac{p^{n+1}\!-1}{p-1}-1}t,  
\end{eqnarray*}
and from $\tau_K${\rm ,} $\alpha_K${\rm ,} and $d\log\pi_K$ being infinite cycles.
\endproclaim

\demo{Proof} Since $n<v_p(e_K)$, Proposition~\ref{taualpha} and
Addendum~\ref{v_1} show that $\tau_K\alpha_K$ and $\tau_K\alpha_K^p$
are infinite cycles. Hence if $d^r\alpha_K$ is nontrivial then so is
$d^r(\alpha_K^p)$ contradicting the fact that $d^r$ is a derivation. It follows
that both $\alpha_K$ and $\tau_K$ are infinite cycles, and
$d\log\pi_K$ is an infinite cycle by Proposition~\ref{dlogpi}. Hence,
Theorem~\ref{diff} amounts to the statement above. 

Suppose first that $u_K'(0)$ is a unit. We prove the stated formula
for $d^r(\pi_K^{p^v})$ by induction on $0\leq v<n$. The basic case
$v=0$ follows from Proposition~\ref{d2=connes}. Now, assume that the
$d^r$-differential is as stated, for $2\leq r\leq 2(p^v-1)/(p-1)$, and
consider $2(p^v\!-1)/(p-1)<r\leq 2(p^{v+1}\!-1)/(p-1)$. We note that
$\pi_K^i(t\kappa^p)^jd\log\pi_K$ is equal to zero in $\hat
E^r(C_{p^n},T(A|K))$, if $v_p(i)=s<v$ and $j\geq (p^s\!-1)/(p-1)$.

By definition, $\tau_K=u_K^{(-n)}(\pi_K)^p\,t$, so that
$$t^{p^{v-1}} = u_K^{(v-n)}(\pi_K^{p^v})^{-1}\tau_K^{p^{v-1}},$$
and since $\tau_K$ is an infinite cycle, we find
$$d^r(t^{p^{v-1}})  = 
-\frac{u_K^{(v-n)}{}'(\pi_K^{p^v})}{u_K^{(v-n)}(\pi_K^{p^v})}t^{p^{v-1}}
\cdot d^r(\pi_K^{p^v}).$$
The first factor on the right is a unit in $\hat E^r$, for $r$ in the
stated range, and\break hence, we can evaluate $d^r(\pi_K^{p^v})$ from the
value of $d^r(t^{p^{v-1}})$, which is\break known by
Proposition~\ref{KeqKsub0}. It follows that $d^r(\pi_K^{p^v})$ is equal
to zero, if $r<\break 2(p^{v+1}\!-1)/(p-1)$. If $r=2(p^{v+1}\!-1)/(p-1)$, 
\begin{eqnarray*}
d^r(t^{p^{v-1}}) & = &\lambda_v \cdot
(t\kappa)^{\frac{p^{v+1}\!-1}{p-1}-1}td\log(-p)\cdot t^{p^{v-1}} \\
{} & =& \lambda_v\cdot 
\frac{u_K^{(-n)}{}'(\pi_K)\pi_K}{u_K^{(-n)}(\pi_K)}
(t\kappa)^{\frac{p^{v+1}\!-1}{p-1}-1}td\log\pi_K\cdot t^{p^{v-1}} \\
{} & =& \lambda_v\cdot
\frac{u_K^{(v-n)}{}'(\pi_K^{p^v})\pi_K^{p^v}}{u_K^{(v-n)}(\pi_K^{p^v})}
(t\kappa)^{\frac{p^{v+1}\!-1}{p-1}-1}td\log\pi_K\cdot t^{p^{v-1}}.  
\end{eqnarray*}
The second equation uses $-p=\pi_K^{e_K}\theta_K(\pi_K)^{-1}$ and
$\theta_K(\pi_K)=u_K(\pi_K)^{p-1}$, and the third follows from
Lemma~\ref{cartierlemma} since, as noted above,
$\pi_K^i(t\kappa^p)^jd\log\pi_K$ is equal to zero in $\hat
E^r(C_{p^n},T(A|K))$, if $v_p(i)=s<v$ and $j\geq
(p^{s+1}\!-1)/(p-1)$. The stated formula for $d^r(\pi_K^{p^v})$
follows. Similarly, we see that $d^r(\pi_K^{p^n})$ is equal to zero,
if $r<2(p^{n+1}\!-1)/(p-1)$, and the differential on $u_n$ follows
from Proposition~\ref{KeqKsub0}. For degree reasons, all further
differentials are zero.

To treat the general case, let $\Pi$ be the pointed monoid
$\{0,1,\pi,\pi^2,\dots\}$ with base point $0$. The choice of
uniformizer $\pi_K$ determines a map of ${\Bbb T}$-spectra
$$\rho_K\:T(W(k)|K_0)\wedge |N\s^\cy(\Pi)| \to T(A|K),$$
which is multiplicative with component-wise multiplication on the
left; compare Section~A.1 below. As a differential graded
$k$-algebra,
$$\bar\pi_*(T(W(k)|K_0)\wedge |N\s^\cy(\Pi)|)=
\Lambda\{d\log(-p),d\pi\}\otimes S\{\kappa,\pi\},$$
and the map of homotopy groups with $\Z/p$-coefficients induced from
$\rho_K$ is the unique map of differential graded $k$-algebras that is
$\bar\pi_*T(W(k)|K_0)$-linear and takes $\pi$ to $\pi_K$. We claim
that in the spectral sequence
$$\hat E^*(\Pi) = \hat E^*(C_{p^n},T(W(k)|K_0)\wedge |N\s^\cy(\Pi)|),$$ 
the nonzero differentials are generated multiplicatively from
$$d^{2(\frac{p^{v+1}-1}{p-1})}(\pi^{p^v})  =
-\lambda_v\cdot(t\kappa)^{\frac{p^{v+1}-1}{p-1}-1}t\cdot
\pi^{p^v-1}d\pi,\hskip10.3mm0\leq v<n,$$
from the differentials on the $t^{p^{v-1}}$, $1\leq v<n$, and the
differential on $u_n$ given by Proposition~\ref{KeqKsub0}, and from
$t\kappa^p$, $d\log(-p)$, $\pi^{p^n}$ and $\pi^{p^n-1}d\pi$ being
infinite cycles. This proves the proposition since
$\hat E^*(C_{p^n},T(A|K))$ is a module spectral sequence over
$\hat E^*(\Pi)$. 

To prove the claim, we choose a totally ramified extension $K/K_0$
with $\mu_p\subset K$ such that $n<v_p(e_K)$ and $u_K'(0)$ is a
unit. The proposition already has been established for $\hat
E^*(C_{p^n},T(A|K))$. As cyclic sets
$$N\s^\cy(\Pi)=\bigvee_{s\geq 0}N\s^\cy(\Pi;s),$$
where the $s$-th summand has $n$-simplices
$(\pi^{i_0},\dots,\pi^{i_n})$ with $i_0+\dots i_n=s$, and the 
spectral sequence $\hat E^*(\Pi)$ decomposes accordingly. It will
suffice to show that for $0\leq v\leq n$, the differentials in the
$p^v$-th summand spectral sequence,
\begin{eqnarray*}
\hat E^2(\Pi,p^v) & = \Lambda\{u_n,d\log(-p)\} \otimes
S\{t^{\pm1},\kappa\} \otimes k\{\pi^{p^v},\pi^{p^v-1}d\pi\} \\
{} & \Rightarrow\bar\pi_*(\tate(C_{p^n},T(W(k)|K_0)\wedge
|N\s^\cy(\Pi,p^v)|), \\
\end{eqnarray*}
are multiplicatively generated from the stated differentials on
$\pi^{p^v}$, the differentials on the $p$-powers of $t$, and the
differential on $u_n$, and from $d\log(-p)$ and $\pi^{p^v-1}d\pi$
being infinite cycles. We note that the map $\rho_{K*}\:\hat
E^2(\Pi,p^v) \to \hat E^2(C_{p^n},T(A|K))$ is a monomorphism. Indeed,
$\pi_K^{p^v}$ and
$$d\log(-p) =
\frac{u_K^{(-n)}{}'(\pi_K)\pi_K}{u_K^{(-n)}(\pi_K)}d\log\pi_K$$
are nonzero, since $p^v<e_K$ and since $u_K'(0)$ is a unit,
respectively. It follows, by induction on $r$, that $\rho_{K*}\:\hat
E^r(\Pi,p^v)\to \hat E^r(C_{p^n},T(A|K))$ is a monomorphism and that
the $d^r$-differential is as stated. For instance, $\pi^{p^v-1}d\pi$
is an infinite cycle because
$\rho_{K*}(\pi^{p^v-1}d\pi)=\pi_K^{p^v}d\log\pi_K$ is.
\enddemo

\demo{Proof  of Theorem~{\rm \ref{diff}}} Let $n\geq 1$ and $K$ be
given. We prove by induction on $q$ that the $d^q$-differential
in $\hat E^*(C_{p^n},T(A|K))$ is as stated. The basic case $q=2$
follows from Propositions~\ref{d2=connes} and~\ref{pi_*alpha}. So
assume the statement for $q-1$ and suppose first that
$2(p^s\!-1)/(p-1)<q\leq 2(p^{s+1}\!-1)/(p-1)$ with $s<n$. We recall
from Lemma~\ref{Er} that $\hat E^q=\hat E^q(C_{p^n},T(A|K))$ is given
by
\begin{eqnarray*}
\hat E^q &=  & \bigoplus_{v=1}^{s-1}
\Lambda\{u_n\}\otimes k\big\{\tau_K^a\pi_K^r
\alpha_K^dd\log\pi_K\mid v_p\{a,r,d\}_K\!=v,
d<\textstyle{\frac{p^{v+1}\!-1}{p-1}-1}\big\}\\
 &&\oplus\;\Lambda\{u_n,d\log\pi_K\}\otimes 
k\big\{\tau_K^a\pi_K^r\alpha_K^d\mid
v_p\{a,r,d\}_K\!\geq s\big\}. 
\end{eqnarray*}
Since the elements $\tau_K^a\pi_K^r\alpha^dd\log\pi_K$ are infinite
cycles, and since $d^q(u_n)$ is zero by Proposition~\ref{KeqKsub0}, it
suffices to evaluate $d^q(\tau_K^a\pi_K^r\alpha_K^d)$ with
$v_p\{a,r,d\}_K\geq s$. To this end, we find a totally ramified
extension
$$L=K[\pi_L]/(\pi_L^{e_{L/K}}+\pi_K\theta_{L/K}(\pi_L))$$
such that $n<v_p(e_L)$ and such that the map
$$\iota_*\:\hat E_{*,t}^q(C_{p^n},T(A|K))
\to \hat E_{*,t}^q(C_{p^n},T(B|L))$$
is a monomorphism, for $t\geq q-1$. Since the differential structure
of the right-hand spectral sequence is known from
Proposition~\ref{nlessv}, this allows us to evaluate the
$d^q$-differential in the spectral sequence on the left. We consider
the extension $L/K$ with $e_{L/K}=p^{n+1}$ and $\theta_{L/K}(x)=x-1$,
and recall from~(\ref{E2iota}) that the map of $E^2$-terms is given by
\begin{eqnarray*}
\iota_*(\tau_K^a\pi_K^r\alpha_K^d) & =&
(1-\pi_L)^{-\{a,r,d\}_K}\tau_L^a\pi_L^{e_{L/K}r}\alpha_L^d, \\
\iota_*(d\log\pi_K) & =&\pi_L(1-\pi_L)^{-1}d\log\pi_L.  
\end{eqnarray*}
Hence, the induced map of $E^q$-terms takes
$\tau_K^a\pi_K^r\alpha_K^d$ with $v_p\{a,r,d\}_K\geq s$ to
$$(1-\pi_L^{p^s})^{-p^{-s}\{a,r,d\}_K} \cdot
\tau_L^a\pi_L^{e_{L/K}r}\alpha_L^d,$$
and $\tau_K^a\pi_K^r\alpha_K^dd\log\pi_K$ with $v_p\{a,r,d\}_K\geq s$
and $d\geq(p^s\!-1)/(p-1)-1$ to
$$(1-\pi_L^{p^s})^{-p^{-s}\{a,r,d\}_K}
\pi_L^{p^s}(1-\pi_L^{p^s})^{-1} \cdot 
\tau_L^a\pi_L^{e_{L/K}r}\alpha_L^d d\log\pi_L,$$
where the latter statement uses Lemma~\ref{cartierlemma} and that
$\pi_L^i\cdot\tau_L^a\pi_L^{e_{L/K}r}\alpha_L^dd\log\pi_L$ is equal to 
zero in $\hat E^q(C_{p^n},T(B|L))$, if $v_p(i)<s$. It is clear that
this map is a monomorphism in the stated range. Indeed, $r\leq e_K-1$
and $e_L=e_{L/K}e_K$, and therefore, $e_{L/K}r+p^s\leq
e_L-p^{n+1}+p^s\leq e_L-1$. 
It follows immediately from Proposition~\ref{nlessv} that
$d^q(\iota_*(\tau_K^a\pi_K^r\alpha_K^d))$ vanishes, if
$q<2(p^{s+1}\!-1)/(p-1)$, and a straightforward calculation shows that
$$d^q(\iota_*(\tau_K^a\pi_K^r\alpha_K^d)) =
\iota_*(-\lambda_s \cdot p^{-s}\{a,r,d\}_K \cdot 
(t\kappa)^{\frac{p^{s+1}\!-1}{p-1}-1}td\log\pi_K \cdot
\tau_K^a\pi_K^r\alpha_K^d),$$
if $q=2(p^{s+1}\!-1)/(p-1)$. Since $\iota_*$ is a monomorphism, we
conclude that 
\begin{eqnarray*}
  d^q(\tau_K^a\pi_K^r\alpha_K^d)  &=&
-\lambda_s \cdot p^{-s}\{a,r,d\}_K \cdot 
(t\kappa)^{\frac{p^{s+1}\!-1}{p-1}-1}td\log\pi_K \cdot
\tau_K^a\pi_K^r\alpha_K^d \\
 & =& -\lambda_s \cdot p^{-s}\{a,r,d\}_K \cdot 
u_K^{(s-n)}(\pi_K^{p^s})^{-p}\\&&
\cdot\ (\tau_K\alpha_K)^{\frac{p^{s+1}\!-1}{p-1}-1}\tau_Kd\log\pi_K
\cdot \tau_K^a\pi_K^r\alpha_K^d  
\end{eqnarray*}
as desired. Finally, an analogous argument shows that
$d^q(\tau_K^a\pi_K^r\alpha_K^d)$ is equal to zero, if
$2(p^n\!-1)/(p-1)<q<2(p^{n+1}\!-1)/(p-1)$, and the stated
differential on $u_n$ follows from Proposition~\ref{KeqKsub0}. All
further differentials vanish for degree reasons.
\enddemo

5.6.\quad We conclude this section with a proof of the following
result, which was used in the proof of Proposition~\ref{pi_1exact} above
for $n>3$.

\specialnumber{5.6.1} \proclaim{Lemma}\label{frobsurj}
For all $i\geq 0${\rm ,} the Frobenius induces a
surjection{\rm ,}
$$F\:\TR^n_{2i+1}(A|K;p) \twoheadrightarrow \TR^{n-1}_{2i+1}(A|K;p).$$
\endproclaim

\demo{Proof} For $i>0$, the group $\TR^n_i(A|K;p)$ is a sum of a
uniquely divisible group and a $p$-torsion group of bounded height.
Indeed, this is true when $n=1$, and the general case follows
inductively from the cofibration sequence
$${}_h\!\TR^n(A|K;p)\xto{N}\TR^n(A|K;p)\xto{R}\TR^{n-1}(A|K;p)$$
and the spectral sequence~\eqref{homotopyorbitss}. Since $FV=p$, the
Frobenius induces a surjection of uniquely divisible summands. Hence,
it suffices to prove that the statement of the lemma holds after
$p$-completion. And, by Addendum~\ref{gammahatadd}, we may instead
show that the canonical map 
$$\gamma_n\:\pi_{2i+1}(\coborel({\Bbb T},T(A|K)),\Zp)
\to\pi_{2i+1}(\coborel(C_{p^n},T(A|K)),\Zp)$$
is surjective. To this end, we consider the spectral sequences 
\begin{eqnarray*}
E^2_{s,t}({\Bbb T})&=&H^{-s}(BS^1,\pi_t(T(A|K),\Zp))
\Rightarrow\pi_{s+t}(\coborel({\Bbb T},T(A|K)),\Zp),\\
E^2_{s,t}(C_{p^n})&=&H^{-s}(BC_{p^n},\pi_t(T(A|K),\Zp))
\Rightarrow\pi_{s+t}(\coborel(C_{p^n},T(A|K)),\Zp), 
\end{eqnarray*}
both of which are strongly convergent second quadrant spectral
sequences. This means that the filtration of
$\pi_*(\coborel(G,T(A|K)),\Zp)$ associated with the spectral sequence
$E^*(G)$ is complete and separated and that there is a canonical
isomorphism
$$\gr^s\pi_{s+t}(\coborel(G,T(A|K)),\Zp)
\cong E^\infty_{s,t}(G).$$
It will therefore suffice to show that
$$\gr^s(\gamma_n)\:\gr^s\pi_{2i+1}(\coborel({\Bbb T},T(A|K)),\Zp)
\to\gr^s\pi_{2i+1}(\coborel(C_{p^n},T(A|K)),\Zp)$$
is a surjection for all $s\leq 0$ and $i\geq 0$. The induced map of
$E^2$-terms is given by the map on cohomology induced from the
inclusion $C_{p^n}\to{\Bbb T}$, and hence, is surjective for $s$
even. Moreover, by Remark~\ref{concodd}, $\pi_*(T(A|K),\Zp)$ is
concentrated in odd degrees with the exception of $\pi_0(T(A|K),\Zp)$,
and hence, the nonzero differentials in the spectral sequence
$E^r({\Bbb T})$ must originate on the line $t=0$. It follows that
for $s$ even and $t>0$, the map
$$\gamma_{n*}\:E^r_{s,t}({\Bbb T})\to E^r_{s,t}(C_{p^n})$$
is surjective for all $2\leq r\leq\infty$. (Since these groups do not
support nonzero differentials, they are stable for $r>s$.) Since only
the groups $E^r_{s,t}(C_{p^n})$ with $s$ even and $t>0$ can contribute
to $\pi_{2i+1}(\coborel(C_{p^n},T(A|K)),\Zp)$, this shows that the
map $\gr^s(\gamma_n)$ is indeed surjective.
\enddemo

\section{The pro-system $\TR\cs_*(A|K;p,\Z/p^v)$} \label{main}

6.1.\quad In this  section, we prove the main theorem of this
work. Suppose that $\mu_{p^v}\subset K$ such that we have the maps
$$\Sigma^\infty B\mu_{p^v+} \xto{\det} K(K) \xto{\tr} \TR^n(A|K;p).$$
Since $p$ is odd, the Bockstein gives an isomorphism
$$\pi_2(\Sigma^\infty B\mu_{p^v+},\Z/p^v)
\xto{\sim} {}_{p^v}\pi_1(\Sigma^\infty B\mu_{p^v},\Z/p^v)
\stackrel{\sim}{\longleftarrow} \mu_{p^v},$$
and hence, these maps induce
$$\mu_{p^v} \to K_2(K,\Z/p^v) \xto{\tr} \TR_2^n(A|K;p,\Z/p^v)
=\pi_2(\TR^n(A|K;p),\Z/p^v).$$
It follows that there is a canonical map of log Witt complexes
$$W\s\,\omega_{(A,M)}^*\otimes S_{\Z/p^v}(\mu_{p^v})
\to \TR_*\cs(A|K;p,\Z/p^v),$$
where on the second tensor factor on the left, the maps $R$, $F$ and
$V$ act as the identity and the differential $d$ acts as zero. We
recall from Theorem~\ref{pionetr} that this map is an isomorphism in
degrees $0$ and $1$.

By Addendum~\ref{gammahatadd} the map
$$\hat\Gamma_{A|K}\: \TR^n_*(A|K;p,\Z/p^v)
\to \pi_*(\tate(C_{p^n},T(A|K)),\Z/p^v)$$
is an isomorphism in nonnegative degrees. The groups on the right,
for $v=1$, are given by the spectral sequence $\hat E^*=\hat
E^*(C_{p^n},T(A|K))$, which we evaluated in Theorem~\ref{diff} above.
The result is that
\begin{eqnarray*}
\hat E^\infty& = &\,\bigoplus_{v=1}^{n-1}
k\big\{u_n^{\e}\tau_K^a\pi_K^r
\alpha_K^dd\log\pi_K\big|\,v_p\{a,r,d\}_K=v,\,
d<\textstyle{\frac{p^{v+1}-1}{p-1}-1}\big\}\\
{}&&\oplus\;k\big\{\tau_K^a\pi_K^r\alpha_K^d(d\log\pi_K)^\e
\,\big|\,v_p\{a,r,d\}_K\geq n,\,
d<\textstyle{\frac{p^{n+1}-1}{p-1}-1}\big\}, 
\end{eqnarray*}
where $a\in\Z$, $d\in\N_0$, $\e\in\{0,1\}$, and $0\leq r<e_K$, and
where
$$\{a,r,d\}_K=(pa-d)e_K/(p-1)+r.$$
We call the basis of $\hat E^\infty$ as a $k$-vector space exhibited
here the {\it standard} basis.

\specialnumber{6.1.1}\proclaim{Proposition} \label{dimension}
If $\mu_p\subset K$ or if $K=K_0$ then
$\TR_q^n(A|K;p,\Z/p)$ is an $ne_K$\/{\rm -}\/dimensional $k$\/{\rm -}\/vector space{\rm ,} for
all $q\geq 0$.
\endproclaim

\demo{Proof} We fix a total degree $q$ and evaluate the
cardinality of the standard basis of $\hat E^\infty(C_{p^n},T(A|K))$.
An element of the standard basis is in total degree $q=2m+\e$ if and
only if $d-a=m$. We let $v=v_p\{a,r,d\}_K$ and note that
$$\{a,r,d\}_K=de_K+r - pe_Km/(p-1).$$
Hence, the elements of the standard basis of
$\hat E^\infty(C_{p^n},T(A|K))$ in total degree $q$ are
indexed by integers $1\leq v\leq n$, $0\leq r<e_K$ and $d\geq 0$ such
that either $1\leq v<n$ and $v_p(de_K+r-pe_Km/(p-1))=v$ and
$0\leq de_K+r<(\textstyle{\frac{p^{v+1}-1}{p-1}-1})e_K$ or $v=n$ and
$v_p(de_K+r-pe_Km/(p-1))\geq v$ and
$0\leq de_K+r<(\textstyle{\frac{p^{n+1}-1}{p-1}-1})e_K$. But these
requirements are equivalent to the requirement that for all $1\leq
v\leq n$, $de_K+r$ is congruent to $pe_Km/(p-1)$ modulo $p^v$ and
$$(\textstyle{\frac{p^v-1}{p-1}-1})e_K\leq
de_K+r<(\textstyle{\frac{p^{v+1}-1}{p-1}-1})e_K
= (\textstyle{\frac{p^v-1}{p-1}-1})e_K+p^ve_K.$$
It is clear that for each value of $1\leq v\leq n$, there are $e_K$
pairs $(d,r)$ which satisfy this requirement. Hence, the dimension is
equal to $ne_K$ as stated.
\enddemo

\specialnumber{6.1.2} \proclaim{Lemma}\label{bottproduct}
Suppose that the class
$\xi\in\pi_*(\tate(C_{p^n},T(A|K)))$ is represented in $\hat 
E^\infty(C_{p^n},T(A|K))$ by the element
$u_n^{\e}\tau_K^a\pi_K^r\alpha_K^d(d\log\pi_K)^{\delta}$. Then 
the product $b_n\cdot\xi$ is represented by
$\pm u_n^{\e}\tau_K^{a+a'}\pi_K^{r'}\alpha_K^{d+a'+1}
(d\log\pi_K)^{\delta}${\rm ,} where\break $r+e_K/(p-1)=a'e_K+r'$ and $0\leq r'<e_K$.
\endproclaim

\demo{Proof} We show that the map induced from multiplication by $b_n$,
$$b_n\:\hat E^3(C_{p^n},T(A|K)) \to \hat E^3(C_{p^n},T(A|K)),$$
is given by the stated formula. It suffices to consider the case
$n=1$. Indeed,
\begin{eqnarray*}
 & &F^{n-1}\: \hat E_{s,t}^3(C_{p^n},T(A|K)) 
\to \hat E_{s,t}^3(C_p,T(A|K)), \\
&&V^{n-1}\:  \hat E_{s,t}^3(C_p,T(A|K)) 
\to \hat E_{s,t}^3(C_{p^n},T(A|K)),  
\end{eqnarray*}
are isomorphisms for $s$ even and odd, respectively, and commute with
multiplication by the Bott element, since $F^{n-1}(b_n)=b_1$. Suppose
first that $v_p(e_K)>1$ such that
$$\hat E^3(C_p,T(A|K)) = \Lambda\{u_1,d\log\pi_K\}\otimes
S\{\tau_K^{\pm1},\pi_K^p,\alpha_K\}/(\pi_K^{e_K}).$$
It will suffice to prove that $b_1\cdot\pi_K^r$ is equal to
$\pm\tau_K^{a'}\pi_K^{r'}\alpha_K^{a'+1}$. This follows from the
``multiplicative extension'' $\pi_K^{e_K}=-\tau_K\alpha_K$. More
precisely, Proposition~\ref{taualpha} shows that the elements
$\pi_K^p$ and $-\tau_K\alpha_K$ represent the classes
$\hat\Gamma_{A|K}(\pi_K)$ and $\hat\Gamma_{A|K}(\pi_K^{e_K/p})$,
respectively. We also recall from~\eqref{bottelement} that the element
$-\pi_K^{e_K/(p-1)}\alpha_K$ represents the Bott element $b_1$. But
$\alpha_K$ survives the spectral sequence and represents a homotopy
class, say, $\tilde\alpha_K$. Hence, $-\pi_K^{e_K/(p-1)}\alpha_K$ also
represents the class
$-\hat\Gamma_{A|K}(\pi_K^{e_K/p(p-1)})\tilde\alpha_K$. The claim
follows, if $v_p(e_K)>1$. In general, we pick a totally ramified 
extension $L/K$ such that $v_p(e_L)>1$ and such that the map
$$\iota_*\:\hat E^3(C_p,T(A|K)) \to \hat E^3(C_p,T(B|L))$$
is a monomorphism.
\enddemo

We note that multiplication by $b_n$ preserves the symbol
$$\{a,r,d\}_K=\{a+a',r',d+a'+1\}_K,$$
and that the class $b_n^q$ is represented by
$\pm\tau_K^{q_1}\pi_K^{q_0e_K/(p-1)}\alpha_K^{q_1+q}$ with
$q=\break q_1(p-1)+q_0$ and $0\leq q_0<p-1$.

\specialnumber{6.1.3} \proclaim{Lemma}\label{image}
An element of the standard basis of
$\hat E^\infty(C_{p^n},T(A|K))$ represents a homotopy class in the image of
the composite
$$W_n\,\omega_{(A,M)}^*\otimes S_{\Z/p}(\mu_p)
\to \TR^n_*(A|K;p,\Z/p)
\to \bar\pi_*\tate(C_{p^n},T(A|K))$$
if and only if $\{a,r,d\}_K\geq 0$.
\endproclaim

\demo{Proof} The map of the statement is an isomorphism in degrees $0$
and $1$ by Theorem~\ref{pionetr} and Addendum~\ref{gammahatadd}. 
Indeed, in these dimensions $\{a,r,d\}_K$ is automatically
nonnegative since $a=d$. We must thus show that for all $q\geq 0$ and
$\e=0,1$, the map
$$\bigoplus_{s\leq 0}\hat E_{s,\e-s}^\infty(C_{p^n},T(A|K))\to
\bigoplus_{s\leq 0}\hat E_{s,2q+\e-s}^\infty(C_{p^n},T(A|K))$$
induced by multiplication by the $q$-th power of the Bott element is a
surjection onto the stated subspace. Suppose for example that a
homotopy class is represented in the spectral sequence by the element
$\tau_K^a\pi_K^r\alpha_K^{a+q}$ and write $r-qe_K/(p-1)=-a_0e_K+r_0$
with $0\leq r_0<e_K$. The \pagebreak requirement $\{a,r,a+q\}_K\geq 0$ is then
equivalent to $a_0\leq a$, and by Lemma~\ref{bottproduct}
$$b^q\cdot\tau_K^{a-a_0}\pi_K^{r_0}\alpha_K^{a-a_0}
=\pm\tau_K^a\pi_K^r\alpha_K^{a+q}.$$
The other elements of the standard basis are treated similarly.
\enddemo

\specialnumber{6.1.4}\proclaim{Theorem} \label{mainthm}
Suppose $K$ contains the $p$\/{\rm -}\/th roots of
unity. Then the canonical map is a pro\/{\rm -}\/isomorphism\/{\rm :}\/
$$W\s\,\omega_{(A,M)}^*\otimes S_{\Z/p}(\mu_p)
\xto{\sim} \TR\cs_*(A|K;p,\Z/p).$$
\endproclaim

\demo{Proof} Let $E\s^*$ denote the pro-system on either side of the
map in the statement. The standard filtration, given by
$$\Fil^sE_n^*=V^sE^*_{n-1}+dV^sE_{n-1}^*,$$
is a descending filtration with $s\geq 0$. The filtration has length
$n$ in level $n$, i.e. $\Fil^nE_n^*$ is trivial. The map of the
statement clearly preserves the filtration. We show that for all
$q\geq 0$, there exists $N\geq 1$ such that for all $n\geq 1$ and
$0\leq s<n-N$, the canonical map
$$\gr^s(W_n\,\omega_{(A,M)}^*\otimes S_{\Z/p}(\mu_p))_i
\to\gr^s\TR^n_i(A|K;p,\Z/p)$$
is an isomorphism when $0\leq s<n-N$. Since the structure maps in the
pro-systems preserve the standard filtration, the theorem follows. 

We have already proved that the map of the statement is an isomorphism
in degrees $0$ and $1$. Hence, it suffices to show that for all $q\geq
0$, there exists $N\geq 1$ such that for all $n\geq 1$, $0\leq s<n-N$
and $\e=0,1$, multiplication by the $q$-th power of the Bott element
induces an isomorphism
$$\gr^s\TR^n_{\e}(A|K;p,\Z/p)
\xto{\sim}\gr^s\TR^n_{2q+\e}(A|K;p,\Z/p).$$
We claim that any $N\geq 1$ with $p(q+1)e_K/(p-1)<p^N$ will do.

For surjectivity we use   Lemma~\ref{image}. Consider an element of
the standard basis in degree $2q+\e$ with symbol $\{a,r,d\}_K$. Since
$d\geq 0$ and $d=a+q$, we have $a\geq -q$, and hence
\begin{eqnarray*}
\{a,r,d\}_K&=&ae_K-qe_K/(p-1)+r\\
{}&\geq& -pqe_K/(p-1)+r>-p^N.\\
\end{eqnarray*}
Therefore, if $v_p\{a,r,d\}_K\geq N$ we have $\{a,r,d\}_K\geq 0$. It
follows that multiplication by the $q$-th power of the Bott element
induces a surjection of all summands in $\hat E^\infty(C_{p^n},T(A|K))$
except for the summands with $v<N$. But these summands all represent
homotopy classes of filtration greater than or equal to $n-N$. Indeed,
by Proposition~\ref{connesandtate}
\begin{eqnarray*}
V^s(u_{n-s}\tau_K^a\pi_K^r\alpha_K^dd\log\pi_K)&
=&u_n\tau_K^a\pi_K^r\alpha_K^dd\log\pi_K,\\
d(u_n\tau_K^a\pi_K^r\alpha_K^dd\log\pi_K)&
=&\tau_K^a\pi_K^r\alpha_K^dd\log\pi_K. 
\end{eqnarray*}
Thus elements of the standard basis with $\{a,r,d\}_K<N$ are either in
the image of $V^{n-N}$ or  of $dV^{n-N}$.

To prove injectivity, we first note that for an element of the
standard basis of $\hat E^\infty(C_{p^n},T(A|K))$ in total degree $2q+\e$,
the requirement that
$$0\leq d<\frac{p^{v+1}-1}{p-1}-1$$
is equivalent to the requirement that
$$r-\frac{pqe_K}{p-1}\leq\{a,r,d\}_K<-\frac{pqe_K}{p-1}
+e_K\frac{p^{v+1}-1}{p-1}+r-e_K.$$
We show that $v_p\{a,r,d\}_K=v\geq N$ and
$\{a,r,d\}_K<e_K(p^{v+1}-1)/(p-1)$ implies that
$$\{a,r,d\}_K<-\frac{pqe_K}{p-1}
+e_K\frac{p^{v+1}-1}{p-1}+r-e_K.$$
Indeed, the largest integer which is both congruent to zero modulo
$p^v$ and smaller that $e_K(p^{v+1}-1)/(p-1)$ is
$e_Kp^{v+1}/(p-1)-p^v$. Thus $\{a,r,d\}_K\leq e_Kp^{v+1}/(p-1)-p^v$,
and it suffices to check that
$$e_Kp^{v+1}/(p-1)-p^v<-\frac{pqe_K}{p-1}+e_K\frac{p^{v+1}-1}{p-1}+r-e_K.$$
But this is equivalent to the inequality
$$p^v>\frac{p(q+1)e_K}{p-1}-r,$$
which is satisfied for $n<N$. This shows that the map induced by
multiplication by the $q$-th power of the Bott element induces a
monomorphism of all summands in $\hat E^\infty(C_{p^n},T(A|K))$ except for
the summands with $v< N$. The theorem follows.
\enddemo 

\demo{Proof of Theorem~{\rm \ref{thmC}}}
The proof is by induction on
$v$; the basic case $v=1$ is Theorem~\ref{mainthm}. In the induction
step, we write $q=2s+\e$ with $0\leq\e\leq 1$ and consider the diagram
of pro-abelian groups
$$\begin{array}{ccccc}
{W\s\,\omega_{(A,M)}^{\e}\otimes\mu_{p^{v-1}}^{\otimes s}}
&\longrightarrow& 
{W\s\,\omega_{(A,M)}^{\e}\otimes\mu_{p^v}^{\otimes s}}
&\longrightarrow& 
{W\s\,\omega_{(A,M)}^{\e}\otimes\mu_p^{\otimes s}}
\\[4pt]
\scs{\sim} &&\big\downarrow &&\scs{\sim} \\[4pt]
{\TR_q\cs(A|K;p,\Z/p^{v-1})} &\longrightarrow&  
{\TR_q\cs(A|K;p,\Z/p^v)} &\longrightarrow& 
{\TR_q\cs(A|K;p,\Z/p),}\end{array}$$
where, inductively, the right- and left-hand vertical maps are
pro-isomorphisms. The lower sequence is exact at the middle. Hence, it
will suffice to show that the upper horizontal sequence is a
short-exact sequence of pro-abelian groups. Clearly, we can assume
that $s=0$. If $\e=0$, the sequence is exact since $W_n(A)$ is torsion
free, for all $n\geq 1$. (This does not use the fact that $\mu_{p^v}\subset
K$.) If $\e=1$, only the injectivity of the left-hand map requires
proof. To this end, we consider the diagram
$$\begin{array}{ccccc}
{W\s(A)\otimes\mu_p}  &\longrightarrow&
{W\s\,\omega_{(A,M)}^1\otimes\Z/p^{v-1}}&  \longrightarrow &
{W\s\,\omega_{(A,M)}^1\otimes\Z/p^v}  \\[4pt]
 \scs{\sim}&& \scs{\sim}&&\big\downarrow\\[4pt]
{\TR_2\cs(A|K;p,\Z/p)} &\shtck{\beta} &
{\TR_1\cs(A|K;p,\Z/p^{v-1})} &\longrightarrow &
{\TR_1\cs(A|K;p,\Z/p^v),} \end{array}
 $$
where the left-hand and middle vertical maps are pro-isomorphisms by
induction, and where the lower sequence is exact. It will suffice to
show that the upper left-hand horizontal map is zero. But this map
takes $x\otimes\zeta$ to $xd\log\s\zeta$, and since $\zeta$ has a
$p^{v-1}$st   root, $d\log\s\zeta$ is divisible by $p^{v-1}$.
\enddemo  

\specialnumber{6.1.5}\numbereddemo{{R}emark} It follows from Theorem~\ref{thmC} that if
$\mu_{p^v}\subset K$, the map
$$W\s(A)\otimes\mu_{p^v}\xto{\sim}{}_{p^v}W\s\,\omega_{(A,M)}^1,$$
which takes $x\otimes\zeta$ to $xd\log\s{\zeta}$, is a
pro-isomorphism. It would be desirable to have an algebraic proof
of this fact.
\enddemo

\specialnumber{6.1.6}\proclaim{Theorem} \label{maintc}
There are natural isomorphisms{\rm ,} for
$s\geq 0$\/{\rm :}\/
\begin{eqnarray*}
\TC_{2s}(A|K;p,\Z/p) & \cong& H^0(K,\mu_p^{\otimes s})
\oplus H^2(K,\mu_p^{\otimes(s+1)}),\\
\TC_{2s+1}(A|K;p,\Z/p) & \cong& H^1(K,\mu_p^{\otimes(s+1)}). 
\end{eqnarray*}
\endproclaim

\demo{Proof} Since the extension $K(\mu_p)/K$ is tamely ramified, we
may assume that $\mu_p\subset K$. Indeed, Theorem~\ref{tamedescent}
shows that the canonical map
$$\TC_*(A|K;p,\Z/p) \xto{\sim}
\TC_*(A(\mu_p)|K(\mu_p);p,\Z/p)^{\Gal(K(\mu_p)/K)}$$
is an isomorphism, and the analogous statement holds for
$H^*(K,\mu_p^{\otimes s})$. If $\mu_p\subset K$, Theorem~\ref{mainthm}
shows that for $s\geq 0$ and $0\leq\e\leq 1$, the canonical map
$$\TC_{\e}(A|K;p,\Z/p)\otimes\mu_p^{\otimes s} \xto{\sim}
\TC_{2s+\e}(A|K;p,\Z/p)$$
is an isomorphism, and hence, it suffices to prove the statement in
degrees $0$ and $1$. 

In degree one, the cyclotomic trace induces an isomorphism
$$K^\times/K^{\times p}=K_1(K,\Z/p)
\xto{\sim} \TC_1(A|K;p,\Z/p),$$
and by Kummer theory, the left-hand side is $H^1(K,\mu_p)$,
\cite[p.~155]{serre}. In degree zero, we use the fact that Addendum
\ref{localizationtc} gives an exact sequence
$$0\to\TC_0(A;p,\Z/p) \to \TC_0(A|K;p,\Z/p) \to \TC_{-1}(k;p,\Z/p)
\to 0.$$
The left-hand term is naturally isomorphic to $\Z/p=K_0(A,\Z/p)$ by\break
\cite[Th.\ D]{hm}, and the left-hand map has a natural retraction
given by
$$\TC_0(A|K;p,\Z/p) \to \TR_0(A|K;p,\Z/p)^F = \Z/p.$$
It remains to show that the right-hand term in the sequence is
naturally isomorphic to $H^2(K,\mu_p)$. We recall from
\cite[p.~186]{serre} the natural short exact sequence
$$0\to H^2(k,\mu_p)\to H^2(K,\mu_p)\to H^1(k,\Z/p)\to 0.$$
Since $k$ is perfect, the left-hand term vanishes,
\cite[p.~157]{serre}. Let $\bar k$ be an algebraic closure of $k$.
The normal basis theorem shows that $H^i(k,\bar k)$ vanishes for
$i>0$, and hence the cohomology sequence associated with the sequence
$$0\to\Z/p\to\bar k\xto{\varphi-1}\bar k\to 0$$
gives a natural isomorphism $k_\varphi\xto{\sim}H^1(k,\Z/p)$. Finally,
since $k$ is perfect, the restriction induces a natural isomorphism
\vglue12pt
\hfill ${\displaystyle \TC_{-1}(k;p,\Z/p)=W(k)_F/pW(k)_F\xto{\sim}k_\varphi.
}$
\enddemo
\vglue9pt

\specialnumber{6.1.7}\numbereddemo{{R}emark} If $\mu_p\subset K$, we can also give the following
noncanonical description of the groups $\TC_*(A|K;p,\Z/p)$. Let
$\zeta\in\mu_p$ be a generator, let $b=b_{\zeta}$ be the corresponding
Bott element, and let $\pi=\pi_K\in A$ be a uniformizer. Then for
$s\geq 0$,
\begin{eqnarray*}
\TC_{2s}(A|K;p,\Z/p)  & =& \Z/p\cdot b^s \oplus
k_\varphi\cdot\partial(d\log\pi\cdot b^s),\\
\TC_{2s+1}(A|K;p,\Z/p) & = &\Z/p\cdot b^sd\log\s\pi \oplus
k_\varphi\cdot\partial(b^{s+1})\oplus k^{e_K}, 
\end{eqnarray*} 
where $k_\varphi$ is the cokernel of $1-\varphi\:k\to k$, $e_K$ is the
ramification index, and $\partial$ is the boundary homomorphism in the
long-exact sequence
$$\cdots \xto{\partial} \TC_q(A|K;p,\Z/p) \to
\TR_q(A|K;p,\Z/p) \xto{1-F}
\TR_q(A|K;p,\Z/p) \xto{\partial} \dots.$$
The summand $k^{e_K}$ in the second line maps injectively to the
kernel of $1-F$, the inclusion
$$\eta\:k^{e_K}=\bigoplus_{i=0}^{e_K-1}k\to\TR_{2s+1}(A|K;p,\Z/p)$$
given, on the $i$-th summand, by
$$\eta_i(a)=\sum_{v\geq 0}a^{p^{-v}(\frac{p^{v+1}-1}{p-1})}
u_K(\ul{\pi})^{-p}dV_\pi^v(\ul{\pi}^i)\cdot b^s
+\sum_{v>0}F^v(au_K(\ul{\pi})^{-p}d(\ul{\pi}^i))\cdot b^s.$$
The sum on the right is finite and the sum on the left converges.
\enddemo  

We shall need a special case of the Thomason-Godement construction of
the hyper-cohomology spectrum associated with a presheaf of spectra on
a site, \cite[\S3]{gh}. Suppose that $F$ is a functor which to every
finite subextension $L/K$ in an algebraic closure $\bar K/K$ assigns a
spectrum $F(L)$. For the purpose of this paper, we   write
\begin{equation}\label{Fet}
F^{\hbox{{\eightpoint \'et}}}(K)={\displaystyle\mathop{\rm holim}_{\lrar\atop{L/K}}}{\coborel(G_{L/K},F(L))}.
\speqnu{6.1.8}
\end{equation}
There is a natural strongly convergent spectral sequence
\begin{equation}\label{dss}
E^2_{s,t}=H^{-s}(K,{\displaystyle\mathop{\rm  lim}_{\longrightarrow\atop{L/K}}}{\pi_tF(L)})
\Rightarrow\pi_{s+t}F^{\hbox{{\eightpoint \'et}}}(K), \speqnu{6.1.9}
\end{equation}
which is obtained by passing to the limit from the spectral sequences
for the group cohomology spectra
$$E^2_{s,t}=H^{-s}(G_{L/K},\pi_tF(L))
\Rightarrow\pi_{s+t}\coborel(G_{L/K},F(L)).$$
Indeed, filtered colimits are exact so we get a spectral sequence with
abutment
$${\displaystyle\mathop{\rm lim}_{\longrightarrow\atop{L/K}}}{\pi_*\coborel(G_{L/K},F(L))}
\xto{\sim}\pi_*F^{\hbox{{\eightpoint \'et}}}(K),$$
and the identification of the $E^2$-term follows from the isomorphism
\begin{eqnarray*}
{\displaystyle\mathop{\rm lim}_{\longrightarrow\atop{L/K}}}{H^*(G_{L/K},\pi_*F(L))}
&\xto{\sim}&{\displaystyle\mathop{\rm
lim}_{\longrightarrow\atop{L/K}}}\, {H^*(G_{L/K},({\displaystyle\mathop{\rm
lim}_{\longrightarrow\atop{N/L}}}\, {\pi_*F(N)})^{G_L})} \\ {}& =& H^*(K,{\displaystyle\mathop{\rm
lim}_{\longrightarrow\atop{N/K}}}\, {\pi_*F(N)}). 
\end{eqnarray*}
This isomorphism, which can be found in \cite[\S2 Prop.\
8]{serre1}, is a special case of the general fact that on a site with
enough points, the Godement construction of a presheaf calculates the
sheaf cohomology of the associated sheaf.

\specialnumber{6.1.10}\proclaim{Theorem} The canonical map is an isomorphism in degrees $\geq
1$\/{\rm :}\/
$$\gamma_K\:K_*(K,\Z/p^v)\to K_*^{\hbox{{\eightpoint \'et}}}(K,\Z/p^v).$$
\endproclaim

\demo{Proof} It suffices to consider the case $v=1$. In the diagram
$$\begin{array}{ccc}
{K(K)}&\stck{\gamma_K}&
{K^{\hbox{{\eightpoint \'et}}}(K)}  \\[4pt]
 \scs{\tr} && \scs{\tr} \\[4pt]
{\TC(A|K;p)} &\stck{\gamma_K} &
{\TC^{\hbox{{\eightpoint \'et}}}(A|K;p)}, \end{array}
$$
the left-hand vertical map induces an isomorphism on homotopy groups
with $\Z/p$-coefficients in degrees $\geq 1$. This follows from
Addendum~\ref{localizationtc} and \cite[Th.~D]{hm}. We use
Theorem~\ref{maintc} to prove that the right-hand vertical map induces
an isomorphism on homotopy groups with $\Z/p$-coefficients and that the
lower horizontal map induces an isomorphism on homotopy groups with
$\Z/p$-coefficients in degrees $\geq 0$.

We first prove the statement for the map induced from the cyclotomic
trace
$$K^{\hbox{{\eightpoint \'et}}}(K)\to\TC^{\hbox{{\eightpoint \'et}}}(A|K;p).$$
The spectral sequence \eqref{dss} for $K$-theory with
$\Z/p$-coefficients takes the form
$$E^2_{s,t}=H^{-s}(K,\mu_p^{\otimes(t/2)})\Rightarrow
K_{s+t}^{\hbox{{\eightpoint \'et}}}(K,\Z/p).$$
Indeed, since $K$-theory commutes with filtered colimits, this follows
from
$$K_t(\bar K,\Z/p)=\mu_p^{\otimes (t/2)},$$
which is proved in Suslin's celebrated paper~\cite{suslin} or follows
from Theorem~\ref{maintc} above. Similarly, it follows also
from Theorem \ref{maintc} that   the spectral sequence~\eqref{dss}
for topological cyclic homology takes the form
$$E^2_{s,t}=H^{-s}(K,\mu_p^{\otimes(t/2)})
\Rightarrow \TC_{s+t}^{\hbox{{\eightpoint \'et}}}(A|K;p,\Z/p).$$
Finally, it is clear that the cyclotomic trace induces an isomorphism
of\break $E^2$-terms.

It remains to show that the map
$$\gamma_K\:\TC_i(A|K;p,\Z/p)\to\TC_i^{\hbox{{\eightpoint \'et}}}(A|K;p,\Z/p)$$
is an isomorphism for $i\geq 0$. The domain and range of $\gamma_K$
are abstractly isomorphic by Theorem 6.1.6.  We must   show that
$\gamma_K$ is an isomorphism for $i\geq 0$. By theorem
\ref{tamedescent} we may assume that $\mu_p\subset K$ and that the
residue field $k$ is algebraically closed. When $\mu_p\subset K$, we
have a commutative square
$$\begin{array}{ccc}
{\TC_{\e}(A|K;p,\Z/p)\otimes\mu_p^{\otimes s}} 
&\stck{\gamma_K\otimes\id}  &
{\TC_{\e}^{\hbox{{\eightpoint \'et}}}(A|K;p,\Z/p)\otimes\mu_p^{\otimes s}}  \\[4pt]
\scs{\sim}&&\scs{\sim}\\[4pt]
{\TC_{2s+\e}(A|K;p,\Z/p)} &\stck{\gamma_K\otimes\id} &
{\TC_{2s+\e}^{\hbox{{\eightpoint \'et}}}(A|K;p,\Z/p),} \end{array}
$$
and the vertical maps are isomorphisms for $s\geq 0$ and $0\leq\e\leq
1$. Hence, it suffices to show that $\gamma_K$ is an isomorphism in
degrees $0$ and $1$. And for $k$ algebraically closed, the term
$H^2(K,\mu_p)\xto{\sim}H^1(k,\Z/p)$ in degree zero vanishes. Thus the
edge homomorphism of the spectral sequence~\eqref{dss},
$$\e_K\:\TC_0^{\hbox{{\eightpoint \'et}}}(A|K;p,\Z/p)\to H^0(K,\Z/p),$$
is an isomorphism, and since the composite
$$\TC_0(A|K;p,\Z/p) \xto{\gamma_K} \TC_0^{\hbox{{\eightpoint \'et}}}(A|K;p,\Z/p)
\xto{\e_K} H^0(K,\Z/p\Z)$$
is an isomorphism, then so is $\gamma_K$. In degree one, we use the
spectral sequence~\eqref{dss} for topological cyclic homology with
$\Qp/\Zp$-coefficients. As a $G_K$-module
$${\displaystyle\mathop{\rm
lim}_{\longrightarrow\atop{L/K}}} {\TC_1(B|\,L;p,\Qp/\Zp)}\stackrel{\sim}{\longleftarrow}
{\displaystyle\mathop{\rm lim}_{\longrightarrow\atop{L/K}}}{K_1(L,\Qp/\Zp)}\xto{\sim}
K_1(\bar K,\Qp/\Zp)=\mu_{p^\infty},$$
and the composite
$$\TC_1(A|K;p,\Qp/\Zp)\xto{\gamma_K}
\TC_1^{\hbox{{\eightpoint \'et}}}(A|K;p,\Qp/\Zp)\xto{\e_K}H^0(K,\mu_{p^\infty})$$
is an isomorphism. It follows that $\gamma_K$ is an isomorphism in
degree one.
\enddemo 
 
\centerline{\bf Appendix A. Truncated polynomial algebras}
\vglue18pt
A.1.\quad  Let $\pi=\pi_K\in A$ be a uniformizer
and let $e=e_K$ be the ramification index. Then $A/pA=k[\pi]/(\pi^e)$.
The structure of the topological Hochschild spectrum of this
$k$-algebra was examined in \cite{hm1}. We recall the result.

Let $\Pi=\Pi_e$ be the pointed monoid $\{0,1,\pi,\dots,\pi^{e-1}\}$
with base-point $0$ and with $\pi^e=0$ such that $A/p$ is the pointed
monoid algebra $k(\Pi)=k[\Pi]/k\{0\}$. Then we have from \cite[Th.\ 7.1]{hm} a natural ${\cal F}$-equivalence
of\break ${\Bbb T}$-spectra
$$T(k)\wedge |N\s^\cy(\Pi)| \xto{\sim} T(k(\Pi))$$
defined as follows: Let $C^b({\cal P}_{k(\Pi)})$ be the
category of bounded complexes of finitely generated projective
$k(\Pi)$-modules and consider $\Pi$ as a category with a single object
and endomorphisms $\Pi$. The functor $\Pi\to
C^b({\cal P}_{k(\Pi)})$, which takes the unique object to $k(\Pi)$
viewed as a complex concentrated in degree zero and which takes
$\pi^i\in\Pi$ (resp.\ $0\in\Pi$) to multiplication by $\pi^i\in k(\Pi)$
(resp.\ $0\in k(\Pi)$), induces
$$|N\s^\cy(\Pi)| \to |N\s^\cy(C^b({\cal P}_{k(\Pi)}))|
= T(k(\Pi))_{0,0},$$
and then the desired map is given as the composite
$$T(k)\wedge |N\s^\cy(\Pi)| \to T(k(\Pi))\wedge T(k(\Pi))
\xto{\mu} T(k(\Pi)).$$
Since $k$ and $\Pi$ are commutative, the equivalence is
multiplicative with com\-ponent-wise multiplication on the left. In
particular, the induced map on homotopy groups is an isomorphism of
differential graded $k$-algebras
$$\pi_*(T(k)\wedge |N\s^\cy(\Pi)|) \xto{\sim} \pi_*T(k(\Pi)),$$
where the differential is given by Connes' operator~\eqref{connes}. We 
give the realization $|N\s^\cy(\Pi)|$ the usual CW-structure,
\cite[Th.\ 14.1]{may1} (with the simplices $\Delta^n$ and the
disks $D^n$ identified through a compatible family of orientation-preserving homeomorphisms). Then the
skeleton filtration gives a spectral sequence of differential graded $k$-algebras
$$E_{s,t}^2=\pi_tT(k)\otimes\tilde H_s(|N\s^\cy(\Pi)|;k)
\Rightarrow \pi_{s+t}(T(k)\wedge |N\s^\cy(\Pi)|).$$
The same statements are true for ordinary Hochschild homology. If $k$
is a perfect field of characteristic $p>0$, $\pi_*\HH(k)=k$
concentrated in degree zero (see e.g. \cite[Lemma  5.5]{hm}). Hence,
the spectral sequence collapses and the edge homomorphism gives an
isomorphism of differential graded $k$-algebras
\begin{equation}\label{hhiso}
\pi_*(\HH(k)\wedge |N\s^\cy(\Pi)|) \xto{\sim}
\tilde H_*(|N\s^\cy(\Pi)|;k). \speqnua{1.1}
\end{equation}
The spectral sequence also collapses for $T(k)$. Indeed, the inclusion
of the zero-skeleton gives a map of ring spectra $H(k)\to T(k)$ from
the Eilenberg-MacLane spectrum for $k$, so we have a multiplicative map
\begin{equation}\label{multiso}
\pi_*T(k)\otimes\tilde H_*(|N\s^\cy(\Pi)|;k) \xto{\sim}
\pi_*(T(k)\wedge |N\s^\cy(\Pi)|)\hskip.5in \speqnua{1.2}
\end{equation}
given as the composite of the external product
$$\pi_*T(k)\otimes \pi_*(H(k) \wedge |N\s^\cy(\Pi)|) \xto{\wedge}
\pi_*(T(k)\wedge H(k)\wedge |N\s^\cy(\Pi)|)$$
and the map induced from $\mu\:T(k)\wedge H(k)\to T(k)$. It follows
that the spectral sequence collapses and that this map is an
isomorphism of graded $k$-algebras. However, the map $H(k)\to T(k)$
is not equivariant, so this isomorphism does not preserve the
differential.

Let $N_*(k(\Pi))$ be the normalized standard complex, \cite[Chap.\ IX,
\S7]{ce}. The K\"unneth isomorphism determines an isomorphism of
complexes
$$k(\Pi)\otimes_{k(\Pi)^e}N_*(k(\Pi)) \xto{\sim}
\tilde C_*(|N\s^\cy(\Pi)|;k),$$
and since $N_*(k(\Pi))\xto{\mu}k(\Pi)$ is a resolution of $k(\Pi)$ by
free $k(\Pi)^e$-modules, we have a canonical isomorphism of graded
$k$-algebras
$$\Tor_*^{k(\Pi)^e}(k(\Pi),k(\Pi)) \xto{\sim}
\tilde H_*(|N\s^\cy(\Pi)|;k).$$
To evaluate this, we consider instead the resolution
$R_*(k(\Pi)) \xto{\e} k(\Pi)$ of \cite{bag},
\begin{eqnarray*}
R_*(k(\Pi)) & =& k(\Pi)^e\otimes\Lambda\{c_1\}\otimes\Gamma\{c_2\}, \\
\delta(c_1) & = &\pi\otimes 1-1\otimes\pi, \hskip6mm
\delta(c_2^{[d]}) = \frac{\pi^e\otimes 1-1\otimes\pi^e}{\pi\otimes
1-1\otimes\pi} \cdot c_1c_2^{[d-1]},  
\end{eqnarray*}
where $\Gamma\{c_2\}$ is a divided power algebra and $c_2^{[d]}$
the $d$-th divided power of $c_2$. An augmentation-preserving chain map
$g\:R_*(k(\Pi))\to N_*(k(\Pi))$ is given by
\begin{eqnarray*}
g(c_2^{[d]}) & = &\sum 1\otimes x^{k_0}\otimes x \otimes x^{k_1}\otimes
\dots\otimes x\otimes x^{k_d}, \\
g(c_1c_2^{[d]}) & =& \sum 1\otimes x\otimes x^{k_0}\otimes 
\dots \otimes x\otimes x^{k_d}, 
\end{eqnarray*}
where both sums run over tuples $(k_0,\dots,k_d)$ with
$k_0+\cdots+k_d=d(e-1)$ and $0\leq k_i<e$. (The summands with some
$k_i=0$, for $0\leq i<d$, are zero.) Hence, if $e$ annihilates $k$,
we have an isomorphism of differential graded $k$-algebras
\begin{equation}\label{dgaiso}
k(\Pi)\otimes\Lambda\{c_1\}\otimes\Gamma\{c_2\}
\xto{\sim}\tilde H_*(|N\s^\cy(\Pi)|;k), \speqnua{1.3}
\end{equation}
where $d\pi=c_1$ and $dc_2^{[d]}=0$. The value of the differential is
readily verified using the standard formula, \cite[Prop.\ 1.4.6]{h}. 

\specialnumber{A.1.4}\proclaim{Proposition} \label{TparkparPi}
Let $k$ be a perfect field of
characteristic $p>0$ and suppose $p$ divides $e$. Then there is a
canonical isomorphism of differential graded $k$\/{\rm -}\/algebras
$$S\{\sigma\}\otimes k(\Pi)\otimes\Lambda\{c_1\}\otimes\Gamma\{c_2\}
\xto{\sim}\pi_*T(k(\Pi)),$$
where $d\pi=c_1$ and $d(c_2^{[d+1]})=-(e/p)\pi^{e-1}c_1c_2^{[d]}\sigma$.
\endproclaim

\demo{Proof} The map of the statement is given by the
maps~(\ref{multiso}) and~(\ref{dgaiso}). Since both are isomorphisms
of graded $k$-algebras, it remains only to verify the differential
structure. The formula for $d\pi$ is clear since the edge homomorphism
$$\pi_q(T(k)\wedge |N\s^\cy(\Pi)|) \to \tilde H_q(|N\s^\cy(\Pi)|;k)$$
is an isomorphism for $q\leq 1$ and commutes with the differential.
But the proof of the formula for $dc_2^{[d]}$ is more involved and
uses the calculation in~\cite[Th.\ B]{hm1} of the homotopy type of
the ${\Bbb T}$-CW-complex $|N\s^\cy(\Pi)|$. As cyclic sets
\begin{equation}\label{decomp}
N\s^\cy(\Pi)=\bigvee_{s\geq 0}N\s^\cy(\Pi;s), \speqnua{1.5}
\end{equation}
where the $s$-th summand has $n$-simplices
$(\pi^{i_0},\dots,\pi^{i_n})$ with $i_0+\dots i_n=s$, and the
realization decomposes accordingly. If we write $s=de+r$ with
$0<r\leq e$ then under the isomorphism of the statement
$$\pi_*(T(k)\wedge|N\s^\cy(\Pi;s)|) \cong
\left\{ \begin{array}{ll}
S\{\sigma\}\otimes k\{\pi^rc_2^{[d]},\pi^{r-1}c_1c_2^{[d]}\}, &
\hbox{if $0<r<e$,} \\
S\{\sigma\}\otimes k\{\pi^{e-1}c_1c_2^{[d]},c_2^{[d+1]}\}, &
\hbox{if $r=e$.} \\
\end{array}\right.$$
The formula we wish to prove involves the case $r=e$. In this case,
\cite[Th.\ B]{hm1} gives a canonical triangle of
${\Bbb T}$-CW-complexes
$${\Bbb T}/C_{(d+1)+}\wedge S^{V_d} \xto{\pr} 
{\Bbb T}/C_{s+}\wedge S^{V_d} \xto{i}
|N\s^\cy(\Pi;s)| \xto{\partial}
\Sigma{\Bbb T}/C_{(d+1)+}\wedge S^{V_d},$$
where $V_d=\C(1)\oplus\dots\oplus\C(d)$. 
If we form the smash product with $T(k)$ and take homotopy groups, the
triangle gives rise to a long-exact sequence, which we now
describe. Let $x_0$ (resp.\ $y_0$) be the class of the $0$-cycle
$C_{d+1}/C_{d+1}$ (resp.\ $C_s/C_s$) and let $x_1$ (resp.\ $y_1$,
resp.\ $z_{2d}$) be the fundamental class of ${\Bbb T}/C_{d+1}$
(resp.\ ${\Bbb T}/C_s$, resp.\ $S^{V_d}$). Then
$$\pi_*(T(k)\wedge{\Bbb T}/C_{n+}\wedge S^{V_d}) \cong
\left\{ \begin{array}{ll}
S\{\sigma\}\otimes k\{x_0z_{2d},x_1z_{2d}\}, &
\hbox{if $n=d+1$,} \\
S\{\sigma\}\otimes k\{y_0z_{2d},y_1z_{2d}\}, &
\hbox{if $n=s$,} \\
\end{array}\right.$$
and the differential is $\pi_*T(k)$-linear and maps
$$\begin{array}{rlrl}
d(y_0z_{2d}) &\hsm = (d+1)y_1z_{2d}, &\qquad
d(y_1z_{2d}) &\hsm= 0, \\
d(x_0z_{2d}) & \hsm= sx_1z_{2d}, &\qquad
d(x_1z_{2d}) &\hsm= 0. \end{array}
$$
The induced maps in the long-exact sequence of homotopy groups
associated with the triangle above all are $\pi_*T(k)$-linear and
$$\begin{array}{rlrl}
\pr_*(y_0z_{2d}) & \hsm= x_0z_{2d}, &\qquad 
\pr_*(y_1z_{2d}) & \hsm= ex_1z_{2d}, \\[3pt]
i_*(x_0z_{2d}) & \hsm = 0, &\qquad
i_*(x_1z_{2d}) & \hsm= \pi^{e-1}c_1c_2^{[d]}, \\[3pt]
\partial_*(\pi^{e-1}d\pi\cdot c_2^{[d]}) & \hsm = 0, &\qquad
\partial_*(c_2^{[d+1]}) & \hsm = -y_1z_{2d}. \end{array}
$$
The statements for the maps $\pr_*$ and $i_*$ are clear from the
construction of the triangle in \cite{hm1}. We verify the statement
for the map $\partial_*$. To this end we first choose a cellular
homotopy equivalence 
$$\alpha\:C_{\pr} \xto{\sim} |N\s^\cy(\Pi;s)|$$
such that we have a map of triangles from the distinguished triangle
given by the map $\pr$ to the triangle above. Since the cellular
chain functor carries distinguished triangles of CW-complexes to
distinguished triangles of chain complexes, we have
\begin{eqnarray*}
\partial_*(\alpha_*((0,y_1z_{2d}))) & =& y_1z_{2d}, \\
\alpha_*((x_1z_{2d},0)) & =& \pi^{e-1}c_1c_2^{[d]}.  
\end{eqnarray*}
Hence, it suffices to show that $\alpha_*((0,y_1z_{2d}))$ is homologous to
$-c_2^{[d+1]}$. To do this, we consider the diagram
$$\begin{array}{ccc}
{\tilde H_{2d+2}(|N\s^\cy(\Pi;s)|;\Z/p)\;}&\stackrel{\beta}{\hookrightarrow} &
{\tilde H_{2d+1}(|N\s^\cy(\Pi;s)|;\Z)} \\[4pt]
\scr{\alpha_\ast}\big\uparrow\scr{\sim}&&\scr{\alpha_\ast}\big\uparrow\scr{\sim}\\[4pt]
{\tilde H_{2d+2}(C_{\pr};\Z/p)\;} 
&\stackrel{\beta}{\hookrightarrow}  &
{\tilde H_{2d+1}(C_{\pr};\Z)} \end{array}
$$
with injective horizontal maps. A straightforward calculation shows
that (on the level of chains) the top Bockstein takes $c_2^{[d+1]}$ to
$(e/p)\pi^{e-1}c_1c_2^{[d]}$ and the bottom Bockstein takes
$(0,y_1z_{2d})$ to $-(e/p)x_1z_{2d}$. We have already noted that the
right-hand vertical map takes $(x_1z_{2d},0)$ to
$\pi^{e-1}c_1c_2^{[d]}$. This completes the proof of the stated
formula for $\partial_*$.

We now prove the formula for $d(c_2^{[d]})$. First note that we can
write
$$d(c_2^{[d]})=d_1(c_2^{[d]})+d_2(c_2^{[d]}),$$
where $d_1$ (resp.\ $d_2$) is defined in same way as $d$ but with
${\Bbb T}$ acting in the first (resp.\ second) smash factor of
$T(k)\wedge |N\s^\cy(\Pi;s)|$ only. Since the differential $d_2$
commutes with the isomorphism
$$\pi_*T(k)\otimes\tilde H_*(|N\s^\cy(\Pi;s)|;k)
\xto{\sim} \pi_*(T(k)\wedge |N\s^\cy(\Pi;s)|),$$
we find that $d_2(c_2^{[d]})=0$. Hence, we can ignore the
${\Bbb T}$-action on $|N\s^\cy(\Pi;s)|$. We have a map of
  triangles of (nonequivariant) CW-complexes
$$\begin{array}{ccccccc}
{{\Bbb T}/C_{(d+1)+}\wedge S^{V_d}} &\hsm\shtck{\pr}  \hsm&
{{\Bbb T}/C_{s+}\wedge S^{V_d}}& \shtck{i} \hsm&
{C_{\pr}} &\hskip-.25in\hsm\shtck{\partial}\hsm&\hskip-.2in
{\Sigma{\Bbb T}/C_{(d+1)+}\wedge S^{V_d}}\\[4pt]
\scs{f}&\hsm\hsm &\scs{g} &\hsm\hsm&\scr{\sim}\scs{h} &\hsm\hsm&\hskip-.25in\scs{\Sigma
f}\\[4pt] {S^{2d+1}}&\hsm\stck{\Sigma^{2d+1}e}\hsm &
{S^{2d+1}} &\hskip-.25in\stck{\Sigma^{2d+1}i} \hsm&
{\Sigma^{2d+1}M_e} &\hsm\stck{-\Sigma^{2d+1}\beta} \hsm&
\hskip-.3in{S^{2d+2},} \end{array}
$$
such that $f_*$ (resp.\ $g_*$) maps $x_1z_{2d}$ (resp.\ $y_1z_{2d}$) to
the fundamental class of $S^{2d+1}$. Hence, it suffices to show that
the image of $-h_*((0,y_1z_{2d}))=1\cdot\susp(\e)$ under
$$d\:\pi_{2q+2}(T(k)\wedge\Sigma^{2d+1}M_e) \to
\pi_{2q+3}(T(k)\wedge\Sigma^{2d+1}M_e)$$
is equal to $-(e/p)\sigma\cdot\susp(1)=-(e/p)h_*((x_1z_{2d},0))$. To this
end, we consider the diagram
$$\begin{array}{ccccc}
{\pi_1(M_e\wedge T(k))} &\hsm\stck{\susp}\hsm &
{\pi_{2d+2}(\Sigma^{2d+1}M_e\wedge T(k))}&\hsm\stck{\tw_*}\hsm &
{\pi_{2d+2}(T(k)\wedge\Sigma^{2d+1}M_e)}   \\[4pt]
\scs{d} &\hsm\scr{(-1)}\hsm&\hsm\scs{d}\hsm &&\hsm\scs{d} 
 \\[4pt]
{\pi_2(M_e\wedge T(k))} &\hsm\stck{\susp}\hsm &
{\pi_{2d+3}(\Sigma^{2d+1}M_e\wedge T(k))}&\hsm\stck{\tw_*} \hsm&
{\pi_{2d+3}(T(k)\wedge\Sigma^{2d+1}M_e),} \end{array}$$
which commutes up to the indicated sign. By the definition of the class
$\sigma$, the left-hand vertical map takes $\e\cdot 1$ to
$(e/p)1\cdot\sigma$. Hence, the right-hand vertical map takes
$1\cdot\susp(\e)$ to $-(e/p)\sigma\cdot\susp(1)$. The stated formula
for $d(c_2^{[d+1]})$ follows.
\enddemo

\specialnumber{A.1.6}\proclaim{Addendum} \label{tatepAfpp}
The nonzero differentials in the
spectral sequence
\begin{eqnarray*}
\hat E^2(C_{p^n},T(k(\Pi))) & =&
\Lambda\{u_n,c_1,\e\}\otimes S\{t^{\pm1},\sigma,\pi\}/(\pi^e)
\otimes\Gamma\{c_2\} \\
{} & \Rightarrow&\bar\pi_*(\tate(C_{p^n},T(k(\Pi))))  
\end{eqnarray*}
are generated from $d^2\e=t\sigma${\rm ,} $d^2\pi=tc_1${\rm ,}
and $d^2c_2^{[d+1]}=-(e/p)t\pi^{e-1}c_1c_2^{[d]}\sigma$.
\endproclaim

\demo{Proof} The $d^2$-differential is given by
Propositions~\ref{d2=connes} and~\ref{TparkparPi}. It remains only to
show that the higher differentials $d^r$, $r\geq 3$, vanish. 
The decomposition of cyclic sets~(\ref{decomp}) induces one of
spectral sequences. And if we write $s=de+r$ with $0<r\leq e$, then
the $E^3$-term of the $s$-th summand is concentrated on the lines
$E_{*,d}^3$ and $E_{*,d+1}^3$, if $0<r<e$, and on the lines
$E_{*,d+1}^3$ and $E_{*,d+2}^3$, if $r=e$. In either case, all further
differentials must be zero for degree reasons.
\enddemo

\specialnumber{A.1.7}\proclaim{Proposition} \label{gammahatpilemma}
Let $n\leq v_p(e)$. The images of
$\ul{\pi}_n$ and $\ul{\pi}_n^{e/p^n}$ by the map
$$\hat\Gamma\:\bar\pi_*(T(k(\Pi))^{C_{p^{n-1}}}) \to 
\bar\pi_*(\tate(C_{p^n},T(k(\Pi)))) $$
are represented in the spectral sequence $\hat E^*(C_{p^n},T(k(\Pi)))$
by the infinite cycles $\pi^{p^n}$ and $tc_2${\rm ,} if $v_p(e)>n${\rm ,} and
by $\pi^{p^n}$ and $-(e/p^n)u_1\pi^{e-1}c_1${\rm ,} if $v_p(e)=n$.
\endproclaim

\demo{Proof} The statement only involves the summand
$|N\s^\cy(\Pi,e)|$. We consider the map of spectral sequences induced
from the linearization map,
$$l_*\:\hat E^*(C_{p^n},T(k)\wedge |N\s^\cy(\Pi,e)|) \to
\hat E^*(C_{p^n},\HH(k)\wedge |N\s^\cy(\Pi,e)|).$$
In the left-hand spectral sequence, $E^3=E^\infty$, and in the right-hand spectral sequence, $E^2=E^\infty$. The induced map of
$E^\infty$-terms may be identified with the canonical inclusion
\begin{eqnarray*}
{} &&\hskip-24pt \Lambda\{u_n\}\otimes S\{t^{\pm1}\}\otimes
k\{\pi^{e-1}c_1,c_2+\e\cdot(e/p)\pi^{e-1}c_1\} \\
{} &&\qquad  \hookrightarrow \Lambda\{u_n,\e\}\otimes S\{t^{\pm1}\}\otimes
k\{\pi^{e-1}c_1,c_2\}.  
\end{eqnarray*}
Since the map is injective, it suffices to show that
$l_*(\hat\Gamma(\ul{\pi}_n^{e/p^n}))$ is represented in the sequence
on the right by $-u_n\pi^{e-1}c_1$, if $v_p(e)=n$, and by $tc_2$, if
$v_p(e)>n$. In the proof of this, we shall use the notation and
results of Sections~4.2 and~4.3 above.

We have from~\cite[\S1]{bhm} the ${\Bbb T}$-equivariant
homeomorphism
$$D\:|\sd_{p^n}N\s^\cy(\Pi,e)| \xto{\sim} |N^\cy(\Pi,e)|,$$
where on the left, the action by the subgroup
$C_{p^n}\subset{\Bbb T}$ is induced from a simplicial
$C_{p^n}$-action. It follows that this space has a canonical
$C_{p^n}$-CW-structure, and the homeomorphism $D$ then defines a
$C_{p^n}$-CW-structure on $|N\s^\cy(\Pi,e)|$. We fix, as in the
proof of Proposition~\ref{TparkparPi}, a cellular homotopy equivalence
$$\alpha\:C_{\pr} \xto{\sim} |N\s^\cy(\Pi,e)|$$
with the $C_{p^n}$-CW-structure on $C_{\pr}$ induced from the
$C_{p^n}$-CW-structure of ${\Bbb T}=S(\C)=E_1$ given in
Section~4.4 above. The cellular complex $C_*=\tilde
C_*(C_{\pr};k)$ is canonically identified with the complex
$$k[C_{p^n}]\cdot (0,x_1) \xto{\delta}
k\cdot (x_1,0)\oplus k[C_{p^n}]\cdot (0,x_0) \xto{\delta}
k\cdot(x_0,0),$$
where $\delta((0,x_1))=-(e/p^n)(x_1,0)-(g-1)(0,x_0)$,
$\delta((x_1,0))=0$, and $\delta((0,x_0))=-(x_0,0)$. One shows as in
the proof of Proposition~\ref{TparkparPi} that the cycles
$\alpha_*((x_1,0))$ and $\alpha_*(N(0,x_1))$ represent the classes
$\pi^{e-1}c_1$ and $-c_2$, respectively.
\vglue9pt

We now turn to the spectral sequence
$\hat E^*=\hat E^*(C_{p^n},\HH(k)\wedge C_{\pr})$. There are canonical
isomorphisms of complexes
$$\hat E_{*,t}^1 \cong
(\tilde P\otimes\Hom(P,\bar\pi_t(\HH(k)\wedge C_{\pr})))^{C_{p^n}}
\cong (\tilde P\otimes\Hom(P,\bar H_t(C_*)))^{C_{p^n}}$$
with the left-hand isomorphism given by Lemma~\ref{E^1} and the right-hand isomorphism by~(\ref{hhiso}). We claim that in fact
\begin{equation}\label{totalcx}
\bar\pi_*(\tate(C_{p^n},\HH(k)\wedge C_{\pr})) 
\cong \bar H_*((\tilde P\otimes\Hom(P,C_*))^{C_{p^n}})\hskip.5in \speqnua{1.8}
\end{equation}
and that the spectral sequence $\hat E^*$ is canonically isomorphic to
the one associated with the double complex on the right. To see this,
we filter $M_p$, $\tilde E$, $E$, and $C_{\pr}$ by the skeletons. We
get, as in Section~4.3, a conditionally convergent spectral
sequence
$$
E_{s,t}^2 = \bar H_s((\tilde P\otimes\Hom(P,\pi_t\HH(k)\otimes
C_*))^{C_{p^n}}) 
 \Rightarrow \bar\pi_{s+t}(\tate(C_{p^n},\HH(k)\wedge C_{\pr})), $$
which collapses since $\pi_t\HH(k)$ vanishes for $t>0$. The
edge homomorphism gives the desired isomorphism. Moreover, under this
isomorphism, the filtration of $\tilde E$ and $E$, which gives rise to
the spectral sequence $\hat E^*$, corresponds to the filtration of the
complexes $\tilde P$ and $P$. Tracing through the definitions, one
readily sees that the class $l_*(\hat\Gamma(\ul{\pi}_n^{e/p^n}))$ is
represented by the element $y_0\otimes Nx_0^*\otimes (x_0,0)\in\hat
E_{0,0}^1$. To finish the proof, we note that in the total
complex~(\ref{totalcx}),
\begin{eqnarray*}
&& \hskip-22pt \delta(N(y_0\otimes x_1^*\otimes (0,x_1)-
y_0\otimes x_0^*\otimes (0,x_0))) 
\\&&  = y_0\otimes Nx_0^*\otimes (x_0,0) +y_0\otimes Nx_2^*\otimes
N(0,x_1)+(e/p^n)y_0\otimes Nx_1^*\otimes (x_1,0),  
\end{eqnarray*}
and in the lower line, the first summand represents
$l_*(\hat\Gamma(\ul{\pi}_n^{e/p^n}))$, the second $-tc_2$, and the
third $(e/p^n)u_n\pi^{e-1}c_1$. The statement follows, since 
$-tc_2$ and $u_n\pi^{e-1}c_1$ are not boundaries.
\enddemo


\begin{references}

\bibitem{boardman}
\name{J.~M.\ Boardman}, Conditionally convergent spectral sequences, preprint,
  available at  {hopf.math.purdue.edu},
  1981.

\bibitem{bokstedt}
\name{M.~B\"okstedt}, Topological {H}ochschild homology, preprint, Bielefeld
  1985. 

\bibitem{bhm}
\name{M.~B\"okstedt, W.-C.\ Hsiang,} and \name{I.~Madsen}, The cyclotomic trace and
  algebraic {$K$}-theory of spaces, {\it Invent.\ Math\/}.\  {\bf 111} (1993),
  465--540.

\bibitem{bm}
\name{M.~B\"okstedt} and \name{I.~Madsen}, Topological cyclic homology of the
  integers, in {$K$}-{\it theory\/} (Strasbourg, 1992), {\it
  Ast\'erisque\/} {\bf 226\/} (1994),
  57--143.


\bibitem{ce}
\name{H.~Cartan} and \name{S.~Eilenberg}, {\it Homological Algebra\/}, Princeton
  Univ.\ Press, Princeton, NJ, 1956.

\bibitem{dundas}
B.~I.\ Dundas, {$K$}-theory theorems in topological cyclic homology, 
{\it J.\
  Pure Appl.\ Algebra\/} {\bf 129} (1998), 23--33.

\bibitem{dm}
\name{B.~I. Dundas} and \name{R.~McCarthy}, Topological {H}ochschild homology of ring
  functors and exact categories, {\it J.\ Pure Appl.\ Alg\/}.\  
  {\bf 109} (1996),
  231--294.

\bibitem{dwyermitchell}
\name{W.~G.\ Dwyer} and \name{S.~A.\ Mitchell}, On the {$K$}-theory spectrum of a ring of
  algebraic integers, {$K$}-{\it theory\/} {\bf 14} (1998), 201--263.

\bibitem{frohlich}
\name{A.~Fr\"ohlich}, {\it Galois Module Structure of Algebraic Integers},
{\it Ergebnisse
  der Mathematik und ihrer Grenzgebiete\/} {\bf 1}, Springer-Verlag,
  New York, 1983.

\bibitem{gh}
\name{T.~Geisser} and \name{L.~Hesselholt}, Topological cyclic homology of schemes,
  in {$K$}-{\it theory\/} (Seattle, WA, 1997), {\it Proc.\ Sympos.\  Pure 
  Math\/}.\ {\bf 67} (1999), 41--87.

\bibitem{greenlees}
\name{J.~P.~C.\ Greenlees}, Representing {T}ate cohomology of {$G$}-spaces,
  {\it Proc.\ Edinburgh Math.\ Soc\/}.\  {\bf 30} (1987), 435--443.

\bibitem{greenleesmay}
\name{J.~P.~C.\ Greenlees} and \name{J.~P.\ May}, {\it Generalized {T}ate Cohomology}, 
{\it Mem.\
  Amer.\ Math.\ Soc\/}.\  {\bf 113} (1995), A.\ M.\ S., Providence,
  RI.

\bibitem{ega4}
\name{A.~Grothendieck}, {\it {\rm \'{\it E}}l{\rm \'{\it e}}ments de g{\rm \'{\it e}}om{\rm \'{\it e}}trie
alg{\rm \'{\it e}}brique.
  {\rm IV}. {{\rm \'{\it E}}}tude locale des sch{\rm \'{\it e}}mas et des morphismes de 
sch{\rm \'{\it e}}mas{\rm ,}
{\it {I}nst.\ {H}autes {\'E}tudes {S}ci.\ {P}ubl.\ {M}ath\/}}.\  
  {\bf 32} (1967).

\bibitem{bag}
\name{J.~A. Guccione, M.~J.~Redondo J.~J.~Guccione, A.~Solotar}, and \name{O.~E.\ Villamayor},
   Cyclic homology of algebras with one generator, $K$-{\it
  theory\/} {\bf 5}
  (1991), 51--68.

\bibitem{harrissegal}
\name{B.~Harris} and \name{G.~Segal}, {$K_i$} groups of rings of algebraic integers,
  {\it Ann.\ of Math\/}.\  {\bf 101} (1975), 20--33.

\bibitem{h}
\name{L.~Hesselholt}, On the $p$-typical curves in {Q}uillen's {$K$}-theory,
  {\it Acta Math\/}.\  {\bf 177} (1997), 1--53.

\bibitem{hm3}
\name{L.~Hesselholt} and \name{I.~Madsen}, On the de {R}ham-{W}itt complex in mixed
  characteristic, {\it Ann.\ Sci.\ {\rm \'{\it E}}cole Norm.\ Sup.}, to appear.

\bibitem{hm1}
\bibline, Cyclic polytopes and the {$K$}-theory of truncated polynomial
  algebras, {\it Invent.\ Math\/}.\  {\bf 130} (1997), 73--97.

\bibitem{hm}
\bibline, On the {$K$}-theory of finite algebras over {W}itt vectors of
  perfect fields, {\it Topology\/} {\bf 36} (1997), 29--102.
 
\bibitem{hkr}
\name{G.~Hochschild, B.~Kostant}, and \name{A.~Rosenberg}, Differential forms on
  regular affine algebras, {\it Trans.\ Amer.\ Math.\ Soc\/}.\  
  {\bf 102} (1962), 383--408.

\bibitem{hovey}
\name{M.~Hovey}, Spectra and symmetric spectra in general model categories, 
  {\it J. Pure Appl.\ Algebra\/} {\bf 156} (2001), 63--127.

\bibitem{hoveypalmieristrickland}
\name{M.~Hovey, J.~H. Palmieri}, and \name{N.~P.\ Strickland}, {\it Axiomatic Stable
  Homotopy Theory}, {\it Mem.\ Amer.\ Math.\ Soc\/}.\  {\bf 128} (1997).

\bibitem{hyodokato}
\name{O.~Hyodo} and \name{K.~Kato}, Semi-stable reduction and crystalline cohomology
  with logarithmic poles, in {\it P\'eriodes $p$-adiques\/}
  (Bures-sur-Yvette, 1988),  {\it Ast\'erisque\/} {\bf 223},  1994,
  221--268.

\bibitem{k}
\name{K.~Kato}, Logarithmic structures of {F}ontaine-{I}llusie, {\it Algebraic
  Analysis, Geometry, and Number Theory\/}, {\it Proc.\  JAMI Inaugural
  Conference\/} (Baltimore, 1988), Johns Hopkins Univ.\ Press, 
  Baltimore, MD, 1989,
 191--224.

\bibitem{lms}
\name{L.~G.\ Lewis, J.~P.\ May}, and \name{M.~Steinberger}, {\it Equivariant Stable Homotopy
  Theory}, {\it Lecture Notes in Math\/}.\ {\bf 1213}, Springer-Verlag,
  New York, 1986.

\bibitem{lichtenbaum}
\name{S.~Lichtenbaum}, Values of zeta functions, {\'e}tale cohomology, and
  algebraic {$K$}-theory, in {\it Algebraic {$K$}-theory\/}, 
  {II} (Battelle Memorial
  Inst., Seattle, Washington, 1972), {\it Lecture Notes in Math\/}.\
  {\bf 342}, 
  Springer-Verlag, New York, 1973, 489--501.

\bibitem{lm}
\name{A.~Lindenstrauss} and \name{I.~Madsen}, {T}opological {H}ochschild homology of
  number rings, {\it Trans.\ Amer.\ Math.\ Soc\/}.\  {\bf 352} 
  (2000), 2179--2204.

\bibitem{loday}
\name{J.-L.\  Loday}, {\it Cyclic Homology\/} ({\it {A}ppendix {\rm E}\/} 
by M.\ O.\ Ronco),
  {\it Grundlehren der Mathematischen Wissenschaften\/}  {\bf 301}, 
  Springer-Verlag, New York, 
  1992.

\bibitem{maclane}
\name{S.~MacLane}, {\it Categories for the Working Mathematician\/}, {\it
Grad.\ Texts in Math\/}.\ 
  {\bf 5}, Springer-Verlag, New York, 1971.

\bibitem{m}
\name{I.~Madsen}, Algebraic $K$-theory and Traces, in {\it Current Developments in
  Mathematics 1995\/},  Internat.\  Press, Cambridge, MA, 1996, 191--321.

\bibitem{mandellmay}
\name{M.~A.\ Mandell} and \name{J.~P.\ May}, {\it Equivariant Orthogonal Spectra and
  {$S$}-modules\/}, {\it Mem.\ Amer.\ Math.\ Soc\/}.\ {\bf 159} (2002).

\bibitem{mat}
\name{H.~Matsumura}, {\it Commutative Ring Theory}, {\it Cambridge Studies in
Adv.\ 
  Math\/}.\  {\bf 8},  Cambridge Univ.\  Press, Cambridge,  U.K., 1986.

\bibitem{may1}
\name{J.~P. May}, Simplicial objects in algebraic topology, {R}eprint of the
  1967 original, {\it Chicago Lectures in Math\/}., Univ.\  of Chicago Press,
  Chicago, IL, 1992.

\bibitem{mc}
\name{R.~McCarthy}, The cyclic homology of an exact category, {\it J.\ Pure
Appl.\
  Alg\/}.\  {\bf 93} (1994), 251--296.

\bibitem{mumford1}
\name{D.~Mumford}, {\it Lectures on Curves on an Algebraic Surface\/}, {\it
Ann.\ of 
  Math.\  Studies\/} {\bf 59},  Princeton Univ.\  Press, Princeton, NJ,
  1966.

\bibitem{Quillen1}
\name{D.~Quillen}, On the cohomology and {$K$}-theory of the general linear
  group over a finite field, {\it Ann.\ of Math\/}.\  {\bf 96} (1972), 
  552--586.

\bibitem{quillen}
\bibline, Higher algebraic {$K$}-theory {I}, in {\it Algebraic 
{$K$}-theory I\/}:
  {\it Higher {$K$}\/}-{\it theories\/} 
  (Battelle Memorial Inst., Seattle, Washington, 1972),
  {\it Lecture Notes in Math\/}.\ {\bf 341}, Springer-Verlag, 
  New York, 1973.

\bibitem{quillen2}
\bibline, Higher algebraic {$K$}-theory, in {\it Proc.\ 
  {I}nternational {C}ongress of {M}athematicians\/} (Vancouver, B. C., 1974),
  vol.~1, Canad.\ Math.\ Congress, Montreal, Que., 1975, 171--176.

\bibitem{cf}
\name{J-P.\ Serre}, {Local class field theory}, in 
{\it Algebraic Number Theory\/}, Thompson, Washington, D.C., 
  1967, 128--161.

\bibitem{serre}
\bibline, {\it Local Fields}, {\it Grad.\  Texts in Math\/}.\ {\bf 67},
  Springer-Verlag, New York, 1979.

\bibitem{serre1}
\bibline, {\it Galois Cohomology}, Springer-Verlag, New  York, 1997.

\bibitem{soule}
\name{C.~Soul\'e}, On higher {$p$}-adic regulators, in {\it Algebraic
{$K$}-theory\/}
  ({E}vanston 1980), {\it Lecture Notes in Math\/}.\ {\bf 854\/},
  Springer-Verlag,
  New York, 1981, 372--401.

\bibitem{suslin}
\name{A.~A.\ Suslin}, On the {$K$}-theory of algebraically closed fields,
  {\it Invent.\ Math\/}.\  {\bf 73} (1983), 241--245.

\bibitem{suslin1}
\bibline, On the {$K$}-theory of local fields, {\it J.\ Pure Appl.\
Alg\/}.\ 
  {\bf 34} (1984), 304--318.

\bibitem{thomason}
\name{R.~W.\ Thomason}, Algebraic {$K$}-theory and \'etale cohomology, {\it
Ann.\
  Sci.\ \'Ecole Norm.\ Sup\/}.\  {\bf 13} (1985), 437--552.

\bibitem{tt}
\name{R.~W. Thomason} and \name{T.~Trobaugh}, Higher algebraic {$K$}-theory of schemes
  and of derived categories, {\it Grothendieck Festschrift\/}, Volume III, 
  {\it Progr.\  in
  Math\/}.\ {\bf 88} (1990), 247--435.

\bibitem{tsalidis}
\name{S.~Tsalidis}, Topological {H}ochschild homology and the homotopy descent
  problem, {\it Topology\/} {\bf 37} (1998), 913--934.

\bibitem{w}
\name{F.~Waldhausen}, Algebraic {$K$}-theory of spaces, in {\it Algebraic and 
Geometric
  Topology\/}, {\it Lecture Notes in Math\/}.\ {\bf  1126}, 
  Springer-Verlag, New  York, 1985,
  318--419.

\bibitem{weibel}
\name{C.~A.\ Weibel}, {\it An Introduction to Homological Algebra}, {\it Cambridge 
Studies
  in Adv.\ Math\/}.\ {\bf 38}, Cambridge Univ.\  Press, Cambridge, U.K.,
  1994.

\end{references}
\end{document}